\documentclass{SIAMbook2023} 
\usepackage{amsmath}  
\usepackage{epsfig}
\usepackage{graphicx}
\usepackage{makeidx}
\usepackage{multicol} 
\usepackage{hyperref}
\usepackage{mathrsfs} 
\usepackage{algorithmic} 
\usepackage{hyperref}
\usepackage{url}
\usepackage{float}
\usepackage{epstopdf}
\usepackage{subcaption}
\let\oldmarginpar\marginpar
\renewcommand\marginpar[1]{\-\oldmarginpar[\raggedleft\footnotesize #1]%
{\raggedright\footnotesize #1}}
\floatname{algorithm}{Algorithm}

 \usepackage[normalem]{ulem}
 \usepackage{epstopdf} 
\usepackage{amssymb}
\usepackage{xcolor} 

\usepackage{longtable,booktabs,array}
\usepackage{calc} 
\usepackage{etoolbox}
\makeatletter
\patchcmd\longtable{\par}{\if@noskipsec\mbox{}\fi\par}{}{}
\makeatother
\IfFileExists{footnotehyper.sty}{\usepackage{footnotehyper}}{\usepackage{footnote}}
\makesavenoteenv{longtable}
\usepackage{graphicx}
\makeatletter
\def\maxwidth{\ifdim\Gin@nat@width>\linewidth\linewidth\else\Gin@nat@width\fi}
\def\maxheight{\ifdim\Gin@nat@height>\textheight\textheight\else\Gin@nat@height\fi}
\makeatother

\newcommand{\cedp}{\newpage{\pagestyle{empty}\cleardoublepage}}

\makeindex

\usepackage[pageref]{backref}
\renewcommand*{\backrefalt}[4]{%
\ifcase #1 %
(Not cited)%
\or
(Cited on p.~#2)%
\else
(Cited on pp.~#2)%
\fi
}

\usepackage{times}
\usepackage[T1]{fontenc}

\providecommand{\tightlist}{\setlength{\itemsep}{0pt} \setlength{\parskip}{0pt}}
\newcommand \be {\begin{equation}}
\newcommand \ee {\end{equation}}
\usepackage{subcaption}

\newcommand \RR    {\mathbb{R}}

\newcommand \Lcal    {\mathcal{L}}

\newcommand \del   {\partial} 
\newcommand \RD {{\mathbb R}^D}
\newcommand \bei {\begin{itemize}}
\newcommand \eei {\end{itemize}}
\usepackage{xstring}
\usepackage{xparse}

\NewDocumentCommand{\ensureFinalDot}{m}{%
  \IfEndWith{#1}{.}{#1}{%
    \IfEndWith{#1}{?}{#1}{%
      \IfEndWith{#1}{!}{#1}{%
        \IfEndWith{#1}{:}{#1}{%
          \IfEndWith{#1}{;}{#1}{#1.}%
        }%
      }%
    }%
  }%
}

\let\ParagraphOriginal\paragraph

\RenewDocumentCommand{\paragraph}{s o m}{%
  \IfBooleanTF{#1}
    {%
      \ParagraphOriginal*{\ensureFinalDot{#3}}%
    }{%
      \IfNoValueTF{#2}
        {\ParagraphOriginal{\ensureFinalDot{#3}}}%
        {\ParagraphOriginal[\ensureFinalDot{#2}]{\ensureFinalDot{#3}}}%
    }%
}

\title{\bf \huge Reproducing kernel methods  
\\
for machine learning, PDEs, and statistics}

\author{\Large Philippe G. LeFloch\footnote{Laboratoire Jacques-Louis Lions, Sorbonne Universit\'e and Centre National de la Recherche Scientifique, 4 Place Jussieu, 75258 Paris, France. Email: \tt contact@philippelefloch.org}, Jean-Marc Mercier$^2$,  and Shohruh Miryusupov\footnote{MPG-Partners, 136 Boulevard Haussmann, 75008 Paris, France. 
\newline
Email: \tt jean-marc.mercier@mpg-partners.com, shohruh.miryusupov@mpg-partners.com.}
}

\date{September 2025}

\begin{document}

\maketitle
\cedp
\frontmatter
\tableofcontents
\cedp

%


\mainmatter
 
\chapter[Introduction]{Introduction} 
\index{Sample!file} 

\section{Main objective}
\label{main-objective}

This monograph offers a modern presentation of kernel-based algorithms with a focus on applications and use cases, through the prism of reproducible numerical tests. We primarily intend to tackle industrial use cases coming from computational physics and mathematical finance, and seek widespread applications across various areas, such as statistics, or artificial intelligence (physics-informed systems, reinforcement learning, machine learning, generative methods, etc.). Our algorithms are built upon a strategy based on the theory of reproducing kernel Hilbert spaces (RKHS) and the theory of optimal transport. On the other hand, the proposed algorithms have been implemented using an open-source library, CodPy\footnote{Our acronym stands for the curse of dimensionality in Python.}, and all figures of this book can be found on a companion website\footnote{\url{https://codpybook-read-the-docs.readthedocs.io/en/latest/}} with commented, reproducible Python code, in order to facilitate the diffusion of RKHS approaches to students, teachers, or practitioners.

To present the theoretical principles and the techniques employed in their applications, we have structured this monograph into two main parts. In Chapters 2 to 5, we focus on the fundamental principles of kernel-based representations, where the numerical part is supplemented with illustrative examples only. Next, in Chapters 6 to 9 we discuss the application of these principles to many classes of concrete problems. Chapter 6 deals with machine learning (supervised, unsupervised, generative methods). Chapter 7 covers some aspects of the numerical approximation of partial differential equations, which forms the foundation of physics-informed approaches. Estimating a model, to optimize a feed-back action, is the heart of reinforcement learning, which approaches and algorithms are adapted to an RKHS framework in chapter eight. Finally, we describe how to adapt RKHS to time-series predictions and market generators in Chapter 9, dedicated to mathematical finance applications.

We have aimed to make this monograph as self-contained as possible, primarily targeting engineers. We have intentionally omitted theoretical aspects of functional analysis and statistics which can be found elsewhere in the existing literature, and we chose to emphasize the operational applications of kernel-based methods. We solely assume that the reader has a basic knowledge of linear algebra, probability theory, and differential calculus. Our core objective is to provide a framework for applications, enabling the reader to apply the proposed techniques.

Obviously, this text cannot cover all possible directions on the vast subject that we touch upon here. Yet, we hope that this monograph can put in light the particularly robust strengths of kernel methods, and contribute to bridge, on the one hand, basic ideas of functional analysis and optimal transport theory and, on the other hand, a robust framework for machine learning and related topics. With this emphasis in mind, we have designed here novel numerical strategies, while demonstrating the versatility and competitiveness of these kernel methods for dealing with artificial intelligence problems, which is supported in this book by several benchmarks.


\section{Outline of this monograph}\label{outline-of-this-monograph}

We provide here a comprehensive study of kernel-based machine learning methods and applies them across a diverse range of topics in applied mathematics, finance, and engineering. It is organized as follows.
\begin{itemize}

\item Chapter~2 presents the core aspects of kernel techniques, starting from the basic concepts of reproducing kernels, moving on to kernel engineering, and then discussing interpolation/extrapolation operators, known as \textit{kernel ridge regression.} This chapter also introduces the notions of  discrepancy error and kernel-based norms, paving the way for designing effective performance indicators to assess the relevance of kernel ridge regression in applications.

\item In Chapter~3, we define and investigate the properties of kernel-based operators in greater depth. These operators play a key role in the discretization of partial differential equations, making them particularly useful in physics and engineering. Interestingly, they also find major applications in machine learning, especially for predicting deterministic, non-stochastic functions of the unknown variables.

\item Chapter~4 presents some clustering methods, adapted to an RKHS framework. Clustering is an important field for data exploration. This chapter focuses on clustering methods which lower the numerical burden of RKHS approaches. It allows us to design numerically efficient, large-scale dataset strategies. Clustering is presented here as a general combinatorial approach, which is intrinsically linked to an optimal transport viewpoint. It provides a rich combinatorial layout to tackle optimal transport-based algorithms next.

\item Chapter~5 is dedicated to statistical RKHS methods, and begins reviewing key concepts from optimal transport theory, from both continuous and discrete standpoints. We include a comparative analysis of two legacy kernel methods, namely the Nadaraya-Watson and kernel ridge regression, used for conditional density and conditional distribution modeling, respectively. This legacy analysis aims to introduce the need for optimal transport techniques to design RKHS generative methods, which is one of the central topics of this monograph.

\item  Chapter~6 focuses on supervised, unsupervised, and generative machine learning methods. We compare our framework against various machine learning methods, benchmarking across multiple scenarios and performance indicators, while analyzing their suitability for different types of learning problems.

\item Chapter~7 showcases the efficiency of the kernel techniques in solving partial differential equations on unstructured meshes. We consider a range of academic problems, starting from the Laplace equation to specific problems from fluid dynamics, together with the Lagrangian methods employed in particle, mesh-free methods. This chapter also highlights the power of the proposed framework in enhancing the convergence of Monte Carlo methods and briefly discusses automatic differentiation ---an essential yet intrusive tool.

\item Chapter~8 proposes a detailed description of a scalable standardized kernel approach to popular reinforcement learning algorithms, where agents interact with environments having continuous states and discrete action spaces, dealing with possibly unstructured data. These algorithms are namely Q-Learning, Actor Critic, Policy Gradient, Hamilton-Jacobi-Bellman (HJB) and Heuristic Controls. We analyzed them through the RKHS prism, and implemented them with default CodPy settings, to emphasize the sample efficiency, accuracy, and robustness of this approach, benchmarking our algorithms on simple games.

\item Chapter~9 is dedicated to mathematical finance, with a focus on \textit{econometrics}, whose purpose is to model time-series forecasts. The RKHS approach amounts to turning any grounded quantitative model into market generators, used for risk management or investment strategic purposes. We then present a use case, which reduces to studying the robustness of extrapolation methods, focusing on \textit{real-time} hedging or risk management methods. Market generators provide straightforward Monte Carlo pricing methods for simple financial products. For more complex derivatives or portfolios, recall that payoffs and risk factors in finance are called rewards and states for reinforcement learning. In particular, Hamilton-Jacobi-Bellman methods to approximate value functions with RKHS methods, described in the previous chapter, can be straightforwardly adapted to pricing, forming a natural bridge between mathematical finance and reinforcement learning.
\end{itemize}

By presenting a fresh perspective on kernel-based methods and offering a broad overview of their applications, we hope that this text will serve as a resource for researchers, students, and professionals in the fields of scientific computation, statistics, mathematical finance, and engineering sciences. In our endeavor to make the theoretical framework and algorithms proposed in this monograph accessible and user-friendly, our companion library integrates Python codes, documented in LaTeX as Jupyter notebooks, with a high-performance C++ core. This library provides a robust and versatile toolset for tackling a wide range of practical challenges. This open-source code aims to help readers learn and test the proposed algorithms, while also offering a foundation for the techniques that can be tailored to specific applications. Additionally, we expect this monograph to be updated later on, when new applications become available. For the convenience of the reader a selection of the acronyms used in this monograph is provided in Table~\ref{tab:acronyms}. 


\section{About this work}
\label{main-contribution}

This text relies on the joint research by the first authors~\cite{LeFloch-Mercier:2002,LeFlochMercier-2015, LeFloch-Mercier:2017,LeFloch-Mercier:2020b,LeFloch-Mercier:2020a,LeFloch-Mercier:2021,LeFlochMercier-technical1,LeFlochMercier-2024JCP} over the past fifteen years, published in several academic and industrial journals. This text was initially conceived as documentation for a non-regression tool for the CodPy Library, needed for its industrial use. In collaboration with the third author, it gradually evolved over the years toward a user manual for internal purposes and has shifted from a user manual to the present form, while striving to serve the first purposes. This evolution is motivated by the numerical efficiency of kernel approaches, supported in this monograph by numerous benchmarks in chapters devoted to applications. RKHS-based methods are particularly suited to design energy-efficient, fast-learning, and data-efficient algorithms. These properties are relevant for industrial competitiveness as well as climate concerns, but also for the development of artificial intelligence itself, which today reaches limits due to enormous energy and data consumption. Despite their many advantages, kernel methods remain underutilized due to computational complexity, scalability limitations, and a lack of dedicated software frameworks and analytical techniques. 

One of our purposes here is to popularize RKHS methods, filling some of these gaps: we provide a prototype of a dedicated RKHS, CPU-parallel, C++ library, as well as new approaches and algorithms built on top of it. Indeed, we hope, doing so, to pave the way toward more energy- and data-efficient algorithms for the artificial intelligence community.

There is a vast literature available on kernel methods and reproducing kernel Hilbert spaces which we do not attempt to review here. Our focus is on providing a practical framework for the application of such methods. For a comprehensive review of the theory we refer to textbooks by Berlinet and Thomas-Agnan \cite{BerlinetThomasAgnan:2004}, Fasshauer \cite{Fasshauer:2006,Fasshauer:2007,Fasshauer}, and Suzuki \cite{Suzuki}. 

The methodology in Chapters 1 to 5 and the kernel-based mesh-free algorithms presented in Chapter 7 are based on research papers by P.G. LeFloch and J.-M. Mercier. We refer the reader 
to~\cite{LeFloch-Mercier:2002} (a class of fully discrete, entropy-conservative schemes), 
\cite{LeFlochMercier-2015} (the convex hull algorithm), 
\cite{LeFloch-Mercier:2017} (a new method for solving the Kolmogorov equation), 
\cite{LeFloch-Mercier:2020b} (study the integration error via the kernel discrepancy),
\cite{LeFloch-Mercier:2020a} (a class of mesh-free algorithms),
\cite{LeFloch-Mercier:2021} (the transport-based mesh-free method),
\cite{LeFlochMercier-technical1} (predictive machines with uncertainty quantification), 
\cite{LeFlochMercier-2024JCP} (mesh-free algorithms for finance and machine learning). 
The other algorithms were developed in unpublished reports in collaboration with 
S. Miryusupov~\cite{LeFlochMercierMiryusupov:2021-t3,LeFlochMercierMiryusupov:2021a,LeFlochMercierMiryusupov:2021-t1,LMM-SSRN1,LeFlochMercierMiryusupov:2024sept}. 
For additional information on mesh-free methods in fluid dynamics and material science, the reader is referred to  \cite{BabuskaBanerjeeOsborn:2003,BessaFosterBelytschkoLiu:2014,GuntherLiu:1998,HaghighatRaissibMoureGomezJuanes:2021, KoesterChen:2019, KorzeniowskiWeinberg:2021,LiLiu:2004,Liu:2003,Liu:2016, Nakano:2017,NarcowichWardWendland:2005, OhDavisJeong:2012,SalehiDehghan:2013,SirignanoSpiliopoulos:2018,ZhouLi:2006}. For various developments in mathematical finance and applications, we refer  
to~\cite{AntonovKonikovSpector:2015,bollerslev1986,BraceGatarekMusiela:1997,Buehler:2010,cox1985,EFO:2021,engle1982,Harrison,heston1993,HugeSavine:2020,Mania2018,Mercier:2009,MercierMiryusupov:2019,Nadaraya:1964,TDGJ,vasicek1977}, and for machine learning to \cite{ChengHC2021,Chowdhury2017proceedings,Hastie,HSS:2008}. 

\begin{table}
\begin{tabular}{ll}
\textbf{Acronym} & \textbf{Meaning} \\
\hline
RKHS & Reproducing Kernel Hilbert Space \\
MMD & Maximum Mean Discrepancy \\
RL & Reinforcement Learning \\
OT & Optimal Transport \\
AD & Automatic Differentiation \\
PDE & Partial Differential Equation \\
PCA & Principal Component Analysis \\
RBF & Radial Basis Function \\  
RMSE & Root Mean Square Error \\
SVM & Support Vector Machine \\
GAN & Generative Adversarial Network \\ 
CDF & Cumulative Distribution Function \\ 
PDF & Probability Density Function \\ 
SDE & Stochastic Differential Equation \\
\end{tabular}
\caption{Selected acronyms cited in this monograph.}
\label{tab:acronyms}
\end{table}


\section{Acknowledgments}

\hskip.5cm P.G. LeFloch is a research professor (Directeur de recherche) at the Laboratoire Jacques-Louis Lions at Sorbonne University, supported by the Centre National de la Recherche Scientifique (CNRS), Agence Nationale de la Recherche (ANR), and MSCA Staff Exchange Grant Project 101131233 funded by the European Research Council (ERC). Part of this research was conducted while he was a visiting research fellow at the Courant Institute of Mathematical Sciences, New York University.  

J.-M. Mercier and S. Miryusupov are researchers at MPG-Partners. They are deeply indebted to their colleagues at MPG-Partners for their constant support, precious advices, and team spirit. We cannot mention them all, but we wish to express special thanks to our esteemed supervisors, G. Mathieu and M. Poirier.  

The authors also wish to thank all those who generously contributed to this project, especially Aymen Abidi, Max Aguirre, Ana\"{\i}s Barbiche, Barry Eich, and Amine El Marraki.  Last but not least, the authors would like to thank all of their colleagues for their scientific support and their close relatives for their love and patience! 


\clearpage 

\part{A framework based on reproducing kernels and optimal transport}

\chapter{Fundamental notions on reproducing kernels}
\label{basic-notions-about-reproducing-kernels}

\section{Discrete reproducing kernel Hilbert spaces}
\label{purpose-of-this-chapter}

\subsection{A definition of positive-definite kernels} 

We begin the presentation of our methods with the notion of reproducing kernels, which plays a pivotal role in building representations and approximations of both data and solutions, in combination with other features at the core of RKHS algorithms. Among these are other important concepts, such as transformation maps. These maps are mainly viewed in this chapter as simple tools to adapt kernels to data, but they are also introduced as a powerful technique for developing deep kernel architectures, which offer the flexibility to tailor basic kernels to address specific challenges, tackled in Chapter~\ref{optimal-transport-and-statistical-kernel-methods}, devoted to RKHS approaches to statistical methods. In the present chapter, we focus our attention on the notion of kernels from a discrete point of view. The well-established RKHS theory, available for the more general continuous case, is presented to support and clarify this viewpoint.

A \emph{kernel}, denoted by \( k: \mathbb{R}^D \times \mathbb{R}^D \to \mathbb{R} \) (respectively $\mathbb{C}$), is a symmetric real-valued function (respectively complex-valued), meaning that it satisfies
\be 
k(x, y) = k(y, x) \ (\text{respectively } \overline{k(y, x)}), \qquad x, y \in \mathbb{R}^D.
\ee
Consider a set of points (\textit{features}) in \( \mathbb{R}^D \), that is, a distribution
\(X \in \RR^{N_x,D}\), with $X=(x^1,\ldots,x^{N_x})$, and $x^n = (x_1^n,\ldots,x^{n}_D)$. Similarly, consider \(Y \in \RR^{N_y,D}\) and the following rectangular matrix
\be \label{KMATRIX}
    K(X, Y) = 
\begin{bmatrix}
        k(x^1,y^1) & k(x^1,y^2) & \ldots & k(x^1,y^{N_y}) \\
        k(x^2,y^1) & k(x^2,y^2) & \ldots & k(x^2,y^{N_y}) \\
        \vdots & \vdots & \ddots & \vdots \\
        k(x^{N_x},y^1) & k(x^{N_x},y^2) & \ldots & k(x^{N_x},y^{N_y})
\end{bmatrix}.
\ee
When \( X = Y \), the square matrix \( K(X, X) \) is often referred to as the \emph{Gram matrix} of the kernel function with respect to the chosen data points.

The function \( k \) is said to be a \emph{positive kernel} if, for any finite collection of \textit{distinct} points \( X = (x^1, \dots, x^{N_x}) \subset \mathbb{R}^D \) and for any set of coefficients \( c^1, \dots, c^{N_x} \in \mathbb{R} \), not all zero, the quadratic form satisfies
\be \label{PDK}
\sum_{i=1}^{N_x} \sum_{j=1}^{N_x} c^i c^j k(x^i, x^j) \ge 0.
\ee
A kernel that is strictly positive-definite (respectively positive) defines a strictly positive-definite (respectively positive) Gram matrix \( K(X, X) \).

More generally, a kernel \( k \) is said to be \emph{conditionally positive-definite} if it is positive-definite only on a certain sub-manifold of \( \mathbb{R}^D \times \mathbb{R}^D \). In other words, the positivity condition holds only when \( X \) and \( Y \) are restricted to belong to this sub-manifold, which may be referred to as the \emph{positivity domain}. By definition, this domain is a subset of \( \mathbb{R}^D \times \mathbb{R}^D \) on which \( k \) is ensured to be positive-definite. Finally, positive-definite kernels are not necessarily positive-valued. For example, the sinus-cardinal $k(x,y)=\frac{\sin(x-y)}{x-y}$ defines a positive-definite kernel.

An important property of positive-definite kernels is thus to define symmetric, positive-definite Gram matrices \(K(X,X)\), which are invertible provided \(X\) consists of \textit{distinct} points in \(\RR^{N_x,D}\).

The \textit{features} set $X$ can contain discrete values, for instance $x_d^i \in \{0,1\}$, or continuous ones $x_d^i \in \RR$, or a combination of both. It can represent any data: of course, point coordinates, but also labels, images, videos, sounds, text, DNA bases, etc.


\subsection{Kernel-based approximations}

Given a set \(Y\), a kernel defines a discrete space of vector-valued functions (regressors), in which functions are parametrized by matrices \(\theta \in \mathbb{R}^{N_y,D_f}\), as follows: 
\be
\label{Hk}
\mathcal{H}_{k,Y} = \Bigl\{
f_{k,\theta}(\cdot) = \sum_{n=1}^{N_y} \theta^n k(\cdot,y^n) = K(\cdot,Y) \theta, \quad \theta \in \mathbb{R}^{N_y,D_f} \Bigr\}.
\ee
The notation \(\cdot\) is here a placeholder, and we will also use the notation \(f_{k, \theta} (Z) \in \mathbb{R}^{N_z, D_f}\), to describe the values of a given function on a set \(Z\). The parameters \(\theta\) are computed, or \textit{fitted}, to a given continuous function \(f(\cdot) \in \mathbb{R}^{D_f}\), known from its discrete values \(X,f(X)\), with possibly \(Y \neq X\), according to the relation
\be\label{FIT}
 f_{k, \theta} (\cdot) = K(\cdot,Y) \theta, \quad \theta = \left( K(X,Y) + \epsilon R(X,Y) \right)^{-1} f(X),
\ee
in which \(\epsilon \ge 0\), and \(R(Y,X)\) is a matrix defining an
optional regularizing term. Throughout, matrix inversions such as in \eqref{FIT} are performed using a least-square method, corresponding to \textit{kernel ridge regression}, although other methods could also be considered.

Consider now \eqref{FIT} as a matrix product between \(\mathcal{P}_k(\cdot,X)\) and \(f(X)\):
\be\label{OP}
 f_{k,X,Y}(\cdot) = \mathcal{P}_{k,Y}(\cdot,X) f(X), \qquad \mathcal{P}_{k,Y}(\cdot,X) =K(\cdot,Y)\left( K(X,Y) + \epsilon R(X,Y) \right)^{-1}.
\ee 
If evaluated pointwise, \(\mathcal{P}_{k,Y}(x,X)\) is a vector
of size \(N_x\). If the evaluation is performed on a set \(Z\), \(\mathcal{P}_{k,Y}(Z,X)\) is a matrix of size \(N_z,N_x\). The operator \(\mathcal{P}_{k,Y}(\cdot,X)\) is called \textit{projection operator}, as it computes the projection of any function onto the space \(\mathcal{H}_{k,Y}\), known through its values \(X,f(X)\).

In the same way, we can derive from \eqref{FIT} the following gradient formula
\be\label{GRAD}
 \nabla f_{k, \theta} (\cdot) = (\nabla K)(\cdot,X) \theta, \quad \theta = \left( K(Y,X) + \epsilon R(Y,X) \right)^{-1} f(Y),
\ee
where the gradient operator is defined by partial derivatives \(\nabla= \left[\frac{\partial}{\partial x_1}, \ldots, \frac{\partial}{\partial x_D}\right]\), and \((\nabla k)(x,y)\) denotes the gradient of the kernel with respect to \(x\) or \(y\), since kernel functions are symmetric. We also use the short-hand notation \(f_{k}\), or \(\mathcal{P}_{k}\), but depending on context, we may disambiguate using \(f_{k, \theta}\), \(f_{k, \theta, X, Y}(Z)\), \(\mathcal{P}_{k,Y}(\cdot,X)\), etc. 

The discrete function space \(\mathcal{H}_{k,X}\) is equipped with the scalar product 
\be
\label{kprod}
 \langle f_{k,\theta_f},g_{k,\theta_g} \rangle_{\mathcal{H}_{k,X}} = \sum_{i}\sum_{j} \theta_f^i \theta_g^j k(x^i,x^j) = \theta_f^T K(X,X) \theta_g.
\ee
Consider a strictly positive-definite kernel and the \textit{extrapolation mode} \(X=Y, \epsilon = 0\) in \eqref{FIT}. It satisfies \(\langle k(x^i,\cdot),k(x^j,\cdot) \rangle_{\mathcal{H}_{k,X}} = k(x^i,x^j)\), or more generally the reproducing property
\be
\label{REPRO}
 \langle f_{k}(\cdot),k(x^i,\cdot) \rangle_{\mathcal{H}_{k,X}} = f(x^i), \quad f_k \in  \mathcal{H}_{k,X}.
\ee
For this mode, we can further reduce the scalar product expression to \(\langle f_{k}(\cdot),g_{k}(\cdot) \rangle_{\mathcal{H}_{k,X}} = \langle \theta_f(X), g(X) \rangle = \langle f(X), \theta_g(X) \rangle\). Unless specified, the default scalar product is the Euclidean one, \(\langle X, Y \rangle = \sum_n x^n y^n\), or the Frobenius one \(\langle X, Y \rangle = \sum_{n,d} x_d^n y_d^n\).

The couple \(X,f(X)\) represents observations of inputs / outputs, called features/labels, usually referred to as the \textit{training set}. The set \(f(X)\) can be discrete, or continuous, or a random distribution, and there is not a clear separation between the role of \(X\) and \(f(X)\).


\subsection{Error estimates based on kernel discrepancy}
\label{error-estimates-based-on-the-kernel-based-discrepancy}

Consider the reproducing mode \(\epsilon=0\) in the fitting formula \eqref{FIT}. On any \textit{test set} \(Z \in \mathbb{R}^{N_z,D}\), the following error estimate holds:
\be\label{estim}
   \Bigg| \frac{1}{N_z} \sum_{n=1}^{N_z}  f(z^n) - f_{k,\theta}(z^n) \Bigg| \le \Big( d_k\big(X,Y\big) + d_k\big(Y,Z\big) \Big) \, \|f\|_{\mathcal{H}_{k}}
\ee
for any vector-valued function \(f:\mathbb{R}^D \to \mathbb{R}^{D_f}\), where \(d_k(X,Y)^2\) is a distance between sets, defined as 
\be\label{dk}
\frac{1}{N_x^2}\sum_{n=1,m=1}^{N_x,N_x} k\big(x^n,x^m\big) + \frac{1}{N_y^2}\sum_{n=1,m=1}^{N_y,N_y} k\big(y^n,y^m\big) - \frac{2}{N_x N_y}\sum_{n=1,m=1}^{N_x,N_y} k\big(x^n,y^m\big).
\ee
In our presentation, we use the terms \emph{maximum mean discrepancy} (MMD)\footnote{first introduced in \cite{GR:2006}} and \emph{kernel discrepancy} interchangeably.

The key term \(d_k\big(X,Y\big)\) is a kernel-related distance between a set of points that we call the \emph{discrepancy functional}. It is a rather natural quantity, and we expect that the accuracy of an extrapolation diminishes when the extrapolation set \(Z\) moves away from the sampling set \(X\). Indeed, considering \(X=Y\) in \eqref{estim} leads to a pointwise, \(L^\infty\) estimate at any point \(z\):
\be
   | f(z) - f_{k,X}(z) | \le d_k\big(X,z\big) \, \|f\|_{\mathcal{H}_{k}}
\ee
We emphasize here that the two terms on the right-hand side of \eqref{estim} are computationally realistic and can be systematically applied to assess the validity of an extrapolation. To be specific, in order to estimate the norm on the right-hand side, we must approximate \(\|f\|^2_{\mathcal{H}_{k}} \ge \|f\|^2_{\mathcal{H}_{k,X}} = \langle \theta, f(X) \rangle\). This estimation is reasonable, as it uses all the available information on \(f\). 

Finally, note that this error measurement can be used to build other estimates, for instance, considering \(g(\cdot) = | f(\cdot) - f_k(\cdot) |^2\) in \eqref{estim}, we get an \(\ell^2\) estimate:
\be\label{err}
  \big\| f(Z) - f_k(Z) \big\|_{\ell^2} \le \Big( d_k\big(X,Y\big) + d_k\big(Y,Z\big) \Big) \|f\|_{\mathcal{H}_{k}}. 
\ee
Observe that the important reproducing property of RKHS methods \eqref{REPRO}, corresponding to the \textit{extrapolation mode} \(X=Y, \epsilon = 0\) in \eqref{FIT}, can be formulated in several equivalent ways:
\be\label{REPRO2}
  f_{k,X,X}(X) = f(X), \quad \text{or } d_k(X,X) = 0, \quad \text{or } \mathcal{P}_{k,X}(X,X) = I_d.
\ee


\section{Continuous reproducing kernel Hilbert spaces}
\label{sec-conti}

\subsection{Discrete versus continuous}

We anchor the previous numerical, discrete section to the well-established RKHS theory, developed since the beginning of the twentieth century, which sheds light on the discrete standpoint. We discuss here the discrete framework, which may be skipped by readers interested in the algorithmic aspects of RKHS methods.

The links between the discrete and the continuous standpoints are covered by well-known theoretical results.
\begin{itemize}
\item The existence of positive-definite kernels is ensured by the Bochner theorem. Universal kernels can describe any continuous function. The Mercer theorem provides a general separable form for positive-definite kernels.

\item Starting from a positive-definite kernel \(k\), the Moore-Aronszajn result gives a characterization of the Hilbert function space \(\mathcal{H}_{k,X}\), as the size of the sequence \(X \in \mathbb{R}^{N_x,D}\) goes to infinity (\(N_x \to + \infty\)) in \eqref{Hk}.

\item The representer theorem allows us to characterize the projection operator \eqref{FIT} as the minimum of a functional but also shows that a broader class of functionals attains its minimum in \(\mathcal{H}_{k,X}\).

\item We introduced the kernel discrepancy (or MMD) \(d_k\big(X,Y\big)\). This functional provides a distance metric among measures for most of our kernels, and we can use it as an alternative statistical test for rejection in the two-sample problem, for instance.
\end{itemize}

Let us conclude here by discussing the limiting case in which the number of samples or features \(X\), i.e., \(N_x\), tends to \(+\infty\). The set \(X\) usually consists of samples from a distribution \(\mu(\cdot)\), and, for our applications, we assume that the distribution of \(X\) converges in the sense of measures. 
To define this convergence, we consider the equi-weighted measure \(\delta_X=\frac{1}{N} \sum_{n=1}^N \delta_{x^n}\), where \(\delta_{x}\) is the Dirac measure at a point \(x\). Then we assume the convergence \(\delta_X \rightharpoonup \mu(\cdot)\) in the sense that 
\be\label{mu}
 \langle \delta_X,\varphi \rangle_{\mathcal{D}',\mathcal{D}} = \frac{1}{N} \sum_{n=0}^N \varphi(x^n) \mapsto 
  \int_{\mathbb{R}^D} \varphi(x) \, d \mu(x), \quad N \to + \infty,
\ee
for any continuous function \(\varphi\). This definition coincides with the standard convergence in the Schwartz space of distributions \(\mathcal{D}'\).


\subsection{Bochner theorem and universal kernels}
\label{bochner-theorem}

One of the most used classes of positive-definite kernels is \textit{translation invariant}, that is, of the form \(k(x,y)=g(x-y)\), for which an application of Bochner theorem states that this class of kernels is positive-definite (respectively non-negative) if and only if their Fourier transform is a positive (respectively non-negative) measure: 
\be\label{Bochner}
k(x,y)= \int_{\mathbb{R}^D} e^{2i\pi \langle x-y,\omega\rangle} \, d\mu(\omega) 
= \widehat{\mu}(x-y), \quad \mu > 0 \ (\text{respectively } \mu \ge 0),
\ee
where \(\widehat{\mu}\) is the (inverse) Fourier transform. Translation-invariant kernels with strictly positive Fourier transform, i.e., \(\mu > 0\) in \eqref{Bochner}, are \textit{universal}, meaning that the generated space \(\mathcal{H}_k\) induced by the kernel is dense in \(\mathcal{C}_0(\mathbb{R}^D)\), the space of continuous functions vanishing at infinity.

Formally, the formula \eqref{Bochner} can be interpreted as a scalar product, that is,
\be
k(x,y) = \langle \widehat{\mu}^{1/2}(x),\widehat{\mu}^{1/2}(-y) \rangle,
\ee
with 
\be
\widehat{\mu}^{1/2}(\cdot) = \int_{\mathbb{R}^D} e^{\pi \langle \cdot,\omega\rangle} \, d\mu^{1/2}(\omega).
\ee
Mercer theorem (discussed next) states that this inner scalar product structure is universal for positive-definite kernels.


 \subsection{Mercer theorem}

Indeed, the construction of the space \(\mathcal{H}_{k,X}\) in \eqref{Hk}, consisting of all linear combinations of the \textit{basis functions} \(k(\cdot, x^n)\), starts from a kernel and defines a Hilbert space of functions equipped with the scalar product \eqref{kprod}. Let us now consider the converse, namely, starting with any Hilbert space of functions, denoted \(f \in \mathcal{H}_k\), for which we assume an inner scalar product of the kind \eqref{kprod}, that is, there exists a symmetric, square positive-definite matrix \(Q\) such that \(\langle f, g \rangle_{\mathcal{H}_k} = f^T Q g\).

Let us assume that the point evaluation \(x \mapsto f(x)\) is continuous and linear for any \(x \in X\). Hence, applying the Riesz representation theorem, there exists a unique vector \(k(\cdot,x^n) \in \mathcal{H}_k\) such that
\be
f(\cdot) = \langle k(\cdot,x^n), f(\cdot) \rangle_{\mathcal{H}_k}.
\ee
In the discrete case, we can compute explicitly these vectors, defined as the solution of the equation
\be
\langle k(\cdot,x^n), k(\cdot,x^m) \rangle_{\mathcal{H}_k} = k(x^n,x^m), \quad n,m = 1,\dots,N_x.
\ee
This last equation induces a relation between the Gram matrix \(K = \left(k(x^n,x^m)\right)_{n,m=1}^{N_x}\) and the inner scalar product matrix given by \(K Q K = K\), with a formal solution \(K = Q^{-1}\).

This construction amounts to stating that any symmetric positive-definite matrix \(K\) can be written as a Gram matrix representing inner products between a set of vectors: using the eigenvector decomposition \(Q = V D V^T\), where \(D\) is diagonal, and \(V\) is an orthogonal matrix of eigenvectors, then the Gram matrix is defined as the matrix product
\be
K = [V D^{-1/2}][D^{-1/2} V^T],
\ee
each element of \(K\) being computed by an inner scalar product.

This elementary linear algebra result is the Mercer theorem in the finite-dimensional discrete setting, which generalizes in the continuous case as follows.

\begin{theorem}
Let \(\mathcal{H}\) be a Hilbert space of functions. If \(k(\cdot,\cdot)\) is a positive-definite kernel on \(\mathcal{H}\), then
\be
k(x,y) = \sum_{j=1}^{+\infty} \lambda_j \, e_j(x) \, e_j(y),
\ee
where \(e_j(x)\) is an orthonormal basis of \(\mathcal{H}\) defined through the relation
\be
 e_n(y) = \int_{\mathbb{R}^D} k(x,y) \, e_n(x) \, d\mu(x). 
\ee
\end{theorem}


 \subsection{Moore-Aronszajn theorem}

The discrete space \(\mathcal{H}_{k,X}\) in \eqref{Hk}, consisting of all linear combinations of the \textit{basis functions} \(k(\cdot, x^n)\), is built by starting from a kernel and then defining a Hilbert space of functions, equipped with the scalar product \eqref{kprod}. The Moore-Aronszajn theorem states that this construction still holds in the limit \(N_x \to + \infty\).

\begin{theorem}
Suppose \(k\) is a symmetric, positive-definite kernel on a set \(X\). Then there is a unique Hilbert space of functions on \(X\) for which \(k\) is a reproducing kernel.
\end{theorem}

Let us make some comments on this theorem: it states that for every positive-definite kernel \(k : X \times X \to \mathbb{R}\), there is a unique associated function space \(\mathcal{H}_{k,X}\) for which \(k\) has the reproducing property
\be
f(x) = \langle k(\cdot,x), f(\cdot) \rangle_{\mathcal{H}_{k,X}}.
\ee
The map
\be
x \mapsto k(x,\cdot) \in \mathcal{H}_{k,X},
\ee
is called the feature map. Due to the properties of scalar products, the kernel function satisfies the following properties:
\be
k(x,y) = k(y,x),\quad k(x,x) = \|k(\cdot,x)\|^2_{\mathcal{H}_{k,X}} \ge 0, \quad k(x,y)^2 \le k(x,x) k(y,y).
\ee

In particular, most of our kernels are defined on the whole space \(k:\mathbb{R}^D \times \mathbb{R}^D \to \mathbb{R}\). The Moore-Aronszajn theorem states that the following function space
\be\label{HK}
\mathcal{H}_k = \text{Span} \left\{k(\cdot, x) \, : \, x \in \mathbb{R}^D \right\},
\ee
equipped with the scalar product and norm \eqref{kprod}, is complete, and hence defines a well-defined Hilbert space.

In applications, we usually consider a smaller set, given by
\be
\mathcal{H}_{k,X} = \text{Span} \left\{k(\cdot, x) \, : \, x \in X \subset \mathbb{R}^D \right\},
\ee
clearly satisfying \(\mathcal{H}_{k,X} \subset \mathcal{H}_{k}\). Assuming the weak convergence \eqref{mu} as \(N_x \to + \infty\), the space generated \(\mathcal{H}_{k,X} \to \mathcal{H}_{k,\mu}\), as \(N_x \to + \infty\), where \(\mathcal{H}_{k,\mu}\) is a localized version of the space \(\mathcal{H}_k\). For any \(f, g \in \mathcal{H}_{k,X}\), denote \(\theta_f(\cdot), \theta_g(\cdot)\) their components, then \(\mathcal{H}_{k,X}\) is endowed with the scalar product:
\be
\label{kprodc}
\langle f_{k,\theta_f}, g_{k,\theta_g} \rangle_{\mathcal{H}_{k,X}} = \int_{\mathbb{R}^D}\int_{\mathbb{R}^D} \theta_f(x) \theta_g(y) k(x,y) \, d\mu(x) \, d\mu(y).
\ee


\subsection{The representer theorem}
\label{representer_theorem}

The representer theorem provides the theoretical basis {for approximating general regularized minimization problems on reproducing kernel Hilbert spaces. Let us first motivate this theorem before stating it. Consider a positive-definite kernel \(k\), and define the associated space \(\mathcal{H}_k\) as in \eqref{HK}, and let \(f \in \mathcal{H}_k\) and \(X\) a finite set of distinct sample points. Consider the functional:
\be
J(\varphi) = \sum_n \|\varphi(x^n) - f(x^n)\|^2_{\ell^2} + \epsilon \varphi^T(X) R \varphi(X).
\ee
Consider the minimization problem \(f_k = \arg \inf_{\varphi \in \mathcal{H}_k} J(\varphi)\), whose solution coincides with the fitting formula \eqref{FIT}. The representer theorem ensures that this minimum, attained in the finite-dimensional subspace \(\mathcal{H}_{k,X} \subset \mathcal{H}_k\), is indeed a global minimum in \(\mathcal{H}_k\), and extends to a broader class of functionals, as now stated. 

\begin{theorem}\label{repth}
Let \(\mathcal{X}\) be a nonempty set, \(k\) a positive-definite real-valued kernel on \(\mathcal{X} \times \mathcal{X}\) with associated reproducing kernel Hilbert space \(H_k\), and let \(R : H_k \to \mathbb{R}\) be a differentiable regularization function.  Then, given a training sample \((X, f(X)) \in \mathcal{X} \times \mathbb{R}\) and an arbitrary error function
\be
E : (\mathcal{X} \times \mathbb{R}^2)^m \to \mathbb{R} \cup \lbrace +\infty \rbrace,
\ee
a minimizer
\be
f^{*} = \underset{f\in H_k}{\operatorname{arg\,min}} \left\lbrace E\left( (x^1, f(x^1), f_k(x^1)), \ldots,  (x^n, y^n, f(x^n)) \right) + R(f) \right\rbrace
\ee
of the regularized empirical risk admits a representation of the form
\be
f^{*}(\cdot) = \sum_{i = 1}^n \alpha_i k(\cdot, x_i),
\ee
where \(\alpha_i \in \mathbb{R}\) for all \(1 \le i \le n\), if and only if there exists a nondecreasing function \(h : [0, +\infty) \to \mathbb{R}\) such that \(R(f) = h(\|f\|)\).
\end{theorem}

\subsection{Two-sample problem and characteristic kernels}
 
Consider two measures \(d\mu, d\nu\), and a kernel \(k\). The discrete formula \eqref{dk} corresponds to the following, more general, kernel discrepancy:
\be \label{eq:TSP}
 d_k(\mu,\nu) = \int \int d_k(x,y) \, d\mu(x) \, d\nu(y), \quad d_k(x,y) = k(x,x) + k(y,y) - 2k(x,y).
\ee
The two-sample problem consists in identifying those kernels for which this formula provides a distance between measures, in which case the kernel is called \emph{characteristic}. Specifically, a sufficient condition for a kernel to be characteristic is to be \textit{universal}; see Section~\ref{bochner-theorem}.


\section{Examples and properties of reproducing kernels}

\subsection{List of positive-definite kernels}

Throughout, we work with both positive-definite and conditionally positive-definite kernels. A list of kernels\footnote{available in the CodPy library} is provided in Table~\ref{tab:302} and visualized in Figure~\ref{fig:303}. We emphasize that in practical applications, a \textit{scaling of these basic kernels} is often required to suitably handle input data. Such a transformation is made to guarantee that the kernel captures the relevant features at the right scale. The details of these transformations are discussed later on, in the section on transformation maps.

\begin{longtable}[]{@{}
  >{\raggedright\arraybackslash}p{0.22\linewidth}
  >{\raggedright\arraybackslash}p{0.78\linewidth}@{}}
\caption{\label{tab:302}List of available kernels}\tabularnewline
\toprule
\textbf{Kernel} & \textbf{Expression} \( k(x, y) \) \\
\midrule
\endfirsthead
\toprule
\textbf{Kernel} & \textbf{Expression} \( k(x, y) \) \\
\midrule
\endhead
\midrule
\multicolumn{2}{r}{\textit{Continued on next page}} \\
\midrule
\endfoot
\bottomrule
\endlastfoot
Dot product & \( k(x, y) = x^T y \) \\
Periodic Gaussian & \( k(x, y) = \prod_{d} \theta_3(x_d - y_d) \) \\
Matérn & \( k(x, y) = \exp(-|x - y|) \) \\
Matérn tensorial & \( k(x, y) = \exp(-\prod_{d} |x_d - y_d|) \) \\
Matérn periodic & \( k(x, y) = \prod_{d} \frac{\exp(|x_d - y_d|) + \exp(1 - |x_d - y_d|)}{1 + \exp(1)} \) \\
Multiquadric & \( k(x, y) = \sqrt{1 + \frac{|x - y|^2}{c^2}} \) \\
Multiquadric tensorial & \( k(x, y) = \prod_{d} \sqrt{1 + \frac{(x_d - y_d)^2}{c^2}} \) \\
Sinc square tensorial & \( k(x, y) = \prod_{d} \left( \frac{\sin(\pi (x_d - y_d))}{\pi (x_d - y_d)} \right)^2 \) \\
Sinc tensorial & \( k(x, y) = \prod_{d} \frac{\sin(\pi (x_d - y_d))}{\pi (x_d - y_d)} \) \\
Tensor product kernel & \( k(x, y) = \prod_{d} \max(1 - |x_d - y_d|, 0) \) \\
Truncated kernel & \( k(x, y) = \max(1 - |x - y|, 0) \) \\
Polynomial kernel & $k(x,y) = \big(1+\frac{<x,y>}{D} \big)^p$ \\
Polynomial convolutional kernel & $k(x,y) = \| 1+\frac{x \ast y}{D} \|_{\ell^p}^p$ \\
\end{longtable}


Certain kernels are particularly useful in specific applications due to their mathematical properties and the structures they capture in data. We outline some key applications. 
\begin{itemize}
    \item \textit{ReLU kernel}. The rectified linear unit (ReLU) is widely used in neural networks, particularly as an activation function in deep learning. Its applications include image recognition, natural language processing, and various tasks in artificial intelligence. The positive-definite version of the ReLU-activation function is the tensor-product in Table~\ref{tab:302}, given by
\be
k(x, y) = \max(1 - |x-y|, 0).
\ee
It provides a useful and widely adopted choice in machine learning applications. This kernel is also quite interesting for its link to finite-difference approaches.

    \item \textit{Gaussian kernel}. The Gaussian kernel assigns higher weights to points that are closer to one another, making it ideal for tasks where local similarity matters, such as image recognition. It is also a key component in clustering algorithms, kernel density estimation, and dimensionality reduction techniques such as kernel PCA (Principal Component Analysis).

    \item \textit{Multiquadric kernel}. The multiquadric kernel, along with its tensorized variants, is based on radial basis functions (RBF). It is particularly useful for smoothing and interpolating scattered data, making it a valuable tool in applications such as weather forecasting, seismic data analysis, and computer graphics.

    \item \textit{Sinc and sinc square kernels}. The Sinc kernel and its squared tensorial form play a crucial role in signal processing and image analysis. These kernels accurately model periodicity in signals and images, making them well-suited for applications such as speech recognition, image denoising, and pattern recognition.
    
   \item \textit{Linear and polynomial regression kernels.} Given a mapping \( S:\mathbb{R}^D \to \mathbb{R}^P \) and a function \( g: \mathbb{R} \to \mathbb{R} \), the following construction provides a positive-definite kernel for any scalar-valued function $g(\cdot):$
\be
k(x,y) = g(\langle S(x), S(y) \rangle_{\mathbb{R}^P}), \qquad x,y \in \mathbb{R}^D,
\ee
where \( g \) is called the \emph{activation function}, and \( \langle \cdot, \cdot \rangle \) denotes the standard inner product in \( \mathbb{R}^P \). A useful example considers \( S(x) \) as the successive powers of the coordinate functions \( x_d \), yielding
\be
k(x,y) = \langle (1, x, x^T x, \ldots), (1, y, y^T y, \ldots) \rangle.
\ee
This corresponds to the classical kernel associated with \emph{linear regression} on a polynomial basis. This kernel is positive-definite, but the null space of the associated kernel matrix might be non-trivial, typically if the dimension $D$ exceeds the number of monomials used for the regression.

   \item \textit{Convolutional kernels.} In some applications we need kernels presenting some invariance properties. For instance, image detection often requires kernels that are invariant by translation of data, as the same object can be in different locations in a picture. A simple way to build invariant kernels is to consider any kernel that is based on the scalar product $k(x,y)=\varphi(<x,y>)$ and instead define $\varphi( x \ast y )$, where $x \ast y$ is a convolution. We illustrate this construction with the polynomial kernel (see Table~\ref{tab:302}), built upon the following definition for a discrete convolution
\be \label{conv_kernel}
  x \ast y = \Big[ \sum_{d=1,\ldots,D}x_d y_{(m+d)\%D} \Big]_{m=1}^D,
\ee
that is, a one-dimensional periodic convolution. 
\end{itemize}

\begin{figure}
    \centering
    \includegraphics[width=0.7\textwidth]{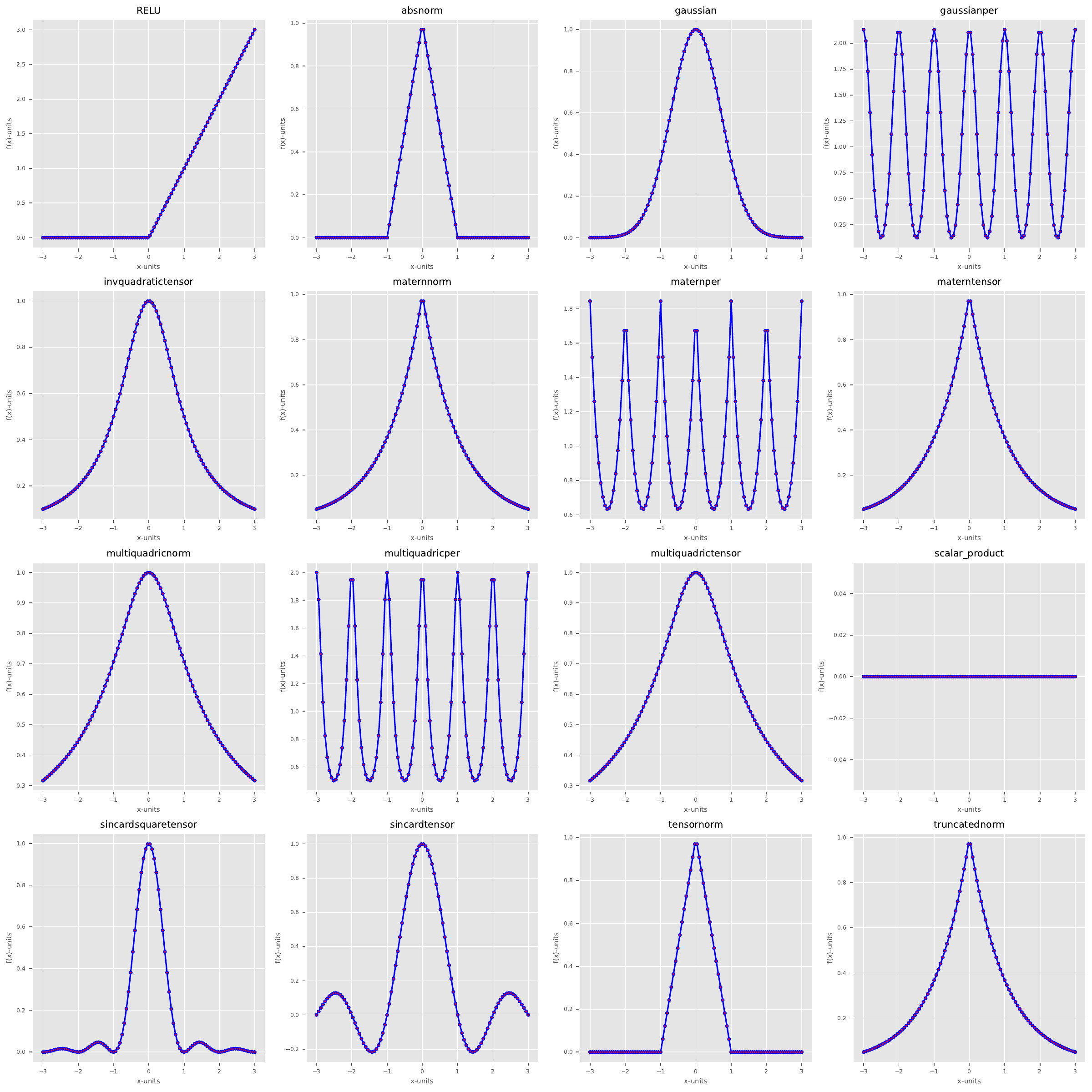}
    \caption{\label{fig:303}A list of kernels}
\end{figure}

We conclude by providing an example of a kernel matrix, considering the so-called \emph{tensornorm kernel}, described below. Typical values for this matrix are presented in Table~\ref{tab:304}, which displays the first four rows and columns.
  
\begin{table}[htbp]
\caption{\label{tab:304}First four rows and columns of a kernel matrix $K(X,Y)$}
\centering
\begin{tabular}[t]{r|r|r|r}
\hline
4.000000 & 3.873043 & 3.746087 & 3.619130\\
\hline
3.873043 & 3.833648 & 3.714253 & 3.594858\\
\hline
3.746087 & 3.714253 & 3.682420 & 3.570586\\
\hline
3.619130 & 3.594858 & 3.570586 & 3.546314\\
\hline
\end{tabular}
\end{table}

 
\subsection{Maps and kernels}  
\label{maps-and-kernels}

In view of the definition of positive-definite kernels in~\eqref{PDK}, the composition of a kernel $k$ and a bijective map $S$, $k\circ S$, remains a positive-definite kernel. Maps are paramount for kernels, and we can roughly categorize them in two different categories. 
\begin{itemize}

\item One category, called \textit{rescaling maps} $S : \RR^D \to \RR ^{D}$ adapts data to kernels, as most require specific and cautious scaling to escape numerical instabilities.

\item The second, \textit{dimensional embedding maps}, embedding the feature space $S : \RR^D \to \RR ^{D_S}$. The kernel trick usually refers to the dimension segmenting case $D_S > D$, since when $D_S < D$, this construction is usually referred to as a mapping to a \textit{latent space}.
\end{itemize}

The list in Table~\ref{tab:307} consists of rescaling maps, available in our framework. Chapter~\ref{optimal-transport-and-statistical-kernel-methods} deals with the construction of dimensional embedding maps.

\begin{longtable}[]{@{}
  >{\raggedright\arraybackslash}p{(\columnwidth - 4\tabcolsep) * \real{0.0142}}
  >{\raggedright\arraybackslash}p{(\columnwidth - 4\tabcolsep) * \real{0.1321}}
  >{\raggedright\arraybackslash}p{(\columnwidth - 4\tabcolsep) * \real{0.8538}}@{}}
\caption{\label{tab:307}List of available maps}\tabularnewline
\toprule\noalign{}
\begin{minipage}[b]{\linewidth}\raggedright
\end{minipage} & \begin{minipage}[b]{\linewidth}\raggedright
Maps
\end{minipage} & \begin{minipage}[b]{\linewidth}\raggedright
Formulas
\end{minipage} \\
\midrule\noalign{}
\endfirsthead
\toprule\noalign{}
\begin{minipage}[b]{\linewidth}\raggedright
\end{minipage} & \begin{minipage}[b]{\linewidth}\raggedright
Maps
\end{minipage} & \begin{minipage}[b]{\linewidth}\raggedright
Formulas
\end{minipage} \\
\midrule\noalign{}
\endhead
\bottomrule\noalign{}
\endlastfoot
1 & Scale to standard deviation & \(S(X) = \frac{x}{\sigma}\), \(\sigma =  \sqrt{\frac{1}{N_x}\sum_{n<N_x} (x^n - \mu)}\), \(\mu =  \frac{1}{N_x}\sum_{n<N_x} x^n.\) \\
2 & Scale to erf & \(S(X) = \text{erf}(x)\), \(\text{erf}\) is the standard error function. \\
3 & Scale to erfinv & \(S(X)=\text{erf}^{-1}(x)\), \(\text{erf}^{-1}\) is the inverse of \(\text{erf}\). \\
4 & Scale to mean distance & \(S(X) = \frac{x}{\sqrt{\alpha}}\), \(\alpha = \sum_{i,k \le N_x}  \frac{\mid x^i-x^k \mid^2}{N_x^2}.\) \\
5 & Scale to min distance & \(S(X) = \frac{x}{\sqrt{\alpha}}\), \(\alpha = \frac{1}{N_x}  \sum_{i\le N_x} \min_{k \neq i} \mid x^i-x^k \mid^2.\) \\
6 & Scale to unit cube & \(S(X) = \frac{x - \min_n{x^n} + \frac{0.5}{N_x}}{\alpha}\),\(\alpha  = \max_n{x^n} - \min_n{x^n}.\) \\
7 & Bandwidth & \(S(X) = hX \), where $h$ is user defined. \\
\label{tab:maplist}
\end{longtable}

Applying a map \(S\) is equivalent to replacing a kernel \(k(x,y)\) by the kernel \(k(S(x),S(y))\). For instance, the use of the ``scale-to-min distance map'' is usually a good choice for Gaussian kernels, as it scales all points to the average minimum distance. As an example, we can transform the given Gaussian kernel using such a map. 

Finally, in Figure~\ref{fig:308} we illustrate the action of maps on our kernels. Here, we should compare the two-dimensional results generated with maps to the one-dimensional results generated without maps and given earlier in Figure~\ref{fig:303}.

\newpage 

\textit{Composition of maps}. Within our framework, we frequently employ maps to preprocess input data prior to the computation based on kernel functions or using model fitting. Each map, with its unique features, can be combined with other maps in order to craft more robust transformations. As an illustrative example, the map used for the defaults kernel is called \textit{standard mean map}, and corresponds to the map
\be
    \label{standardmap}
    S(x) = S_4\circ S_3 \circ S_6(x),
\ee
where $S_i$ denotes the corresponding map in Table~\ref{tab:maplist}. This composite map starts by rescaling all data points to fit within a unit hypercube, followed by the erf inv map, and finally uses a bandwidth type map, with a scaling given by a mean distance adapted to the kernel. This particular transformation has been found to be a good balance for a number of our machine learning algorithms\footnote{For instance, the default kernel in CodPy considers the map \eqref{standardmap} with the Mat\'ern kernel $k(x,y)=\exp(-|x-y|_1)$.} 

\begin{figure}
\centering
\includegraphics[width=0.7\textwidth, keepaspectratio]{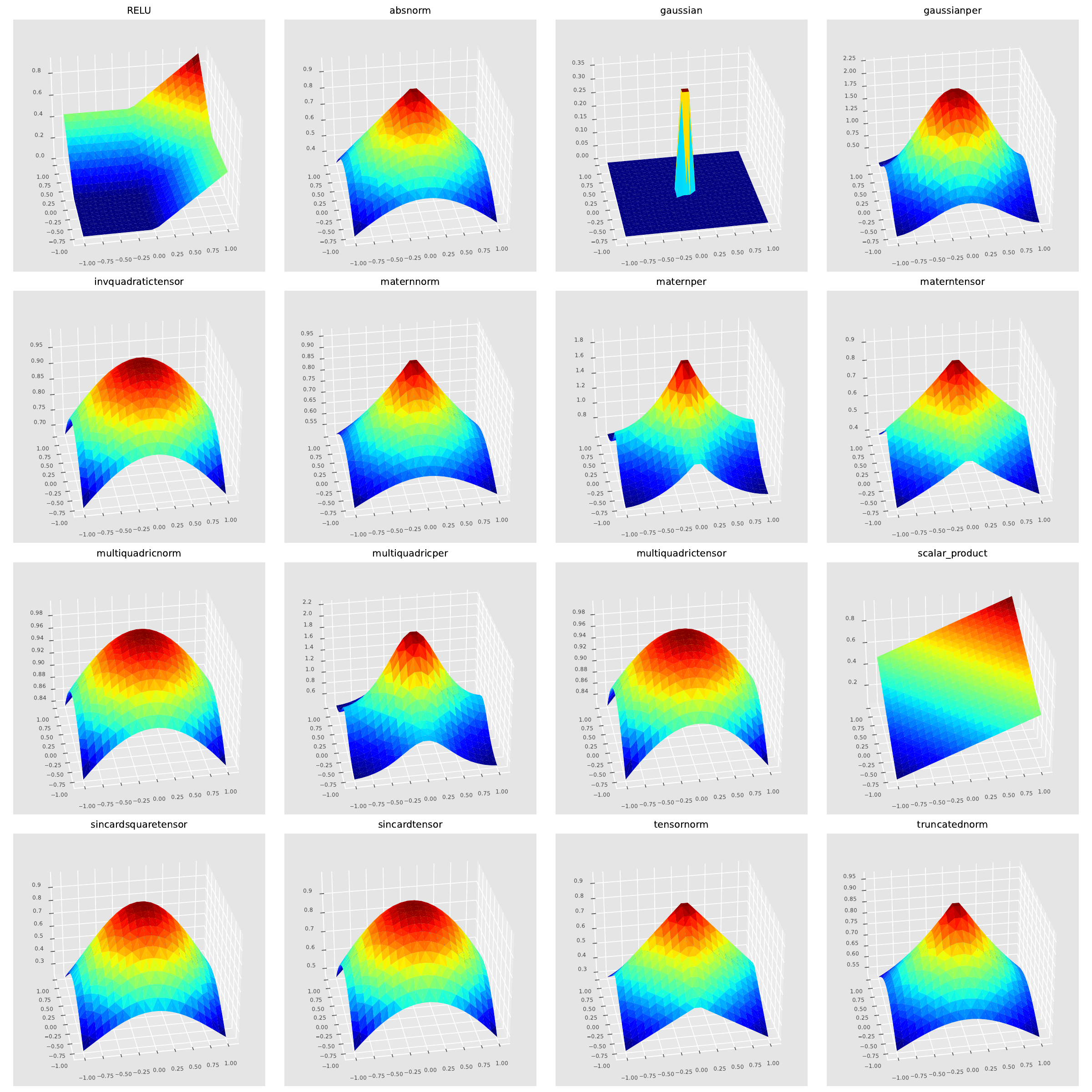}
\caption{\label{fig:308}Kernels transformed with mappings}
\end{figure}



\section{Kernel engineering}\label{kernel-engineering}

\subsection{Perturbative kernel regression}
\label{Perturbative-kernel-regression}

\paragraph{Residual kernel regression}

Perturbative kernel regression allows coupling standard calculus to kernel regression, which can be described as 
$g_{k,\theta}(\cdot) = G(\cdot,f_{k,\theta}(\cdot))$. We describe some useful constructions for further references later on. 
The following approach is a simple, yet useful, construction, while considering regressions that are considered as perturbations of given maps, overloading the kernel ridge regression \eqref{FIT} to output
\be \label{RKR}
 g_{k,\theta}(x) = g(x) + f_{k,\theta}(x), \quad \nabla g_{k,\theta}(x) = (\nabla g)(x)+ \nabla f_{k,\theta}(x).
\ee
In particular, the choice $g(x) = x, (\nabla g)(x) = I_D$, where $I_D$ is the $D$ dimensional identity matrix, is a simple construction allowing to define RKHS regressors satisfying the relation $(\nabla_k x)(X) = I_D$, which is useful in several situations modeling mappings.

\paragraph{Kernel classifiers and their derivatives}

Consider $\pi(\cdot)=(\pi_1(\cdot),\ldots, \pi_{D_\pi}(\cdot))$ be a vector-valued function of probabilities, known from observations $\pi(X)$. We extrapolate it with a kernel regressor using the softmax function; see \eqref{eq:softmax}, as follows: 
\be \label{pi_k}
 \pi_{k}(\cdot) = \text{softmax}\Big( K(\cdot,X) \theta \Big), \qquad \theta =\big( K(X,X) + \epsilon R(X,X) \big)^{-1} \ln \pi(X).
\ee 
We recall that the softmax function, used to any map any vector
\(y=(y^1,\dots,y^{|\mathcal{A}|})\) to a vector of probabilities
\(\pi=(\pi^1,\dots,\pi^{|\mathcal{A}|})\), is characterized as
\be \label{eq:softmax}
 \text{softmax}(y) = \frac{\exp(y^i)}{\sum_{j=1}^{|\mathcal{A}|} \exp(y^j)}=\pi, \quad \frac{\partial}{\partial y^j}\text{softmax}(y^i) = \pi^i(\delta^i(j)-\pi^j).
\ee
The softmax pseudo-inverse is defined as
\(y=\ln \pi\). The gradient of this regressor, modeling $\nabla \pi(\cdot)$,  takes into account the derivative expression of the softmax function \eqref{eq:softmax} as follows

\be
 \nabla \pi_{k, \theta}(\cdot) = \Big( \nabla K (\cdot,X) \theta \Big)\Big( \pi^i_k(\delta^i(j)-\pi^j_k)(\cdot)\Big),
\ee 
where
$\delta^i(j)$ is the Kronecker delta function and the right-hand side is a standard multiplication between matrices of size $(D,D_\pi)$ and $(D_\pi,D_\pi)$.

\subsection{Operations on kernels}
\label{operations-on-kernels}

We now present some operations that can be performed on kernels, and allow us to produce new, relevant kernels. These operations preserve the positivity property which we require for kernels. In this discussion, we are given several kernels denoted by \(k_i(x,y) : \mathbb{R}^D, \mathbb{R}^D \to \mathbb{R}\) (with \(i=1,2, \ldots\)) and their corresponding matrices are denoted by \(K_i\). According to \eqref{OP}, we define the projection operators, considering $\epsilon=0$ for simplicity:
\be 
  \mathcal{P}_{k_i,Y_i}(\cdot,X_i)= K_i(\cdot,Y_i)K_i(X_i,Y_i)^{-1} \in \mathbb{R}^{N_{x_i}}, \qquad i = 1,2,\ldots 
\ee
There are two different possibilities.  
\begin{itemize}

\item Operations on kernels, defined on the same distribution $X_i=X$, $Y_i=Y$, usually intended to tune a kernel to a particular problem.

\item Operations on functional spaces $\mathcal{H}_{k_i,X_i,Y_i}$, as union, resulting in operations that merge the parameters set $\theta_1,\theta_2,\ldots$.
\end{itemize}

We give some elementary examples of basic operations for two kernels, but they can be combined to define more complex combinations.

\paragraph{Adding kernels}\label{adding-kernels}

The operation \(k_1 + k_2\) is defined from any two kernels and consists of adding the two kernels straightforwardly. If \(K_1\) and \(K_2\) are the kernel matrices associated with the kernels \(k_1\) and \(k_2\), then we define the sum as \(K(X,Y)\in \mathbb{R}^{N_x, N_y}\) with corresponding projection \(\mathcal{P}_{k}(\cdot) \in \mathbb{R}^{N_x}\), as follows:
\be\label{add}
K(X,Y)=K_1(X,Y)+K_2(X,Y), \quad \mathcal{P}_{k}(\cdot) = K(\cdot, X)K(X,Y)^{-1}.
\ee
The functional space generated by \(k_1+k_2\) is then\\
\be
\mathcal{H}_k = \Big\{\sum_{1 \leq m \leq N_x} a^m (k_1(\cdot, x^m)+k_2(\cdot, x^m))\Big\}. 
\ee


\paragraph{Multiplying kernels}\label{multiplying-kernels}

A second operation \(k_1 \cdot k_2\) is also defined from any two kernels and consists in multiplying them together. A kernel matrix \(K(X,Y)\in \mathbb{R}^{N_x, N_y}\) and a projection operator \(\mathcal{P}_{k}(\cdot) \in \mathbb{R}^{N_x}\) corresponding to the product of two kernels are defined as
\be\label{mul}
K(X,Y)=K_1(X,Y) \circ K_2(X,Y), \qquad \mathcal{P}_{k}(\cdot) = K(\cdot,Y)K(X,Y)^{-1}, 
\ee
where \(\circ\) denotes the Hadamard product of two matrices. The functional space generated by \(k_1 \cdot k_2\) is
\be
\mathcal{H}_k = \Big\{\sum_{1 \leq m \leq N_x} a^m k_1(\cdot, x^m) \, k_2(\cdot, x^m) \Big\}.
\ee

\paragraph{Convolution kernels}

Our next operation, denoted by \(k_1 \ast_Z k_2\), is defined for any two kernels and consists in multiplying their kernel matrices \(K_1\) and \(K_2\) together, as follows: 
\be\label{ast}
K(X,Y)=K_1(X,Z) K_2(Z,Y), 
\ee
where \(K_1(X,Z) K_2(Z,Y)\) stands for the standard matrix multiplication and $Z$ is a third set, eventually equal to $X$ or $Y$. The projection operator is given by \(\mathcal{P}_{k}(\cdot) = K(\cdot,X)K(X,Y)^{-1}\). This amounts to considering the following kernel $k(x,y) = \int k_1(x,z)k_2(z, y)dZ$. This later formula corresponds to a standard convolution in several cases of interest.


\subsection{Operations on functional spaces}
\label{operations-on-functional-spaces}

\paragraph{Piped kernels}

So far, we considered operations on kernels, but we can also operate on functional spaces.
Consider two kernels $k_1,k_2$, two sets of points or features $X_1,X_2$, with eventually $X_2 \subset X_1$, and denote $X=[X_1,X_2]$ the  joint law.
In addition to kernel summation, we can consider summing kernel functional spaces, that is, by considering the following functional space:
\be
\mathcal{H}_{k,X} = \Big\{\sum_{0 < m \leq N_{X_1}} \theta^m k_1(\cdot, x^m) + \sum_{N_{X_1} < m \leq N_{X_1}+N_{X_2}} \theta^{m+N_x} k_2(\cdot, x^m)  \Big\},
\ee
which is a space having $X$ as a features set, and $N_{X_1}+N_{X_2}$ coefficients or parameters $\theta$. If $\mathcal{H}_{k_1,X_1} \cap \mathcal{H}_{k_2,X_2} \neq \emptyset$, there is no uniqueness of the decomposition $f=f_1+f_2$, with $f_i \in \mathcal{H}_{k_i,X_i}$.

In order to select a unique one, let us introduce a new kernel denoted as \(k = k_1 | k_2\) and we proceed by writing first the projection operator \eqref{OP} as follows:
\be\label{pipe}
  \mathcal{P}_{k}(\cdot) = \mathcal{P}_{k_1}(\cdot) \pi_1 + \mathcal{P}_{k_2}(\cdot)\Big( I_d-\pi_1\Big), 
\ee
where we defined the projection operator on the image (respectively the null space) as $\pi_1$ (respectively $I_d-\pi_1$) as
\be
 \pi_1= \mathcal{P}_{k_1}(X_1).
\ee
Hence, we split the projection operator \(\mathcal{P}_{k}(\cdot)\) into two parts. The first part deals with a single kernel, while the second kernel handles the remaining error. We define its inverse matrix by concatenation:
\be
  K^{-1}(X,Y) = \Big( K_1(X,Y)^{-1}, K_2(X,Y)^{-1}\big( I_{d}-\pi_1\big)\Big) \in \mathbb{R}^{2N_y, N_x}.
\ee
The kernel matrix associated to a ``piped kernel`` pair is then
\be
  K(X,Y) = \Big(K_1(X,Y), K_2(X,Y) \Big) \in \mathbb{R}^{N_x, 2 N_y}. 
\ee
Piping the two kernels \(k_1 | k_2\) is similar to applying a Gram-Schmidt orthogonalization process of the functional spaces \(\mathcal{H}_{k_1,X}\), \(\mathcal{H}_{k_2,X}\).


\paragraph{Piping scalar product kernels: an example with a polynomial regression}

Consider a map \(S : \mathbb{R}^D \to \mathbb{R}^N\) associated with a family of \(N\) basis functions denoted by \(\varphi_n\), namely \(S(x) = \big(\varphi_1(x),\ldots,\varphi_N(x)\big)\). Let us introduce the \emph{dot product kernel}
\be\label{LRK}
    k_1(x,y)= \langle S(x),S(y)\rangle,  
\ee
which can be checked to be conditionally positive-definite. Let us also consider a pipe kernel denoted as \(k_1 | k_2\), where \(k_1\) and \(k_2\) are positive kernels. This construction becomes especially useful in combination with a polynomial basis function \(S(x) = \big(1,x_1,\ldots\big)\). The pipe kernel enables a classical polynomial regression, allowing an exact matching of the moments of a distribution. Namely, any remaining error can be effectively handled by the second kernel \(k_2\). Importantly, this combination of kernels provides a powerful framework for modeling and capturing complex relationships between variables.


\section{Kernel extrapolation}
\label{sec-kern}

\subsection{Inverse of a kernel matrix and reproducibility property}

We now illustrate some aspects of the projection formula \eqref{FIT}. The inverse of a kernel matrix \( K(X, Y)^{-1} \) is computed differently depending on whether \( X = Y \) or \( X \neq Y \).
\bei 

\item[] {\bf Case 1: \( X = Y \).} When \( X = Y \), the inverse is computed using the formula:
\be
K(X, X)^{-1} = (K(X, X) + \epsilon R)^{-1}, 
\ee
where \( \epsilon \geq 0 \) is an optional regularization term, known as the \emph{Tikhonov regularization} parameter. This regularization\footnote{In the CodPy Library, the default value of \( \epsilon \) is \( 10^{-8} \), corresponding to the reproducible modes, but it can be adjusted as needed.} is often necessary to improve numerical stability. The matrix \( R \) is typically chosen as the identity matrix \( I_d \) of size \( N_x \times N_x \); see also Figure~\ref{fig:denoising} for an alternative choice. 

\item[] {\bf Case 2: \( X \neq Y \).} When \( X \neq Y \), the inverse is computed using the least-squares method:
\be\label{pseudo-inverse}
  K(X,Y)^{-1} = (K(Y,X)K(X,Y) + \epsilon R)^{-1} K(Y,X),
\ee
where \( R \) now has dimensions \( N_y \times N_y \).

\eei

The matrix product \( K(X, X) K(X, X)^{-1} \) may not always coincide with the identity matrix. This discrepancy can arise in the following cases.
\begin{itemize}
    \item If \( \epsilon > 0 \), the Tikhonov regularization parameter modifies the inverse for improved numerical stability.     
    \item If the chosen kernel is not \textit{strictly} positive-definite and leads to a kernel matrix \( K(X,X) \), which is not of full rank. For instance, when using a linear regression kernel (see Section~\ref{kernel-engineering}), the projection operator \( \mathcal{P}_k(\cdot)\) behaves as a projection onto the image space of \( k(\cdot,\cdot) \), which only captures some moments of functions.
\end{itemize}


\subsection{Computational complexity of kernel methods}
\label{proposed-methodology}

Our algorithms provide us with general functions in order to make predictions, once a kernel is chosen. That is, considering the projection operator \eqref{OP}, considered here without regularization ($\epsilon=0$) for simplicity
\be\label{Pk}
f_k(\cdot)  = \mathcal{P}_{k,Y}(\cdot,X)f(X) = K(\cdot,Y)K(X,Y)^{-1}f(X), 
\ee
defines a supervised learning machine, which we call a \emph{feed-forward operator}, and \(\mathcal{P}_{k}(\cdot) \in \mathbb{R}^{N_y}\) is the  \emph{projection operator} \eqref{OP}, as it realizes the projection of any function on the discrete space \(\mathcal{H}_{k,Y}\) from observed values $f(X)$. Observe that \eqref{Pk} includes two contributions, namely the kernel matrix \(K(X,Y)\) and the \emph{projection set of variables} denoted by \(Y \in \mathbb{R}^{N_y, D}\).

To motivate the role of the argument \(Y\), let us consider the particular choices of the reproducible mode in \eqref{FIT}, which \emph{do not depend} upon \(Y\).
\begin{align}\label{EI}
 &\text{Extrapolation operator: } \mathcal{P}_{k,X}(\cdot,X) = K(\cdot,X)K(X, X)^{-1}. 
\end{align}
In some applications, these operators may lead to certain computational issues, due to the fact that the kernel matrix \(K(X, X) \in \mathbb{R}^{N_x, N_x}\) must be inverted, as is clear from \eqref{EI},
and the overall algorithmic complexity of \eqref{Pk} is of the order
\be
D \, \big( N_x \big)^3. 
\ee
This is a rather costly computational process when faced with a large set of input data. Specifically, this is our motivation for introducing the additional variable \(Y\) which has the effect of lowering the computational cost. When computing $f_k(Z)$ on a distribution $Z$, the overall algorithmic complexity of \eqref{Pk} is of the order
\be
D \, \big( (N_y)^3 + (N_y)^2 N_x + (N_y)^2 N_z \big). 
\ee
Importantly, the projection operator \(\mathcal{P}_{k,Y}(Z,X)\) is \textit{linear} in terms of, both, input and output data $X$ and $Z$. Hence, while keeping the set \(Y\) to a reasonable size, we can consider a large set of data, as input or output, at the expense of losing the reproducibility property. This approach led to reproducible RKHS methods for large datasets, with similar algorithmic complexity; see Section~\ref{large-scale-dataset}.

Choosing a well-adapted set \(Y\) is often a major source of optimization. We are going to use this idea intensively in several applications. For example, kernel clustering methods (which we will describe later) aim at minimizing the error implied by kernel ridge regression with respect to the set \(Y\). This technique also connects with the idea of \emph{sharp discrepancy sequences} to be defined later. 


\subsection{Deep kernel architecture}
\label{deep-kernel-architectures}

Neural networks (NNs) are also kernel-based methods, sharing close similarities with RKHS methods. The key distinction is that NNs consider non-symmetrical, hence non-positive-definite kernels . Let us describe deep-learning construction for supervised learning, considering a distribution $X \in \mathbb{R}^{D_x}$ and a function $f(X) \in \mathbb{R}^{D_{f}}$. Borrowing vocabulary from the neural network community, given $N$ layers, each of them containing $D^n$ neurons, we define a matrix having prescribed size parameters $\theta^n \in \mathbb{R}^{ D_n,D_{n-1}}$ and a function $\sigma^n$, called an \textit{activation function}. A deep-learning method of depth $N$ is defined as the following recurrent construction:
\be\label{DNN}
f_{\sigma,\theta}(x) = \theta^{N} \sigma^N\Big(\theta^{N-1} \sigma^{N-1} \big( \ldots \sigma^1(\theta^{0} \ x ) \big)\Big), 
\ee
with the convention $\theta^{0} \in \mathbb{R}^{ D_0,D_x}$, $\theta^{N} \in \mathbb{R}^{D_N,D_f}$. With one layer, the function $x \mapsto \sigma(\theta \ x)$ is called a perceptron, and $x \mapsto \theta^1 \sigma(\theta^{0} \ x)$ adds a linear layer to the perceptron, etc. Let us denote $\theta = \theta^1,\ldots,\theta^N$. A general fitting procedure for the constructions \eqref{DNN} considers the following problem:
\be\label{FITDNN}
\overline{\theta} = \arg \inf_{\theta}\| f(X)- f_{k,\theta}(X)\|, \quad \theta = (\theta^0, \ldots, \theta^N)
\ee
where $\| \cdot \|$ is a general \textit{loss function}, usually tackled with a descent approach, also called \textit{back--propagation}, typically implemented using the stochastic gradient algorithm Adam. 

The kernel construction \eqref{FIT}, that is, $f_{k,\theta}(x) = k(x,Y)\theta$, can be interpreted as a double layer, one being a linear layer. Considering at each level a kernel $k^n$, a deep kernel architecture can be described in an RKHS setting as follows:
\be \label{DKN}
f_{k,\theta}(x) = k^N\Big( k^{N-1}\big( \ldots ,Y^{N-1}\big)\theta^{N-1},Y^N\Big)\theta^N. 
\ee
For these deep architectures, each level consists of transforming a feature distribution $X^{n-1}$ to the following one $X^n = k^{n-1}(X^{n-1},Y^{n-1})\theta^{n-1}$ (like Russian dolls matryoshkas). 

Stochastic gradient algorithms might be used to fit the parameters $\theta$ in \eqref{DKN}. However, if we consider the standard mean square error loss function, deep kernel architectures can be fit using \eqref{FIT}, which is an efficient computational approach. 

For kernel architectures, it is not clear whether several layers are beneficial or not. Indeed, in practice, two-layer architectures are often enough:
\be 
f_{k,\theta}(x) = k^1\Big( k^{0}\big( x ,Y^{0}\big)\theta^{0},Y^1\Big)\theta^1, 
\ee
in which we call the distribution $Y^1$ \textit{latent distribution}, lying in a space $\mathbb{R}^{D_1}$ called \textit{latent space}. To define a smooth, invertible mapping $x \mapsto  k^{0}\big(x ,Y^{0}\big)\theta^{0}$, is the main topic of Section~\ref{Maps-and-Generative-methods}, devoted to mappings and generative methods.


\subsection{Basic numerical examples}
\label{basic-numerical-examples}

\paragraph{Test description} 

In most applications, we are given $X,f(X)$, also called \textit{training set} in the machine-learning vocabulary, and seek from this set to infer the value $f(Z)$ on another set $Z$, also called \textit{test set} in the machine-learning vocabulary. This is illustrated now with some simple function extrapolation problems, using the formula \eqref{FIT}. 

In our first test, we use a generator that selects \(X\) (respectively \(Y, Z\)) as \(N_x\) (respectively \(N_y, N_z\)) points generated regularly (respectively randomly, regularly) on a cube (a segment if $D=1$) $[-1,+1]^D$, and define the function \(f\) as a sum of a periodic and a linear, polynomial function:
\begin{equation} \label{2D}
f(X) = \Pi_{d=1,\ldots,D} \cos (4\pi x_d) + \sum_{d=1..D} x_d.
\end{equation}
To observe extrapolation and interpolation effects, a validation set \(Z\) is distributed over a larger cube $[-1.5,+1.5]^D$. As an illustration, in Figure~\ref{fig:xfxzfz} we show both graphs \((X, f(X))\) (left, training set),\((Z, f(Z))\) (right, test set) for $D=1$.

\begin{figure}
\centering
\includegraphics[width=0.7\textwidth, keepaspectratio]{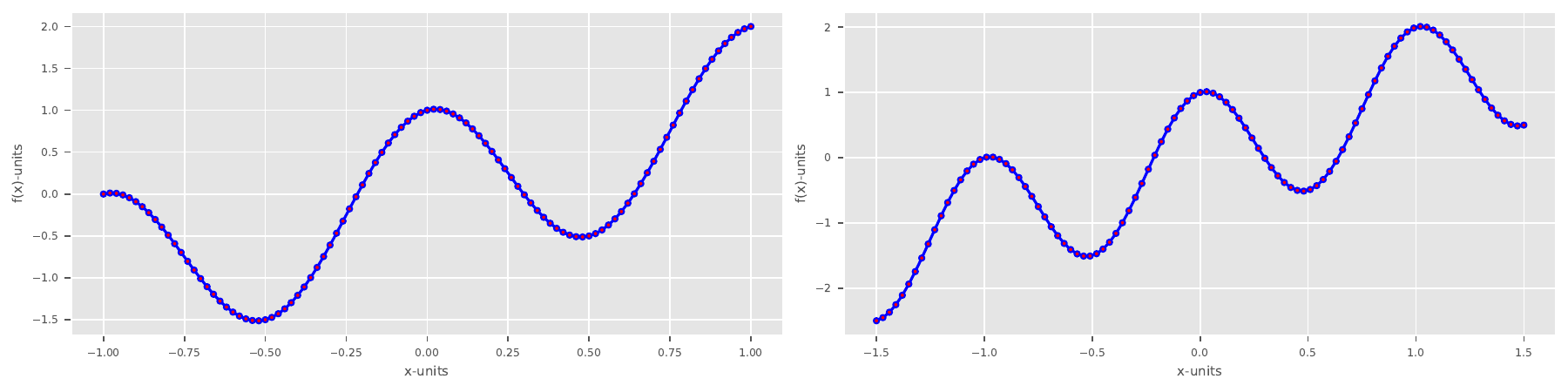}
\caption{\label{fig:unnamed-chunk-10}\label{fig:xfxzfz}Training set $(x,f(x))$ (left) and test set $(z,f(z))$ (right)}
\end{figure}

\paragraph{Guidelines to benchmark methodology}

We profit from this section to provide some guidelines for a general benchmark methodology of predictive machines. Quantitative benchmark methods are usually based on two criteria.
\begin{itemize}
\item
  The first one is a score, which is a metric quantifying the quality of the prediction. In this simple illustration, we will be using the mean-squared error (MSE) $\frac{1}{N_z}\|f_k(Z)-f(Z)\|_{\ell^2}$.
\item
  The second one is the execution time. A prediction machine should be always thought in terms of computational efficiency: with a given computational budget or electricity consumption, a given algorithm can reach a given accuracy. 
\end{itemize}

\paragraph{Guidelines to kernel parameters}

As shown in the previous sections, the external parameters of a kernel-based prediction machine typically consist of \emph{a positive-definite kernel} function and \emph{a map}. Most of the maps that we use for our kernels are \emph{stateful}, in the sense that they depend on external parameters, however these parameters do not require to be input, as they are fit automatically to arbitrary set $X$; see Section~\ref{maps-and-kernels}.

However, in the formula \eqref{FIT}, the set $Y$ is usually understood as an optional degree of freedom, which is an external parameter set. We distinguish between several options.
\begin{itemize}
\item
  First, we can choose \(Y=X\), which corresponds to the \emph{extrapolation}, satisfying the reproducing property, typically resulting in highest accuracy at the expense of heavier computational effort. 
\item
  Alternatively, we can select another, usually smaller set for \(Y\). The purpose here is usually to trade accuracy for execution time and is better suited for larger training sets. Strategies to select such a set are discussed in Section~\ref{clustering-strategies}
\end{itemize}

\paragraph{A qualitative comparison between kernels}

To illustrate the impact of different kernels and maps on our extrapolation machine, we consider a 
one-dimensional test and compare the predictions achieved by using various kernels; see Figure~\ref{fig:311-5}. As one can see, gluing together a periodic kernel with a polynomial kernel through piping can capture both components of the function.

\begin{figure}
\centering
\includegraphics[width=0.7\textwidth, keepaspectratio]{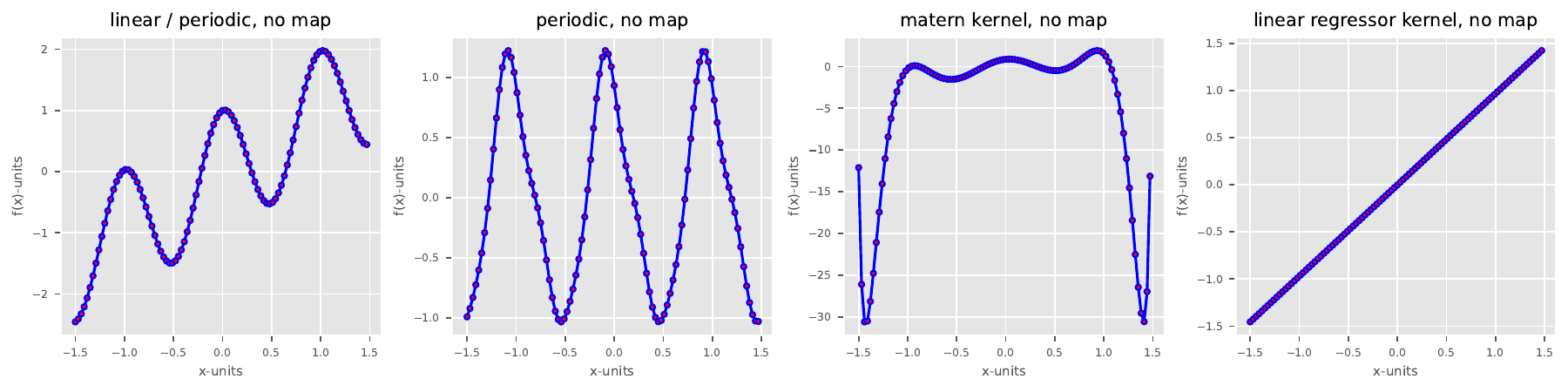}
\caption{\label{fig:311-5}A qualitative comparison between kernels}
\end{figure}

\paragraph{A qualitative comparison between methods} 

Reproducing kernel Hilbert space regressors equipped with \emph{universal} kernels are universal approximators and can, in principle, model broad classes of functions in arbitrary dimension. In practice, strong non-kernel baselines include neural networks and tree-based ensembles. We also consider support vector regression (SVR) with an RBF kernel. Choosing among these methods requires a careful, reproducible benchmarking protocol. 

In this test we evaluate kernel ridge regression (KRR)~\eqref{OP}\footnote{Implemented with CodPy with a periodic kernel \url{https://codpy.readthedocs.io/en/dev/}}, and compare it with two standard regression models: a feed-forward neural network (FFN)\footnote{Implemented with PyTorch \url{https://pytorch.org}}, support vector regression {SVR) with RBF kernel, a decision tree (DT), a random forest (RF), Adaboost\footnote{Implemented with scikit-learn \url{https://scikit-learn.org/stable/}} and XGBoost\footnote{Implemented with the XGBoost Library \url{https://xgboost.readthedocs.io/en/stable/}}.

\begin{figure}
\centering
\includegraphics[width=0.7\textwidth, keepaspectratio]{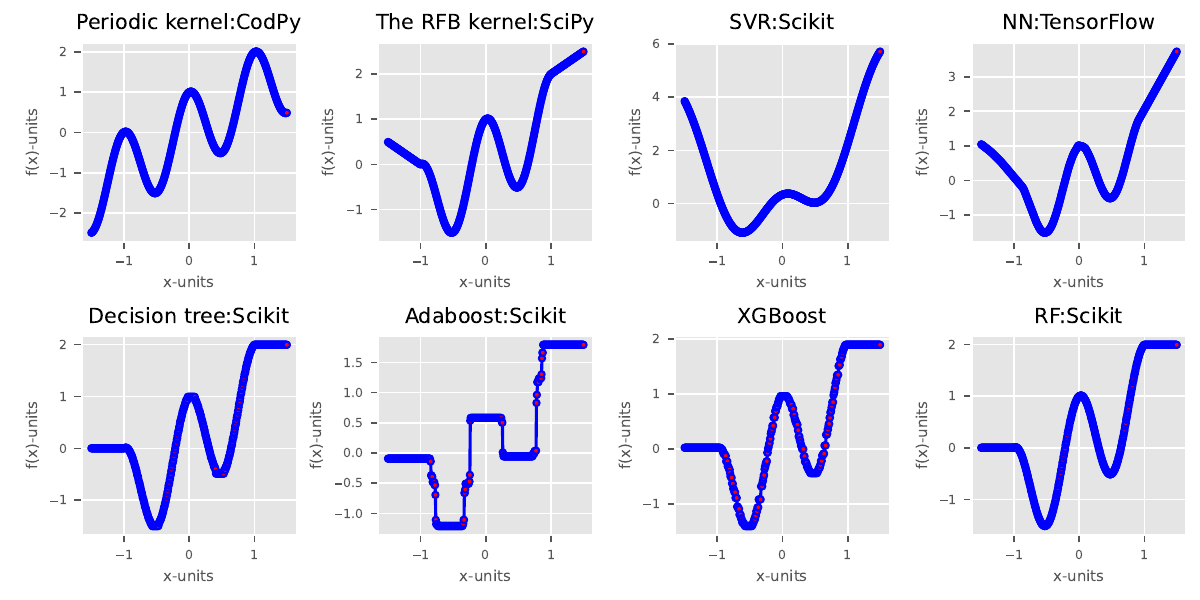}
\caption{Periodic function extrapolation test with KRR, SVR, FFN, DT, Adaboost, XGBoost, RF
\label{fig:a1a2a3}}
\end{figure}

\begin{figure}
\centering
\includegraphics[width=0.7\textwidth, keepaspectratio]{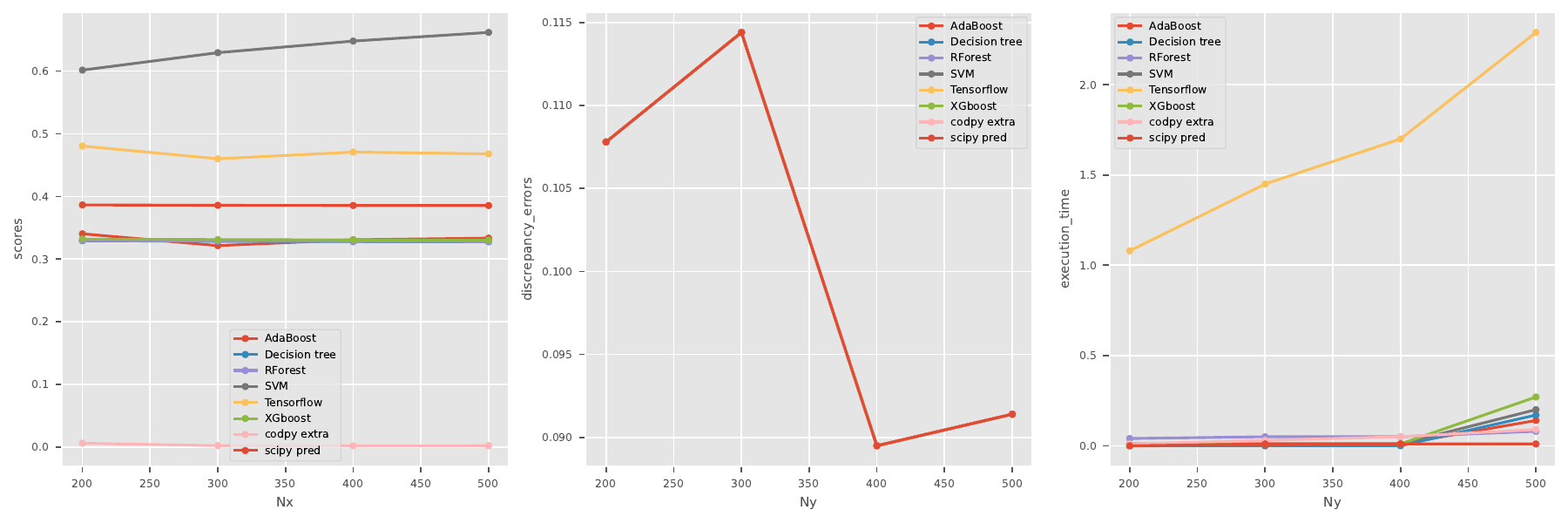}
\caption{\label{fig:unnamed-chunk-12}\label{fig:a1a2a4}RMSE, MMD, and execution time}
\end{figure}

In Figure~\ref{fig:a1a2a3}, we can observe the extrapolation performance of each method. It is evident that the periodic kernel-based method outperforms the other methods in the extrapolation range between \([-1.5,-1]\) and \([1,1.5]\). This finding is also supported by Figure~\ref{fig:a1a2a4}, which shows the RMSE error for different sample sizes \(N_x\).

Observe that the choice of method does not affect the function norms and the discrepancy errors. Although the periodic kernel-based method performs better in this example, our goal is not to establish its superiority, but to present a benchmark methodology, especially when extrapolating test set data that are far from the training set.

\paragraph{A quantitative comparison between methods}

We demonstrate that the dimensionality of the problem, denoted by $D$, does not affect neither the code interface to our extrapolation method, nor its performance. In term of code, the dimension $D$ is a simple input parameter for this test.

To illustrate this point, we repeat the same steps as in the previous section but simply prescribing \(D=2\) (i.e., a two-dimensional case). For the two-dimensional case, we still can easily plot the training set and the test set; see Figure~\ref{fig:xfxzfz2}. 

The reader can easily test other different values of dimensions \(D\) in our code site. If the dimensionality is greater than two, visualization of input and output data are more difficult. A simple choice is to use a two-dimensional visualization by plotting \(\widetilde{X},f(X)\), where \(\widetilde{X}\) is obtained either by choosing two indices $i_1,i_2$ and plot \(\widetilde{X}=X[i_1,i_2]\) or by performing a principal component analysis (PCA) over \(X\) and setting \(\widetilde{X}=\text{PCA}(X)[i_1,i_2]\).

\begin{figure}
\centering
\includegraphics[width=0.7\textwidth, keepaspectratio]{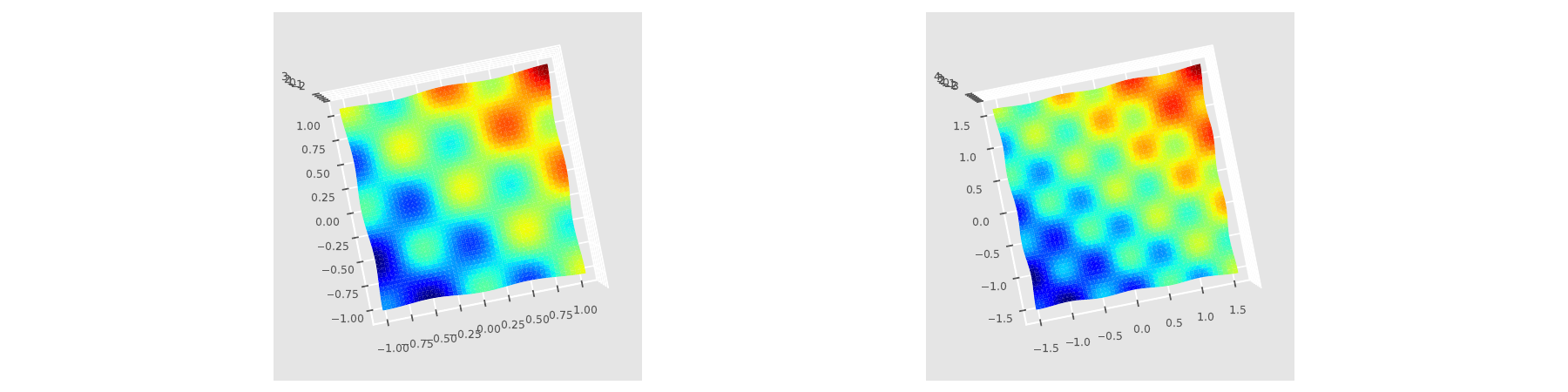}
\caption{\label{fig:xfxzfz2}Train set vs. test set}
\end{figure}

We generate data using some scenarios and visualize the results using Figure~\ref{fig:a1a2a5}. The left and right plots show the training set \((X,f(X))\) and the test set \((Z,f(Z))\), respectively. Observe that \(f\) is the two-dimensional periodic function defined at \eqref{2D}.

\begin{figure}
\centering
\includegraphics[width=0.7\textwidth, keepaspectratio]{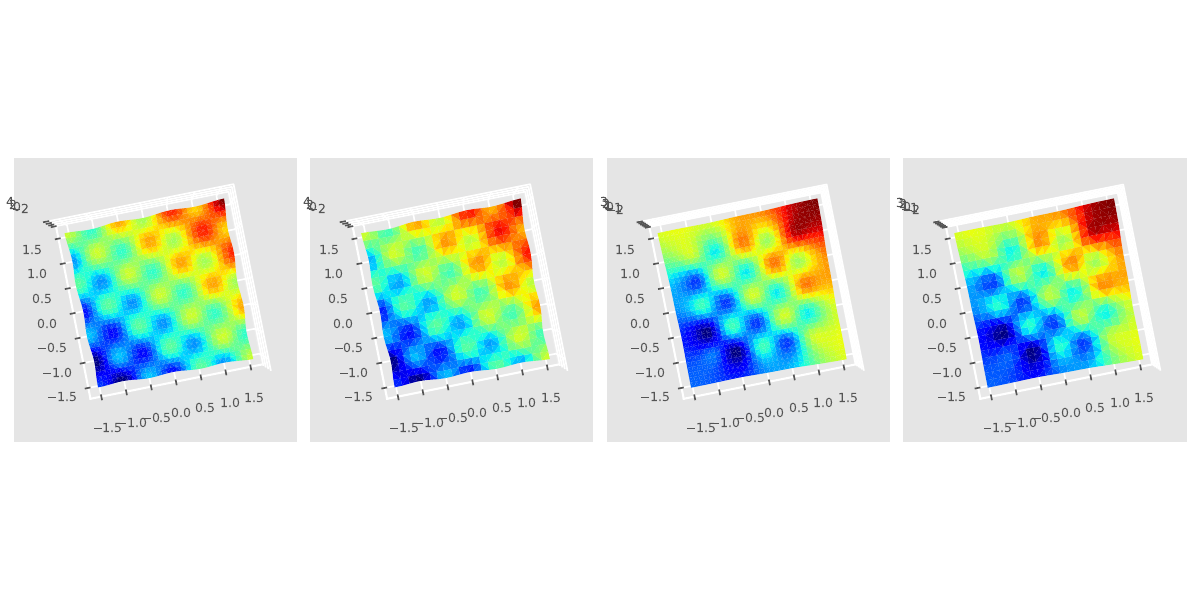}
\caption{\label{fig:a1a2a5}RBF (first and second) and periodic Gaussian kernel (third and fourth)}
\end{figure}

\paragraph{Maps can ruin your prediction}

Drawing upon the notation introduced in the preceding chapter, we examine the comparison between the ground truth values \((Z,f(Z)) \in \mathbb{R}^{N_z,D} \times \mathbb{R}^{N_z,D_f}\) and the corresponding predicted values \((Z,f_k(Z)) \). To further clarify the role of distinct maps in computation, we rely on a particular map referred to as the mean distance map. This map scales all points to the average distance associated with a Gaussian kernel. The resulting plot, presented in Figure~\ref{fig:310}, underscores the substantial influence of maps on computational results.

It is crucial to observe that the effectiveness of a specific map can differ significantly depending upon the choice of kernel. This variability is further illustrated further in Figure~\ref{fig:310}.

\begin{figure}
\centering
\includegraphics[width=0.7\textwidth, keepaspectratio]{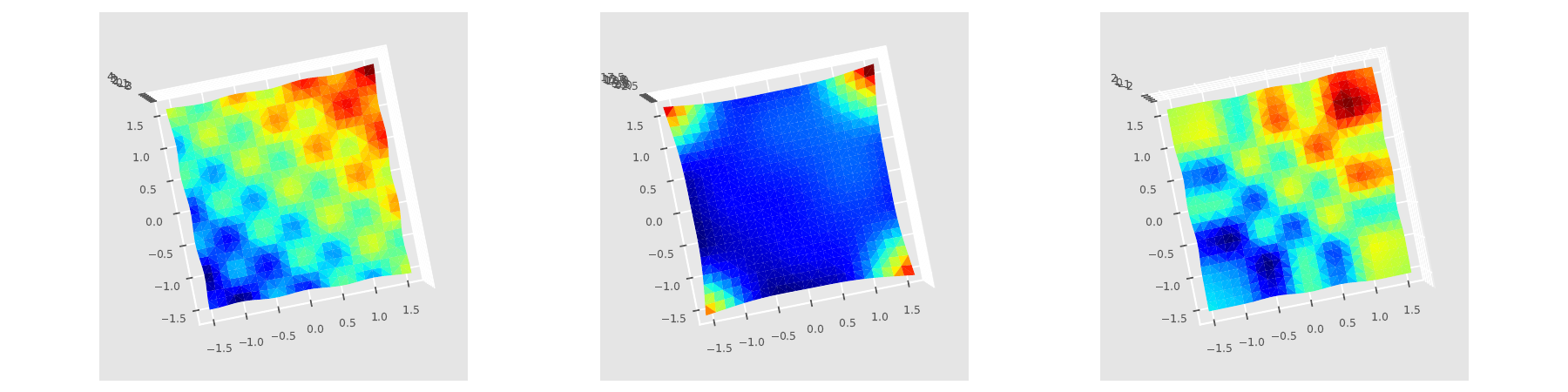}
\caption{\label{fig:unnamed-chunk-24}\label{fig:310}A ground truth value (first), Gaussian (second), and Mat\'ern kernels (third) with mean distance map}
\end{figure}


\section{Error measurements with discrepancy}
\label{sec-err}

\subsection{Distance matrices}

Distance matrices provide a valuable tool for assessing the accuracy of computations. Given a positive kernel \( k : \mathbb{R}^D \times \mathbb{R}^D \to \mathbb{R} \), we associate the \emph{distance function} \( d_k(x, y) \), defined for \( x, y \in \mathbb{R}^D \) as
\be\label{ES}
d_k(x,y) = k(x,x) + k(y,y) - 2k(x,y).
\ee
For positive kernels, \( d_k(\cdot, \cdot) \) is continuous, non-negative, and satisfies the condition \( d_k(x, x) = 0 \) for all relevant \( x \).

Given two collections of points, \( X = (x^1, \dots, x^{N_x}) \) and \( Y = (y^1, \dots, y^{N_y}) \) in \( \mathbb{R}^D \), we define the associated \emph{distance matrix} \( D(X,Y) \in \mathbb{R}^{N_x \times N_y} \) as
\be
    D(X,Y) = 
\begin{bmatrix}
        d_k(x^1,y^1) & d_k(x^1,y^2) & \ldots & d_k(x^1,y^{N_y}) \\
        d_k(x^2,y^1) & d_k(x^2,y^2) & \ldots & d_k(x^2,y^{N_y}) \\
        \vdots & \vdots & \ddots & \vdots \\
        d_k(x^{N_x},y^1) & d_k(x^{N_x},y^2) & \ldots & d_k(x^{N_x},y^{N_y})
\end{bmatrix}.
\ee
Distance matrices play a crucial role in numerous applications, particularly in clustering and classification tasks.
Table~\ref{tab:306} presents the first four rows and columns of a kernel-based distance matrix \( D(X,X) \). As expected, all diagonal values are zero.

\begin{table}[htbp]
\caption{\label{tab:306} First four rows and columns of a kernel-based distance matrix \( D(X,X) \)}
\centering
\begin{tabular}{r|r|r|r}
\hline
0.00 & 0.08 & 0.16 & 0.24 \\
\hline
0.08 & 0.00 & 0.08 & 0.16 \\
\hline
0.16 & 0.08 & 0.00 & 0.08 \\
\hline
0.24 & 0.16 & 0.08 & 0.00 \\
\hline
\end{tabular}
\end{table}


\subsection{Kernel maximum mean discrepancy functional}

We now deal with an interesting aspect of the discrepancy functional, $d_k(\cdot,X)$, plotted in Figure~\ref{fig:unnamed-chunk-38} (respectively Figure~\ref{fig:unnamed-chunk-39}) for three kernels $k$, and a simple random one-dimensional (respectively two-dimensional) distribution \(X \in \mathbb{R}^{N_x}\). 

\textit{An example of smooth kernel}. First, consider the discrepancy induced by a Gaussian kernel 
\be
k(x,y) = \exp(-(x-y)^2),
\ee
which generates functional spaces made of smooth functions.
In Figure~\ref{fig:MMD1}, we show the function $D_x(\cdot,X)$ in red color. Additionally, we display in blue  a linear interpolation of the function \(d_k(x,x^n)\), \(n=1 \ldots N_x\) in Figure~\ref{fig:MMD1} to demonstrate that this functional is smooth but neither convex nor concave. Notably, the minimum of this functional is achieved by a point which is not part of the original distribution \(X\). For a two-dimensional example, we refer to Figure~\ref{fig:MMD3} (left-hand) for a display of this functional.

\textit{An example of Lipschitz continuous kernels: ReLU}. Let us now consider a kernel that generates a functional space with less regularity. The ReLU kernel is the following family of kernels which essentially generates the space of functions with bounded variation: 
\be
k(x, y) = \max(1 -|x-y|,0). 
\ee
As shown in Figure~\ref{fig:MMD1} (middle), the function \(y \mapsto d_k(x,y)\) is only piecewise differentiable. Hence, in some cases, the functional \(d_k(x,y)\) might have an infinite number of solutions (if a ``flat'' segment occurs), but a minimum is attained on the set \(X\). Figure~\ref{fig:MMD3} (middle) displays the two-dimensional example.

\textit{An example of continuous kernel: Mat\'ern}. The Mat\'ern family generates a space of continuous functions, and is defined by the kernel 
\be
k(x,y) = \exp(-|x-y|). 
\ee
In Figure~\ref{fig:MMD1}, we observe that the function \(y \mapsto d_k(x,y)\) has concave regions almost everywhere, and the minimum of the functional is a point of $X$.

This form of the discrepancy implies that a global minimum of the functional $d_k(\cdot,X)$ should be looked first as an element of $X$, involving usually combinatorial algorithms. Then, depending on the kernel, a descent algorithm is necessary to reach a global minimum. 

\begin{figure}
\centering
\includegraphics[width=0.99\textwidth, keepaspectratio]{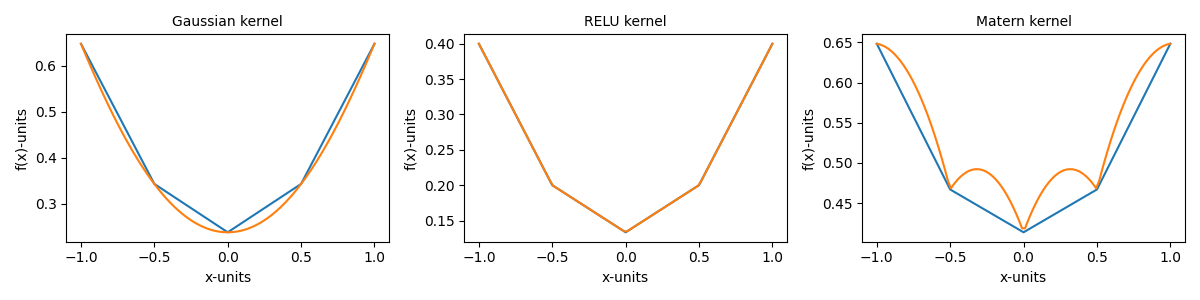}
\caption{\label{fig:unnamed-chunk-38}\label{fig:MMD1} Distance functional for the Gaussian, Mat\'ern, and  ReLU kernels (1D)}
\end{figure}

\begin{figure}
\centering
\includegraphics[width=0.99\textwidth, keepaspectratio]{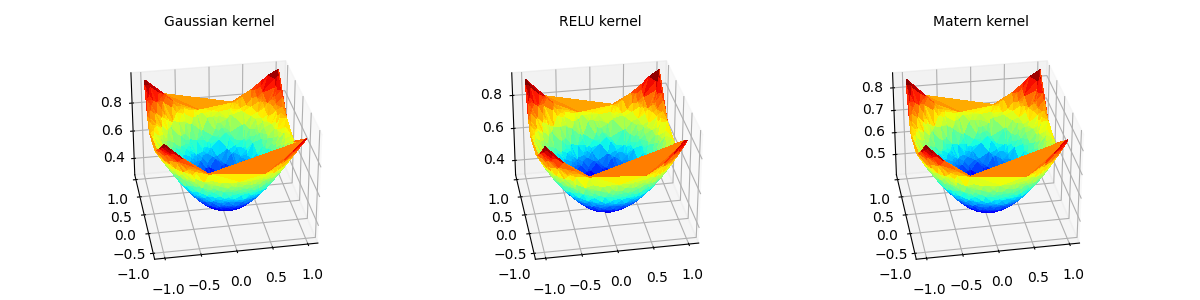}
\caption{\label{fig:unnamed-chunk-39}\label{fig:MMD3} Distance functional for the Gaussian, Mat\'ern, and  ReLU kernels (2D)}
\end{figure}


\clearpage 

\chapter{Discrete operators based on reproducing kernels}\label{kernel-based-operators}

\section{Objective of this chapter}
\label{introduction-1}

We now define and study classes of operators constructed from a reproducing kernel. We begin by introducing interpolation and extrapolation operators, which play a central role in machine learning. Interpreting these constructions as operators enables the development of kernel-based approximations of discrete differential operators, including the gradient and divergence. 
These discrete operators prove useful in various contexts, particularly for modeling
physical phenomena governed by partial differential equations (PDEs).


\section{Discrete kernel operators}
 
\subsection{Standpoint} 

Consider the fitting formula \eqref{FIT}, with $\epsilon = 0$, i.e. without any regularization terms for simplicity. It motivated the earlier introduction of the operator $\mathcal{P}_{k,Y}(\cdot,X)$ in \eqref{OP}, and we call $\mathcal{P}_{k,Y}$ the \textit{projection operator}: this operator projects
any function $f \in \mathcal{H}_k$ onto its discrete representation $f_{k}(\cdot) \in \mathcal{H}_{k,Y}$ as follows: 
\be\label{nablak}
	f_k(\cdot) = \mathcal{P}_{k,Y}(\cdot,X) f(X), \quad \mathcal{P}_{k,Y} : \mathcal{H}_k \mapsto \mathcal{H}_{k,Y}
\ee
This operator-based viewpoint allows us to define various other kernel-based operators. For instance, considering the gradient formula \eqref{GRAD}, we similarly define the \textit{gradient operator}
\be
\label{nabla}
 \nabla f_{k}(\cdot) = \nabla_k(\cdot,X) f(X), \quad \nabla_k(\cdot,X) = (\nabla K)(\cdot,Y)K(X,Y)^{-1}.
\ee
The projection operator $\mathcal{P}_{k,Y}(\cdot,X) \in \mathbb{R}^{N_x}$ is a vector field, 
and the gradient $\nabla_k(\cdot) \in \mathbb{R}^{D,N_x}$ is a matrix field. When these operators are evaluated on a set $Z$, $\mathcal{P}_{k,Y}(Z,X)$ becomes a matrix of size $N_z \times N_x$, while $\nabla_k(Z,X)$ is of size $N_z \times D \times N_x$. From now on, we will use interchangeably $\nabla_k(Z,X)$, $\nabla_k(Z)$, or even $\nabla_k$, whenever a shorter notation is not ambiguous. 


\subsection{Transpose of operators and Laplace-Beltrami operator} 
\label{laplace-beltrami-operator}

Let us consider a set $Z=(z^1,\ldots,z^{N_z})$, any scalar-valued function $\phi$ (respectively  vector-valued $\varphi$) belonging to $\mathcal{H}_k$,  (respectively  $\varphi \in (\mathcal{H}_k)^D$), and consider the transpose of the gradient operator, defined as
\be
 <\nabla_k(Z) \phi(X),\varphi(Z)> = <\phi(X),\nabla_k(Z)^T \varphi(Z)>.
\ee
where $\nabla_k(Z)$ is a matrix of size $N_x,D\times N_z$. This last formula is the discrete formulation of an integration by part involving the gradient and the divergence, written formally as $\int <\nabla \phi(z),\varphi(z)> dZ = \int  \phi(x) (\nabla \cdot \varphi)(x) dX$. So we introduce the transpose operator of the gradient, homogeneous to a discrete \textit{divergence operator}, denoted $\nabla_k \cdot$ as follows:
\be
\label{div}
 (\nabla_k \cdot)(\cdot,X) = K(Y,X)^{-1}(\nabla K)(Y, \cdot).
\ee
This operator defines a matrix field $(\nabla_k \cdot)(\cdot) \in \mathbb{R}^{N_x,D}$, and operates on functions in $\mathcal{H}_{k,Z}$ consistently with a divergence operator $(\nabla_k \cdot)(Z)\varphi(X) \sim -\nabla \cdot(\varphi dZ)$, as we estimate it on a distribution $Z$.
We can now introduce the kernel \textit{Laplace-Beltrami operator}, which is constructed from these two operators as follows:
\be\label{Delta}
\Delta_k(\cdot) = \nabla_k^T(\cdot) \nabla_k(\cdot),
\ee
defining a field of vectors having size $N_x$. While evaluated on the set $X$, $\Delta_k(X)$ is a matrix of size $N_x, N_x$. This operator thus provides an
approximation of the Laplace-Beltrami operator for any function in $\mathcal{H}_{k,X,Y}$. 

The Laplace-Beltrami operator is a central notion in many areas, including fluid mechanics, image analysis and signal processing. In particular, the Laplacian arises for solving PDE boundary value problems (e.g. Poisson, Helmholtz), and is involved in many time evolution problems involving diffusion or propagation, such as the heat equation or the wave equation, as well as stochastic processes of martingale type, as we will illustrate later.


\subsection{Inverse of operators and variational formulation} 
\label{inverse-of-operators-and-variational-formulation}
\label{the-helmholtz-hodge-decomposition}


\paragraph{Gateaux derivative.}

We fix $\mu$ and seek to characterize the minimum of this functional as $u = \arg \min_{v \in H^1_\mu } J_\mu(v)$.
 This functional is associated with the equation $\mathcal{L}$ as $J_\mu'(f) = \mathcal{L}$, in the sense of Gateaux\footnote{see \url{https://en.wikipedia.org/wiki/Gateaux_derivative}} : for any functions \(\varphi \in \mathcal{D}\), the space of smooth functions, this minimum is characterized as
\be
  J_\mu(u+\varphi) - J_\mu(u) =  \int \left(<A \nabla u,\nabla \varphi> + u\varphi  - \varphi \ f \right) d\mu + \mathcal{O}(\varphi^2) \ge 0.
\ee
Taking $\psi = -\varphi$ in the last inequality shows that $\lim_{\epsilon \to 0} \frac{J(u+\epsilon\varphi) - J(u)}{\epsilon}=0$. Hence, the minimizer $u \in H_\mu^1$ satisfies the weighted weak problem:
\be
\label{WFPE}
   \int <A \nabla u,\nabla \varphi>d\mu + \int u\varphi d\mu + \int (\varphi \ f) d\mu  = 0, \quad \forall \varphi \in H_\mu^1.
\ee
When the functions $A, \mu$ and $f$ are sufficiently smooth, Green identity leads to the strong form of the equation:
\be
  -\nabla \cdot (\mu A \nabla u ) + u\mu = f \mu \quad \text{in} \ \Omega, \quad u=0 \ \text{on} \ \partial\Omega.
\ee
We next present the associated RKHS-based numerical method to solve this elliptic system. 

\paragraph{Representer theorem}

Consider a kernel $k$, its generated RKHS $\mathcal{H}_k$, and a \textit{mesh} $X=(x^1,\ldots,x^N)$ such that $\delta_X = \frac{1}{N} \sum_n \delta_{x^n} \sim d\mu$. Denote $\mathcal{H}_{k,X} \subset \mathcal{H}_k$ the discrete functional space. We approximate
\be
    J_\mu(u) \sim J_{\delta_X}(u)= \frac{1}{2N} \sum_{n=1}^N  <A(x^n) \nabla_k u(x^n), \nabla_k u(x^n)> + \frac{1}{N} \sum_{n=1}^N   (u f)(x^n). 
\ee
We can characterize the minimum of this functional, i.e. $u = \arg \min_{v \in \mathcal{H}_{X}}J_{\delta_X}(v)$, leading to the discrete system
\be
  \Big[\nabla_k \cdot ( A \nabla_k u ) \Big](X) = f(X)
\ee
Section~\ref{representer_theorem} proposed an example of application of the representer theorem, to characterize the projection operator $\mathcal{P}_k(\cdot)$ through a variational problem. We now study two other fundamental variational problems for applications, while Chapter~\ref{application-to-partial-differential-equations} below provides also some examples.

\paragraph{Helmholtz-Hodge decomposition}\label{the-helmholtz-hodge-decomposition}

In many areas of fluid mechanics, for example to analyze turbulence problems, study flow past obstacles, and develop numerical methods for simulating fluid flows, one uses the so-called Helmholtz-Hodge decomposition. One important component of this decomposition is the Leray operator, which can be used to orthogonally decompose any field and is key to understanding some important structures of fluid flows.

Consider the following \textit{Helmholtz-Hodge} decomposition of a vector field $u$ into a potential scalar field $h$ and a divergence-free vector $\zeta$.
\be\label{HH}
u = \nabla h + \zeta, \quad \nabla \cdot \zeta = 0.
\ee
Provided $u$ is smooth enough, this decomposition is unique. The Helmholtz-Hodge decomposition is an essential technique in functional analysis and is used in many applications, ranging from optimal transport to fluid or electromagnetic field motions.

Consider a probability measure $d\mu = \mu(\cdot) dx$ and the following functional associated to the Helmholtz-Hodge decomposition:
\be
    J_\mu(h) = \frac{1}{2} \int | \nabla h - u|^2 d \mu, 
\ee
which is defined for the weighted Sobolev space $h \in H^1_\mu = \{h : \int |\nabla h|_2^2 d \mu < +\infty\}$. Using the Gateaux derivatives to characterize a minimum of this functional shows that a minimum is characterized in a weak sense by the equation
\be
 \nabla\cdot \nabla (h\mu) = \Delta (h\mu) = \nabla \cdot (u\mu), \quad \zeta = (u - \nabla h)\mu,
\ee
which satisfies $\nabla \cdot \zeta = 0$.

Consider $X,Y$ two sets of distinct points, the subspace $\mathcal{H}_{k,Y} \subset \mathcal{H}_k$, and instead of considering the integral $J(h)$, approximate it as a point-wise sum on the set $X$. We solve the discrete problem:
\be\label{HHRKHS}
h = \arg \inf_{v \in \mathcal{H}_{k,Y}} J_X(v), \quad J_X(v) =  \|(\nabla_k v - u)(X)\|_{\ell^2}^2, \quad \zeta = u - \nabla h. 
\ee
The weak formulation of this variational problem can be derived considering that  $J(h) \le J(h+\varphi)$, for any $\varphi \in \mathcal{H}_{k,Y}$, leading to the equation $<(\nabla_k h - u)(X), \nabla_k \varphi(X)>=0$. Using the operator $\nabla_k(\cdot)$ and its transpose $\nabla_k(\cdot)^T$, we get the following equation
\be
 (\nabla_k^T \nabla_k)(X) h(X) = (\nabla_k^T u)(X). 
\ee
In the left-hand side, one recognizes the Laplace-Beltrami operator, providing a formal solution
\be\label{Delta-inv}
 \overline{h}(X) = \Delta_k^{-1} (\nabla_k \cdot u)(X), \quad \Delta_k = \nabla_k^T \cdot \nabla_k,
\ee
where the inverse of the Laplace-Beltrami operator $\Delta_k^{-1}$ is obtained considering a least-square inversion of the matrix $\Delta_k$.

The operator $\mathcal{L}_k = I_d - \nabla_k \Delta_k^{-1} \nabla_k \cdot$ computes the divergence-free component $\zeta$ in \eqref{HHRKHS}, and is called the \textit{Leray} projection.

\paragraph{Trace operators and boundary conditions}

Many problems in mathematical physics require functions $u \in \mathcal{H}$, a Hilbert space of continuous functions, to have prescribed values or relations at a physical boundary $\Gamma$. We consider the example of a Dirichlet condition $u(\Gamma) = \varphi$, prescribing values of a function at the boundary $\Gamma$, but more complex boundary conditions as Neumann $(\nabla u)(\Gamma) = \psi$ or mixed can be treated in a similar manner. This amounts to consider formally the following subspace of $\mathcal{H}$:
\be\label{BC}
    u \in \mathcal{H}^{\varphi(\Gamma)} = \Big\{ u \in \mathcal{H} : u(\Gamma) = \varphi(\Gamma) \Big\}
\ee
In order to model numerically such subspaces, consider a set of points $X$ and a positive-definite kernel $k$ generating $\mathcal{H}_k$, the RKHS $\mathcal{H}_{k,X}$, a set of points $Z$ modeling the boundary $\Gamma$, and denote $\mathcal{H}_{k,X}^{\varphi}$, a closed subspace of  $\mathcal{H}_{k,X}$ modeling $\mathcal{H}_k^{\varphi}$ as follows:
\be\label{KBC}
    v \in \mathcal{H}_{k,X}^{\varphi(Z)} = \Big\{ u_k \in \mathcal{H}_{k,X} : u_k(Z) = \varphi(Z) \Big\},
\ee
In order to characterize this discrete subspace, meaning computing the decomposition $v = u_{k,\varphi} + u_{k,\varphi}^\perp$, satisfying $u_{k,\varphi}(Z) = \varphi(Z)$, let us introduce the following minimization problem
\be
 u_{k,\varphi} = \arg \inf_{v \in \mathcal{H}_{k,X}} J(v), \quad J(v) = \|(u_k - v)(X)\|_{\ell^2}^2  + \|(u_k - v)(Z)\|_{\ell^2}^2 ,
\ee
where values of any functions $u_k \in \mathcal{H}_{k,X}$ are extrapolated on the boundary points $Z$ according to $u_k(Z) = \mathcal{P}_{k,X}(Z)u_k(X)$. Solving this minimization problem, the decomposition is given by
\be
    u_{k,\varphi(X)} = \Big(I_d- \mathcal{P}_{k,X}(Z)^T\mathcal{P}_{k,X}(Z) \Big)^{-1} \Big( u_k(X) - \mathcal{P}_{k,X}(Z)^T \varphi(Z) \Big).
\ee
This formula defines thus the projection of any function $u_k \in \mathcal{H}_{k,X}$ into $\mathcal{H}_{k,X}^{\varphi(Z)}$.


\section{A zoo of kernel operators}

\subsection{Interpolations and extrapolation operators}
\label{interpolations-and-extrapolation-operators}

\paragraph{Coefficient operator}
\label{coefficient-operator}

We now provide several numerical illustrations of discrete RKHS operators. Recall that the projection operator \(\mathcal{P}_k(\cdot)\) \eqref{OP} maps any continuous function onto a basis of functions. With the notation in the previous chapter, given a kernel \(k\) and two sets \((X,Y)\), let us consider the fitting procedure
\be \label{coefy}
 f_k(\cdot) = K(\cdot,Y) \theta, \quad \theta = K(X,Y)^{-1}f(X) \in \mathbb{R}^{N_x, D_f},
\ee
where, \(\theta\) represents the coefficients of the decomposition of a function \(f\). In other words, \(f\) can be written as a linear combination of the basis functions \(K(\cdot,y^n)\), where \(n\) ranges from \(1\) to \(N_y\). We refer to the matrix $K(X,Y)^{-1}$ as the coefficient operator.

\paragraph{Partition of unity}\label{partition-of-unity}

The notion of partition of unity is, both, a standard and a very useful concept. Consider the projection operator $\mathcal{P}_{k,Y}(\cdot,X)$, which is a vector having  $N_x$ components. This operator can also be viewed as $N_x$ real valued functions $\psi^n(x,X)$ defined as
\be \label{PU}
 \psi^N_x(\cdot) = \mathcal{P}_{k,Y}(\cdot,X)^n = K(\cdot,Y)K(X,Y)^{-1} 1^n(X), \quad 1_n(X) = \delta_n(m),
\ee
where \(\delta_{n,m}\) denotes the Kronecker delta symbol, which we refer to as the partition of unity. The reproducing property can be written on the dataset as
\be 
 \psi^n_{k,X}(x^m)= \delta_{n}(m), 
\ee
Figure~\ref{fig:312} illustrates this notion with an example of four partition functions.

\begin{figure}
\centering
\includegraphics[width=0.7\textwidth, keepaspectratio]{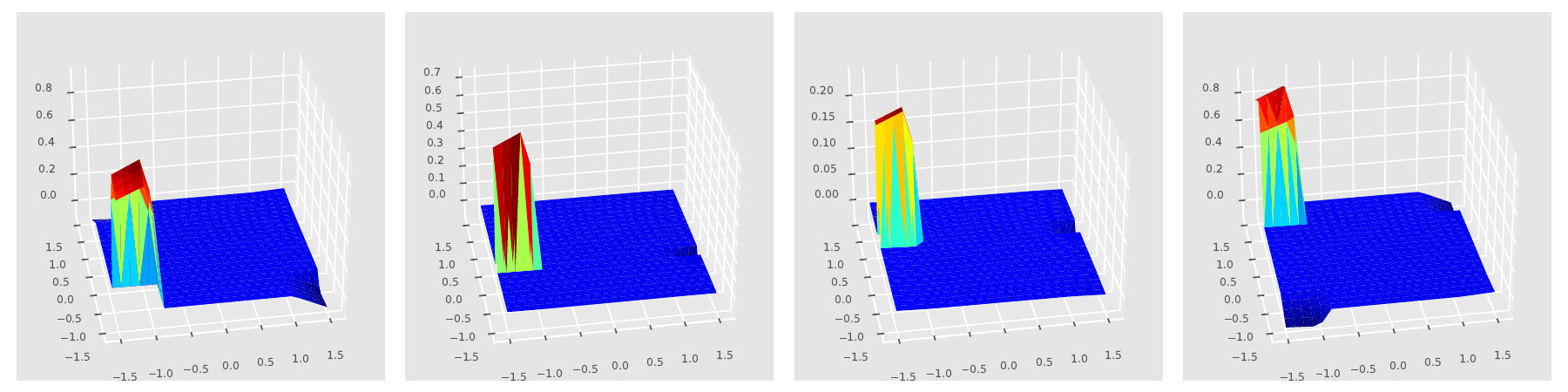}
\caption{\label{fig:312}Four partition of unity functions}
\end{figure}


\subsection{Discrete differential operators}
\label{discrete-differential-operators}

\paragraph{Gradient operator}\label{gradient-operator}

To illustrate numerically the discrete differential operators, we will be using the two-dimensional function $f$ defined at \eqref{2D}, as well as the three sets $X,Y,Z$ defined in Section~\ref{basic-numerical-examples}.
In Figure~\ref{fig:nabla1}, we begin with a gradient computation of the vector-valued function \(f\), using the expression
\be
( \nabla_k f)(Z) = (\nabla_{k})(Z) f(X) \in \mathbb{R}^{N_z, D_f}. 
\ee
This figure plots a comparison between the exact gradient of the original function and their corresponding values computed using the operator \eqref{nablak}, thus for the two dimensions. The left-hand plot corresponds to the original function, while the right-hand plot shows the computed values.

\begin{figure}
\centering
\includegraphics[width=0.7\textwidth, keepaspectratio]{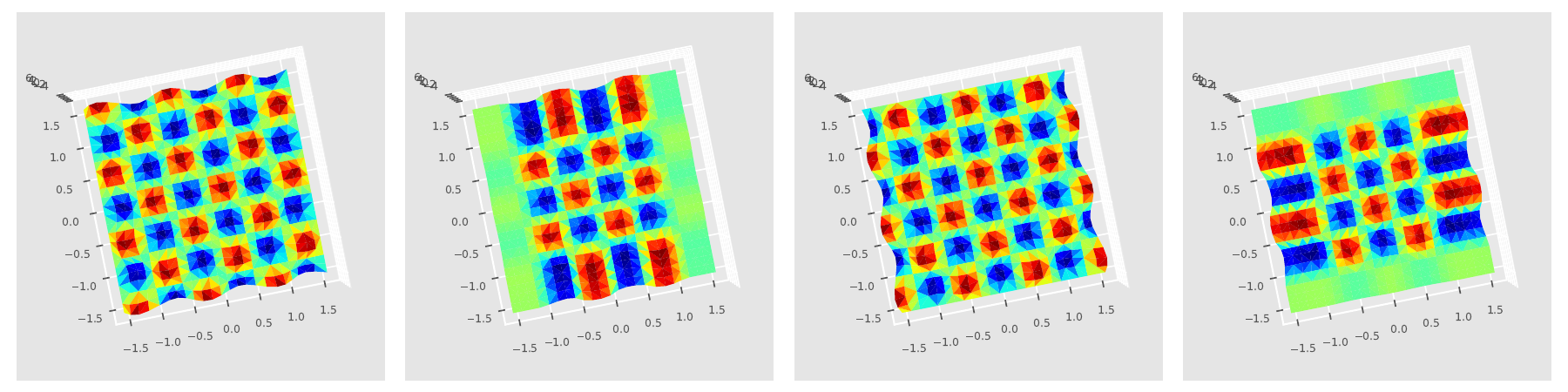}
\caption{\label{fig:unnamed-chunk-28}\label{fig:nabla1}The first two graphs correspond to the first dimension (original on the left-hand, computed on the right-hand). The next two graphs correspond to the second dimension (original on the left-hand, computed on the right-hand).}
\end{figure}

\paragraph{Divergence and Laplace-Beltrami operator}\label{divergence-operator}

We illustrate a divergence computation of a vector-valued function \(g(Z)\), coming from the expression~\eqref{div}
\be
 (\nabla_k \cdot g)(X) = K(Y,X)^{-1}(\nabla K)(Y, Z) g(Z).
\ee
To test the consistency of our operators, we consider $g(Z) = (\nabla_k f)(Z)$ in the previous expression, and thus should be equal to $(\nabla_k \cdot g)(X) = \Delta_k f(X)$, $\Delta_k$ being the Laplace-Beltrami operator of Section~\ref{laplace-beltrami-operator}. Figure~\ref{fig:nablaTnabla} compares thus this expression to the Laplace operator (see the discussion below).

\begin{figure}
\centering
\includegraphics[width=0.7\textwidth, keepaspectratio]{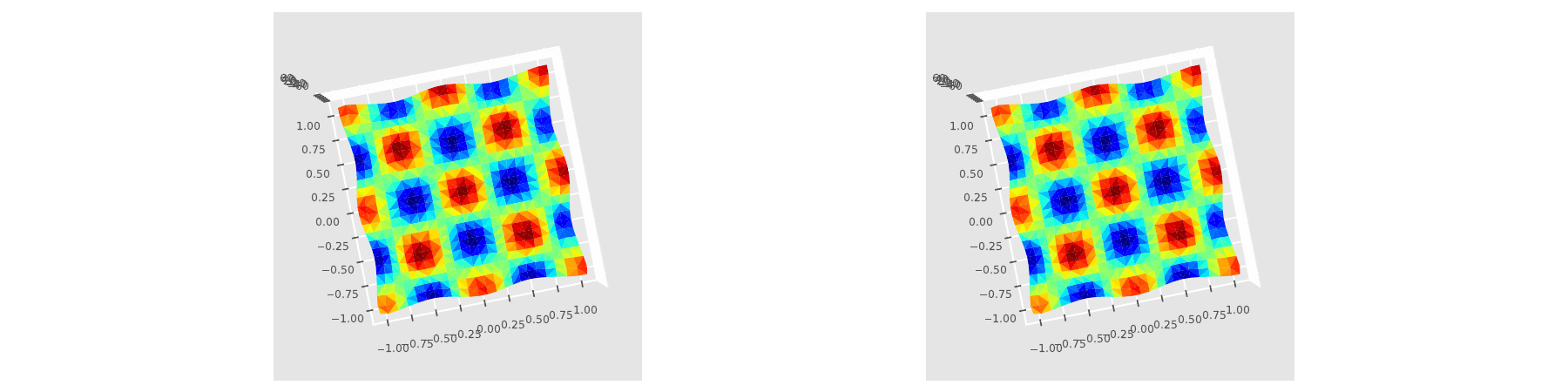}
\caption{\label{fig:unnamed-chunk-29}\label{fig:nablaTnabla}Comparison of the divergence of a gradient with the Laplace operator}
\end{figure}

\paragraph{Leray operators}
\label{leray-orthogonal-operator}

The discrete modeling of the Leray operator is described, as well as its orthogonal projector, in Section~\ref{the-helmholtz-hodge-decomposition} (see the Helmholtz-Hodge decomposition), as
\be
  \mathcal{L}_{k,X}^\perp = \nabla_{k}\Delta_{k}^{-1} \nabla_{k}^T, \quad \mathcal{L}_{k,X} = I_d - \mathcal{L}_{k}^\perp
\ee
This operator acts on any vector field \(f(Z) \in \mathbb{R}^{ D \times N_z, D_f}\), and produces two vector fields that are orthogonal in the following, discrete, sense:
\be
  f(z)  = L_{k}^\perp f(Z) + L_{k}f(Z), \quad \nabla_{k}^T L_{k}f(Z) = 0.
\ee
The construction of this discrete Leray decomposition enjoys the same orthogonality properties as the continuous ones of the original Helmholtz-Hodge decomposition.

In particular, the Leray decomposition of a vector field having form $f(\cdot)=\nabla \varphi(\cdot)$, where $\varphi$ is a scalar function, should be trivial. In Figures~\ref{fig:Leray}--\ref{fig:LerayT}, we test this idea, comparing the action of the Leray decomposition on a vector field with the function \((\nabla f)(Z)\), $f$ defined in~\eqref{2D}, showing that this conclusion might be tempered due to extrapolation errors.

\begin{figure}
\centering
\includegraphics[width=0.7\textwidth, keepaspectratio]{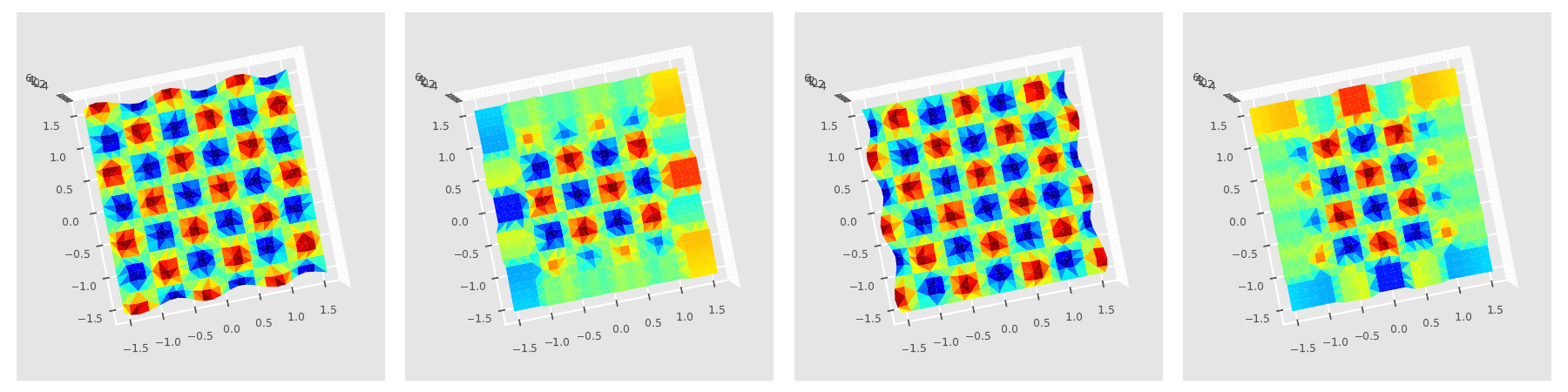}
\caption{\label{fig:unnamed-chunk-35}\label{fig:LerayT}Comparing $\nabla f$ and its orthogonal Leray projection on each two directions}
\end{figure}

\begin{figure}
\centering
\includegraphics[width=0.7\textwidth, keepaspectratio]{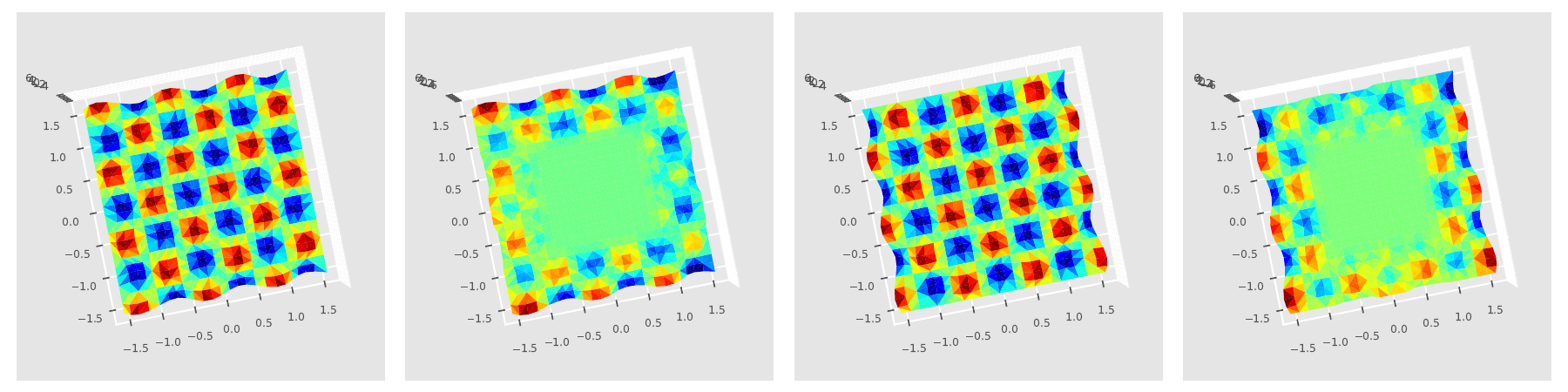}
\caption{\label{fig:unnamed-chunk-36}\label{fig:Leray}Comparing $\nabla f$ and its Leray projection in each direction}
\end{figure}


\subsection{Discrete integral operators}
\label{inverse-laplace-operator}
\label{integral-operator---inverse-gradient-operator}

\paragraph{Inverse Laplace operator}

The inverse Laplace operator can be formally defined as the pseudo-inverse of the Laplacian operator \(\Delta_{k}(X) \in \mathbb{R}^{N_x, N_x}\), defined in \eqref{Delta}. This operator corresponds to a solution to the following functional, where $f$ is any continuous function:
\be
u = \arg \inf_{v \in \mathcal{H}_k} J(v), \quad J(v) = \int |\nabla v|^2 -\int vf. 
\ee
The discrete, RKHS representation of this minimization problem can be expressed as
\be
u = \arg \inf_{v \in \mathcal{H}_{k,Y}} J(v), \quad J(v) = \frac{1}{2}\| (\nabla_k v)(X) \|_{\ell^2}^2 - <v(X),f(X)>.
\ee
Thus it can be computed formally as follows:
\be\label{Deltainvk}
  u = \Delta_{k}^{-1} f,  
\ee
$\Delta_{k}^{-1} \in \mathbb{R}^{ N_x , N_x}$ being the pseudo-inverse of the matrix $\Delta_{k}^{-1}$. Figure~\ref{fig:Delta1} compares \(f(X)\) with \(\Delta_{k}^{-1}\Delta_{k}f(X)\). This latter operator is a projection operator (hence is stable).

\begin{figure}
\centering
\includegraphics[width=0.7\textwidth, keepaspectratio]{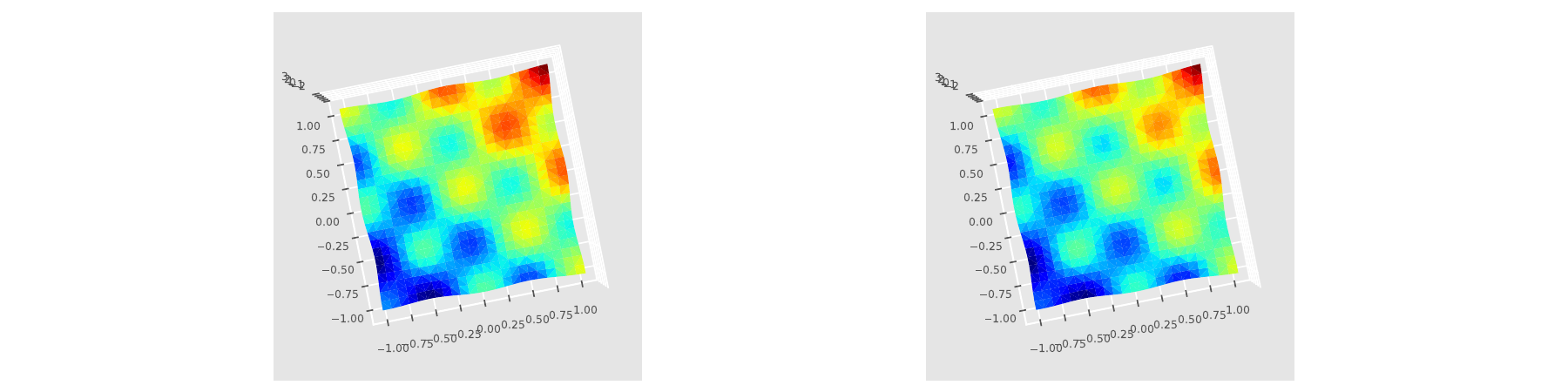}
\caption{\label{fig:unnamed-chunk-30}\label{fig:Delta1}Comparison of the original function with the product of Laplace and its inverse}
\end{figure}

In Figure~\ref{fig:Delta2}, we compute the operator \(\Delta_{k}\Delta_{k}^{-1}f(X)\) to check that the pseudo-inverse commutes, i.e., applying the Laplacian operator and its pseudo-inverse in any order produces the same result. This property should hold for strictly positive-definite kernels, and we perform it to check the consistency of the framework.

\begin{figure}
\centering
\includegraphics[width=0.7\textwidth, keepaspectratio]{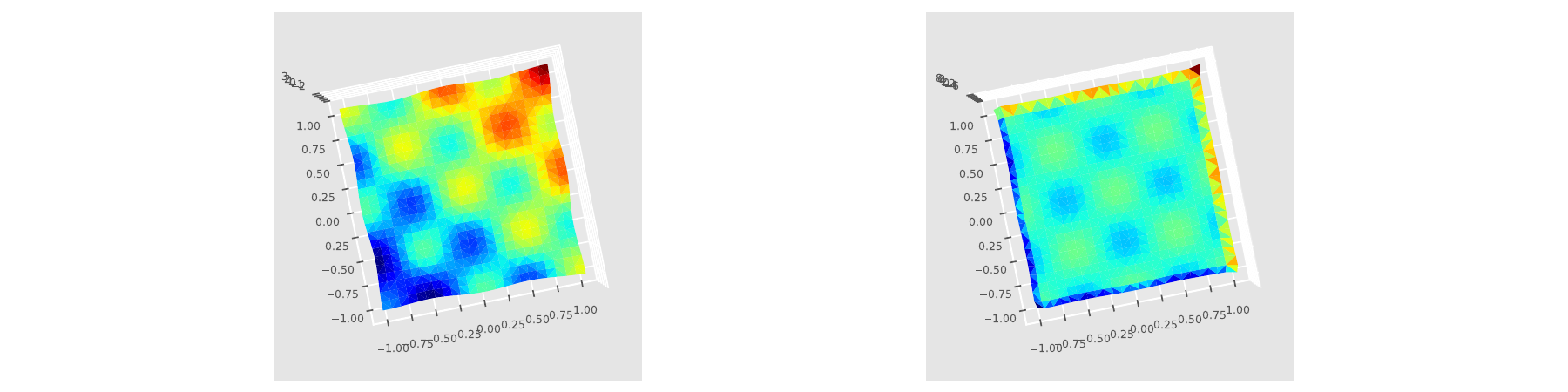}
\caption{\label{fig:unnamed-chunk-31}\label{fig:Delta2}Comparison of the original function with the product of the inverse of the Laplace operator and the Laplace operator}
\end{figure}

\paragraph{Integral operator - inverse gradient operator}

The operator \(\nabla^{-1}_k\) is defined as the integral-type operator
\be \label{nablainv}
  \nabla^{-1}_{k} =  \Delta_{k}^{-1}\nabla_{k}^T \in \mathbb{R}^{ N_x, DN_z}. 
\ee
This operator, already introduced considering the minimization problem \eqref{HHRKHS}, acts on any field of vectors \(v(Z) \in \mathbb{R}^{ D \times N_z, D_{v_z}}\) and produces a matrix
\be
  \nabla^{-1}_{k}v(z) \in \in \mathbb{R}^{ D \times N_z, D_{v_z}}.
\ee
In Figure~\ref{fig:NablainvNabla} we test whether
\be
(\nabla_k )^{-1}(\nabla_k) f(X)
\ee
coincides or at least is a good approximation of \(f(X)\), which is a property of strictly positive-definite kernels $k$. 

\begin{figure}[H]
\centering
\includegraphics[width=0.7\textwidth, keepaspectratio]{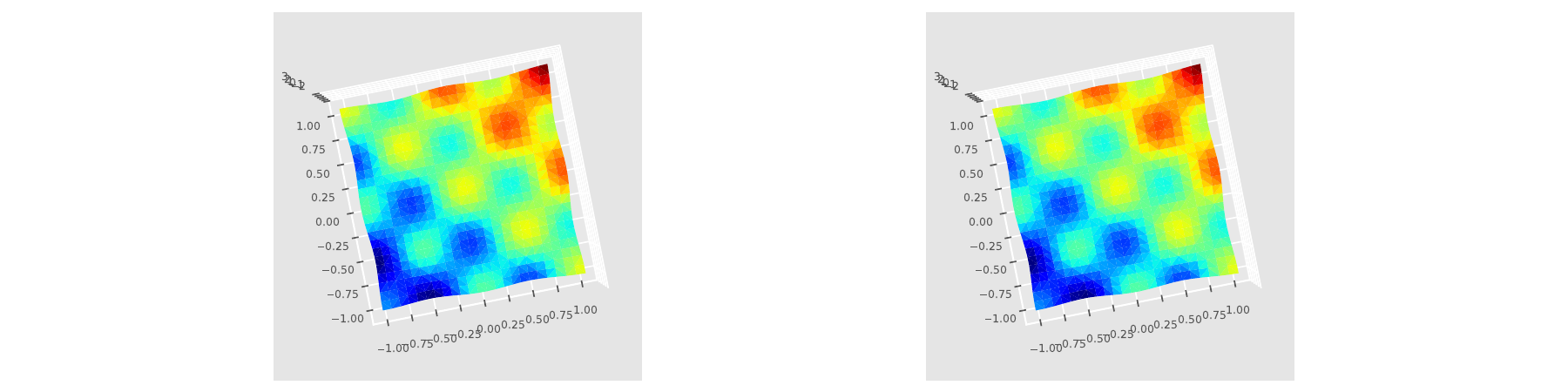}
\caption{\label{fig:unnamed-chunk-32}\label{fig:NablainvNabla}Comparison of the original function with the product of the gradient operator and its inverse}
\end{figure}

\begin{figure}[H]
\centering
\includegraphics[width=0.7\textwidth, keepaspectratio]{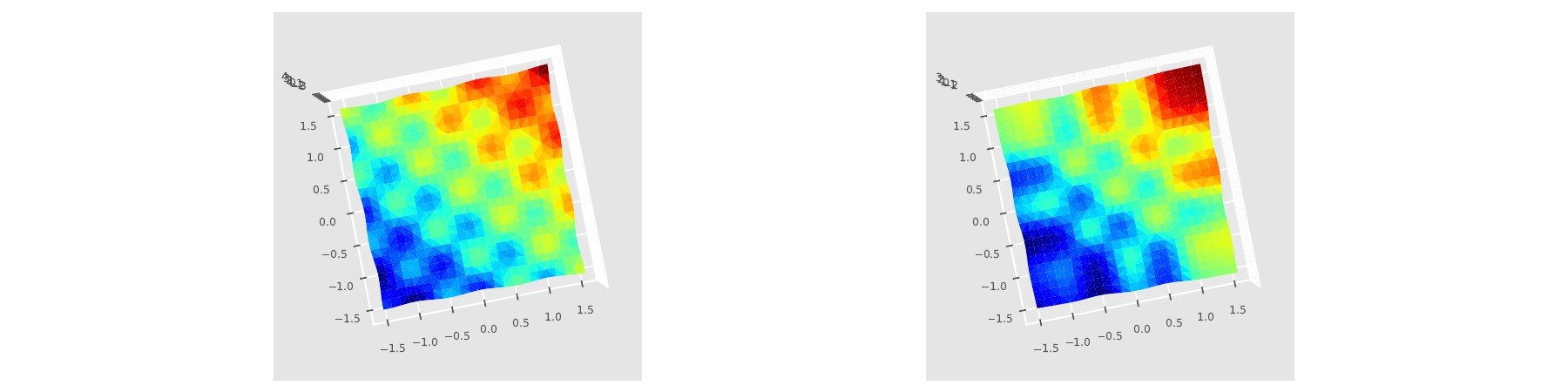}
\caption{\label{fig:unnamed-chunk-33}\label{fig:NablainvNabla2}Comparison of the original function with the product of the inverse of the gradient operator and the gradient operator}
\end{figure}

\paragraph{Integral operator - inverse divergence operator}\label{integral-operator---inverse-divergence-operator}

The following operator \((\nabla_k^T)^{-1}\) is another integral-type operator of interest. We define it as the pseudo-inverse of the divergence operator \(\nabla^T\) as follows:
\be
  (\nabla^T_{k})^{-1}= \nabla_{k}\Delta_{k}^{-1}.
\ee
It corresponds to the following continuous, minimization problem
\be
h = \arg \inf_{v \in (\mathcal{H}_k)^D} J(v), \quad J(v) = \int |\nabla \cdot v - u|^2
\ee
According to Theorem~\ref{repth} (the representer theorem), this problem can be discretized as
\be
h = \arg \inf_{v \in \mathcal{H}_{k,Y}} J_X(v), \quad J_X(v) =  \|(\nabla_k^T v - u)(X)\|_{\ell^2}^2.
\ee
To check the consistency of this operator,  in Figure~\ref{fig:unnamed-chunk-34} we compute \(\nabla_{k}^T(\nabla^T_{k})^{-1}= \Delta_{k}\Delta_{k}^{-1}\). Thus, the following computation should give comparable results as those obtained in our study of the inverse Laplace operator in Section~\ref{inverse-laplace-operator}.

\begin{figure}[H]
\centering
\includegraphics[width=0.7\textwidth, keepaspectratio]{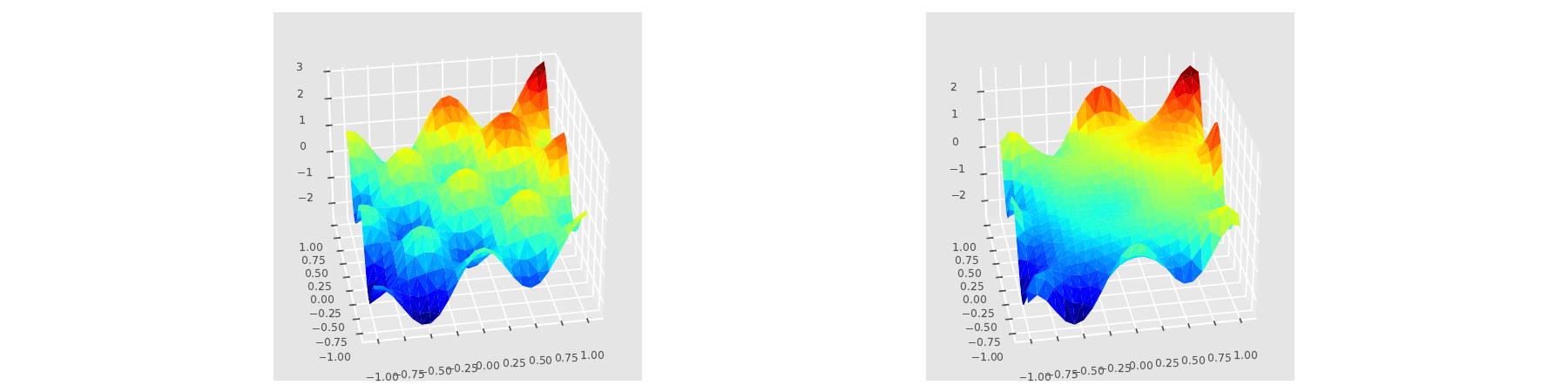}
\caption{\label{fig:unnamed-chunk-34}\label{fig:nablaTinv}Comparison of the product of the divergence operator and its inverse with the product of Laplace operator and its inverse}
\end{figure}


\clearpage 

\chapter{Clustering strategies}\label{clustering-strategies}

\section{Introduction}\label{introduction}

\paragraph{Main purpose.}

In this chapter, we examine various clustering techniques from the perspective of reproducing kernel Hilbert spaces. Clustering is a fundamental task in exploratory data analysis, aiming to group similar data points based on a given similarity or cost criterion. In the context of kernel methods, clustering also serves a critical computational role, as it can substantially reduce the computational complexity associated with kernel operations. Furthermore, clustering methods are closely connected to optimal transport theory, a powerful framework with numerous applications in statistical kernel methods, which we discuss in Section~\ref{optimal-transport-and-statistical-kernel-methods}.

We now introduce the main concepts that underlie our treatment of clustering. In general, clustering methods rely on a notion of \textit{distance}, or more broadly, a \textit{cost} function $d(X,Y)$, where $X,Y$ represent finite point distributions in $ \mathbb{R}^{N_x,D_x}$ and $\mathbb{R}^{N_y,D_y}$, respectively, with $N_y < N_x$. This distance allows us to define the concept of \textit{centroids} formally as
\be
\label{centroids}
\overline{Y} = \mathop{\arg\inf}_{Y \in \mathbb{R}^{N_y,D_y}} d\big(Y,X\big).
\ee
Given an arbitrary distribution $X$, distances define also \textit{assignment maps} $\sigma$, that are functions taking any features as inputs, and which outputs a label 
$\sigma_X(\cdot) : \mathbb{R}^{D_x} \to \{ 0, \ldots, N_x \}$. Assuming that the distance $d$ is convex, this function realizes a partition of $\mathbb{R}^{D_x}$ into $N_x$ connected components as follows: 
\be
\label{assignment}
\sigma_X(\cdot) = \mathop{\arg\inf}_{n=1,\ldots,N_x} d\big(\cdot,x^n\big).
\ee
This partition allows for \textit{semi-supervised methods}, that are predictive methods commonly based on a nearest-neighbour methods, an example of which being given by the formula
\be
\label{semisupervised}
f_{d}(\cdot) = f(x^{\sigma_Y(\cdot)}), 
\ee
where $X,f(X)$ is a labeled dataset, $Y$ is the set of centroids in \eqref{centroids} and $\sigma_Y$ associated assignment map defined in~\eqref{assignment}. 

\paragraph{Numerical strategy.} 

We now discuss our strategy to compute the centroids set \eqref{centroids}, adapted to different distances $d(\cdot,\cdot)$.  While gradient descent methods are straightforward to implement for solving~\eqref{centroids} (see Section~\ref{explicit-descent-algorithms}), they face several limitations: high computational costs in large-scale or high-dimensional settings, and a tendency to get trapped in local minima due to the nonconvex nature of the optimization landscape.

To address these issues, we first consider a discrete relaxation of the centroid problem 
\be\label{Ycentroids}
\overline{Y} = \mathop{\arg\inf}_{Y \subset X} d\big(Y,X\big),
\ee
which restricts the set of centroids $Y$ to subsets of the dataset $X$. This formulation admits fast combinatorial algorithms, which can yield good initializations for subsequent continuous optimization via gradient descent. However, in many practical scenarios, the optimal centroids are indeed contained within the original dataset.

Section \ref{general-purpose-algorithms}, below, introduces a family of general-purpose algorithms that implement this strategy. These algorithms are highly versatile, and offer robustness and ease of use across a wide range of problems. Nevertheless, they may sacrifice some performance in domain-specific contexts, where more specialized algorithms can be advantageous.


\section{General purpose algorithms}
\label{general-purpose-algorithms} 

\subsection{Greedy search algorithm}
\label{greedy-search-algorithms}

The following algorithm provides us with a family of efficient and
versatile algorithms for the approximation of the following (quite
broad) class of problems: \be
\label{GS}
\inf_{Y \subset X} D\big(Y,X\big), \quad D(Y,X) = \sum_n d(Y,x^n),
\ee in which \(D(Y,X)\) denotes a (user defined) distance
measure between sets. We consider here a greedy search algorithm, approximating \eqref{GS} recursively as follows:
\be
\label{GS=2}
Y^{n+1} = Y^n \cup \arg \sup_{x \in X} d\big(Y^n,x\big). 
\ee 
We will present some useful examples of this algorithm in
the context of kernel interpolation. 

This approach leads us to Algorithm~\ref{cl0} as stated next, where the
parameter \(M\) (which may be taken to be \(1\) by default) is introduced as a simple
optimization mechanism in order to trade accuracy versus time. Other
strategies can be implemented, more adapted to the problem at hand. Such
an algorithm usually relies on optimization techniques within the main loop
for faster evaluation of \(d\big(Y^n, \cdot\big)\), using
pre-computations or \(d\big(Y^{n-1}, \cdot\big)\).

\begin{algorithm}\label{cl0}
\begin{algorithmic}
\REQUIRE{a training set $X$, a distance measure $d(\cdot,\cdot)$, two integers: $1 \le N_y \le N_x$, where $N_y$ denotes the number of clusters, $M$ is an optional \textit{batch} number (taken to $1$ by default), and $Y^0 \subset X$ is  an optional set of initial points.}
    \ENSURE{A set of indices $\sigma : [1, \ldots, N_y]  \to [1, \ldots, N_x]$, defining $X\circ\sigma \subset X$ of size $N_y$ as approximate clusters.}
    \STATE{for $n=0,\ldots,N_y/M$ 
             find a new numbering, say $X^n$, according to the decreasing order of $d\big(Y^n, x^p\big)$, $p=1,\ldots,N_x$ \\
              $Y^{n+1}=Y^{n}\cup X[1,\ldots,M]$.}
\end{algorithmic}
\end{algorithm}


\subsection{Permutation algorithm}
\label{linear-sum-assignment-problems}

\paragraph{Purpose.}

Next, let us focus on the Linear Sum Assignment Problem (LSAP), a foundational problem in combinatorial optimization with numerous applications across engineering, statistics, and computer science. LSAP has been extensively studied and is well-documented in both academic and applied literature\footnote{\url{https://en.wikipedia.org/wiki/Assignment_problem}}. We also present a generalization of LSAP, for which we design a class of efficient descent-type algorithm that are amenable to parallelization, thereby enhancing performance on large-scale problems. 

\paragraph{Linear sum assignment problem.}

Let $c(x,y)$ be a given distance or cost function between two points, define the matrix $C(X,Y)$ having entries $c(x^i,y^j)$, and consider the following distance between set
\be
  d(X,Y) = \text{Tr}\big(C(X,Y) \big).
\ee
Solving \eqref{Ycentroids} with such a distance amounts to find a permutation of the indices, characterized as follows: 
\be
\label{LSAP}
  \overline{\sigma} = \mathop{\arg\inf}_{\sigma \in \Sigma} \sum_{i=1}^M c(i,\sigma^i), 
\ee where \(\Sigma\) is the set of all injective re-numberings
\(\sigma:[1,\ldots,M]  \to [1,\ldots, N]\) with \(M \le N\). Of
course, when \(M=N\), this is nothing but the set of permutations. These
algorithms \textit{put together} two distributions according to a
similarity criterion given by a general, rectangular \textit{cost} or
\textit{affinity} matrix \(C \in \mathbb{R}^{M,N}\).

\paragraph{Discrete descent algorithm.}

We now consider a generalization of LSAP, where the cost functional $C(\sigma)$ is not necessarily linear. We define the problem as
\be\label{GLSAP}
\overline{\sigma} = \mathop{\arg\inf}_{\sigma \in \Sigma} C(\sigma), 
\ee 
where \(\Sigma\) denotes the set of all possible
permutations, and \(C(\sigma)\) is a general cost functional. A special
case of this problem is indeed the LSAP problem above, which corresponds to a
linear cost function \(C(\sigma) = \sum_i c(i, \sigma^i)\). However, we do consider other forms of cost functionals.

For problems where the permutation gain, defined as
\(s(i,j,\sigma) = C(\sigma)-C(\sigma_{i,j})\) can be efficiently
computed, we use a simple descent algorithm. Here \(\sigma_{i,j}\)
represents the permutation obtained by swapping the indices \(\sigma^i\)
and \(\sigma^j\). The LSAP problem \eqref{LSAP} corresponds to a permutation gain
given by
\(s(i,j,\sigma) = c(i,\sigma^i)+c(j,\sigma^j)-c(i,\sigma^j)- c(j,\sigma^i)\).

This approach leverages the fact that any permutation
\(\sigma\) can be decomposed into a sequence of elementary two-element
swaps. 
In its simplest form, this algorithm is summarized in Algorithm~\ref{PBA}. There exist more performing algorithms than this
discrete descent approach but adapted to specific situations. 

\begin{algorithm}
\label{PBA}
\begin{algorithmic}[1]
\REQUIRE A permutation-gain function $s(i,j,\sigma)$, where $\sigma:[1,\ldots,N_J]  \to [1,\ldots,N_I]$ is any injective mapping (that is, a permutation if $N_I=N_J$).
\vskip.15cm
\ENSURE An injective mapping $\sigma:[1,\ldots,N_J]  \to [1,\ldots,N_I]$ achieving a local minima to \eqref{GLSAP}.
\WHILE{FLAG = True} 
\STATE FLAG $\leftarrow$ False
\FOR{$i=1 \ldots, N_I$,$j=1 \ldots, N_J$} 
\IF{$s(i,j,\sigma) < 0$}
\STATE $\sigma \leftarrow \sigma^{i,j}$
\STATE FLAG $\leftarrow$ True
\ENDIF
\ENDFOR
\ENDWHILE
\end{algorithmic}
\end{algorithm}

The \texttt{For} loop in this algorithm can be further adapted to
specific scenarios or strategies, for instance a simple adaptation to
symmetric permutation gain function \(s(i,j,\sigma)\) (as met with the
LSAP problem) is given as follows:
\texttt{for\ i\ in\ {[}1,\ N{]},\ for\ j\ \textgreater{}\ i} in order to
minimize unnecessary computation.

While these algorithms generally produce sub-optimal solutions for
nonconvex problems, they are robust and tend to converge in finite
time, i.e.~within a few iterations of the main loop, and careful
choice of the initial permutation \(\sigma\) can escape sub-optimality.
They are especially useful as auxiliary methods in place of more complex
global optimization techniques or for providing an initial solution to a
problem. Additionally, they are advantageous in finding a local minimum
that remains \emph{close} to the original ordering, preserving the
inherent structure of the input data.

However, these algorithms have some limitations. Depending on the
formulation of the permutation gain function \(s(i,j,\sigma)\), it is
possible to parallelize the \texttt{For} loop. Nonetheless,
parallelization often requires careful management of concurrent
read/write access to memory, which can complicate implementation and
alter performances. Furthermore, theoretical bounds on the algorithm
complexity are typically not available, making performances prediction
difficult in all cases.

\subsection{Explicit descent algorithm}
\label{explicit-descent-algorithms}

We now present a generalized gradient-based algorithm, specifically
designed for the minimization of functionals of the form
\(\inf_X J(X)\), where the gradient \(\nabla J(X)\) is locally convex
and is explicitly known. In this scenario, we can apply the simplest
form of a gradient descent scheme, often referred to as an Euler-type
method. In its continuous form, it is written as
\(\frac{d}{dt} X(t) = -\nabla J(X(t))\), \(X(0)=X^0\) and, in its
numerical formn, it takes the form \be
\label{DA}
  X^{n+1} = X^{n} - \lambda^n\nabla J(X^{n}).
\ee The term \(\lambda^n\) is known as the
\textit{learning rate}. In this situation, as the gradient of the
functional is explicitly given, we can compute sharp bounds over
\(\lambda^n\), allowing to apply root-finding methods, such as the Brent
algorithm\footnote{For instance \href{https://en.wikipedia.org/wiki/Brent's_method}{Brent's method on Wikipedia}},
and efficiently locate the minimum while avoiding \textsl{instability}
issues often met with Euler schemes.

\begin{algorithm}
\label{DA=algo}
\begin{algorithmic}[1]
\REQUIRE A function $J(\cdot)$, its gradient $\nabla J(\cdot)$, a first iterate $X^0$, tolerance $\epsilon >0$ or number of maximum iterations $N$.
\vskip.15cm
\ENSURE A solution $X$ achieving a local minimum of $J(X)$.

\WHILE{$\|\nabla J(X^n)\| > \epsilon$ or $n < N$} 
\STATE compute $\lambda^{n+1} = \mathop{\arg\inf}_{\lambda} \|\nabla J(X^\lambda)\|$, where $X^{\lambda}= X^{n} - \lambda \nabla J(X^{n})$ with a root-finding algorithm.
\STATE $ X^{n+1} = X^{n} - \lambda^{n+1}\nabla J(X^{n})$
\ENDWHILE
\end{algorithmic}
\end{algorithm}

\subsection{Illustration with the LSAP problem}

We now illustrate the LSAP problem through a concrete numerical example. Consider a cost matrix \(A=a(n,m) \in \mathbb{R}^{N,M}\), with randomly generated entries shown in Table~\ref{tab:401}.  The goal is to compute a matching \(\sigma \) minimizing the total cost, given by 
\be
  \sigma = \mathop{\arg\inf}_{\sigma \in \Sigma} \sum_{n=1}^M a(n,\sigma(m)), 
\ee
where \(\Sigma\) is the set of all matchings, that is, the set of injective couplings $\sigma :[1, \ldots, M]  \to [1, \ldots, N]$, which is the set of permutations if $N=M$. Let us start with a simple illustration in the case $N=M$, considering the following four by four random matrix in Table~\ref{tab:401}. The total cost before permutation is simply \(Tr(A)\), given in Table~\ref{tab:unnamed-chunk-43}.  After solving the LSAP, we obtain the optimal row permutation \(\sigma\) in Table~\ref{tab:unnamed-chunk-45}Permutation}. 

\begin{table}[htbp]
\caption{\label{tab:401}A 4x4 random matrix}
\centering
\begin{tabular}[t]{r|r|r|r}
\hline
0.2617057 & 0.2469788 & 0.9062546 & 0.2495462\\
\hline
0.2719497 & 0.7593983 & 0.4497398 & 0.7767106\\
\hline
0.0653662 & 0.4875712 & 0.0336136 & 0.0626532\\
\hline
0.9064375 & 0.1392454 & 0.5324207 & 0.4110956\\
\hline
\end{tabular}
\end{table}
\begin{longtable}[t]{r}
\caption{\label{tab:unnamed-chunk-43}Total cost before permutation}\\
\hline
1.465813\\
\hline
\end{longtable}
\begin{longtable}[t]{r|r|r|r}
\caption{\label{tab:unnamed-chunk-45}Permutation}\\
\hline
1 & 3 & 2 & 0\\
\hline
\end{longtable}
\begin{longtable}[t]{r}
\caption{\label{tab:unnamed-chunk-47}Total cost after ordering}\\
\hline
0.6943549\\
\hline
\end{longtable}

Applying this permutation, we derive the reordered matrix \(A^\sigma = A[\sigma]\)  and compute the new total cost: \(Tr(A^\sigma)\), as given in \eqref{tab:unnamed-chunk-47}. 
This simple example demonstrates the purpose of LSAP-type algorithms, which is to find a permutation minimizing the assignment cost.

The standard LSAP assumes square matrices (equal input sizes), but practical applications often involve rectangular matrices where \(M \le N\). Our framework supports this case, as illustrated in Figure~\ref{fig:408}, where the LSAP is applied to sets of unequal sizes. These cases arise, for example, when clustering a large dataset using a smaller set of prototype centroids.

\begin{figure}
\centering
\includegraphics[width=0.7\textwidth, keepaspectratio]{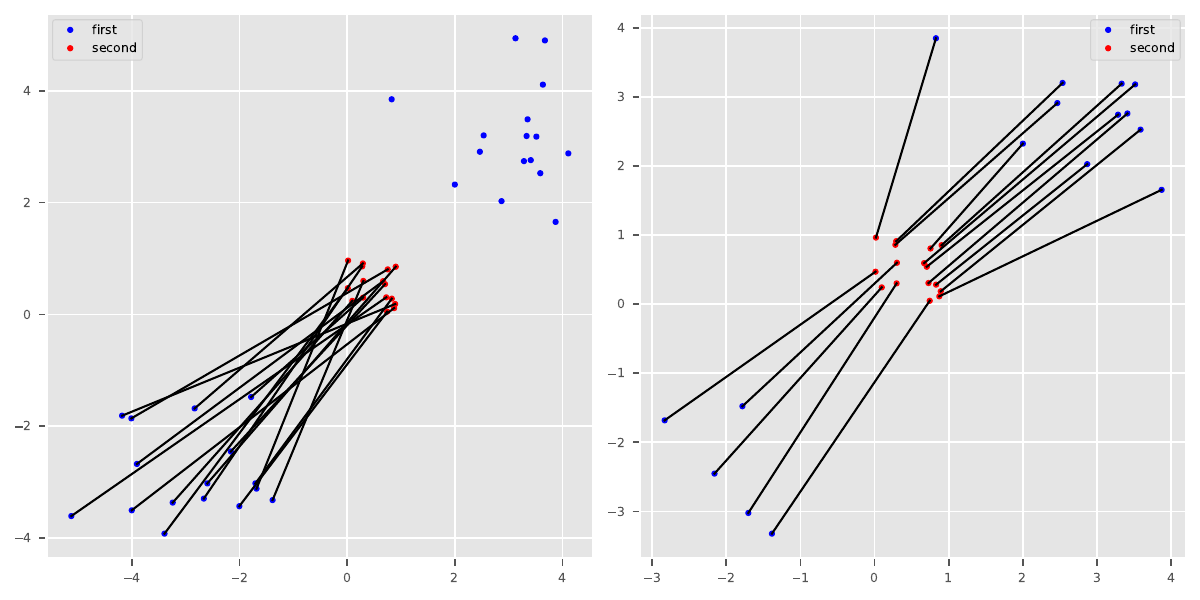}
\caption{\label{fig:unnamed-chunk-54}\label{fig:408}LSAP with different input sizes: $M < N$ (left) and $M=N$ (right)}
\end{figure}


\section{Clustering algorithms for kernels}
\label{clustering-algorithms}


\subsection{Proposed strategy}\label{proposed-strategy}

In light of the error formula~\eqref{dk}, the optimal choice of cluster centres \(Y\)  for approximating a kernel-induced function space
\(\mathcal{H}_{k,X}\)  is determined by solving the following minimization problem:
\be \label{SDS}
    \mathop{\arg\inf}_{Y \subset \RR^{N_y,D}} d_k\big(Y,X\big)^2.
\ee 
We refer to sequences that minimize this criterion as \textit{sharp discrepancy}
sequences, as they yield the most suitable meshes for kernel-based methods, particularly when $X$ represents a continuous distribution. While such sequences can be computed analytically for specific kernels--e.g., Fourier-based kernels  
--most practical settings require numerical approximations of \eqref{SDS}. 

We emphasize that this clustering formulation, based on minimizing the full discrepancy \eqref{SDS} considers
the full discrepancy functional \(d_k\big(Y,X\big)^2\),differs fundamentally from kernel k-means, a more heuristic method typically minimizing intra-cluster variance in the feature space.

Given the computational cost of solving~\eqref{SDS} exactly, we adopt a multi-stage strategy combining speed and accuracy. Specifically, we employ a sequence of three algorithms--ordered from fastest to most computationally intensive--to approximate optimal clusters. The output of one algorithm is used to initialize the next, forming a hierarchical refinement scheme.

The central idea is to first approximate the centroids as a subset \(Y \subset X\) using fast combinatorial algorithms, and then optionally refine them outside of $X$ via descent methods. In what follows, $N_y$ denotes the number of desired clusters, and
\(N_x\) is the size of the training set.

\subsection{Greedy clustering method}
\label{greedy-clustering-methods}

Our initial approximation replaces~\eqref{SDS} with the discrete optimization problem: \be 
\label{SDSD}
    \inf_{Y \subset X} d_k\big(Y,X\big),
\ee 
which we solve using the Greedy Search Algorithm described in Section~\ref{greedy-search-algorithms} (Algorithm~\ref{cl0}). The kernel-based distance functional used at each iteration is given by
\be
    d(Y,\cdot) = \frac{1}{(N_y+1)^2}\Big(2 \sum_{m=1}^{N_y} k(y^m,\cdot)+k(\cdot,\cdot)\Big) - \frac{2}{(N_y+1)N_x}\sum_{n=1} ^{N_y} k(x^n,\cdot).
\ee 
The cross-terms \(\sum_{n=1}^{N_x} k(x^n,x^m)\) can be precomputed at cost \(\mathcal{O}(N_x^2)\), with a memory
complexity \(\mathcal{O}(N_x)\) assuming on-the-fly evaluation.  At each iteration \(n\) of Algorithm~\ref{cl0}, evaluating \(d(Y^n,\cdot)\)  requires approximately
\((n \times N_x -n)\) kernel evaluations. These can be parallelized or accelerated using a precomputed Gram matrix, provided the matrix of size $N_x^2$ fits into available memory.

The overall computational complexity is \(\mathcal{O}(N_y^2 N_x + N_x^2)\).
This algorithm can handle medium-to-large datasets, depending on
settings. We consider here a version that necessitate to
precompute the full Gram matrix, requiring \(\mathcal{O}(N_x^2)\) calls
to the kernel, hence, is adapted to medium datasets. For large-scale data, we may consider kernel k-means, which relies on a similar distance but is less computationally demanding and generally provides lower-quality clusters.

Solving the sharp discrepancy problem~\eqref{SDS} is particularly valuable for mesh generation tasks, where the goal is to approximate the entire function space \(\mathcal{H}_{k,X}\), as in unsupervised learning. Alternatively, when approximating a specific function \(f\)
via interpolation (see equation~\eqref{FIT}), the setting becomes supervised. In that case, a modified discrepancy functional provides an effective selection criterion.

We define the
following distance for any \(1 \le p\le +\infty\): \be 
\label{DF}
 D(Y,X) = \sum_{n=1}^{N_x} \big|f(x^n) - f_{k, \theta_Y}(x^n)\big|^p, \quad \theta_Y=K(Y,Y)^{-1}f(Y),
\ee where \(f\) is a given, vector-valued function. Minimizing this discrepancy over subsets $Y\subset X$ selects training points that best represent the target function \(f\). This approach is closely related to kernel adaptive mesh and control variate techniques, which are valuable for both numerical simulations and statistical estimation.

To accelerate this greedy selection process, we apply block matrix inversion techniques, enabling efficient updates of the inverse Gram matrix \(K(Y^n, Y^n)^{-1}\) from \(K(Y^{n-1}, Y^{n-1})^{-1}\).  With this optimization, the overall complexity remains
\(\mathcal{O}(N_y^3 + N_x N_y^2)\) with moderate memory requirements. This makes the method particularly effective for selecting a large number of representative clusters in high-dimensional or data-intensive settings.

\subsection{Subset clustering method}
\label{subset-clustering-methods}

Our next clustering algorithm allows us to refine the previous
approximation of \eqref{SDSD}, using the discrete permutation algorithm
\eqref{GLSAP}. This method aims to improve the quality of the centroids $Y\subset X$ by locally optimizing the assignment of points using permutation-based descent.

Let $\sigma:[0\ldots N_x]  \to [0\ldots N_x]$ be a permutation and consider a number of clusters $N_y < N_x$. In view of the descent algorithm ~\eqref{PBA}, we define the permutation gain function:
\be
  s(i,j,\sigma) = K_{XY}^{\sigma}-K_{XY}^{\sigma_{i,j}}+K_{YY}^{\sigma}-K_{YY}^{\sigma_{i,j}},
\ee
where $\sigma_{i,j}$ denotes the permutation obtained by swapping indices \(i\) and \(j\), and where the vectors \(K_{XY}^{\sigma},K_{YY}^{\sigma}\) are computed as \be
    K_{YY}^{\sigma} = \frac{1}{N_y^2} \sum_{n,m}^{N_y,N_y} k(x^{\sigma^n},x^{\sigma^m}), \quad K_{XY}^{\sigma}= - \frac{2}{N_x N_y} \sum_{n,m}^{N_x,N_y} k(x^{n},x^{\sigma^m}),
\ee 
in which the computation can be accelerated using pre-computational techniques.
We define the active cluster set as the permuted set \(Y = X\circ \sigma[1,\ldots,N_y]\). To ensure that swapping indices only involve transferring points between the selected centroids and the remaining dataset, we set \(s(i,j,\sigma)=0\), if \(i\ge N_y\) and \(j\ge N_y\), or
\(i\le N_y\) and \(j\le N_y\). Thus, we restrict the for-loop in Algorithm~\ref{PBA} to: $\texttt{for\ i\ in\ {[}1,\ N\_Y{[},\ for\ j\ in\ {[}N\_Y,\ N\_X{[}}$.

This approach yields a sub-optimal solution to~\eqref{SDSD}, refining the centroid selection \(Y\subset X\). As the method relies on local optimization, the quality of the final solution is heavily influenced by the initialization. For this reason, it is recommended to initialize $\sigma$  using the output of the greedy algorithm described in Section~\ref{greedy-clustering-methods}.

Although a specific complexity bound is not provided, empirical evidence suggests a polynomial runtime of approximately \(\mathcal{O}(N_x^2 N_y)\). The method is amenable to parallelization. However, due to its quadratic dependency on $N_x$, it may become computationally prohibitive for very large datasets, in which case scalable alternatives such as kernel k-means may be preferred.

\hypertarget{sharp-discrepancy-sequences}{%
\subsection{Sharp discrepancy
sequences}\label{Discrepancy functional and sharp-discrepancy-sequences}}

We now revisit the full optimization problem~\eqref{SDS} to construct sharp discrepancy sequences, which correspond to optimal centroid configurations minimizing the full kernel discrepancy. These can be computed via the general descent method presented in Algorithm~\ref{DA=algo}, using the functional: \(J(Y)=d_k(X,Y)\), and employing the explicit expression for the gradient \(\nabla J(Y)\), derived in equation~\eqref{nabla}.

To avoid convergence to poor local minima, we initialize the descent using the output of the subset permutation method described above, i.e., by setting the initial iterate \(X^0\) to the centroids $Y\subset X$ obtained from the discrete optimization step.

Two important considerations apply to this method. 
\begin{itemize}
\item Computational complexity. While this gradient-based algorithm produces high-quality centroids, it is computationally expensive. It may be impractical for very large datasets or high-dimensional feature spaces. In such cases, the faster combinatorial algorithms provide viable approximations at the expense of some quality loss.
\item
  Functional concavity. For certain classes of non-smooth kernels \(k\), the discrepancy functional may exhibit concavity on large subsets of the domain. In these cases, continuous descent algorithms like Algorithm~\ref{DA=algo} are ill-suited, and combinatorial optimization strategies should be preferred for global search.
\end{itemize}


\hypertarget{balanced-clustering}{%
\subsection{Balanced clustering}\label{balanced-clustering}}

Balanced clustering is a method to define clusters of
comparable size, an important property for large scale dataset approach. There exist several algorithms capable to handle
balanced clusters. 
In our context, we
work with cluster algorithms based upon an induced distance
\(d(\cdot,\cdot)\), chosen to be the Euclidean distance for k-means and the
kernel discrepancy for kernel-based clustering algorithms. 

We propose a general approach to balanced clustering based on a distance matrix 
\(D \in \mathbb{R}^{N_y,N_x}\), computed through a given distance \(\Big(d(y^i,x^j)\Big)_{i,j}\), which encodes the relation between data point $x^j$ and centroid \(y^i\).  The centroids \(y^i\) are assumed to be determined via an external method, as the ones described in this chapter. Our goal is to solve the following discrete optimal transport
problem, where \(\%\) holds for modulo \be\label{balanced}
  \inf_{\sigma \in \Sigma} \sum_{n=1}^{N_x} d(y^{(\sigma^n\ \% \ N_y)},x^n).
\ee 

This assignment ensures that each cluster receives approximately the same number of data points. The objective~\eqref{balanced} can be optimized via the discrete permutation descent algorithm (Algorithm~\ref{PBA}), using the gain function: 
\be
  \sigma(i,j,\sigma)=D(i,\sigma^{i \% N_y})-D(j,\sigma^{j \% N_y}) - D(i,\sigma^{j \% N_y}) +D(j,\sigma^{i \% N_y}).
\ee 
By design, this algorithm produces balanced clusters, assigning each point
\(x^n\) to a cluster \(y^{\sigma^n \ \% \ N_y}\). An associated allocation function is also naturally defined: 
\(l(\cdot) = l^{\mathop{\arg\inf}_n d(x^n,\cdot)) \ \% \ N_y}\), which maps each data point to its corresponding cluster index.

This approach also allows for flexibility in initialization. For instance, we may initialize $Y$ with \(N_y\) randomly selected points from \(X\), and use either Euclidean or kernel-based distances. Despite its simplicity, this method is numerically efficient, and our tests show that it yields high-quality, balanced clusters with minimal computational overhead.

\subsection{Numerical illustration} 

\paragraph{Illustration of balanced clustering methods} 

We begin by illustrating the behavior of various clustering strategies, including our proposed methods, in comparison to the standard k-means algorithm as implemented in the scikit-learn Library.  This implementation serves as a benchmark due to its high efficiency, scalability, and widespread use.

All clustering algorithms are tested using the kernel defined in \eqref{standardmap} (made standard in the CodPy Library). We generate a synthetic dataset in $\mathbb{R}^2$ consisting of five well-separated Gaussian blobs (i.e., a mixture of five equally weighted Gaussian components). Each clustering algorithm is applied to partition the dataset into $N_y$ clusters.

Figure~\ref{fig:blob_clustering} presents a qualitative comparison of clustering outcomes. Each subplot visualizes the data colored by cluster assignment, with red crosses denoting cluster centroids. The figure is organized in two rows, as follows. 
\begin{itemize}
    \item Top row: unbalanced clustering methods -- greedy discrepancy, sharp discrepancy, standard k-means, and random selection.
    \item Bottom row: their balanced variants, ensuring approximately equal numbers of points per cluster.
\end{itemize}

Balancing is enforced through a permutation-based optimal transport algorithm (see Section~\ref{balanced-clustering}), which assigns data points as evenly as possible across the clusters. This is especially useful in applications involving fairness constraints, sampling quotas, or mitigation of class imbalance.

Visually, balanced clustering methods lead to a more uniform partitioning of the data. While the unbalanced methods tend to produce clusters that vary significantly in size, especially for sharp and random methods, the balanced versions clearly enforce approximately equal cardinality per cluster, at the possible cost of slightly more irregular cluster boundaries.

\begin{figure}
\centering
\includegraphics[width=0.99\textwidth, keepaspectratio]{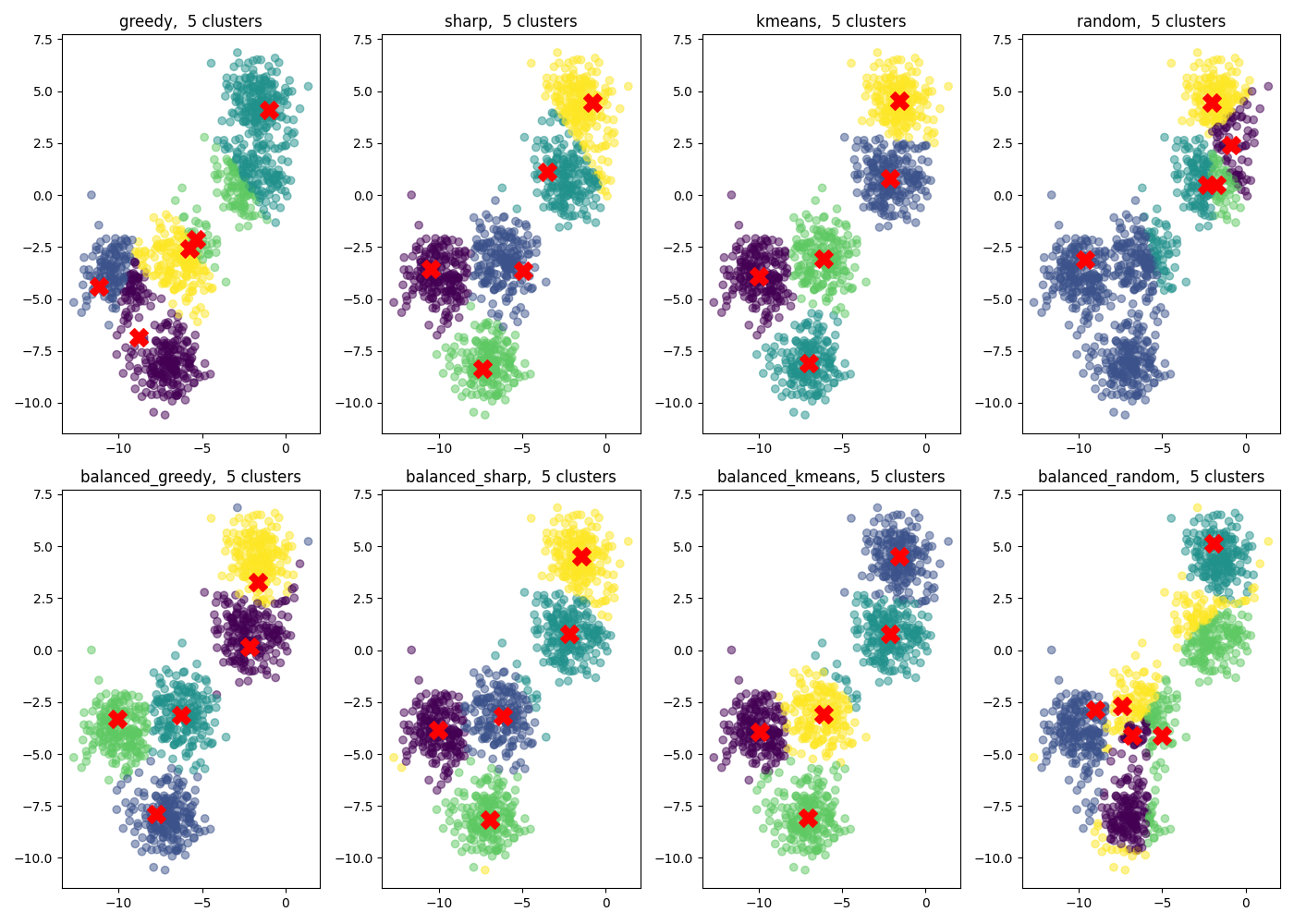}
\caption{\label{fig:blob_clustering}\label{fig:741}Comparison of clustering strategies on a 2D Gaussian mixture with 5 modes. Top row: unbalanced clustering methods (greedy, sharp, k-means, random). Bottom row: balanced variants of each method. Red crosses indicate cluster centroids.}
\end{figure}

\paragraph{Quantitative evaluation of clustering performance} 

To supplement the visual comparison, we also present here quantitative performance metrics for each method, using a typical run with $N_x = 1024$ data points and $N_y = 128$ cluster centres. Results are summarized in Table~\ref{tab:2998}. The metrics include execution time (in seconds), standard k-means inertia (see~the definition in~\eqref{inertia} below), and maximum mean discrepancy (MMD) (see~\eqref{dk} ), which captures how well the distribution of the cluster centroids $Y$ approximates the original distribution $X$ under the kernel $k$.

\begin{longtable}[t]{l|r|r|r|r|r}
\caption{\label{tab:2998}Performance metrics for supervised clustering algorithms}\\
\hline
Method & $N_x$ & $N_y$ & Exec. Time (s) & Inertia & MMD \\
\hline
\endfirsthead
\caption[]{Performance metrics for supervised clustering algorithms \textit{(continued)}}\\
\hline
Method & $N_x$ & $N_y$ & Exec. Time (s) & Inertia & MMD \\
\hline
\endhead
Greedy & 1024 & 128 & 0.0389 & 1028.18 & 0.0000253 \\
Sharp discrepancy & 1024 & 128 & 0.5124 & 1955.15 & 0.0000102 \\
k-means (scikit) & 1024 & 128 & 0.1977 & 408.94 & 0.0007832 \\
Balanced random & 1024 & 128 & 0.0190 & 979.94 & 0.0025103 \\
\hline
\end{longtable}

From this comparison, several observations can be made.
\begin{itemize}
    \item The greedy discrepancy method achieves a strong trade-off between quality and speed, making it particularly suitable for large-scale clustering.
    \item The sharp discrepancy method produces the best MMD (i.e., best approximation of the full data distribution) but is the slowest.
    \item The scikit-learn k-means implementation performs best in terms of inertia, which is expected since this is the criterion it optimizes. It is especially effective when the data are well-clustered and balanced by design.
    \item The balanced random method provides an efficient and simple solution when class balancing is essential, though it sacrifices some accuracy in distributional matching (as shown by its higher MMD).
\end{itemize}


\paragraph{Scalability and convergence behavior} 

We now investigate the behavior of clustering performance as the number of cluster centres $N_y$ increases, focusing on a simplified two-cluster Gaussian mixture (a ``blob dataset'') with $N_x=100$ samples. For  $N_y \in [2, 100]$, we compute
\begin{itemize}
    \item the \textit{discrepancy error}, measuring how well the empirical distribution of cluster centers $Y$ approximates $X$,
    \item and the \textit{inertia} to compare with the classical k-means objective.
\end{itemize}

Figure~\ref{fig:741} presents the evolution of two performance metrics--kernel discrepancy error and inertia--as functions of the number of cluster centres $N_y$, evaluated for three clustering algorithms: the proposed kernel discrepancy minimization method, standard $k$-means, and MiniBatch $k$-means. As the number of clusters increases and approaches the size of the dataset ($N_y \to N_x$), the discrepancy error exhibits a steady and monotonic decline. This behavior reflects the increasing capacity of the centroid set to capture the underlying structure of the data distribution. Notably, the kernel discrepancy method consistently achieves the lowest error across all values of $N_y$, in line with its objective to directly minimize distributional divergence under the kernel.

In terms of inertia, both $k$-means and MiniBatch $k$-means yield lower values, as expected since inertia is the quantity these algorithms are designed to optimize. Nevertheless, the discre\-pancy-based method attains competitive inertia scores, despite not being explicitly tailored for that purpose. This convergence in performance suggests a degree of implicit alignment between minimizing kernel discrepancy and reducing within-cluster variance. The effect becomes more pronounced as $N_y$ increases, indicating that fine-grained approximations of the data distribution tend to naturally support lower intra-cluster dispersion as well.

\begin{figure}
\centering
\includegraphics[width=0.99\textwidth, keepaspectratio]{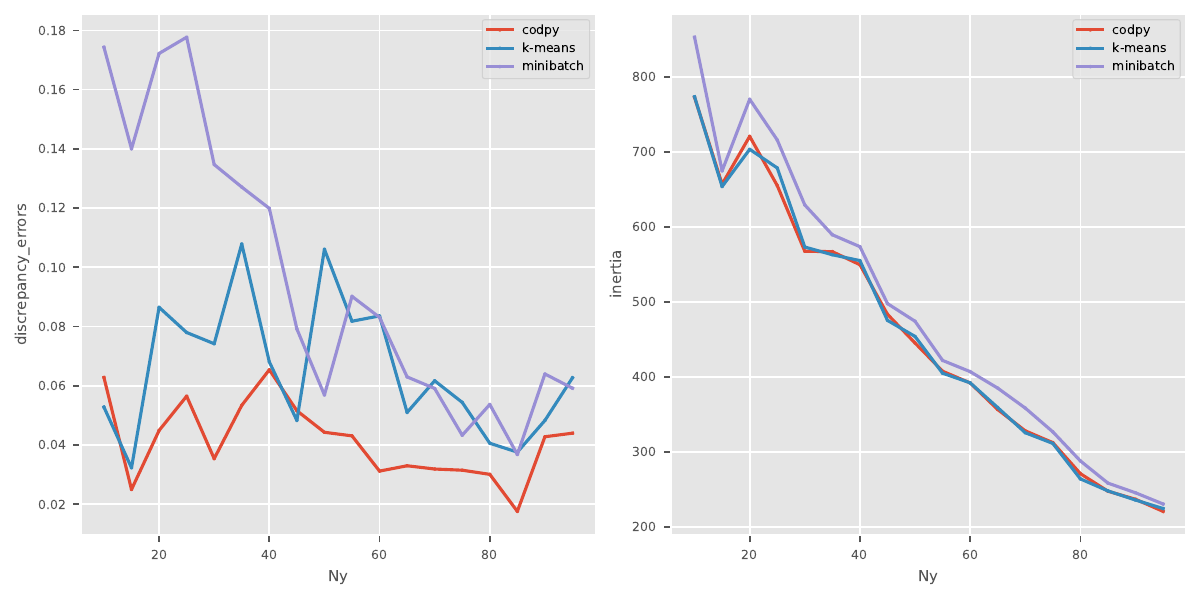}
\caption{\label{fig:741}Comparison of clustering methods across varying numbers of cluster centers $N_y$. Top: kernel discrepancy error. Bottom: k-means inertia. Results are shown for our kernel-based method, standard k-means, and mini-batch k-means.}
\end{figure}


\clearpage

\chapter{Optimal transport and statistical kernel methods}\label{optimal-transport-and-statistical-kernel-methods}

\section{Introduction} 

This chapter focuses on RKHS methods for statistical applications, tackling mainly two related topics. 

On the one hand, we are interested in statistical methods dedicated to (approximate) conditional expectations, conditional densities, conditional distributions, and transition matrix computations. We begin by reviewing the historical Nadaraya-Watson estimator, which relies on a \textit{density model}, adapted to conditional density. This model is compared to a model based on the projection operator~\eqref{FIT}, which is viewed as a \textit{generative model}, adapted to conditional distributions. This comparative study sheds some light on both approaches and motivates the need for \emph{optimal transport techniques} to enhance these historical models; this is done next in the part dedicated to generative methods. Transition probability matrices (also called stochastic matrices) are usually related to Markov-type processes. In this chapter, we consider transition matrices in the context of \textit{martingale stochastic processes.} In such a context, the two models can be used to approximate stochastic matrices, however, optimal transport methods provide a richer layout.

On the other hand, we consider generative methods, which are regarded as the art of defining maps between two arbitrary distributions. There is limited existing work on RKHS-based generative methods, to the best of our knowledge. We provide some arguments to propose a novel method based on optimal transport, which is efficient and easily adapts to a wide class of generative problems.

We emphasize that these methods can be analyzed numerically with optimal transport techniques, as the celebrated Sinkhorn-Knopp algorithm for entropy-regularized optimal transport problems\footnote{see \cite{Cuturi:2013}}. However, instead of this family of algorithms, we have opted\footnote{which, to the best of our knowledge, is a novel contribution to these problems} for an alternative method addressing martingale-related problems with the combinatorial tools in Section~\ref{clustering-strategies}, such as the class of LSAP algorithms, used for instance in \cite{LeFloch-Mercier:2017}. Both approaches (combinatorial vs entropy-regularized) are briefly discussed in Section~\ref{Sinkhorn-vs-Combinatorial-Algorithmic-Approaches}: combinatorial approaches tackle the Monge and Gromov--Monge problems, while entropy-regularized methods seek to solve the Kantorovich and Gromov--Wasserstein problems. The aims of both methods are similar, and some benchmarks can be found in dedicated sections, such as Section~\ref{Maps-and-Generative-methods}
and Section~\ref{application-to-generative-models}.

Due to the importance of the theory of optimal transport, we first provide an overview of its main concepts.


\section{Overview of optimal transport theory}
\label{reminder-on-optimal-transport-theory}

\subsection{Optimal transport on compatible vs. incompatible spaces}
 
Two distributions $X,Y$ being given, lying respectively in $\mathbb{R}^{D_x},\mathbb{R}^{D_y}$, we naturally distinguish between two cases: compatible metric spaces, where $D_x=D_y=D$, which are treated in Section~\ref{Continuous-optimal-transport-on-compatible-spaces}, below; and
 incompatible metric spaces, where $D_x \neq D_y$, which is the topic of~Section~\ref{Continuous-optimal-transport-on-incompatible-spaces}, below. Optimal transport theory is usually formulated between probability distributions on compatible metric spaces, such as $\mathbb{R}^D$, based on a direct point-to-point cost function $c(x, y)$. Such a framework leads to the Monge and Kantorovich formulations and is well-suited for comparing distributions with a common geometry. 

However, modern scenarios ---for instance, comparing shapes, graphs, or latent codes--- often require considering distributions that lie in different metric spaces, with no direct notion of distance between the points, say $x^n$ and $y^m$. In such a case, the Gromov--Wasserstein and Gromov--Monge frameworks provide analogues of optimal transport that align internal structures, rather than pointwise geometry. While Gromov--Wasserstein and Gromov--Monge problems are generally not equivalent to Wasserstein and Monge problems\footnote{However, one-dimensional Gromov--Monge are equivalent to Monge problems, \textit{up to} a mirror symmetry}, they extend the scope of optimal transport to incompatible metric spaces and provide a unified transport framework.


\subsection{Continuous optimal transport on compatible spaces} \label{Continuous-optimal-transport-on-compatible-spaces}

\paragraph{Push-forward maps.} 

We consider the compatible case, that is, $D_x = D_y = D$. We begin by introducing some key concepts from the theory of optimal transport (OT)\footnote{For a comprehensive introduction, see the textbook by Villani \cite{Villani:2009}}. We are going to outline first the relevant notions for the continuous case, as a foundation for our next discussion of the same concepts in a discrete setting.

First of all, push-forward maps are defined (even when $D_x$ and $D_y$ are distinct). We denote by $\mathcal{P}(X)$ the set of all Borel probability measures on a measurable space $X$, and by $\mathcal{C}_b(Y)$ the space of all bounded and continuous functions on a space $Y$. Let $\mu \in \mathcal{P}(\mathbb{R}^{D_x})$ and $\nu \in \mathcal{P}(\mathbb{R}^{D_y})$ be two probability measures. A measurable map $T : \mathbb{R}^{D_x} \to \mathbb{R}^{D_y}$ is said to \emph{transport} $\mu$ to $\nu$, or to \emph{push-forward} $\mu$ to $\nu$, if 
\be\label{transport}
    \int_{\mathbb{R}^{D_x}} (\varphi \circ T)(x) \, d\mu(x) 
    = \int_{\mathbb{R}^{D_y}} \varphi(y) \, d\nu(y) 
    \quad \text{for all } \varphi \in \mathcal{C}_b(\mathbb{R}^{D_y}).
\ee
This is written in a compact form as $T_\# \mu = \nu$. The map $T$ describes how to ``move the mass'' from the source distribution $\mu$ to match the target distribution $\nu$. The formula \eqref{transport} can also be regarded as a change of variables in an integration formula. Indeed, assuming that the map $T$ is sufficiently regular and invertible maps and that $D_x=D_y=D$, the standard change of variable formula specifies the relation between $\mu,\nu$ and $T$ as 
\be \label{COVF}
    \int_{\mathbb{R}^{D}} (\varphi \circ T)(x) \, |\det \nabla T| d\nu(x) = \int_{\mathbb{R}^{D}} \varphi(y) \, d\nu(y) 
    \quad \text{for all } \varphi \in \mathcal{C}_b(\mathbb{R}^{D_y}). 
\ee
Here, $\nabla T$ denotes the Jacobian of the map $T$. 


\paragraph{Monge problem.} 

There exist infinitely many maps satisfying $T_\#\mu = \nu$. To select an \emph{optimal map} for our purpose, we introduce a \emph{cost function}, denoted by $c : \mathbb{R}^{D} \times \mathbb{R}^{D} \to [0, +\infty)$. In this context, the so-called \emph{Monge problem} is defined as
\be\label{MongeProblem} 
   T^* \in \operatorname*{argmin}_{T \, big/ \,  T_\# \mu = \nu} \int_{\mathbb{R}^{D}} c(x, T(x)) \, d\mu(x), 
\ee
and a solution $T^*$ is called an \emph{optimal transport map}. Optimal transport maps do exist between two arbitrary measures, while uniqueness statements require technical assumptions, which are beyond the scope of this monograph.

A commonly used cost function is $c(x, y) = |x - y|_p^p = \sum_d |x_d - y_d|^p$, which leads to the definition of the \emph{$p$-Wasserstein distance} 
\be\label{Wasserstein}
    W_p^p(\mu, \nu) = \inf_{T : T_\# \mu = \nu} \int_{\mathbb{R}^{D}} |x - T(x)|_p^p \, d\mu(x).
\ee
The so-called \emph{polar factorization theorem} then states that the Wasserstein distance $W_2(\mu, \nu)$, corresponding to the Euclidean cost $c(x,y)=|x-y|_2^2$ for the Monge problem~\eqref{MongeProblem}, selects a one-to-one map, expressed as the gradient of a convex function. More precisely, if $T :\mathbb{R}^{D}  \to \mathbb{R}^{D}$ is a mapping and $\mu \in \mathcal{P}(\mathbb{R}^{D})$ is a probability measure, we set $\nu=T_\#\mu$ and assume that $\nu$ is absolutely continuous with respect to the Lebesgue measure. Then the following factorization holds:
\be\label{PF}
    T = (\nabla h)\circ \sigma, \quad \sigma_\# \mu = \mu, \quad h \text{ convex}.
\ee
In this formulation, $\sigma$ is $\mu$-preserving, hence it is interpreted as an arbitrary permutation of equally measurable elements in the support of $\mu$. This decomposition of a mapping $T$ into a map $\nabla h$ and a permutation is unique. In short, the polar factorization for maps states that the Monge problem~\eqref{MongeProblem} with Euclidean cost selects a unique, one-to-one map $T^{*} = \nabla h$ with $h$ convex.


\paragraph{Kantorovich problem.} 

The \emph{Kantorovich relaxation} replaces the map $T$ with a probability measure $\pi \in \mathcal{P}(\mathbb{R}^{D_x} \times \mathbb{R}^{D_y})$ whose marginals are denoted by $\mu$ and $\nu$. The set of such couplings is defined as the collection of all probability measures with fixed marginals $\mu$ and $\nu$ (for all Borel sets $A, B$):
\be
\Pi(\mu, \nu) = \left\{ \pi \in \mathcal{P}(\mathbb{R}^{D} \times \mathbb{R}^{D}) \mid \pi(A \times \mathbb{R}^{D}) = \mu(A),\qquad
 \pi(\mathbb{R}^{D} \times B) = \nu(B) \right\}. 
\ee
The Kantorovich problem seeks to minimize the expected transport cost over all couplings: 
\be\label{KantorovichC}
    \inf_{\pi \in \Pi(\mu, \nu)} \int_{\mathbb{R}^{D_x} \times \mathbb{R}^{D}} c(x,y) \, d\pi(x, y).
\ee
For instance, minimizing the cost of transport between $\mu$ and $\nu$, according to the $p$-Wasserstein cost function $c(x,y)=|x - y|_p^p$, defines the $p$-Wasserstein distance as
\be
    W_p^p(\mu, \nu)  =  \inf_{\pi \in \Pi(\mu, \nu)} \int_{\mathbb{R}^{D} \times \mathbb{R}^{D}} |x - y|_p^p \, d\pi(x, y).
\ee
The two definitions of the $p$-Wasserstein distance in~\eqref{KantorovichC} and~\eqref{Wasserstein} are equivalent. 
 

\paragraph{Dual problem and weak formulation.} 

The dual problem to the Kantorovich problem \eqref{KantorovichC} is given by: 
\be
\label{dualC}
\sup_{\varphi, \psi} \Bigl( \int \varphi d\mu - \int \psi d\nu \Bigr), \quad \text{subject to} \quad \varphi(y) - \psi(x) \le c(x, y),
\ee
where $\varphi : X \to \mathbb{R}, \ \psi : Y \to \mathbb{R}$ represent potential functions. 

Next, consider a coupling $\pi \in \Pi(\mu, \nu)$ and denote by $\pi(y|x)$ the regular conditional distribution of $y$ given $x$ obtained through the formula $d\pi(x, y) = d\pi(y|x) d\mu(x)$, where we assume that $\pi(x, y)$ is dominated by $\mu(x)$. We then reformulate the Kantorovich problem as
\be\label{WOT}
    \inf_{\pi \in \Pi(\mu, \nu)} \int_{\mathbb{R}^{D_x}} \omega(x,\pi(y|x))\, d\mu(x), \quad \omega(x,\pi(y|x)) = \int_{\mathbb{R}^{D_y}} c(x,y) d\pi(y|x).
\ee
The \emph{weak optimal transport} (WOT) problem replaces the local cost function $c=c(x,y)$ with a cost function expressed as $\omega(x,\pi(y|x))$ and involving a conditional distribution. For instance, consider $D_x=D_y=D$, and a coupling $\pi(x,y) \in \Pi(\mu, \nu)$. Then the cost induced by the Euclidean distance defines the \textit{barycentric transport cost} as 
\be\label{BTC}
 V^2_2(\pi)=\int_{\mathbb{R}^{D}} \left|x - \int_{\mathbb{R}^{D}} y \, d\pi(y|x)\right|_2^2 \, d\mu(x), \quad d\pi(x,y) = d\pi(y|x) d\mu(x).
\ee
This formulation penalizes the deviation between each source point $x$ and a \emph{barycenter} computed with the distribution $\pi(y|x)$. 

The set of joint probability $\pi(x,y)$ satisfying  $V^2_2(\pi)=0$, or $\mathbb{E}\big(y|x\big) = x$ is called the set of \textit{martingale couplings}\footnote{In other words, martingale coupling is a transport plan 
$\pi \in \Pi(\mu,\nu)$ such that 
$\mathbb{E}[Y \mid X=x]=x$ for $\mu$-almost every $x$. 
It models displacements with zero conditional drift, 
typical in finance or stochastic processes.} Strassen's Theorem states that this set is nonempty, provided that $\mu$ is dominated by $\nu$ according to the stochastic convex order, notation $ \mu \ll \nu $. 
Let us denote the set of martingale couplings  by $\mathcal{M}(\mu, \nu) = \{ \pi \in \Pi(\mu, \nu) : V_2^2(\pi) = 0\}$ and refer to
\be\label{MOT}
 \inf_{\pi \in \mathcal{M}(\mu, \nu)} \int_{\mathbb{R}^{D} \times \mathbb{R}^{D}} \left|x - y\right|_2^2 d\pi(x,y) 
\ee
as a \textit{Martingale Optimal Transport} (MOT) problem. 


\subsection{Continuous optimal transport on incompatible spaces}
\label{Continuous-optimal-transport-on-incompatible-spaces}

\paragraph{Purpose.}

\paragraph{Purpose.}

However, in many modern applications --such as dimensionality reduction, encoder--decoder frameworks, shape matching, graph alignment, or multimodal translation-- the source and target distributions live in \emph{different} metric spaces, making it impossible to define directly a cost function between points. In these cases, the Gromov--Wasserstein framework extends optimal transport to operate on pairwise structural similarities rather than pointwise distances. These two perspectives ---compatible and incompatible spaces--- can be unified under a broader view of optimal transport, where the key distinction lies in whether a \emph{common ground metric} is available for defining transport costs.


\paragraph{Gromov--Wasserstein and Gromov--Monge problems.}

Given two metric spaces \((\mathcal{X}, d_{\mathcal{X}})\) and \((\mathcal{Y}, d_{\mathcal{Y}})\), with \(\mathcal{X} = \mathbb{R}^{D_x}\) and \(\mathcal{Y} = \mathbb{R}^{D_y}\), each endowed with probability measures \(\mu\) and \(\nu\), and dissimilarity functions \(c_{\mathcal{X}} : \mathcal{X} \times \mathcal{X} \to \mathbb{R}\) and \(c_{\mathcal{Y}} : \mathcal{Y} \times \mathcal{Y} \to \mathbb{R}\), the Gromov--Wasserstein (GW) distance is defined as
\[
\label{GW}
GW_p^p(\mu, \nu) = \inf_{\pi \in \Pi(\mu, \nu)} \int_{\mathcal{X}^2 \times \mathcal{Y}^2} \bigl|c_{\mathcal{X}}(x, x') - c_{\mathcal{Y}}(y, y')\bigr|^p \, d\pi(x, y) \, d\pi(x', y'),
\]
where \(\Pi(\mu, \nu)\) denotes the set of probability couplings with marginals \(\mu\) and \(\nu\).
This formulation aligns the internal structures of the two spaces by comparing all pairs of intra-space distances. (A common choice is \(c_{\mathcal{X}} = d_{\mathcal{X}}\) and \(c_{\mathcal{Y}} = d_{\mathcal{Y}}\).) The minimizer $\pi$ represents an optimal coupling between $\mu$ and $\nu$ in terms of structural preservation.

In the Monge-type variant, referred to as the \emph{Gromov--Monge problem}, we restrict attention to deterministic transport maps \(T : \mathcal{X} \to \mathcal{Y}\) such that \(T_\#\mu = \nu\). In this case, the coupling \(\pi\) is induced by the map, i.e., \(\pi = (\mathrm{id}, T)_\#\mu\), and the GW objective reads
\be\label{eq:GM}
GM_p^p(\mu, \nu)
  = \inf_{\,T:\, T_\#\mu = \nu}
     \int_{\mathcal{X}^2}
       \bigl|c_{\mathcal{X}}(x, x') - c_{\mathcal{Y}}(T(x), T(x'))\bigr|^p
     \, d\mu(x)\, d\mu(x').
\ee
This formulation is the Gromov--Monge problem; when a measurable map \(T\) exists that achieves the infimum, it is called an \emph{optimal Gromov--Monge map} between \(\mu\) and \(\nu\).

Hence, we have introduced first the Monge formulation on compatible spaces to fix notations
and intuition (push-forward maps, Wasserstein distances, etc.).
We then passed to Kantorovich's relaxation in couplings. Finally, we discussed
Gromov--Wasserstein and next Gromov--Monge formulations, which extend the optimal transport theory to incompatible spaces by aligning internal structures rather than pointwise distances.


\subsection{Discrete optimal transport on compatible spaces}
\label{Discrete-Optimal-Transport-on-Compatible-Spaces}

\paragraph{Discrete Monge formulation.}

The discrete standpoint mirrors the continuous presentation. For discrete optimal transport, we consider two equally weighted distributions supported on finite point clouds \(X=\{x^1,\dots,x^N\}\subset\mathbb{R}^{D_x}\) and \(Y=\{y^1,\dots,y^N\}\subset\mathbb{R}^{D_y}\). This setting defines the discrete measures
\be
\delta_X=\frac{1}{N}\sum_{n=1}^N \delta_{x^n}, \qquad
\delta_Y=\frac{1}{N}\sum_{n=1}^N \delta_{y^n},
\ee
where \(\delta_{z}\) denotes the Dirac measure at the point \(z\). For any bijection \(T:X\to Y\) (e.g., \(T(x^n)=y^{\sigma(n)}\) for some permutation \(\sigma\)), the push-forward relation \(T_{\#}\delta_X=\delta_Y\) holds, as in~\eqref{transport}.

Given a cost function \(c:X\times Y\to\mathbb{R}\), let \(C(X,Y)\in\mathbb{R}^{N\times N}\) be the cost matrix with entries \(C_{nm}=c(x^n,y^m)\). The discrete Monge problem then seeks a permutation \(\sigma\in\Sigma\) (the symmetric group on \(\{1,\dots,N\}\)) minimizing the total transport cost:
\be
\label{LSAP-deux}
\overline{\sigma} \in \operatorname*{arg\,min}_{\sigma\in\Sigma}\;\sum_{n=1}^N c\big(x^n,y^{\sigma(n)}\big)
 \;=\; \operatorname*{arg\,min}_{\sigma\in\Sigma}\;\mathrm{Tr}\big(C(X,Y^{\sigma})\big),
\ee
where \(Y^{\sigma}\) denotes \(Y\) with its columns reordered by \(\sigma\). This is the classical \emph{Linear Sum Assignment Problem (LSAP)}, described in Section~\ref{linear-sum-assignment-problems}. Solving~\eqref{LSAP-deux} yields the permuted set \(Y^{\overline{\sigma}}=Y\circ\overline{\sigma}\) and the map \(T^{\overline{\sigma}}\) defined by \(T^{\overline{\sigma}}(x^n)=y^{\overline{\sigma}(n)}\), which realizes the discrete optimal transport and satisfies \(T^{\overline{\sigma}}_{\#}\delta_X=\delta_Y\) in the sense of~\eqref{transport}.

Let us now interpret the polar factorization \eqref{PF} in a discrete setting. Consider $D_x=D_y=D$, the permutation $\overline{\sigma}$ determined by the discrete Monge problem with Euclidean costs and consider the discrete mapping $T^{\overline{\sigma}}(X) = Y^{\overline{\sigma}}$.
This map is the optimal map transporting $\delta_X$ into $\delta_Y$, which is defined only pointwise on the distribution $X$. Applying carelessly\footnote{The assumptions of the polar factorization are not fulfilled here, as we deal with discrete distributions} the polar factorization theorem \eqref{PF} implies that $T^{\overline{\sigma}}(x)$ can be seen as the gradient of a convex function whose values are prescribed. In particular, a numerical method to approximate the convex potential $h$ is, considering a kernel $k$, to compute the Helmholtz-Hodge decomposition of the map $T^{\overline{\sigma}}$ (see~\eqref{Delta-inv}), leading to the numerical approximation
\be\label{PFN}
    h_k(X) = \Delta_k^{-1} (\nabla_k \cdot Y^{\sigma}).
\ee

\paragraph{Discrete Kantorovich relaxation.}

Instead of restricting to deterministic bijections as in the Monge formulation, the Kantorovich formulation allows for probabilistic couplings, represented by doubly stochastic matrices:
\be
\label{KantorovichD}
\overline{\Pi} = \mathop{\arg\min}_{\Pi \in \Gamma} \sum_{n,m=1}^N \pi_{n,m} c(x^n, y^m) = \mathop{\arg\min}_{\Pi \in \Gamma} \langle\Pi,C(X,Y)\rangle,
\ee
where $<\cdot,\cdot>$ is the Frobenius scalar product of matrices, and $\Gamma$ is called the Birkhoff polytope of stochastic matrices:
\be \label{birkhoff_polytope}
    \Gamma  =  \Bigl\{ 
    \gamma \in \mathbb{R}^{N \times N} \, \Big| \, 
    \sum_{m=1}^N \gamma_{n,m} = 1, \;\sum_{n=1}^N \gamma_{n,m} = 1,\; \gamma_{n,m} \ge 0 \Bigr\}.
\ee
This set of squared matrices is closed and convex.

\paragraph{Discrete dual formulation and optimality.}
The discrete version of the continuous dual problem \eqref{dualC} is characterized as
\be
\label{dual}
\sup_{\varphi, \psi} \sum_{n=1}^N \varphi(x^n) - \psi(y^n), \quad \text{subject to} \quad \varphi(x^n) - \psi(y^m) \le c(x^n, y^m),
\ee
where $\varphi : X \to \mathbb{R}, \ \psi : Y \to \mathbb{R}$ are potential functions. Observe that $\psi=\phi$ if the cost function is symmetrical and satisfies $c(x,x)=0$, as is the case for distance cost functions.

\paragraph{Equivalence of discrete optimal transport formulations}

The equivalence of the Monge formulation, the Kantorovich relaxation, and the associated dual problem has been rigorously established\footnote{see Brezis \cite{Brezis:2018}, especially Theorem~1.1 therein} in our context of equally weighted discrete distributions. Notably, the minimum cost obtained in the Monge formulation and the relaxed Kantorovich formulation are the same. This means that  
\be
\min_{\sigma \in \Sigma} \sum_{n=1}^N c(x^n, y^{\sigma(n)}) = \min_{\gamma \in \Gamma} \sum_{n,m=1}^N \gamma_{n,m} c(x^n, y^m).
\ee
Regarding the dual problem \eqref{dual}, the values of the associated dual problem satisfy the so-called complementary slackness condition
\be
\varphi(x^n) - \psi(y^{\overline{\sigma}(n)}) = c(x^n, y^{\overline{\sigma}(n)}).
\ee

\paragraph{Discrete weak optimal transport formulation.}

The weak optimal transport formulation \eqref{MOT} uses the definition of the set of martingale couplings having a null barycentre transport cost \eqref{BTC}. The discrete version of this cost, for a stochastic matrix $\pi \in \Gamma$ is given by
\be
\label{discrete_weak_ot}
V^2_2(\pi) = \sum_{n=1}^N \left | x^n - \sum_{m=1}^N \pi_{n,m} y^m \right |_2^2, \quad \pi \in \Gamma.
\ee
This defines the set of discrete martingales as $\mathcal{M} = \{\pi \in \Gamma : V^2_2(\pi) = 0\}$, and the discrete version of the martingale optimal transport problem \eqref{MOT} is
\be
\label{KantorovichM}
\overline{\Pi} = \mathop{\arg\min}_{\Pi \in \mathcal{M}} \sum_{n,m=1}^N \pi_{n,m} c(x^n, y^m).
\ee
This problem can be seen as a modification of the classical Kantorovich objective through the martingale constraint. Through relaxation, these problems can also be interpreted as standard Kantorovich problems, with cost having an extra penalization term given by the barycentre cost.

\subsection{Discrete optimal transport on incompatible spaces}\label{Discrete-Optimal-Transport-on-Incompatible-Spaces}

\paragraph{Discrete Gromov--Wasserstein problem.} \label{Discrete-GW}

In contrast to the classical optimal transport setting where the cost is defined over a common ground space, the Gromov--Wasserstein (GW) framework allows for comparing probability measures supported on different metric spaces. 

Let $X = \{x^1, \dots, x^N\} \subset \mathbb{R}^{D_x}$ and $Y = \{y^1, \dots, y^M\} \subset \mathbb{R}^{D_y}$ and let $C_X \in \mathbb{R}^{N \times N}$ and $C_Y \in \mathbb{R}^{N \times N}$ be the pairwise distance matrices: $(C_X)_{n,n'} = D_x(x^n, x^{n'})$ and  $(C_Y)_{m,m'} = d_Y(y^m, y^{m'})$.
The discrete Gromov--Wasserstein problem seeks a probabilistic coupling $\Pi \in \mathbb{R}^{N \times N}$, minimizing the distortion between pairwise intra-space distances:
\be
\label{GW-discrete}
GW_2(X,Y) = \mathop{\arg\min}_{\Pi \in \Gamma} \sum_{n,n'=1}^N \sum_{m,m'=1}^N |(C_X)_{n,n'} - (C_Y)_{m,m'}|^2 \, \pi_{n,m} \, \pi_{n',m'} 
\ee
where $\Gamma$ is the set of admissible couplings.

The objective function in \eqref{GW-discrete} defines a quadratic form in $\Pi$ which is nonconvex, NP-hard\footnote{See E.M. Loiola et al., 
A survey for the quadratic assignment problem, European J. Operational Research 176 (2007):657--690.
}, and can be seen as measuring the discrepancy between the geometry of the two spaces under the coupling. 

\paragraph{Discrete Gromov--Monge problem.} \label{Discrete-GM}

The Gromov--Monge (GM) problem is a deterministic counterpart to the Gromov--Wasserstein problem, which seeks a structure-preserving map between discrete metric spaces. Instead of optimizing over couplings, the GM formulation restricts to bijective mappings, similar in spirit to the Monge formulation of optimal transport. 

Gromov--Monge problems and Gromov--Wasserstein approaches coincide in the case of discrete, equi-weighted distributions having the same length $X = \{x^1, \dots, x^N\} \subset \mathbb{R}^{D_x}$ and $Y = \{y^1, \dots, y^N\} \subset \mathbb{R}^{D_y}$. Consider two distance functions $c_\mathcal{X}: X \times X \to \mathbb{R}$ and $c_\mathcal{Y}: Y \times Y \to \mathbb{R}$, which measure intra-space similarities. 
The discrete Gromov--Monge problem seeks a permutation $\sigma \in \Sigma$, where $\Sigma$ is the set of bijections $\{1, \dots, N\} \to \{1, \dots, N\}$, that best preserves pairwise structural relations:
\be
\label{GM-discrete}
GM_2(X, Y) = \mathop{\arg\min}_{\sigma \in \Sigma} \sum_{i,j=1}^N |c_\mathcal{X}(x^i, x^j) - c_\mathcal{Y}(y^{\sigma(i)}, y^{\sigma(j)})|^2.
\ee

This objective aligns the relational structures of $X$ and $Y$ by minimizing the discrepancy between their pairwise costs under a deterministic mapping. The GM problem is a combinatorial quadratic assignment problem (QAP), and as such is generally NP-hard\footnote{See \url{https://en.wikipedia.org/wiki/Quadratic_assignment_problem}}.


\subsection{The class of Sinkhorn-Knopp algorithms} 
\label{Sinkhorn-vs-Combinatorial-Algorithmic-Approaches}

The Sinkhorn\textendash Knopp methodology generates an efficient family of algorithms that allow one to tackle numerically many of the optimal transport problems described above, such as the Monge or Gromov--Wasserstein problems. This approach\footnote{popularized in~\cite{Cuturi:2013} and subsequent works \cite{Bach,MF2011,Cuturi:2016}} is based on the notion of \textit{entropy regularization}.

However, in this monograph, we present and use an alternative approach to the Sinkhorn--Knopp method, based on combinatorial analysis described in Section~\ref{clustering-strategies} (on clustering). We briefly discuss and motivate both approaches below, and we describe the Sinkhorn--Knopp family of algorithms, providing references and links to dedicated libraries\footnote{For a comprehensive introduction to the theory, see Peyré and Cuturi (2019), \emph{Computational Optimal Transport}. For practical Python implementations, the Python Optimal Transport (POT) library is widely used: \url{https://pythonot.github.io/}. For modern and scalable OT implementations, see the \texttt{ott-jax} library: \url{https://ott-jax.readthedocs.io/}.}.


Considering the Kantorovich problem \eqref{KantorovichD}, we observe that the set $\Gamma$ can be characterized as the intersection of two simpler sets $\Gamma=\Gamma^+ \cap \Gamma^-$, where $\Gamma^{+}$ (respectively  $\Gamma^{-}$) is the set of row-stochastic (respectively  column-stochastic) matrices, satisfying $\sum_{m=1}^N \gamma_{n,m} = 1$ (respectively  $\sum_{n=1}^N \gamma_{n,m} = 1$). The scaling $C^+$ of any cost matrix $C$ toward $\Gamma^{+}$ consists of dividing each column of $C$ by its sum: $C_+ = D_+ C$, $D_+ = \text{Diag}(\sum_{n=1}^N C_{n,m})^{-1}$. Similarly, the projection on $\Gamma^{-}$ consists of normalizing by a right multiplication with a diagonal matrix defined by the sum of rows $C_- = C D_-$. This produces a family of algorithms known as \textit{Iterative Proportional Fitting} (IPF) algorithms, to compute factorization $C = D^{+} \Pi D^{-}$, referred also as \textit{Sinkhorn-Knopp}. For example, the following alternate-direction iterative scheme is a classical implementation of these algorithms:
\be\label{IPF}
 \Pi^{2n} = D_+^{2n-1} \Pi^{2n-1}, \quad \Pi^{2n+1} = \Pi^{2n} D_-^{2n}.
\ee
Sufficient conditions for convergence of this algorithm are marginals of the cost matrix $C$ must be strictly positive, and $C$ is not separable, i.e. this matrix does not permute to a block diagonal one. This is the case if all elements of $C$ are strictly positive, for instance\footnote{The IPF algorithm can be unstable or even blow up numerically without these assumptions.}. This algorithm can be seen as an alternated descent algorithm, switching the normalization of rows and columns until a convergence criterion is finally reached. 

The \textit{entropy-based regularization} approach proposes a relaxation to the Kantorovich formulation \eqref{KantorovichD} as follows: 
\be
\Pi^\varepsilon = \mathop{\arg\min}_{\Pi \in \Gamma} \left\{ \langle \Pi, C \rangle - \varepsilon H(\Pi) \right\}, \quad H(\Pi)  =  -\sum_{n,m} \pi_{n,m} \log \pi_{n,m}, 
\ee
where $H(\Pi)$ is called the Shannon entropy, and $\varepsilon > 0$ is a regularization parameter. To solve this problem, an efficient algorithm is the IPF algorithm, but applied to the Gibbs kernel $K=\exp(-C/\epsilon)$ instead of $C$. The algorithm converges since $K$ is now not separable. 

Unlike the exact solutions to the classical LSAP (see~\eqref{LSAP}) solved by the Hungarian algorithm, the entropy formulation yields a \textit{smooth approximation} to the optimal transport plan. The resulting matrix $\Pi^\varepsilon$ is dense and depends critically on the choice of $\varepsilon$: small values recover sharper solutions but may suffer from numerical instability, while large values yield smoother transitions but less accurate approximations. In practice, tuning $\varepsilon$ balances numerical stability against fidelity to the true transport map.

The family of entropy-regularized algorithms is popular, one reason being that they can handle large datasets via parallel computations. However, these algorithms approximate exact transport, unlike combinatorial algorithms. The entropy-regularized approach results in the loss of the reproducible property of kernel projection, which is problematic in some situations. Moreover, a direct combinatorial approach is particularly effective for low to medium-sized datasets, our primary focus in the present work. Combinatorial approaches can also be run in parallel to handle large datasets (see Section~\ref{large-scale-dataset}), usually performing better in terms of accuracy (see Section~\ref{the-bachelier-problem}), and can also tackle both regularized and non-regularized problems. 

The two algorithmic paradigms ---combinatorial and entropy-regularized--- can be viewed as complementary.
While the Monge problem seeks deterministic maps, the Kantorovich relaxation allows probabilistic couplings. Entropy-regularized variants further smooth the solution and aid in optimization.

\subsection{Numerical illustration of optimal transport maps}

In Figure~\ref{fig:407}, we use the LSAP to compute some simple transfer plans, or matching, between two distributions. Consider $X \subset \mathbb{R}^2$ sampled from a bimodal distribution (plot blue in figures) and $Y \subset [0,1]^2$ sampled uniformly (plot red in figures). Using the Euclidean cost $c(x, y) = \|x - y\|_2$, the optimal matching corresponds to Wasserstein transport, resulting in non-crossing assignments, as shown in Figure~\ref{fig:407}.

If a kernel-induced cost is used instead, the resulting matching may differ from the Euclidean match, potentially leading to crossing assignments. For example, a kernel with product structure $k(x, y) = \exp(-\prod_d |x_d - y_d|)$ (up to a rescaling map) results in a different geometry, as illustrated in Figure~\ref{fig:407}.

\begin{figure}
\centering
\includegraphics[width=1.0\textwidth]{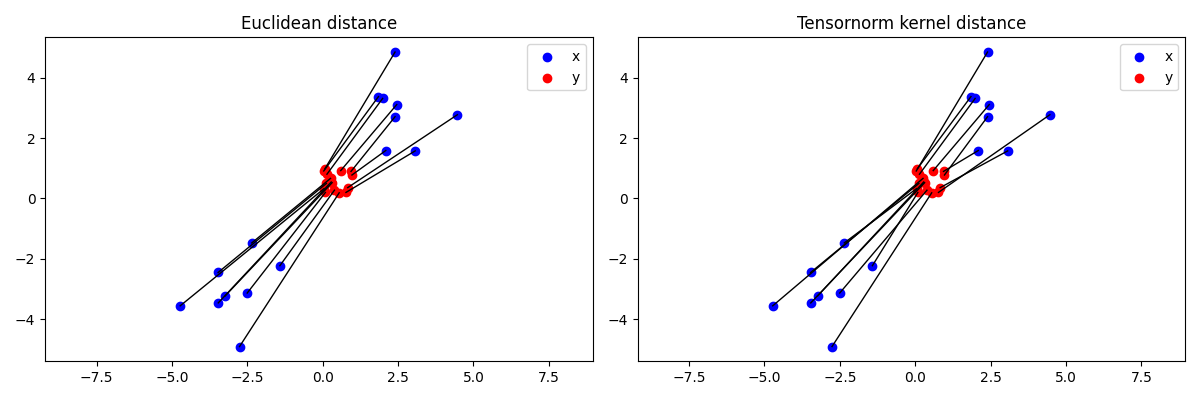}
\caption{\label{fig:407}Same set with two matchings: Euclidean cost (left) and kernel distance MMD cost (right)}
\end{figure}


\section{Conditional expectations and densities,  transition probabilities}
\label{Conditional-Expectations-and-Densities-Transition Probabilities}

\subsection{Purpose}

Let $X, Y \in \mathbb{R}^{D_x} \times \mathbb{R}^{D_y}$ be a pair of random variables defined in a probability space $(\Omega, \mathcal{F}, \mathbb{P})$, with a joint distribution denoted \( (X,Y) \), which consists, from an application point of view, in a discrete joint distribution \( (x^n, y^n) \), \( n=1,\dots,N \). Section~\ref{ACTKBE} presents two methods to tackle two different, although related, tasks: can we estimate the conditional density $\mathbb{E}(Y \mid X=x)$~? Can we estimate the conditional distribution $Y \mid X=x$? We review two different kernel-based approaches, one for each task. 
\begin{itemize}
\item The historical Nadaraya--Watson (NW) estimator, a kernel-based approach, provides a model for conditional density, relying on a density model, called the \textit{kernel density estimate (KDE)}.
\item The kernel ridge projection formula \eqref{FIT} provides a model to approximate conditional distributions, based on the reproducing kernel Hilbert space regression \eqref{FIT}. 
\end{itemize}
Due to the historical importance and prevalence of NW estimators, the next section compares and dissects the links between both approaches, giving an interesting interpretation and viewpoint of the kernel ridge regression formula, which motivates the introduction of optimal transport techniques to enhance these approaches.

We then consider another important recurrent question in statistical analysis: can we estimate the transition probability matrix entries \( \pi_{n,m} = \mathbb{P}(Y = y^m \mid X = x^n) \)? Such problems arise in a wide range of applications, including webpage ranking and risk estimation in mathematical finance. Transition probabilities, which are stochastic matrices, are usually considered with Markov chain models. We focus on a particular Markov-type setting, governed by martingale couplings. This particular structure links to the notion of weak optimal transport, as discussed below. 
Section~\ref{Transition-Probabilities-With-Kernels} presents a combinatorial method, based on the LSAP, to solve the corresponding MOT (Martingale Optimal Transport) problems. As discussed earlier, the Sinkhorn family of algorithms can be adapted to compute MOT problems. One such adaptation is the Entropy-regularized Martingale Optimal Transport (EMOT) algorithm\footnote{introduced in \cite{EMOT:2024}}. Both approaches are discussed in Section~\ref{Sinkhorn-vs-Combinatorial-Algorithmic-Approaches}, and we provide also a benchmark in Section~\ref{the-bachelier-problem}. 

\subsection{Two kernel-based approximations for conditional expectations and densities}
\label{ACTKBE}

\paragraph{Nadaraya--Watson estimator.} 

The Nadaraya--Watson (NW) method, originating in the 1960s, is a classical technique for estimating conditional expectations by \textit{kernel-weighted local averaging}. 

The classical way to introduce NW estimators is to first give a formal representation of a conditional expectation as an integral: consider two random variables $Y,X$. Assume and denote $p(x)$ (respectively  $p(x,y)$), the density of the law $X$ (respectively  joint law $(Y,X$)) with respect to the Lebesgue measure $dx$ (respectively  $dxdy$). Consider the following definition of the conditional expectation of $Y$, knowing $X=x$:
\be\label{CONDYX}
\mathbb{E}[Y \mid X = x] = \int y p(y|x) dy, \quad p(y|x) = \frac{p(x,y)}{p(x)}.
\ee
The term $p(y|x)$ defines the conditional density of $y$ knowing $x$, calculated from $p(x)$, the density of $x$, and $p(x,y)$ the density of the joint law. 

From a numerical point of view, $Y,X$ are known from finite distributions, usually from data pairs $(y^n, x^n)$ issued from the joint law $(Y,X)$. The NW estimator of conditional expectation $\mathbb{E}[Y \mid X = x]$, at a query point \( x \in \mathbb{R}^D \), is the discrete analog of the formula \eqref{CONDYX}:
\be\label{rhonw}
\mathbb{E}_{NW}[Y \mid X = x] = \frac{\sum_n k(x,x^n)y^n}{\sum_n k(x,x^n)}, 
\ee
in which the conditional density is approximated by 
\be
p(y^n|x) \approx \frac{k(x,x^n)}{\sum_n k(x,x^n)}, 
\qquad \quad p(x) \approx \sum_n k(x,x^n).
\ee
The notation $\approx$ means up to a normalizing constant. 
These approximations are called \textit{kernel density estimate (KDE)}.  KDE is a density model, which considers a kernel $k$ (respectively  two kernels $k=(k_1,k_2)$), a distribution $X$ (respectively  joint distribution $(X,Y)$), to approximate a density function $\rho(\cdot)$, denoted $\rho_{k,X}(\cdot)$ (respectively  $\rho_{k,(X,Y)}(\cdot)$), computed from $X$, as follows:
\be \label{nw_kde}
\rho_{k,X}(x) \approx \sum_{n=1}^{N} k(x, x^n), \quad \rho_{k,(X,Y)}(x,y) \approx \sum_{n=1}^{N} k_1(x, x^n)k_2(y, y^n).
\ee
For KDEs, the kernel \( k \) does not need to be positive-definite or associated with RKHS. Indeed, it needs to be non-negative and localized to ensure meaningful weighting. The model simply assigns higher weights to the points \( x^n \) that are closer to the query point \( x \). This results in a \textit{local estimator} which is fast to compute and particularly effective in low-dimensional settings. For the kde expression \eqref{nw_kde} to be consistent with the conditional expectation definition \eqref{CONDYX}, it is usually assumed that the kernel $k_2$ satisfies $\int y k_2(y,y^n)dy = y^n$, advocating for the use of translation invariant kernels $k(x,y) = \varphi(x-y)$. Alternatively, whenever possible, the Nadaraya-Watson estimator should be modified as follows: 
\be\label{NWC}
\mathbb{E}_{NW}[Y \mid X = x] = \frac{\sum_n k_1(x,x^n)\overline{y^n}}{\sum_n k_1(x,x^n)}, \quad \overline{y^n} = \int y k_2(y,y^n)dy.
\ee
To conclude, the NW estimator is a useful tool to estimate conditional density. However, there is no easy way to determine a conditional distribution from its conditional density.\footnote{The rejection algorithm is a natural, but impractical tool to that task.}

\paragraph{Kernel-ridge estimator.} 

In contrast to the local averaging of Nadaraya--Watson \eqref{rhonw}, the kernel regression formula \eqref{FIT} approximates the conditional expectations as follows:
\be \label{rhok}
\mathbb{E}_k[Y \mid X = x] =  \sum_{n=1}^{N_x} \psi_{k,X}^n(x) y^n = Y_k(x),
\ee
where $\psi_{k,X}^n(x) = \mathcal{P}_k(x, X)^n$ are the unity functions defined at \eqref{PU}, which play the role of the kernel function $k(x,x^n)$ in the NW estimator \eqref{rhonw}. However, the summed term $\sum_n \psi_{k,X}^n(x)$ approximates the density of the Lebesgue measure, that is, the constant function one. The density model induced by the kernel ridge regression is the Lebesgue one, adapted to the integration formula.

To better interpret the formula \eqref{rhok}, let us use the definition of a transport \eqref{transport}, resulting in the following version of the conditional expectation \eqref{CONDYX}:
\be\label{transportc}
    \mathbb{E}[Y \mid X = x] = \int y p(y|x) dy = \int T(x,y) dy, \quad T_{\#} [p(y|x) dy] = dy.
\ee
We can now interpret \eqref{rhok} as an approximation of the term $\int T(x,y) dy$, relying on an approximation $T_k$ of the transport map $T$, rather than relying on an estimation of the density $p(y|x)$ as Nadaraya-Watson does. 

Consider the following composite kernel $k(z,z') = k_1(x, x')k_2(y, y')$, as in \eqref{nw_kde}. The map $T_k$ is  determined everywhere as follows:
\be
T_k(x,y) = \mathcal{P}_k(z,Z)Y, \quad z=(x,y), Z=(X,Y).
\ee
This map can be used efficiently to generate samples of a conditioned distribution. 

However, with this generative model, we face the opposite problem than the Nadaraya-Watson model: How can we deduce conditioned density from conditioned distributions? A first method to evaluate the density $p(y|x) dy$ induced by the map $T_k$ is to generate the distribution $\big\{T_k(x,v^n)\big\}_{n=1}^N$, $v^n$ drawn from the uniform distribution and to evaluate the resulting density using the NW estimator. Alternatively, we can estimate directly using the change of variable formula \eqref{COVF}:
\be 
  p(y | x) \approx |\det \nabla T_k(x,y)|,
\ee
$\nabla T_k$ being computed by the gradient formula \eqref{GRAD}. This formula is valid for $D_x=D_y=D$, and assumes strong assumptions, namely $T_k$ is a smooth, invertible map. This observation is at the heart of Section~\ref{Maps-and-Generative-methods}, which uses optimal transport techniques to build such maps.

\paragraph{Illustrative example: estimation of conditional expectations.}

We evaluate the performance of three nonparametric estimators for conditional expectations \( \mathbb{E}[Y \mid X = x] \) in figure ~\ref{fig:Condex}. To assess the ability to capture nonlinear patterns while estimating conditional expectations, this test considers synthetic data from a heteroskedastic model, involving a uniform law $X \sim \mathcal{U}(-1, 1)$, and a Gaussian one for $Y$ with data
\be
Y \mid X = x \sim \mathcal{N}(\mu(x), \sigma^2(x)),  \quad \sigma(x) = 0.1 \cdot \cos\left(\frac{\pi x}{2}\right), \quad \mu(x) =  \sin\left(\pi x\right).
\ee
Figure~\ref{fig:Condex} compares the estimated conditional expectations produced by each method with the ground truth $\mu(x)$ over a dense grid, which is the black curve.

 The green curve is a kernel ridge expectation estimate that uses the standard CodPy kernel, a Mat\'ern kernel with the standard map \eqref{standardmap}, which is, as discussed, an approximation of a constant Lebesgue density. Two estimators are built considering the Nadaraya--Watson method with different kernels: the first, in blue, uses a Gaussian kernel with fixed bandwidth, requiring a careful calibration of the bandwidth parameter for each $x$, and is translation invariant. The second, in red, uses the standard kernel, which is not translation invariant, plot to emphasize the importance of this assumption without the correction \eqref{NWC}.  

\begin{figure}[h!]
\centering
\includegraphics[width=.8\textwidth]{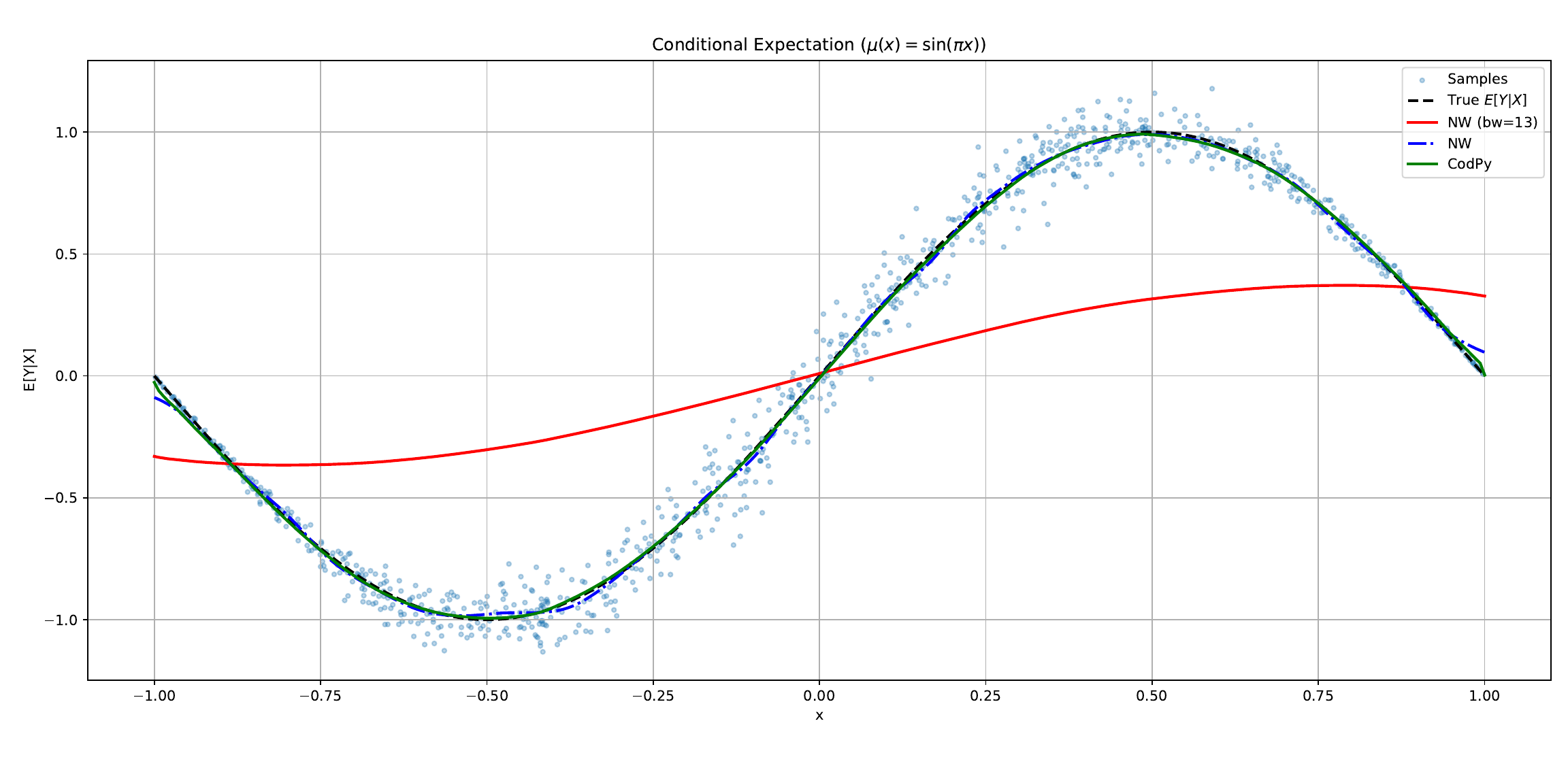}
\caption{Different estimates of conditional expectation}
\label{fig:Condex}
\end{figure}

\paragraph{Illustrative example: different models of conditional density.}

We illustrate in this test the difference in density models of both the kernel ridge and the Nadaraya-Watson approaches. We consider the following distribution, similar to the previous test.  
\be
Y \mid X = x \sim \mathcal{N}(\mu(x), \sigma^2(x)),  \quad \sigma(x) = 0.1 \cdot \cos\left(\frac{\pi x}{2}\right), \quad \mu(x) =  0.1 \cdot \cos\left(\frac{\pi x}{2}\right).
\ee
We sample this distribution, the resulting samples $\{(x^n, y^n)\}_{n=1}^{N}$ are plotted in blue
 in Figure~\ref{fig:DistribAndDensity}, 
 together with a set in red. 
We plot the conditional density $\rho(y \mid X = 0)$ in black in Figure~\ref{fig:DistribAndDensity}, 
 which is the reference one for this test.

We plot in green the density model used by the projection operator, which is, as discussed earlier, an approximation of the Lebesgue measure. 

We plot in red and blue two Nadaraya-Watson estimators. The first one,  in red, using the standard default kernel which is not translation invariant. Observe that this kernel nonetheless captures the correct conditioned distribution, although it does not consider the modified estimator \eqref{NWC}. The reason is that we carefully chose the distribution. The second, in blue, is a Nadaraya-Watson estimator with a translation invariant kernel, which bandwidth has been fit manually. These two curves show that the Nadaraya-Watson method can capture hetero-skedasticity and nonlinear dependencies. 

\begin{figure}[h!]
\centering
\begin{minipage}{.39\textwidth}
\centering
\includegraphics[width=\textwidth]{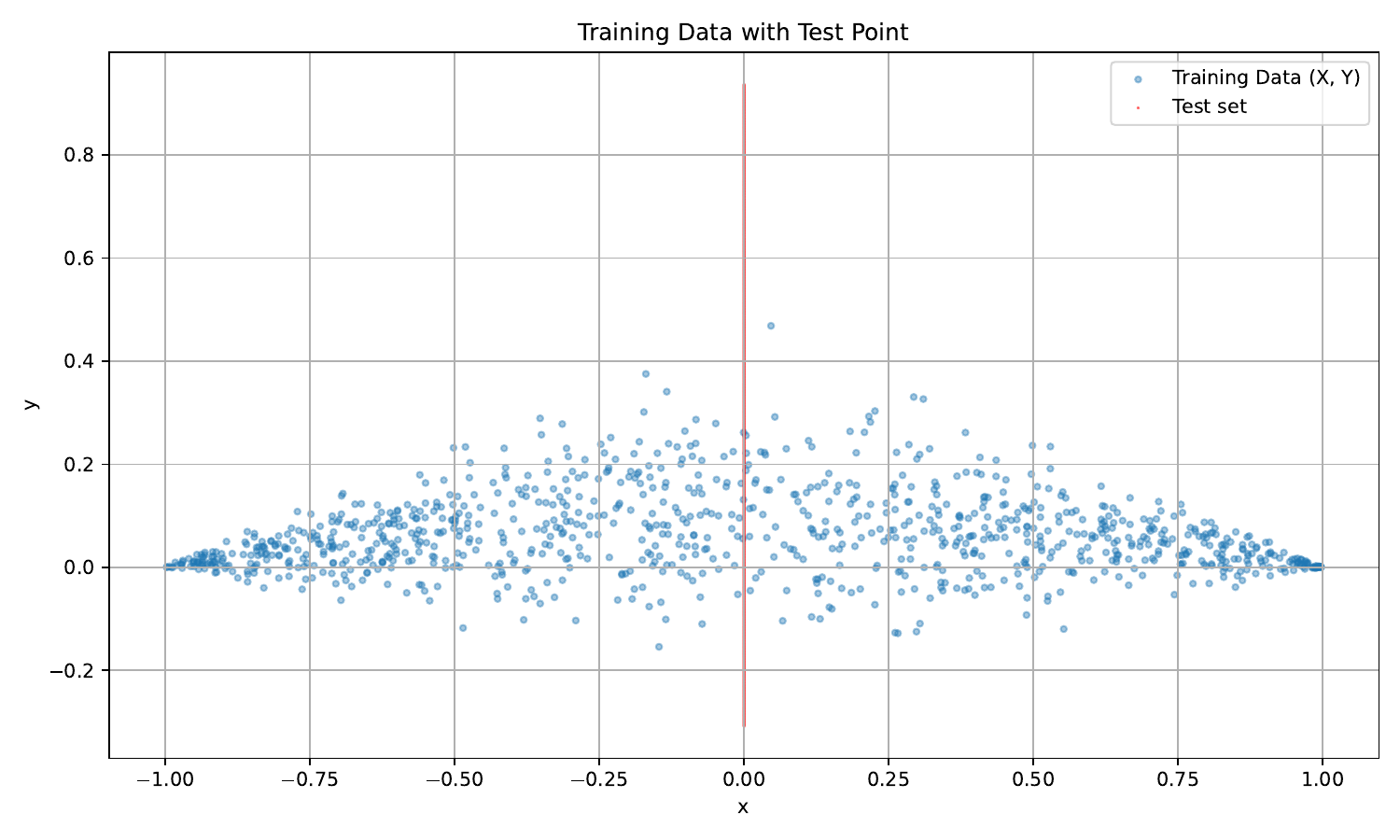}
\caption*{(a) A joint law (blue) with an observation set (red)}
\end{minipage}%
\hfill
\begin{minipage}{.59\textwidth}
\centering
\includegraphics[width=\textwidth]{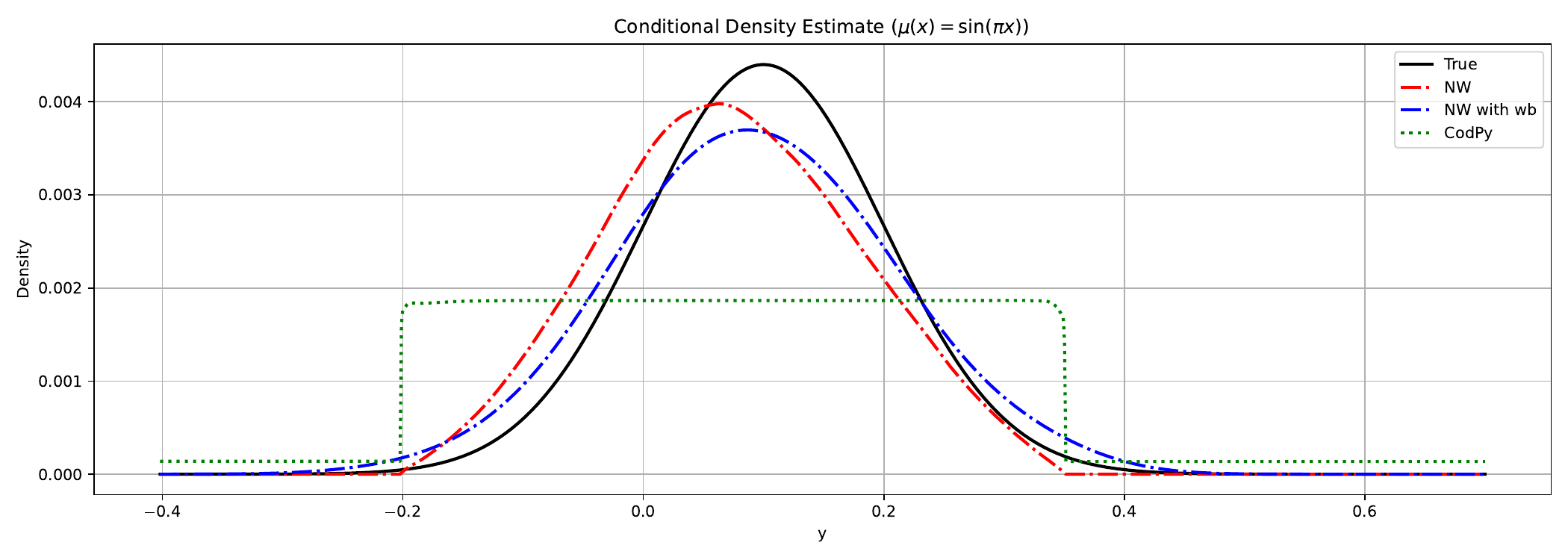}
\caption*{(b) Different density models}
\end{minipage}
\caption{(a) A joint law (blue),with an observation set (red). 
(b) Different density models.}
\label{fig:DistribAndDensity}
\end{figure}

\subsection{Transition probabilities with kernels}
\label{Transition-Probabilities-With-Kernels}

\paragraph{‘Purpose.}

We now present a method to solve the Martingale Optimal Transport problem and display a numerical illustration with the Bachelier problem. 
For these problems, the inputs consist of two distributions $X,Y$, having equal sizes $N_x=N_y = N$, living in the same space $X,Y \in \mathbb{R}^D$, drawn from two laws satisfying the Strassen theorem assumption $dX \ll dY$. Since this problem looks for martingale transfer plans, we can consider w.l.o.g. that both distributions $X$,$Y$ have null means.

To tackle this problem, we observe first that its formulation~\eqref{KantorovichM} is a perturbation of the discrete Kantorovich problem \eqref{KantorovichD}, this last problem being equivalent to the discrete Monge problem \eqref{LSAP}. This motivates the following approach. 
\begin{itemize}
\item Consider two distributions $X,Y$, w.l.o.g. having null mean $\sum_n x^n = \sum_n y^n = 0$, eventually removing their respective means.
\item Compute the permutation $\sigma$ solving the Monge Optimal Transport \eqref{LSAP}, for instance using the LSAP, and relabel $ Y \leftarrow Y^\sigma$.
\item We now build a continuous path of bi-stochastic matrix $\Pi_{t\ge 0} \in \Gamma$, with  $\Pi_0 = I_d$, minimizing the Frobenius norm functional $J(\Pi)=\|X - \Pi_t Y \|_{\ell^2}^2$. The associated continuous-time gradient flow of $J(\Pi)$ is given by $\frac{d}{dt} \Pi_t = (X - \Pi_t Y)Y^T$. A first-order integration yields the approximation \(\Pi_t = I_d + t(X - Y)Y^T\). The matrix \(\Pi_{t}\) remains bistochastic since \(\Pi_t \ 1 = 1 + t(X - Y)Y^T \ 1 = 1\) due to the null-mean assumption. Since our goal is to find
\(\overline{t} = \mathop{\arg\inf}_t J(\Pi_t)\), elementary computations provide the minimum as 
\be \label{PI1}
\overline{t} = \frac{\|(X - Y)Y^T \|_{\ell^2}^2}{\|(X - Y)Y^TY\|_{\ell^2}^2}
\ee 
\end{itemize}
This combinatorial approach to the martingale optimal transport problems \eqref{KantorovichM} is presented as the fixed-point type algorithm \eqref{sumalg}, although the analysis suggests a direct formulation. Indeed, we observed that in some situations a fixed-point approach is more stable, for instance if the assumption $dX \ll dY$ is not fulfilled. It provides a fast and accurate alternative to approximate a stochastic matrix,
although the output matrix might not be positive if $dX \ll dY$, which is an interesting special case. In any case, we can project the resulting matrix on the Birkhoff polytope.

\begin{algorithm}
\begin{algorithmic}[1] \label{sumalg}
\item[\textit{Set-up:}] Kernel $k$, tolerance $\epsilon >0$ or maximum iteration number $M$.
\REQUIRE $X$, $Y$, two distributions of points of equal length.
\ENSURE $\Pi(X, Y)$ a bi-stochastic approximation of the martingale optimal problem \eqref{KantorovichM}.
\STATE Remove the means : $X \leftarrow X - \mathbb{E}(X)$, $Y \leftarrow Y - \mathbb{E}(Y)$.

\STATE Renumber $Y$: $Y \leftarrow Y\circ \sigma$, the permutation $\sigma$ being computed with the LSAP  \eqref{LSAP}, considered with the MMD cost function $d_k(X,Y)$.
\STATE $Y^0=Y$, $\Pi^0 = I_N$, the identity matrix.
\WHILE{ $\|Y^n-Y^{n-1}\| > \epsilon$ or $l < M$} 
\STATE compute $\Pi^{n+1}=I_d+t(X - Y)Y^T$, where $t$ is computed using \eqref{PI1} with $X,Y^n$.
\STATE compute $Y^{n+1} = \Pi^{n+1}Y^{n}$.
\ENDWHILE
\STATE return $\Pi^{n+1}\circ \sigma^{-1}$. 
\end{algorithmic}
\end{algorithm}

\hypertarget{the-bachelier-problem}{%
\paragraph{Numerical illustration with the Bachelier problem}\label{the-bachelier-problem}}

We present the so-called Bachelier test, which is a pertinent test for statistical applications\footnote{particularly for mathematical finance applications, with the following vocabulary correspondences:  $b_t$ are basket values, and $t_2$ is the maturity of a basket option having payoff $P(x) = \max(b(x) - K, 0)$, $K$ being the strike.}, in order to test the transition probability algorithms \eqref{sumalg}. 

We benchmark this algorithm to other algorithms, namely the Sinkhorn method, and the Nadaraya Watson method, both providing alternative approaches to compute conditional expectations. The test is described as follows.
\begin{itemize}
\item Consider a Brownian motion $t \mapsto X_t \in \mathbb{R}^D$, satisfying $dX_t = \sigma dW_t$,
where the matrix $\sigma \in \mathbb{R}^{D,D}$ is randomly generated. The initial condition is $X_0 = 0$ w.l.o.g. Let $\omega \in \mathbb{R}^D$, randomly generated, satisfying $|\omega|_1 = 1$ and denote the scalar values $b_t = <\omega,X_t>$. This last process follows a univariate Brownian motion $d b_t = \theta dW_t$. We normalize $\sigma$ in order to retrieve a constant value for $\theta$, fixed to 0.2 in our tests.

\item Consider two times $1 = t_1 < t_2 = 2$, and a function denoted $P(x) = \max(b(x) - K, 0)$. The goal of this test is to benchmark some methods aiming to compute the conditional expectation $\mathbb{E}^{X_{t_2}}[P(\cdot) \mid X_{t_1}]$, for which the reference value is given by the so-called Bachelier formula
\be\label{BF}
  f(\cdot) = \mathbb{E}^{X_{t_2}}[P(\cdot) \mid X_{t_1}] = \theta \sqrt{t_2-t_1} \ p(d) + (b_{t_1}-K)c(d),\quad d=\frac{b_{t_1}-K}{\theta\sqrt{t_2-t_1}},
\ee
where $p$ (respectively  $c$) holds for the cumulative (respectively  density) of the normal law.

\item Considering two integer parameters $N,D$, in this test we consider three samples of the Brownian motion $X \sim X_{t_1},Y \sim X_{t_2}, Z\sim X_{t_1}$ in $\mathbb{R}^{N,D}$, and produce a transition probability to approximate the transition probability matrix $\Pi(Y | X)$. We finally output the mean square error $\text{err}(Z) = \|f(Z)- \Pi(Y | Z) P(Y) \|_{\ell^2}$.
\end{itemize}

Figure~\ref{fig:BacPerfs} (respectively  figure\ref{fig:ottime-deux}) shows the
score
\(\frac{\text{err}(Z)}{\|f(Z)\|_{\ell^2}+\|\Pi(X_{t_2} | Z) P(X_{t_2})\|_{\ell^2}}\)
for various choice of \(N\) and \(D\) (respectively  the execution time in
seconds) for four methods.
\begin{itemize}
\item COT illustrates our results with Algorithm~\ref{sumalg}.
\item OTT consists in our best trial using Sinkhorn algorithm with the OTT library.
\item Ref is a reference, naive value, computing $P_{k,\theta}(Z)$ with the interpolation formula \eqref{FIT}.
\item Nadaraya Watson implements \eqref{rhonw}.
\end{itemize}

All these methods (labelled with an index \(m\)) output a transition
probability matrix \(\Pi^m(Y|X)\), which implies the estimation
\(f^m(X) = \Pi(Y|X) P(Y)\). We used the extrapolation method \eqref{FIT}
in order to extrapolate \(f^m_{k,\theta}(Z)\) for all methods. In this
test, we observed that both the Sinkhorn algorithm and Algorrthm~\ref{sumalg}
performed well across a wide range of parameters \(N\) and \(D\). Both
approaches reliably computed the transition probability matrix. Notably,
the naive extrapolation method exhibited surprisingly good performance
in the high-dimensional case (\(D=100\)). This unexpected result raises
questions about the relevance of this test in high dimensions, where a
simple linear regression method could achieve similar accuracy when
approximating the Bachelier formula \eqref{BF}.

\begin{figure}
\hypertarget{fig:BacPerfs}{%
\centering
\includegraphics[width=1.\textwidth, keepaspectratio]{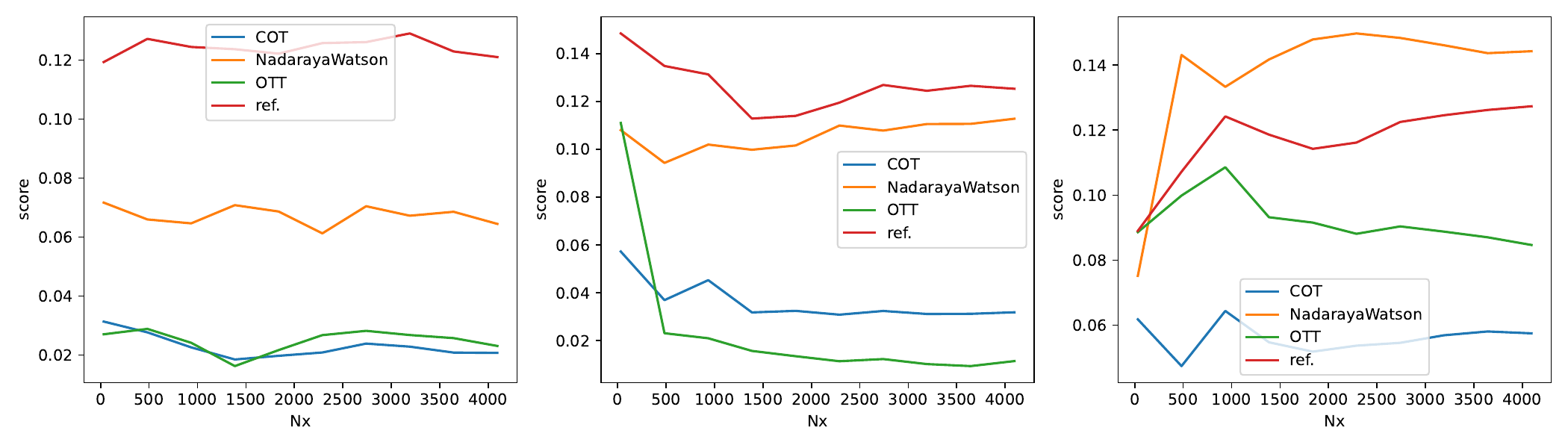}
\caption{Benchmark of scores}\label{fig:BacPerfs}
}
\end{figure}

\begin{figure}
\hypertarget{fig:ottime-deux}{%
\centering
\includegraphics[width=1.\textwidth, keepaspectratio]{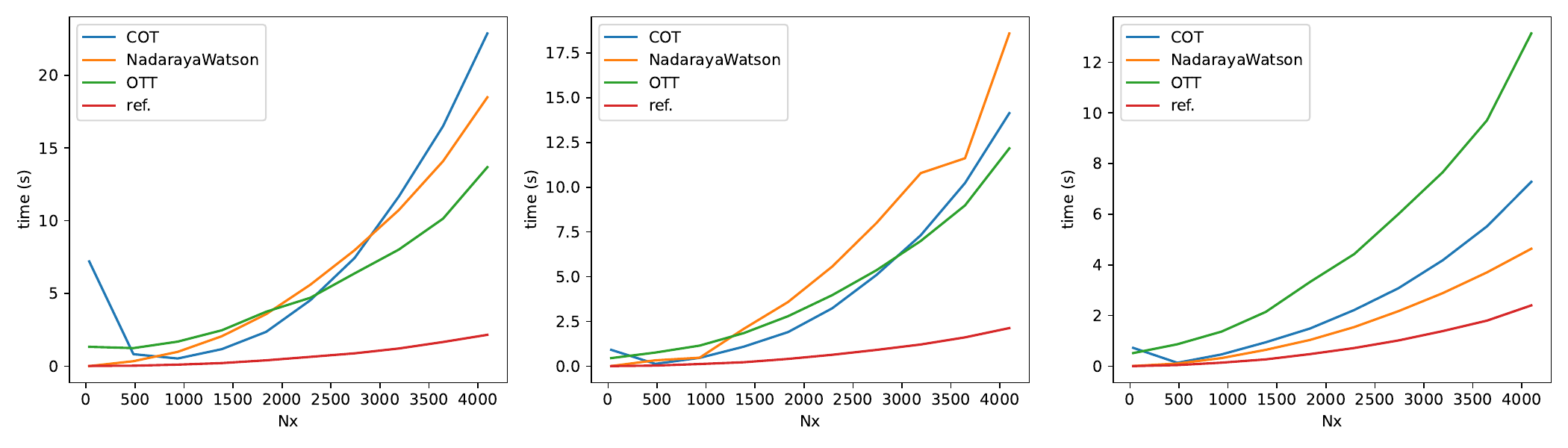}
\caption{Execution time}\label{fig:ottime-deux}
}
\end{figure}


\section{Maps and generative methods: dealing with two distributions}
\label{Maps-and-Generative-methods}

In many applications, it is necessary to define mappings between two
distributions, say \(X\) and \(Y\), with support in
\(\mathbb{R}^{D_x}\) and \(\mathbb{R}^{D_y}\), respectively. These
distributions are typically only known through two discrete distributions of
equal length, say \(X \in \RR^{N,D_x}\) and \(Y \in \RR^{N,D_y}\). We focus now on the construction of such mappings with RKHS methods, which, as discussed in the previous Section~\ref{Conditional-Expectations-and-Densities-Transition Probabilities}, use optimal transport techniques to properly define these mappings.

From a discrete point of view, consider any two distributions of equal length $X \in \mathbb{R}^{N,D_x}$, $Y \in \mathbb{R}^{N,D_y}$. The previous section introduced in \eqref{rhok} the following mapping, which is the projection formula \eqref{FIT} for the reproducible, extrapolation mode (with \(\epsilon=0\)):
\be 
\mathbb{E}_k[Y \mid X = x] =  \mathcal{P}_{k}(x,X) Y = Y_{k,\theta}(x).
\ee
 The corresponding section concluded that optimal transport is needed for this formula to define a one-to-one map, which we discuss now.

This extrapolation mode is the reproducible one,  satisfying here $Y_{k, \theta}(X) = Y$, thus
the mapping $x \mapsto Y_{x,\theta}(x)$ realizes a push-forward, defined at Section~\ref{reminder-on-optimal-transport-theory}, of the discrete distribution \(\delta_X\) to
\(\delta_Y\), where \(\delta_X  =  \frac{1}{N}\sum_n \delta_{x^n}\),
\(\delta_{x}\), \(\delta_Y  =  \frac{1}{N}\sum_n \delta_{y^n}\),
\(\delta_{z}\) being the Dirac measure centered at \(z\), 
that is, we have $Y_{k, \theta \ \#}$ 
$\delta_X = \delta_Y$. However, in this construction, it is important
to observe, that all index permutations
\(\sigma: [1, \ldots, N]  \to [1, \ldots, N]\) of the set \(Y\), denoted $Y \circ \sigma$, would also define
push-forward maps in the sense $(Y\circ \sigma)_{k, \theta \ \#}\delta_X = \delta_Y$. Thus, among all these permutations, the strategy is
to select those permutations which define smooth, invertible transport maps, taking the following form
\be 
\label{ED}
  Y^\sigma_{k,\theta}(\cdot)  =  K(\cdot,X) \theta, \quad \theta = K(X,X)^{-1} (Y\circ \sigma). 
\ee 
To summarize, from the discrete standpoint, finding a smooth and invertible transport map from \(X\) to \(Y\) requires first solving a numerical problem that computes a relevant permutation of \(Y\).

Permutations determining invertible, one-to-one mappings, were discussed in Section~\ref{reminder-on-optimal-transport-theory} (on optimal transport), where a distinction is made according to the compatibility of metric spaces. 
\begin{itemize}
\item
Compatible metric space \(D_x=D_y=D\).  This case,  discussed in the dedicated Section~\ref{Discrete-Optimal-Transport-on-Compatible-Spaces}, compute the permutation as the solution of the Monge problem \eqref{LSAP}, considered with a distance typically defined with the discrepancy distance, say 
\be
   \mathop{\arg\min}_{\sigma \in \Sigma} \text{Tr}\big(M_k(X,Y\circ \sigma)\big). 
\ee 

\item  Not compatible metric spaces \(D_x \neq D_y\). This case, discussed in the corresponding Section~\ref{Discrete-Optimal-Transport-on-Incompatible-Spaces}, computes the permutation as the solution of the Gromov--Monge problem \eqref{GM-discrete}, of the form
\be\label{GMD}
 \mathop{\arg\min}_{\sigma \in \Sigma} \sum_{i,j=1}^N |d_{k_X}(x^i, x^j) - d_{k_Y}(y^{\sigma(i)}, y^{\sigma(j)})|^2.
\ee
In the previous equation, $d_{k_X}(x,x'),d_{k_Y}(y,y')$ are distances defined on each metric space, as for example, the kernel discrepancy \eqref{dk}. 
\end{itemize}

We emphasize that \eqref{GMD} is a numerical efficient approach for the incompatible metric space $D_x \neq D_y$, but it is not the only possible choice. For instance, we designed and tested another interesting alternative approach:

\be 
\label{TC}
  \mathop{\arg\min}_{\sigma \in \Sigma}\|\nabla \Big( Y_{k, \theta}^\sigma(X) \Big) \|_2^2. 
\ee
The term $\nabla Y_{k, \theta}^\sigma(\cdot)$, computed using the gradient approximation \eqref{nabla}, is a matrix of size $D_x,D_y$, thus the norm $\|\cdot\|_2$ holds here for the Frobenius norm for matrix. This approach, reminiscent of the \textit{traveling salesman problem}, fit into the Gromov--Monge formulation framework too. The map is computed similarly to the polar factorization \eqref{PF}-\ref{PFN}, but in the incompatible metric space case. The functional \eqref{TC} gives results quite similar to \eqref{GMD}, although slightly more computationally involved and less stable.

\paragraph{Encoder-decoders and sampling
algorithms}\label{Encoder-decoders-and-sampling-algorithms}

We now describe an algorithm related to this strategy, which has proven useful for our applications. As mentioned earlier, we consider two cases, according whether metric spaces are compatibles or not:
\begin{algorithm}
\begin{algorithmic}[1]
\label{alg1}
\REQUIRE data $X \in \mathbb{R}^{N,D_x}$, $Y \in \mathbb{R}^{N,D_y}$, with $x^i\neq x^j$, $y^i\neq y^j$,  $\forall \ i \neq j$.
\ENSURE A regressor $Y_{k, \theta}(\cdot)=\mathcal{P}_k(\cdot,X) Y^{\sigma}$, modeling an invertible push-forward map, where the permutation $\sigma$ is computed as follows.
\IF{$D_x=D_y$}
  \STATE Compute a permutation $\sigma$ using the LSAP \eqref{LSAP} with cost function $c(i,j)=d_k(x^i,y^j)$
\ELSE
  \STATE Denote  
$$
    s(\sigma) = \sum_{i,j}|d_{k_X}(x^i, x^j) - d_{k_Y}(y^{\sigma(i)}, y^{\sigma(j)})|^2,
$$ 
    corresponding to the Gromov--Monge functional  \eqref{GMD} (or $s(\sigma) = \|\nabla \Big( Y^\sigma_{k, \theta}(X) \Big) \|_2^2$ for \eqref{TC}). 
  \STATE Compute a permutation $\sigma$ using the discrete descent Algorithm~\ref{PBA}.
\ENDIF
\end{algorithmic}
\end{algorithm}

This algorithm is used primarily to generate samples from a discrete distribution \(Y \in \mathbb{R}^{N, D_y}\), assuming that this distribution comes from an i.i.d. sampling of a law having continuous density. The distribution $X$ can be input if the law $Y$ comes from a joint law $(X,Y)$, or be any arbitrary distribution that the user judges pertinent, as for instance drawn from a standard law
\(\mu\) (e.g., a normal distribution). Whatever the distribution \(X \in \mathbb{R}^{N, D_x}\), this procedure computes a regressor
\(Y^\sigma_{k, \theta}(\cdot)\). The map
\(z \mapsto Y^\sigma_k(z)\) provides a generator producing new samples of
\(Y\). This process resembles generative methods, such as those
used in GAN architectures. In this context, the machine learning
community uses specific terminology, for which we determined the following correspondence.
\begin{itemize}
    \item Latent space: The space $\mathbb{R}^{D_x}$ is called the latent space. 
    \item Decoding: The map $x \mapsto (Y \circ \sigma)_{k, \theta}(x)$ is called a decoder.
    \item Encoding: The map $y \mapsto (X \circ \sigma^{-1})_{k, \theta}(y)$, with $\theta = K(Y \circ \sigma, Y \circ \sigma)^{-1} X$, is called an encoder.
    \item Reconstruction: The map $y \mapsto (Y \circ \sigma)_{k, \theta}(X \circ \sigma^{-1})_{k, \theta}(y)$, is called a reconstruction.    
\end{itemize}

The latent space is typically thought of as a lower-dimensional space
that captures the essential features of the data. This provides a
compact and informative representation of the original data \(Y\), allowing efficient encoding and decoding of information while taking advantage of the structural properties of the data.

The choice of the latent space \(\mathbb{R}^{D_x}\) is important for
applications, although it is not well documented. When using a standard
normal distribution and a small latent dimension (e.g. \(D_x = 1\)),
the resulting generator \(x \mapsto Y ^\sigma_{k, \theta}(x)\)
produces samples that closely resemble the original sample \(Y\), which
can pass standard statistical tests such as the Kolmogorov-Smirnov test.
Using larger dimensions (e.g. \(D_x \gg 1\)) results in greater
variability in the samples but may cause the model to fail statistical
tests as the dimension increases. 

One of the main advantages of these kernel generative methods is to
produce continuous generators that can exactly reconstruct the original
variate \(Y\) if needed, as the reproducibility property \eqref{REPRO}
ensures \(Y^\sigma_{k, \theta}(X) = Y\). We now turn our attention to numerical illustrations, which goals are to illustrate, and motivate, the role of the latent space for the generative Algorithm~\ref{alg1}.

\paragraph{One-dimensional numerical illustrations of Monge transport.} 

We start illustrating the encoding/decoding procedure \eqref{alg1} using a simple interface, which we refer to as the \emph{sampling procedure}. This procedure is designed to generate new samples that approximate the distribution of a given dataset $Y \in \mathbb{R}^{N_y \times D_y}$ by constructing a kernel-based regressor and using a latent representation. We begin with a simple one-dimensional Monge problem to demonstrate the generative capabilities of the model. In this test, we consider two types of target distributions: a bimodal Gaussian distribution and a bimodal Student's $t$-distribution\footnote{The word ``Student'' refers to the statistician W. Sealy Gosset, who published under the pseudonym Student.}. 

The bimodal Gaussian is constructed as a mixture of two normal distributions: $\alpha \mathcal{N}(\mu_1, \sigma^2) + (1 - \alpha)\mathcal{N}(\mu_2, \sigma^2)$, with means $\mu_1 = -2$, $\mu_2 = 2$, common variance $\sigma^2 = 0.5^2$, and mixing coefficient $\alpha = 0.5$. Similarly, the bimodal Student's $t$-distribution is defined as a mixture $\alpha t_\nu(\mu_1, \sigma^2) + (1 - \alpha)t_\nu(\mu_2, \sigma^2)$, where $t_\nu$ denotes the Student's $t$-distribution with $\nu = 3$ degrees of freedom.

For each case, we generate a reference dataset $X \in \mathbb{R}^{1000 \times 1}$ sampled from the corresponding true distribution. The sampling procedure described by the algorithm~\ref{alg1} is then applied to produce a dataset $Y \in \mathbb{R}^{1000 \times 1}$, aimed at approximating the structure of the original distribution.

Figure~\ref{fig:409} presents a comparative visualization of the original and generated distributions. For both cases, we display histograms along with kernel density estimates (KDEs) to illustrate the underlying densities. The left panel corresponds to the bimodal Gaussian case, while the right panel shows results for the bimodal Student's $t$-distribution.

The plots show that the generated samples closely approximate the underlying target distributions. Both the Gaussian mixture and the heavier-tailed $t$-mixture are well reproduced by the kernel-based generative model. This intuition is confirmed by two samples statistical tests as Kolmogorov-Smirnov; see Table~\ref{tab:ch5_6_a}.
\begin{figure}
\centering
\includegraphics[width=0.7\textwidth, keepaspectratio]{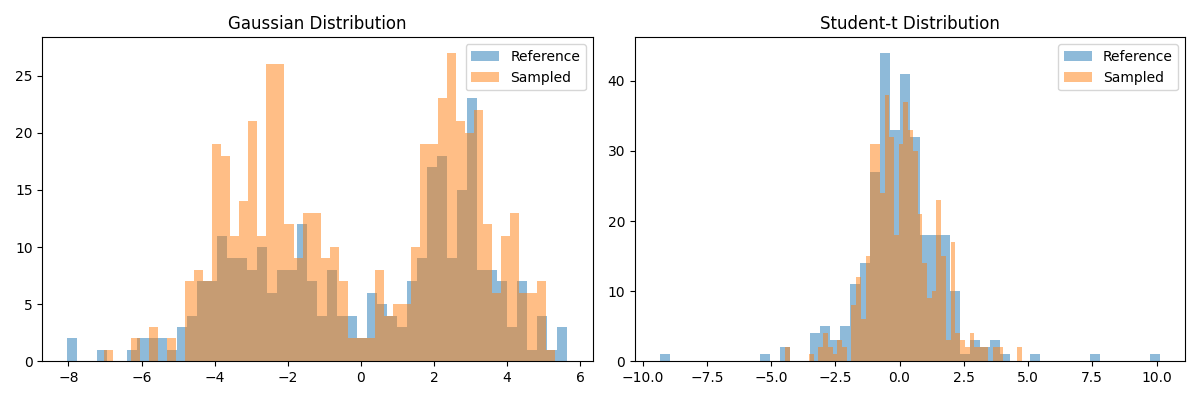}
\caption{\label{fig:409}Comparison of original and sampled one-dimensional distributions. Left: bimodal Gaussian. Right: bimodal Student's $t$-distribution. Histograms and kernel density estimates are shown for both reference and generated datasets.}

\end{figure}

\begin{table}
\caption{Comparison of Gaussian and Student-t distributions using 1D sampling}
\label{tab:ch5_6_a}
\begin{tabular}{llllll}
\toprule
 & Mean & Variance & Skewness & Kurtosis & KS test \\
\midrule
Gaussian:0 & 0.032 (-0.17) & 7.9 (7.6) & 0.052 (0.18) & -1.8 (-1.8) & 0.069 (0.05) \\
Student-t:0 & -0.062 (-0.0084) & 4.8 (4.7) & -0.87 (-0.38) & 4.1 (1.7) & 0.4 (0.05) \\
\bottomrule
\end{tabular}
\end{table}

\paragraph{Two-dimensional numerical illustrations of Monge transport.} 

We repeat the test presented above in one dimension in the Figure~\ref{fig:412} with two-dimensional data, and present scatter plots of the original and generated data for each of the two cases, allowing a visual comparison of the fidelity of the generative process. Usually, two-sample tests as Kolmogorov-Smirnov should deteriorate as the dimension increases; see Tables~\ref{tab:ch5_6_a} and~\ref{tab:ch5_6_b}, which present Kolmogorov-Smirnov tests for each marginal, hence four tests, for both axis and distributions. 

\begin{figure}
\centering
\includegraphics[width=0.7\textwidth, keepaspectratio]{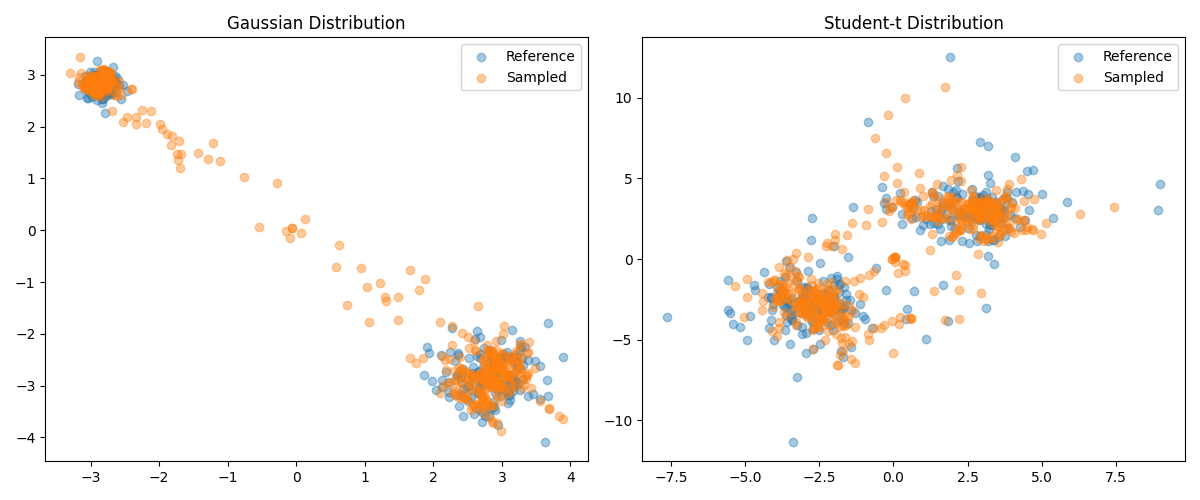}
\caption{\label{fig:412}2D Comparison: bimodal Gaussian (left), and bimodal Student's $t$ (right) versus sampled distributions}
\end{figure}

\begin{table}
\caption{Comparison of Gaussian and Student-t distributions using 1D sampling}
\label{tab:ch5_6_b}
\begin{tabular}{llllll}
\toprule
 & Mean & Variance & Skewness & Kurtosis & KS test \\
\midrule
Gaussian:0 & 0.053 (0.09) & 1 (1.1) & 0.1 (0.17) & -0.58 (-0.71) & 0.89 (0.05) \\
Gaussian:1 & 0.0084 (0.064) & 11 (9.9) & 0.0068 (-0.049) & -1.8 (-1.8) & 0.046 (0.05) \\
Student-t:0 & 0.018 (0.16) & 14 (13) & -0.13 (-0.38) & -0.6 (-0.69) & 0.3 (0.05) \\
Student-t:1 & 0.1 (0.36) & 21 (19) & 0.023 (-0.18) & -1.3 (-1.6) & 0.34 (0.05) \\
\bottomrule
\end{tabular}
\end{table}

\paragraph{High-dimensional numerical illustration of Monge transport and Gromov--Monge transport.}

We now repeat a similar test  Figure~\ref{fig:413} with a bi-modal Gaussian distribution, in fifteen dimensions, comparing Monge and Gromov--Monge methods. The figure plots for each of these two methods the two best and worst axis combinations, according to the Kolmogorov-Smirnov test. As can be seen in the picture, Gromov--Wasserstein-based generative method leads to distributions that are close to the original space, which can pass two samples tests as Kolmogorov Smirnov ones. This property is interesting for industrial applications. 

\begin{figure}
\centering
\includegraphics[width=0.7\textwidth, keepaspectratio]{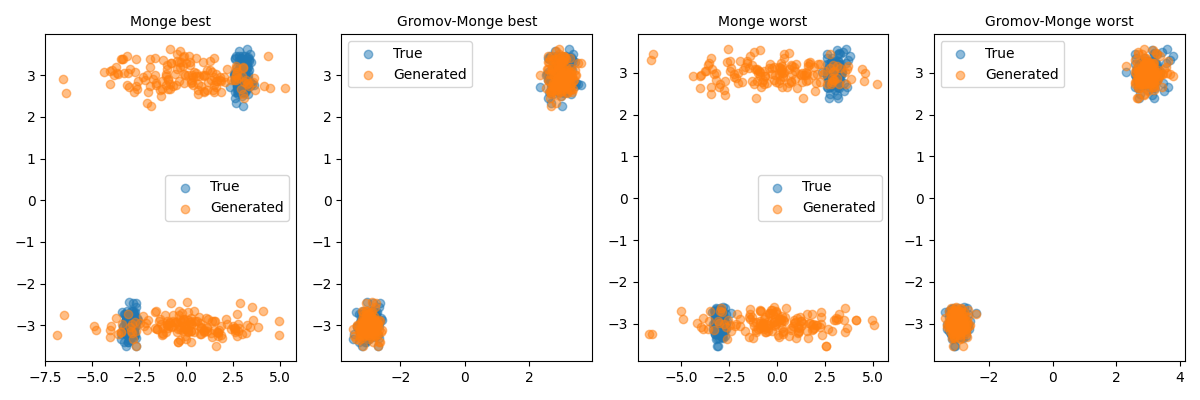}
\caption{\label{fig:413}15 dimensional bimodal Gaussian comparison: Monge transport  (left) and Gromov--Monge transport (right)}
\end{figure}

\paragraph{Gromov--Wasserstein vs. Gromov--Monge for latent parametrization of spherical data}

We investigate the use of optimal transport methods to learn a 1D latent parametrization of a manifold composed of two concentric circular clusters in $\mathbb{R}^2$. Each sample $x^i$ from $X \in \mathbb{R}^{N, 2}$ is generated by $x^i = c^k + r \cdot \frac{\epsilon^i}{\|\epsilon^i\|} + \eta^i$, where $\epsilon^i \sim \mathcal{N}(0, I)$ is a standard Gaussian noise, $\eta^i \sim \mathcal{N}(0, \varepsilon^2 I)$ is small isotropic noise and $c^k$ is the center of cluster $k$ to which $x^i$ belongs.

We align the two-dimensional dataset $X$ to a one-dimensional latent space $Y\subset [0, 1] \in \mathbb{R}^{N,1}$ by determining a permutation $\sigma$ considering a mapping given by the kernel ridge regression $Y_{k,\theta}^\sigma(\cdot) = \mathcal{P}_k(\cdot, X) Y^\sigma$; see Section~\ref{Encoder-decoders-and-sampling-algorithms}. 
 The permutation $\sigma$ is seek minimizing the geometric distortion using two different OT methods.
 \begin{itemize}
     \item Gromov--Wasserstein (GW) approach using the POT Library, where the soft coupling $\overline{\Pi} \in \mathbb{R}^{N \times N}$ solves \eqref{GW-discrete}. Once computed, the permutation is approximated as $\sigma(n) = \mathop{\arg\max}_m \Pi(n,m)$.
      
     \item Gromov--Monge (GM) approach which computes an explicit permutation $\sigma : Y \to X$ by solving the discrete matching problem defined in \eqref{GM-discrete}.
\end{itemize}
 
In both GW and GM settings, the decoder $f : \mathbb{R}^{D_x} \to \mathbb{R}^{D_y}$ is modeled as a smooth kernel-based operator trained on aligned input-output pairs. This decoder enables both faithful reconstruction of the original samples and the generation of new data points by evaluating $f(z)$ for latent codes $z \sim \mathcal{U}(0, 1)$. Label assignment is handled via kernel classification in the data space. Despite the different strategies used to estimate the OT map -- soft coupling in GW versus hard permutation in GM--the remainder of the encoder--decoder pipeline remains identical in structure and implementation.

Figures~\ref{fig:GM-illustration} and~\ref{fig:GW-illustration} illustrate latent encodings, reconstructions, and generated samples. Both models yield sharp latent embeddings, faithful reconstructions, and samples that match the original geometry. This empirical equivalence demonstrates that the proposed permutation-based GM transport is a viable and interpretable alternative to GW in structured generative tasks.

\begin{figure}[ht]
\centering
\includegraphics[width=\textwidth]{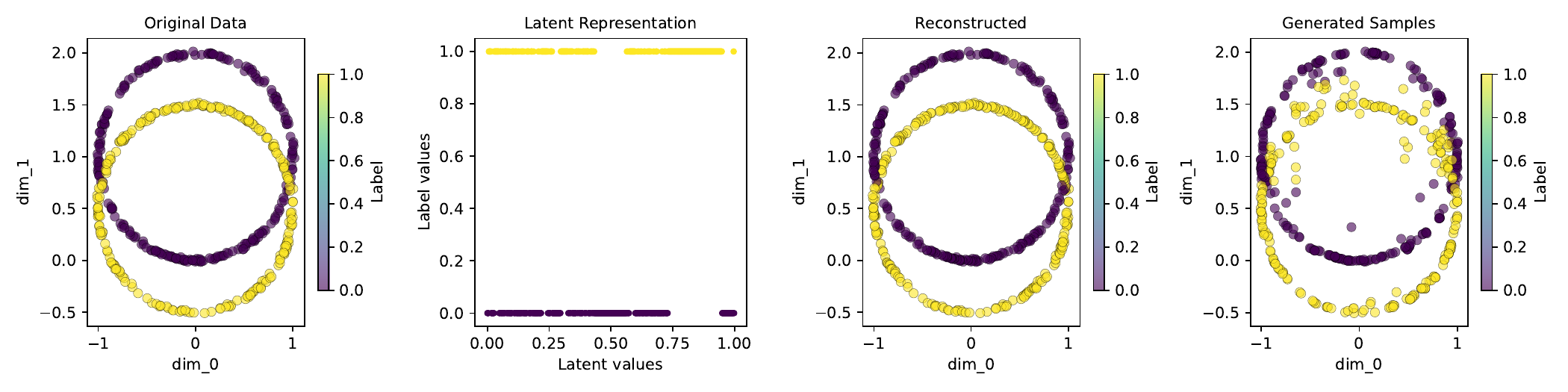}
\caption{\label{fig:GM-illustration}
\textit{Gromov--Monge (GM)}: From left to right -- original 2D data; latent 1D encoding; decoded reconstructions; generated samples. The latent space is binary and structured, consistent with the cluster separation in the data.}
\end{figure}

\begin{figure}[ht]
\centering
\includegraphics[width=\textwidth]{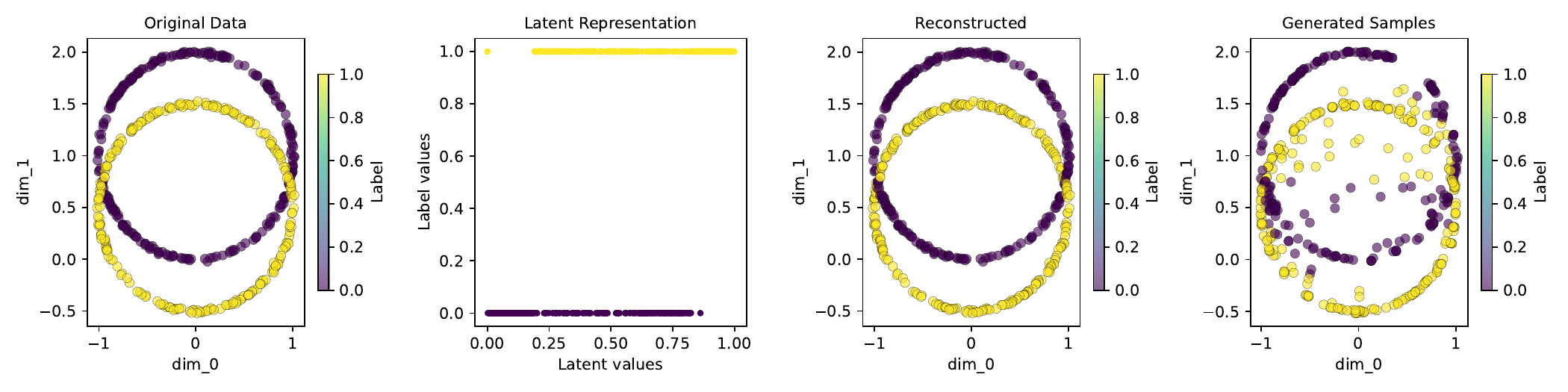}
\caption{\label{fig:GW-illustration}
\textit{Gromov--Wasserstein (GW)}: Same layout as the GM figure. The latent space is more continuous and reveals overlap between clusters. Reconstructions remain accurate, but generated samples are more diffuse.}
\end{figure}

\paragraph{Conditional distribution sampling
model}\label{Conditional-distribution-sampling-model}

The generation of conditioned random variables is an essential part of
the encoder--decoder framework. Given two variates
\(X \in \mathbb{R}^{N, D_x}\) and \(Y \in \mathbb{R}^{N, D_y}\) from two
distributions, our goal is to provide a generator that models the
conditioned distribution \(Y|X=x\).

Considering a latent variable having the form
\(\eta = (\eta_x,\eta_y) \in \mathbb{R}^{N,D_{\eta}}\), with
\(D_{\eta} = D_{\eta_x}+D_{\eta_y}\), the algorithm \eqref{alg1} allows us
to determine these two mappings. 
\begin{itemize}
    \item Encoder (latent variable inference): $x \mapsto (\eta_x\circ \sigma_X)_{k,\theta}(x)$, targeting the latent distribution $\eta_x$.
    \item Decoder (joint distribution generation): $\eta \mapsto ([X,Y] \circ \sigma_{XY})_{k,\theta}(\eta)$ targeting the joint distribution $(X,Y)$.
\end{itemize}

Based on these components, the conditional generator for \(Y|X=x\) is given by
 \be 
\label{condalg}
  \eta_y \mapsto (Y\circ \sigma)_{k,\theta}([\eta_x,\eta_y]), \quad \eta_x = (\eta_x\circ \sigma)_{k,\theta}(x), 
\ee 
where the notation
\(\eta \mapsto (Y \circ \sigma)_{k,\theta}(\eta)\) denotes the $Y$-component of the joint decoder output.  This approach allows us to sample from the conditional distribution by fixing the encoder output $\eta_x$ and varying $\eta_y$.
\begin{itemize}
    \item In some situations, there is no need to define the encoder $\eta_x\circ \sigma_X$, and $x$ is then considered a latent variable.
    \item We can extend this approach to model more elaborate conditioning schemes, such as generating $Y$ conditioned on $X \sim Z$, where $Z$ is a distribution supported in $\mathbb{R}^{D_x}$.

\end{itemize}

\paragraph{Conditioned distributions illustrated with circles}

We now explore some aspects of conditioning on discrete labeled values. In this test, quite similar to the one in Section~\ref{Maps-and-Generative-methods},
%
we consider a low-dimensional feature space \(\mathcal{Y}\) consisting of 2D points \(y = (y_1, y_2)\) lying on two labeled circles \(\{0,1\}\), displayed in Figure~\ref{fig:circlesexample}-(i), with the color code yellow (0), purple (1). Observe that \(\{0,1\}\) are labels in this problem and should not be ordered. Hence we rely on hot encoding, to transform these labels into unordered ones, instead considering conditioning on a two-dimensional label \(x^1=\{1,0\},x^2=\{0,1\}.\). The purpose of this test is to provide a distribution generator $Y |X = x^i$.

To test our mechanism of latent variables, we use a one-dimensional latent variable $\eta_y \in \mathbb{R}$ to encode $Y\in \mathbb{R}^2$. Given a hot-encoded label \(x_i,i=1,2\), we generate samples using the generative conditioned method \eqref{condalg}, hence estimating the conditional distribution \(Y \ |\ X=x_i \).

In doing so, we resample the original distribution, and we test the capability of the generative algorithm to properly identify the conditioned distribution, as well as this choice of latent variable for the kernel generative method.

\begin{figure}
\centering
\includegraphics[width=0.7\textwidth, keepaspectratio]{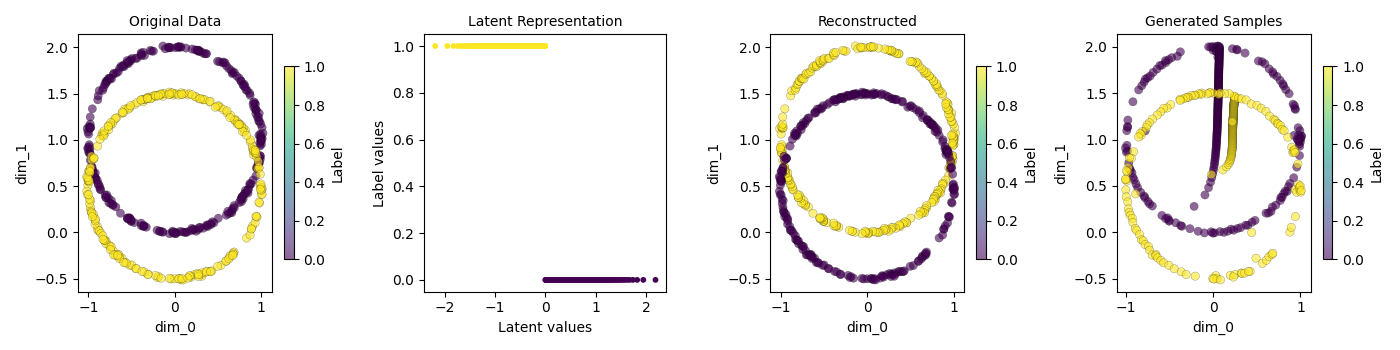}
\caption{\label{fig:unnamed-chunk-98}\label{fig:circlesexample}}
\end{figure}

\paragraph{Inversion of non-invertible mappings.}\label{Inversion-of-non-invertible-mappings.} 

We describe now a general method to properly handle map inversion of non-invertible mappings with the optimal transport techniques described earlier.  Let $y(x) : \mathbb{R}^{D_x}  \to \mathbb{R}^{D_y}$ denote any mapping. There are several situations where we need to compute the inverse mapping $x(y) = y^{-1}(y) : \mathbb{R}^{D_y}  \to \mathbb{R}^{D_x}$, in the non-invertible case, for instance when there might exist multiple solutions to this equation. 
Our main objective at this stage is to perform such inversion in a numerically stable manner.

Let us start discussing this problem from a continuous point of view. Consider $D_x=D_y = D$, where we can use the polar factorization\footnote{see \cite{Brenier:1991,McCann:2001}} (see~\eqref{PF}): we assume and denote $d\mu$, $d\nu$ such that $y_\# d\mu = d\nu$. Provided $y(\cdot) \in L^2(\mathbb{R}^D, d\mu)$, and under certain technical assumptions concerning $d\mu$ and $d\nu$, we can factorize
$y(x) = (\nabla h)(T(x))$, $h$ convex, where $T$ is a permutation, i.e. satisfying $T_\# d \mu = d \mu$. The polar factorization allows us to define the inverse as
\be
 x(y) =  T^{-1} \circ (\nabla h)^{-1}(y),
\ee
which is a safer expression, in the sense that both maps are now considered invertible.
Thus we rely on polar factorizations for inverting, which are more stable to compute, as now both expressions are theoretically invertible. From a discrete point of view,  consider $X=(x^1,\ldots,x^N)$, $Y=(y^1,\ldots,y^N)$, and denote $\delta_X = \frac{1}{N} \sum_n \delta_{x^n}$, $\delta_Y = \frac{1}{N} \sum_n \delta_{y^n}$. Consider the kernel regression~\eqref{FIT} $Y_{k,\theta,X}(\cdot)$, which defines a map transporting $(Y_{k,\theta,X})_\#(\delta_X) = \delta_Y$. The inverse mapping can be defined as $X_{k,\theta,Y}$, satisfying $(X_{k,\theta,Y})_\#(\delta_Y) = \delta_X$, although this inversion is likely to be unstable. The map $X^\sigma_{k,\theta,Y}(\cdot)$, where $\sigma$ is the permutation appearing in~\eqref{ED}, corresponds to the \textit{smooth} part of the polar factorization $\nabla h$. So instead, it is more stable to perform the following computation
\be \label{SI}
 x(y) \sim X_{k,\theta,X^\sigma}(X^\sigma_{k,\theta,Y}(y)).
\ee


\clearpage

\part{Application to machine learning, PDEs, and statistics}

\chapter{Application to machine learning: supervised, unsupervised, and generative methods}\label{overview-on-methods-of-machine-learning}

\section{Purpose of this chapter}

This chapter presents a series of numerical tests that illustrate the application of RKHS methods to a variety of machine learning tasks. Our goal is to evaluate and benchmark kernel-based approaches against standard learning models in supervised and unsupervised settings.

The tests cover regression, classification, clustering, and generative modeling, providing insight into the strengths and limitations of each method under realistic conditions. Special attention is paid to reproducibility and interpretability, with all methods evaluated using standard configurations without hyperparameter tuning.

The following section introduces the core learning paradigms and performance metrics used throughout the chapter, serving as a gentle primer for readers less familiar with classical machine learning evaluation procedures.


\section{Learning models and their evaluation in machine learning}
\label{learning-models-and-their-evaluation-in-machine-learning}

\subsection{Learning paradigms: regression, classification, clustering, and generation}

\paragraph{Aim.}

In the context of machine learning, two major paradigms are commonly distinguished: \textit{supervised learning} and \textit{unsupervised learning}. In supervised learning, the goal is to predict outputs (e.g., labels or values) from paired input--output data. This includes tasks such as \emph{regression}, where the output is continuous, and \emph{classification}, where the output is categorical. Unsupervised learning, on the other hand, deals with data that lacks explicit labels. The aim is to identify underlying structure or representations, such as clusters, low-dimensional manifolds, or generative rules.

\paragraph{Supervised learning: regression} 

Let us first set up the notation that we will use in this chapter.
\begin{itemize}
    \item $D$ is the number of input features (dimensionality of each sample),
    \item $D_f$ is the number of output dimensions (e.g., $D_f = 1$ for scalar regression, $D_f > 1$ for multi-output or vector-valued regression),
    \item $X \in \mathbb{R}^{N_x \times D}$ is the full set of training input points (often used for evaluation),
    \item $Y \in \mathbb{R}^{N_y \times D}$ is a another subset, used for optimization (e.g. Nyström approximation).
    \item $f(X) \in \mathbb{R}^{N_x \times D_f}$ are the corresponding function values on the training set.
\end{itemize}

Given a continuous function $f$, we work with a discrete set $Y$ and the corresponding values $f(Y)$, commonly referred to as a training set. The prediction function is then constructed as in \eqref{FIT}. The set $X$ does not necessarily coincide with the entire sample $Y$, it can be selected using strategies such as the Nyström method or other approximation techniques to reduce computational complexity. 

In ML, computing the coefficients $\theta$ in \eqref{FIT} is usually called as \textit{fitting} the model. Using the fitted model to evaluate $f_{k, \theta}$ at new inputs is referred to as \textit{making predictions}.

The choice of the regularization matrix $R(Y, X)$ in \eqref{FIT} determines the behavior of the model and corresponds to different regression techniques commonly used in machine learning.
\begin{itemize}
  \item \textit{Ridge Regression (Kernel Ridge Regression):} When $R(Y, X) = I$, the identity matrix, the regularization term becomes $\epsilon \| \theta \|_2^2$, penalizing large weights and helping to reduce overfitting. This is the standard choice in kernel ridge regression (KRR).

  \item \textit{Tikhonov regularization:} More generally, $R$ can be a positive semi-definite matrix that encodes prior knowledge about smoothness, scale, or structure (e.g., penalizing certain directions more heavily than others). This includes Ridge as a special case.

  \item \textit{Smoothing splines / Gradient-based smoothing:} If $R$ involves derivatives of the kernel or coordinates (e.g., discretized Laplacians or kernel differential operators), it enforces smoothness of the learned function.
\end{itemize}

The appropriate choice of $R$ depends on the nature of the data and the desired behavior.

\paragraph{Supervised learning: classification}

Classification tasks can be approached as a probabilistic extension of regression models. Rather than predicting continuous outputs, the goal is to assign input samples to discrete classes, using a softmax transformation of the model output. For example, a classifier for kernel ridge regression is described in equation \eqref{pi_k}. Softmax output can be interpreted as a vector of class probabilities. In particular, the $j$-th component of $\pi_{k, \theta}(x)$ represents the probability of the input $x$ belonging to class $j$:
\be
\mathbb{P}(y = j \mid x, \theta) = \pi_{k, \theta}(x)_j.
\ee

\paragraph{Unsupervised learning: clustering}

In this setting, the goal is to extract structure from unlabeled data $X \in \mathbb{R}^{N_x \times D}$. A usual formulation involves selecting a representative subset $Y \subset X$ that minimizes a discrepancy measure or divergence between the full data and the subset:
\be
Y = \mathop{\arg\inf}_{Y \in \mathbb{R}^{N_y, D}} d(Y, X),
\ee
where $d$ may be a classical distance (e.g., Euclidean norm, as in K-means) or a kernel-based discrepancy such as maximum mean discrepancy (MMD) \eqref{dk} for characteristic kernels.

Supervised and unsupervised learning are often interconnected, as follows.
\begin{itemize}
  \item \textit{Semi-supervised learning:} Cluster assignments or representative points $Y$ obtained from unsupervised methods can be used to generate pseudo-labels or guide the training of a supervised learning model. This results in a prediction $f_{\theta}(Z) \in \mathbb{R}^{N_z \times D_f}$; see Section~\ref{greedy-clustering-methods}.
  
  \item \textit{Cluster assignment as prediction:} In clustering, each test point $z^i$ is assigned to its closest prototype or cluster center in the set $Y$ using an \textit{assignment map}:
\be
\sigma(z^i, Y) = \mathop{\arg\min}_j d(z^i, y^j),
\quad \text{where } d(z^i, y^j) \text{ is a pairwise distance}.
\ee
This defines a labeling function
\be
\sigma(Z): [1, \ldots, N_z] \to [1, \ldots, N_y],
\ee
which assigns each point to a cluster index. The distance function $d$ may represent Euclidean distance or a kernel-based discrepancy, and can be expressed as a matrix $d \in \mathbb{R}^{N_x \times N_y}$ encoding all pairwise distances between data points and cluster centers.

\end{itemize}

\paragraph{Supervised learning: generative models} 

Generative models aim to learn a conditional distribution \( \mathbb{P}(y \mid x) \), which allows one to \emph{sample} new outputs \(y\) given an input condition \(x\). This differs from regression or classification, where the goal is to predict a specific output value. Instead, generative models learn the underlying data distribution to produce \emph{new data} which is statistically similar to the training set. We assume the following data. 
\begin{itemize}
  \item \(X \in \mathbb{R}^{N_x \times D_x}\) is a matrix of \emph{conditioning inputs} (e.g., labels, attributes),
  \item \(Y \in \mathbb{R}^{N_y \times D_y}\) is the corresponding set of \emph{high-dimensional outputs} (e.g., images, signals),
  \item \(Z \in \mathbb{R}^{N_z \times D_z}\) is a latent variable sampled from a known prior distribution (typically standard normal or uniform).
\end{itemize}

The generative process involves the following steps. 
\begin{enumerate}
  \item \emph{Learning a map} from latent variables \(z\) and conditioning the input \(x\) to the output space \(y\):  
  \be
  y = g_\theta(z, x),
  \ee
  where \(g_\theta\) is a decoder or generator trained to approximate the true conditional distribution.
  
  \item \emph{Sampling}: Once trained, the model can generate new samples by drawing latent variables \(z \sim \mathbb{P}_Z\) and combining them with a desired conditioning input \(x\).
\end{enumerate}
\noindent An example of generative model frameworks include: conditional kernel methods, i.e. sample \(z \sim \mathbb{P}_Z\), then use kernel-weighted interpolation to map to data space.
%
Usually a generative objective is to ensure that the \emph{samples \(y\)} produced under conditioning \(x\) \emph{match the statistical properties} of the target data distribution \( \mathbb{P}(Y \mid X = x) \), measured via distance metrics such as KL divergence, MMD, or KS test.


\subsection{Performance indicators for machine learning} \label{performance-indicators-for-machine-learning}

We now introduce commonly used performance metrics for machine learning models, divided into two main categories: unsupervised and supervised learning.
\begin{itemize}

\item \textit{Unsupervised learning: distances, divergences, and clustering metrics} 

Unsupervised learning lacks explicit labels, so model evaluation relies on statistical distances, distributional similarity, or clustering consistency rather than accuracy-based metrics.

\textit{$f$-Divergences}. \textit{$f$-divergences} measure the difference between two distributions \( \mathbb{P} \) and \( \mathbb{Q} \) via a convex function \( f : (0, +\infty) \to \mathbb{R} \) satisfying \( f(1) = 0 \). Assuming \( \mathbb{P} \ll \mathbb{Q} \), the general form is 
\be
D_f(\mathbb{P} \| \mathbb{Q}) = \int_{\mathcal{X}} f\left( \frac{d\mathbb{P}}{d\mathbb{Q}} \right) d\mathbb{Q},
\quad \text{or in the discrete case:} \quad
D_f(\mathbb{P} \| \mathbb{Q}) = \sum_x \mathbb{Q}(x) f\left( \frac{\mathbb{P}(x)}{\mathbb{Q}(x)} \right).
\ee
Common examples include: total variation distance (TVD) with \( f(x) = \frac{1}{2}|1 - x| \), Kullback--Leibler (KL) divergence with \( f(x) = x \log x \), and Hellinger distance with \( f(x) = (1 - \sqrt{x})^2 \).

\textit{Integral probability metrics (IPMs)}. IPMs define distances between distributions as the supremum of expectation differences over a function class \( \mathcal{F} \):
\be
\text{IPM}(\mathbb{P}, \mathbb{Q}) = \sup_{f \in \mathcal{F}} \left( \mathbb{E}_\mathbb{P}[f] - \mathbb{E}_\mathbb{Q}[f] \right).
\ee
Examples include the TVD, the Wasserstein-1 distance, which arises from optimal transport theory, and MMD, which uses kernel embeddings in RKHS. (see~ \eqref{dk}).

\textit{Statistical tests}. \textit{Kolmogorov--Smirnov (KS) Test} assesses whether two samples come from the same distribution, using the supremum difference between empirical CDFs:
\be
\label{KS}
\| \widehat{F}_X - \widehat{F}_Y \|_{\ell^\infty} \leq c_\alpha \sqrt{\frac{N_x + N_y}{N_x N_y}},
\ee
where \( c_\alpha \) corresponds to a chosen significance level (e.g., 0.05). For multivariate data, KS is applied marginally.

\textit{Clustering metrics}.
For models like \( k \)-means, internal validation is done using \textit{Inertia} which measures the compactness of clusters:
\be
\label{inertia}
I(X, Y) = \sum_{n=1}^{N_x} \left\| x^n - y^{\sigma(x^n, Y)} \right\|^2,
\quad \text{with} \quad \sigma(x, Y) = \mathop{\arg\min}_j d(x, y^j).
\ee

Lower inertia implies tighter clusters but does not necessarily imply correctness without external labels.

\item \textit{Supervised learning: regression and classification metrics} \label{indicators-for-supervised-learning}

In supervised learning, performance indicators directly compare model outputs to known ground truth labels. They differ by task type: regression or classification.

\textit{Regression metrics}. 
Used when outputs are continuous:

\textit{$\ell^p$ Norms:} The average prediction error is:
\be
\frac{1}{N_x} \| f_{k, \theta}(X) - f(X) \|_{\ell^p}, \quad 1 \le p \le +\infty,
\ee
where \( p = 2 \) gives the root mean square error (RMSE).

\textit{Normalized Error:} A scale-invariant variant is:
\be
\frac{\| f_{k, \theta}(X) - f(X) \|_{\ell^p}}{\| f_{k, \theta}(X) \|_{\ell^p} + \| f(X) \|_{\ell^p}},
\ee
commonly used in finance and other relative-scale domains.

\textit{Classification metrics}. 

In classification tasks, the model outputs a probability distribution over \(C\) classes via the softmax function \eqref{pi_k}. To compute standard classification metrics, the predicted class label is typically taken as the index of the maximum probability:
\be
f_{k,\theta}(X)^n = \mathop{\arg\max}_j \pi_{k, \theta}(x^n)_j.
\ee
This predicted label is then compared to the ground truth class $f(X)^n \in \{1, \ldots, C\}\). The following indicators are commonly used:

\textit{Accuracy:}
\be \label{score}
  \frac{1}{N_x} \#\{ f_{k,\theta}(X)^n = f(X)^n, \quad n = 1, \dots, N_x \}.
  \ee

\textit{Confusion matrix:} Records counts of predicted vs. true classes:
\be
M(i, j) = \#\{ f(x^n) = i \text{ and } f_{k,\theta}(x^n) = j \}
\ee

From the confusion matrix, we define TP (True Positives) are cases where the model correctly predicts a positive outcome. FP (False Positives) are incorrect positive predictions, and FN (False Negatives) are missed positives. TN (True Negatives) are correct negative predictions. PRE refers to Precision, $C$ is the total number of classes, and $i$ indexes over the classes. In statistical terms, a False Positive corresponds to a \textit{Type I error}, while a False Negative corresponds to a \textit{Type II error}. The summary is presented in Table~\ref{tab:class_metr}.

\begin{table}[h!]
\centering
\caption{Summary of classification metrics and their formulas \label{tab:class_metr}}
\renewcommand{\arraystretch}{1.8} 
\begin{tabular}{|l|l|l|}
\hline
\textit{Metric} & \textit{Definition} & \textit{Formula} \\
\hline
Precision (PRE) & Proportion of positive predictions that are correct & $\frac{TP}{TP + FP}$ \\
\hline
Recall (TPR) & Proportion of actual positives correctly predicted & $\frac{TP}{TP + FN}$ \\
\hline
F1 Score & Harmonic mean of Precision and Recall & $ \frac{2\cdot \text{Pre} \cdot  \text{TPR}}{\text{Pre} + \text{TPR}}$ \\
\hline
Micro-average Precision & Aggregated contributions of all classes & $\frac{\sum TP_i}{\sum TP_i + \sum FP_i}$ \\
\hline
Macro-average Precision & Average of class-wise precision scores & $\frac{1}{C} \sum_{i=1}^C PRE_i$ \\
\hline
FPR & False positive rate & $\frac{FP}{FP + TN}$ \\
\hline
\end{tabular}
\end{table}

\end{itemize}


\section{Application to supervised machine learning}
\label{application-to-supervised-machine-learning}

\subsection{Regression and reproducibility with housing price prediction}
\label{regression-problem-housing-price-prediction}

\paragraph{Objective.} In regression tasks with limited data, reproducibility and interpretability are critical--especially when the goal is to understand or audit predictions rather than only optimize performance. This test investigates the behavior of kernel-based regression compared to neural and ensemble models on the well-known Boston Housing dataset, focusing on extrapolation, generalization, and reproducibility. In particular, we test the ability of kernel ridge regression (KRR) to achieve exact interpolation when the data are fully observed, and explore whether statistical discrepancy (MMD) correlates with predictive accuracy.

\paragraph{Dataset.} 

We consider the Boston Housing dataset, which contains information collected by the U.S. Census Service concerning housing in Boston, Massachusetts. The dataset includes 506 observations, each with 13 numerical attributes and a corresponding target value representing median house prices\footnote{see \cite{Harrison}}. Our objective is to assess the extrapolation capabilities of our method. 

\paragraph{Methods and comparison.} 

We evaluate the (KRR)~\eqref{OP}\footnote{Implemented with \url{https://codpy.readthedocs.io/en/dev/}}, 
and we compare it with two standard regression models: a feed-forward neural network (FFN)\footnote{Implemented with PyTorch \url{https://pytorch.org}}, a random forest regressor (RF)\footnote{Implemented with scikit-learn \url{https://scikit-learn.org/stable/}}. 
Given a training dataset \( X \in \mathbb{R}^{N_x \times D} \) and corresponding labels \( f(X) \), we apply each model to predict the labels \( f_z \) of a test set \( Z \), and evaluate the predictions against the true values \( f(Z) \).
We use standard implementations of both FFN (with a typical multi-layer architecture trained using Adam optimizer) and a RF regressor (with 100 trees), without hyper-parameter tuning, to ensure a fair and reproducible comparison.

\paragraph{Results and analysis.} 

Figure~\ref{fig:59898979} presents comparative results in terms of model score, discrepancy error, and execution time as a function of the number of training examples.
\begin{itemize}
    \item The purpose of this test is to illustrate the reproducibility mode of the kernel ridge regressor~\eqref{FIT}. Reproducibility means in this test that zero-error should occur on the training set, which corresponds to the last value ($N_x=506$) of the KRR run of the left in Figure~\ref{fig:59898979}.
    \item The RF regressor demonstrates strong performance and good generalization across all data sizes, with relatively fast training time.
    \item The FFN performs less favorably than Random Forest, particularly on smaller training sets. This is likely due to the limited size of the dataset and the sensitivity of neural networks to hyperparameter tuning and regularization.
    \item $\text{MMD}$ \eqref{dk} correlates closely with prediction performance across all methods, supporting its relevance as an evaluation metric.
\end{itemize}

All methods use the same input data, and the kernel-based models rely on a standard kernel without additional tuning. While further improvements could be made with kernel design or hyperparameter optimization, the purpose here is to provide a fair benchmark comparison using standard configurations.

\begin{figure}
\centering
\includegraphics[width=1.0\textwidth, keepaspectratio]{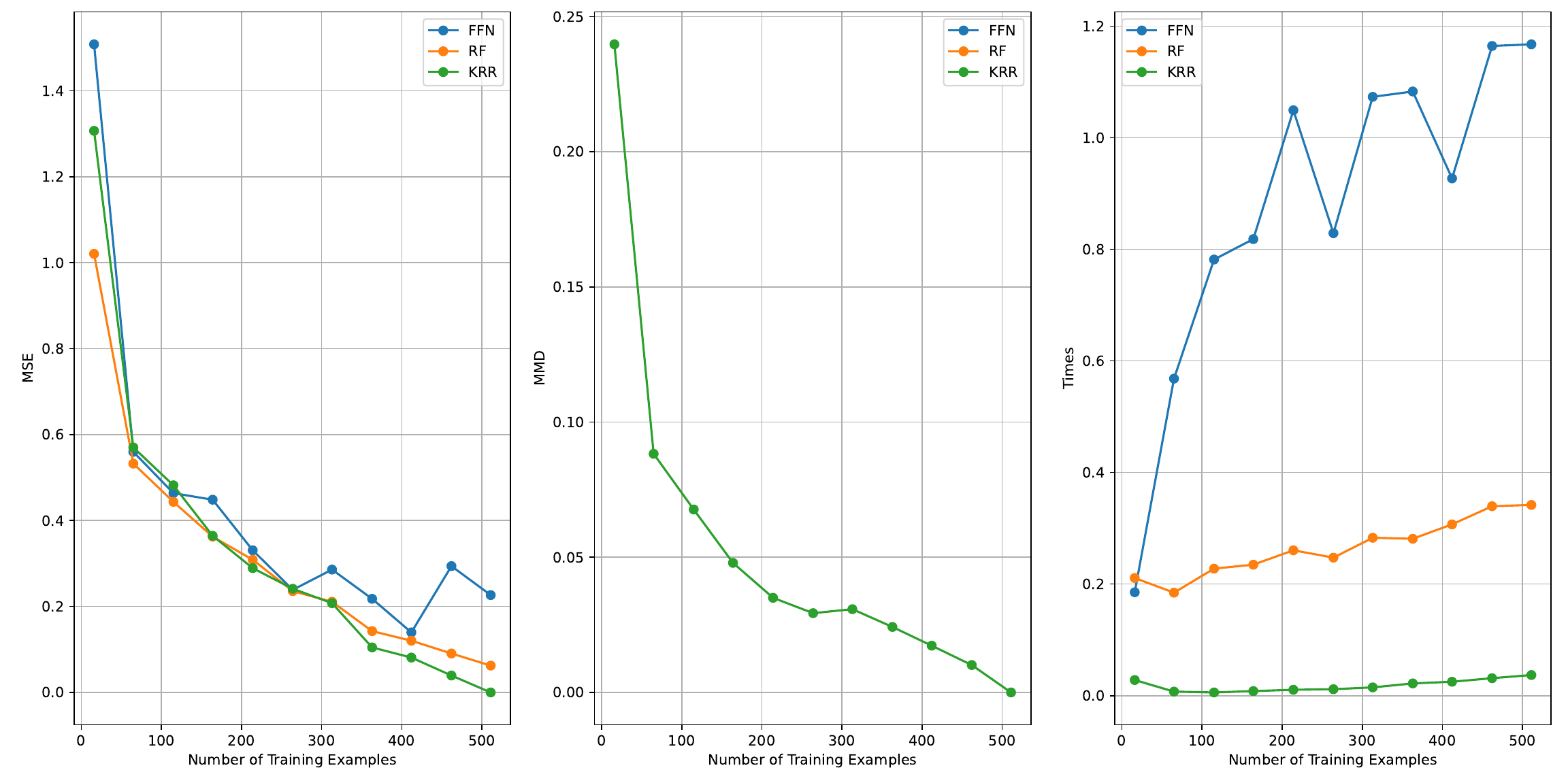}
\caption{\label{fig:59898979} Comparison of regression methods on the Boston Housing dataset in terms of model score (left), MMD, middle, and execution time (right) as a function of the number of training examples. The KRR model exhibits perfect interpolation on the full training set ($N_x = 506$), achieving zero mean squared error--highlighting its reproducibility property. The RF regressor shows strong generalization and fast execution across data sizes, while the {FFN} underperforms, especially on smaller datasets. A consistent correlation is observed between MMD and prediction accuracy, reinforcing MMD as a meaningful error indicator.}
\end{figure}


\subsection{Classification problem: handwritten digits}\label{classification-problem-handwritten-digits}

\paragraph{Objective.}

 Image classification serves as a classical benchmark for testing computational models of learning and generalization. This test aims to evaluate how well different methods--ranging from kernel-based to deep learning approaches--perform on a structured visual recognition task under limited data and tuning constraints. By comparing methods across statistical performance and computational cost, we assess their suitability in low-data or resource-constrained environments.

\paragraph{Dataset.} 

We use the MNIST dataset\footnote{see \cite{LeCunCortesBurges} and Page in Kaggle \url{https://www.kaggle.com/datasets/hojjatk/mnist-dataset}}, a benchmark collection of handwritten digits consisting of 60,000 training images and 10,000 test images. Each image is a grayscale \(28 \times 28\) pixel representation of a digit from 0 to 9, flattened into a vector of dimension \(D = 784\). The goal is to learn a classification function mapping input vectors \(X \in \mathbb{R}^{N_x \times 784}\) to one-hot encoded label vectors \(f(X) \in \mathbb{R}^{N_x \times 10}\), and evaluate its performance on test inputs \(Z \in \mathbb{R}^{N_z \times 784}\), where \(N_z = 10{,}000\).

\paragraph{Models.} We compare multiple classification models ranging from classical machine learning to modern neural architectures. These include the following. 
\begin{itemize}
    \item {K\_RR}: a basic kernel ridge classifier model (see Section~\ref{Perturbative-kernel-regression})\footnote{using the default RKHS setting of our library (see~\eqref{standardmap})}. 
    

    \item {K\_CM}: Same settings, but with some kernel engineering: we use the convolutional filter model $k_c(x,y) = k(x\ast \omega,y\ast \omega)$, where $(x\ast y)^i = \sum_n x^n y^{(i-n)\%N}$ (see~\eqref{conv_kernel}).

    \item {NN\_PM}: a basic deep-learning neural network, deliberately chosen to match the simplicity of the kernel ridge regression model {K\_RR}.

    \item {NN\_VGG}: this is a deep-learning model, inspired by VGGNet\footnote{see VGGNet \url{https://en.wikipedia.org/wiki/VGGNet}}, which employs two convolutional layers followed by max pooling and dropout, culminating in a fully connected classifier. While more expressive than the model {NN\_PM}, the {NN\_VGG} architecture remains compact to match the computational resources consumed by {K\_CM}.

    \item {RF}: a standard tree model, using 100 trees with default parameters. 
\end{itemize}

\paragraph{Results and analysis.} 

Figure~\ref{fig:65898} shows the comparative results in terms of classification accuracy and execution time as a function of the number of training examples.
\begin{itemize}
    \item The {K\_CM} model achieves the highest accuracy, as the cost of the highest computational time, across all training set sizes. With 2,048 training samples, this model can reach 97\% accuracy, which is close to human performance.
    \item The {VGG} model, limited by the computational resources, is at par with {K\_RR}, reaching both 94.5\% accuracy, although {K\_RR} is nearly 80\% more efficient in terms of computational resources for the same accuracy.
    
    \item The deep-learning neural network ({NN\_PM}) lies behind in terms of accuracy, indicating that neural networks architectures need heavy engineering to provide performing methods for image classification, as the VGG ones or ResNet\footnote{see ResNet \url{https://en.wikipedia.org/wiki/Residual_neural_network}}. We observe that these two models arose (in 2014 and 2015, respectively) from convolutional networks (1990\footnote{see ResNet \url{https://en.wikipedia.org/wiki/Convolutional_neural_network)}}).
    
    \item The random forest classifier ({RF}) is the most computationally efficient model but offers limited scalability in terms of accuracy. Its performance plateaus early (around 92\% accuracy with 2,048 samples), reflecting its inherent limitations in modeling smooth, high-dimensional functions like images. This is consistent with the well-known strengths of decision trees on tabular data rather than continuous visual inputs.
    
\end{itemize}

All models use the same raw input data, and kernel methods are applied using standard default kernel, although further gains could be obtained through kernel engineering, our objective here being to provide a consistent and reproducible benchmark using unoptimized RKHS baseline settings.

\begin{figure}
\centering
\includegraphics[width=1.\textwidth, keepaspectratio]{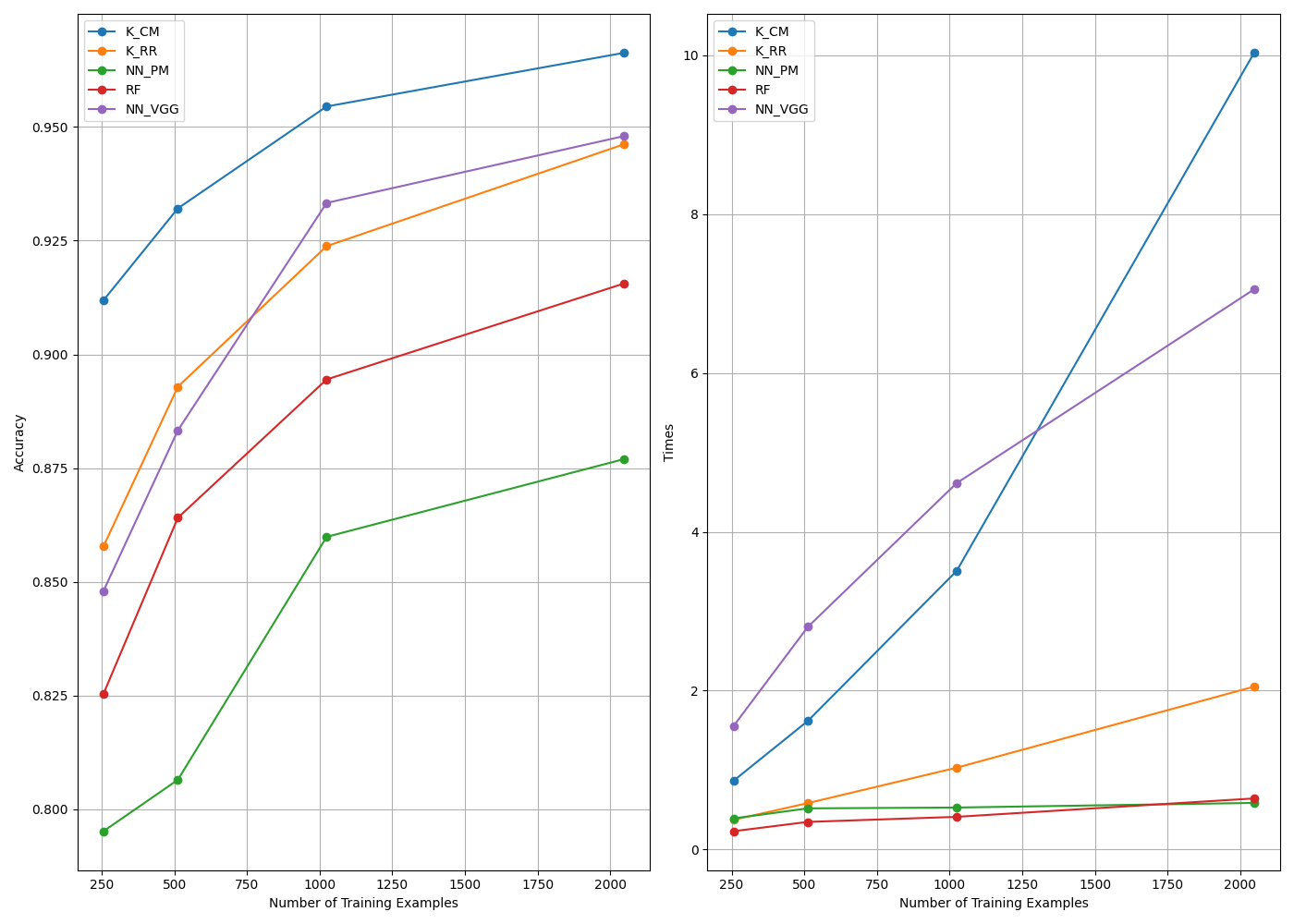}
\caption{\label{fig:65898} Performance comparison of classification models on the MNIST dataset across varying training set sizes. The plots show (left) classification accuracy and (right) execution time. Kernel-based methods ({KRR} and {CKRR}) and the convolutional {VGG} model achieve the highest accuracies, converging above 94\% with sufficient training data. {CKRR} demonstrates similar statistical performance to {KRR} while introducing translation invariance, though at increased computational cost. The regularized {FFN} improves steadily and approaches top-tier performance, while the unregularized {FFN:basic} underperforms throughout. {RF} offers fast execution but saturates earlier in accuracy, reflecting its limitations on image data.}
\end{figure}

 
\subsection{Reconstruction problems: learning from sub-sampled signals in tomography}
\label{reconstruction-problems-learning-from-sub-sampled-signals-in-tomography}

\paragraph{Objective} 

Next, we study the ability of learning-based methods to reconstruct tomographic images from sub-sampled measurement data. Specifically, we focus on SPECT (Single Photon Emission Computed Tomography), where acquiring high-resolution sinograms can be costly or physically constrained. We explore how a kernel-based learning framework can leverage patterns from fully sampled training data to improve reconstructions from sparsely sampled inputs.

This is motivated by various applications where reduced signal acquisition is either necessary or desirable. For example, in nuclear medicine, lowering the concentration of radioactive tracers minimizes radiation exposure. Similarly, faster data acquisition can alleviate the burden on expensive imaging equipment. Beyond medicine, such scenarios occur in fields such as biology, oceanography, and astrophysics.

To assess our method, we compare its performance with a classical iterative reconstruction algorithm -- Simultaneous Algebraic Reconstruction Technique (SART). Our goal is not to outperform SART across the board, but to highlight the potential of pattern-based learning in scenarios where examples recur or share structure.

\paragraph{Dataset and pre-processing.}

We use publicly available high-resolution CT images from Kaggle\footnote{Dataset available at \url{https://www.kaggle.com/vbookshelf/computed-tomography-ct-images}}. The dataset consists of $512\times 512$ grayscale images from $82$ patients, yielding approximately $30$ images per patient. We use images from the first $81$ patients ($2470$ images total) as the training set and reserve all $30$ images of the $82$nd patient for testing.

For each training image, the following steps are applied. 
\begin{enumerate}
\item A high-resolution sinogram is computed using a standard Radon transform with $256$ projection angles, resulting in a ($256\times 256$) matrix.
\item A low-resolution sinogram is generated using only 8 projection angles, resulting in a ($8\times 256$) matrix.
\item The high-resolution image is reconstructed from the high-resolution sinogram using the SART algorithm ($3$ iterations)\footnote{We use the implementation from \texttt{scikit-image} \url{https://scikit-image.org/docs/dev/api/skimage.transform.html\#skimage.transform.iradon_sart}}.
\end{enumerate}
\noindent
Figure~\ref{fig:SPECT} illustrates this process. The left panel shows the reconstructed image, the middle panel shows the high-resolution sinogram, and the right panel shows the corresponding low-resolution sinogram.

\begin{figure}[h]
\centering
\includegraphics[width=0.9\linewidth]{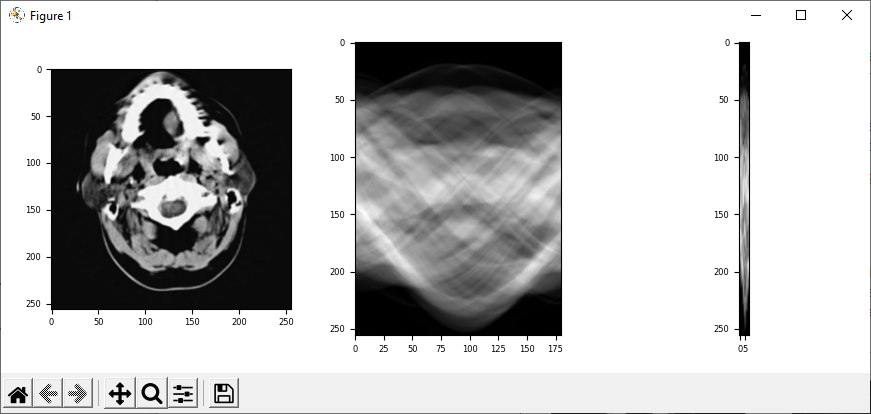}
\caption{\label{fig:SPECT}Example of sinogram generation and reconstruction. Left: reconstructed image from high-resolution sinogram. Middle: high-resolution sinogram. Right: low-resolution sinogram.}
\end{figure}

\paragraph{Learning setup.}

We formulate the reconstruction task as a supervised learning problem, where the input is a low-resolution sinogram and the output is the corresponding high-resolution image. The learning method is based on kernel ridge regression ({KRR}), described previously in Section~\eqref{Pk}.

Let $X \in \mathbb{R}^{2473 \times 2304}$ denote the training inputs, consisting of 2473 flattened sinograms of resolution $8 \times 256$. The associated training outputs are given by $f_x \in \mathbb{R}^{2473 \times 65536}$, representing the corresponding reconstructed images of size $256 \times 256$. For testing, we use $Z \in \mathbb{R}^{29 \times 2304}$, consisting of 29 sinograms from the 82nd patient, and their corresponding high-resolution ground truth images $f(Z) \in \mathbb{R}^{29 \times 65536}$.

The first sinogram of the 82nd patient is deliberately added to the training set to test memorization capabilities.

\paragraph{Results and visual comparison.}

We compare two reconstruction approaches with the ground truth of high-resolution images: SART, which reconstructs from low-resolution sinograms using three iterations, and {KRR}, which applies kernel ridge regression as defined in~\eqref{EI}. Figure~\ref{fig:359} shows the first 8 test images reconstructed using all three methods. Each triplet includes the ground truth (left), SART reconstruction (middle), and {KRR}-based reconstruction (right).

\begin{figure}[h]
\centering
\includegraphics[width=.9\linewidth]{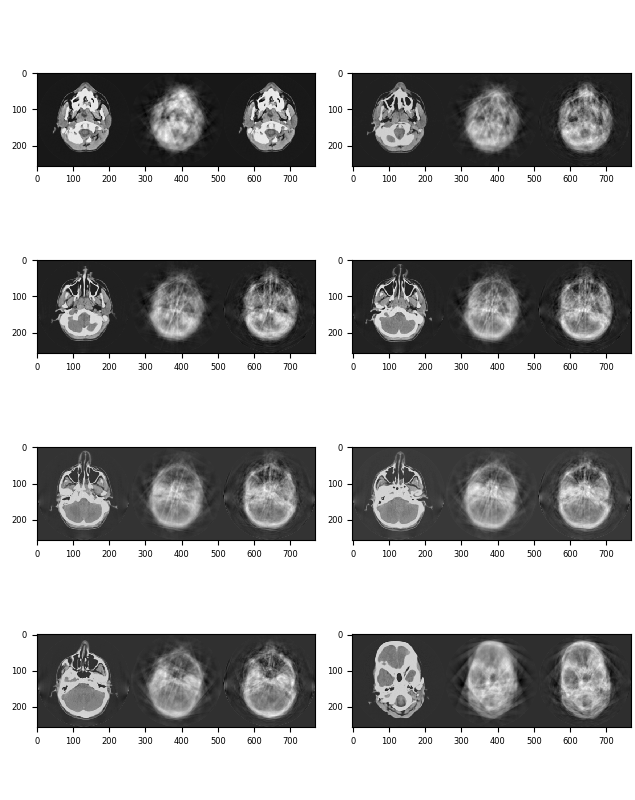}
\caption{\label{fig:359}Qualitative comparison: ground truth (left), SART reconstruction (middle), {KRR}-based reconstruction (right).}
\end{figure}

While the {KRR}-based method shows competitive performance, particularly in capturing familiar structural features, we emphasize that our goal is not to claim superiority over SART in general. Rather, we highlight the  strength of the method in recognizing and reconstructing recurring patterns -- for example, the first image in the test set (included in training) is reconstructed with near-perfect fidelity. This suggests promise for learning-based reconstruction in applications where pattern recurrence is common, such as automated diagnostic support systems.


\section{Application to unsupervised machine learning}
\label{applications-to-unsupervised-machine-learning}

\subsection{Semi-supervised classification with cluster-based interpolation}
\label{classification-problem-handwritten-digits-1}

\paragraph{Objective.} 

In supervised learning, achieving high accuracy with kernel methods such as kernel ridge regression ({KRR}) often requires significant computational resources and memory, particularly for large datasets. In this study, we explore a semi-supervised variant of {KRR}, wherein predictions are made not from the full training set but from a reduced set of cluster centroids. This strategy provides a trade-off between interpolation fidelity and computational efficiency, making it particularly relevant for large-scale problems such as handwritten digit recognition on the MNIST dataset.

\paragraph{Setup.} 

Let $X \in \mathbb{R}^{N_x \times D}$ denote the unlabeled training data, where $D = 784$ for MNIST. A smaller set of representative points $Y \in \mathbb{R}^{N_y \times D}$, with $N_y \ll N_x$, is computed by applying clustering algorithms to $X$. These representatives serve as interpolation centers in the {KRR} framework using formula~\eqref{FIT}. Labels are only required at the $N_y$ cluster centroids, simulating a semi-supervised learning scenario where labeled data are scarce but unlabeled data are abundant.

\paragraph{Cluster construction.} 

The interpolation set $Y$ is obtained by solving the centroid optimization problem:
\be
Y = \mathop{\arg\inf}_{Y \in \mathbb{R}^{N_y \times D}} d(X,Y). 
\ee
where $d(X,Y)$ is a discrepancy measure between the full data $X$ and the representative set $Y$. We evaluate two types of discrepancy. 
\begin{itemize}
\item
  Euclidean k-means:  Here, $d(X,Y)$ is the within-cluster sum of squared distances (inertia) (see~\eqref{inertia}).
\item
  Kernel-based clustering:  Here, $d(X,Y)$  is taken as MMD (see~\eqref{dk}).
\end{itemize}
Due to the nonconvexity of the clustering objective \(d(X,Y)\), solutions to  \eqref{Ycentroids} may not be unique. For example, standard k-means clustering may yield different centroids on different runs. Representative cluster centers obtained by both methods are visualized in Figure~\ref{fig:858}.

\begin{figure}
\centering
\includegraphics[width=0.7\textwidth, keepaspectratio]{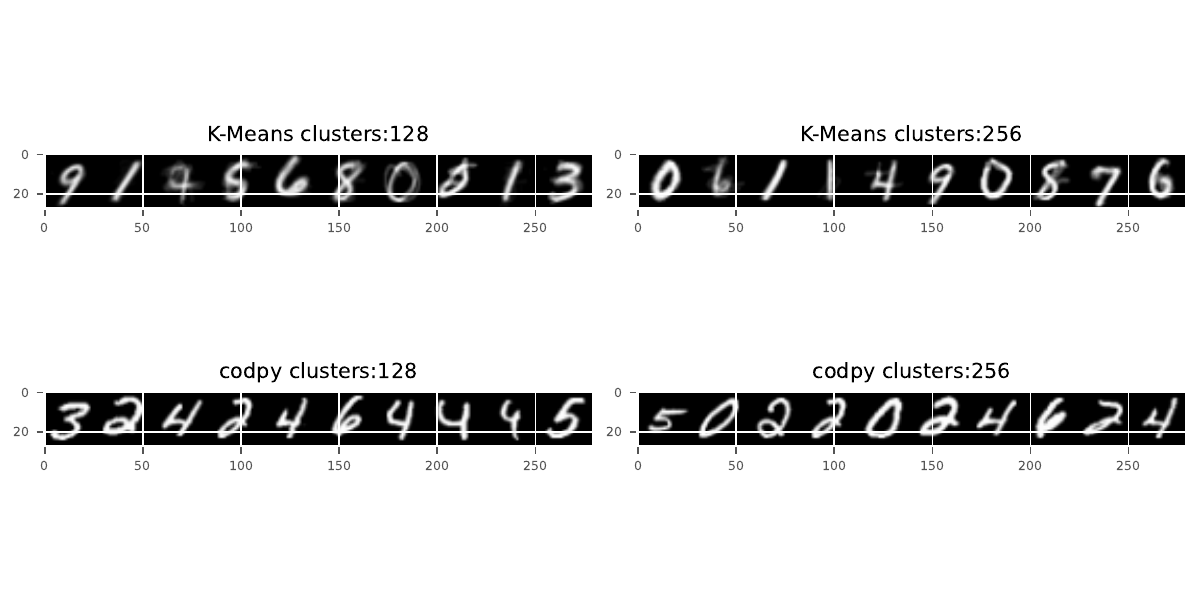}
\caption{\label{fig:858}Scikit (first row) and {KRR} (second row) clusters interpreted as images}
\end{figure}

\paragraph{Methodology.} 

The test considers the MNIST as a semi-supervised learning: we use the training set \(X \in \mathbb{R}^{N_x , D}\), with \(D=784\) dimensions, to compute the cluster centroids \(Y \in \mathbb{R}^{N_y , D}\). Then we use these clusters to predict the test labels \(f(Z)\), corresponding to the test set \(Z \in \mathbb{R}^{N_z , D}\), according to the kernel ridge regression \eqref{FIT}.

\paragraph{Results.} Figure~\ref{fig:859} summarizes the classification accuracy, kernel discrepancy (MMD), inertia, and execution time as a function of the number of clusters $N_Y=10,20, \ldots$. All clustering methods yield comparable classification performance, although kernel-based clustering typically achieves lower MMD. Performance is generally higher than in the baseline configuration of Section~\ref{classification-problem-handwritten-digits}, due to the use of the entire MNIST dataset for training. Notably, MMD serves as a reliable proxy for predicting classification performance in this semi-supervised setup.

\begin{figure}
\centering
\includegraphics[width=1\textwidth, keepaspectratio]{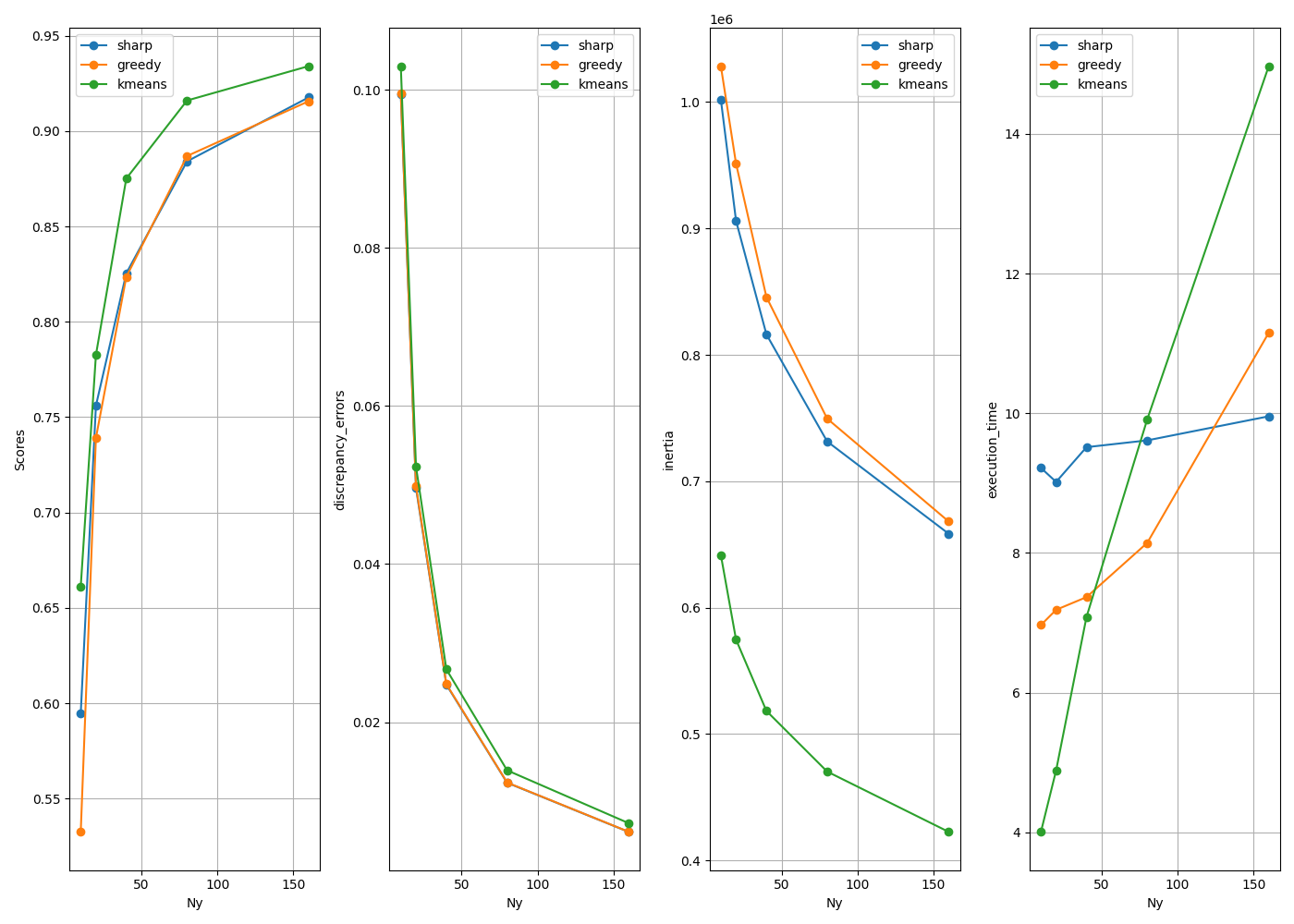}
\caption{\label{fig:859} Scores, MMD, inertia and exec. time for the MNIST dataset}
\end{figure}


\subsection{Credit card fraud detection}\label{credit-card-fraud-detection}

\paragraph{Objective.} 

Credit card fraud detection is a vital task in financial security. Given the rarity of fraudulent events and their often subtle deviations from legitimate behavior, the problem poses a unique challenge. Traditional supervised learning struggles in this domain due to extreme class imbalance and the evolving nature of fraud. In this setting, unsupervised learning methods can be valuable for detecting anomalies without requiring labeled data.

\paragraph{Dataset and methods.} 

The dataset\footnote{Dataset source: \href{https://www.kaggle.com/mlg-ulb/creditcardfraud}{Kaggle Credit Card Fraud Dataset}} contains transactions by European cardholders during two days in September 2013. It includes $284{,}807$ transactions, among which only $492$ are fraudulent (roughly $0.172\%$). Each transaction is represented by 30 numerical features, most of which are anonymized via Principal Component Analysis (PCA). Two features--\texttt{Time} and \texttt{Amount}--remain in their original form. The response variable, \texttt{Class}, is binary, indicating whether a transaction is fraudulent ($1$) or not ($0$).

We apply two unsupervised clustering algorithms.
\begin{itemize}
\item \emph{$k$-means clustering\footnote{{scikit-learn} implementation}}, a classical method that partitions data into clusters by minimizing intra-cluster variance.
\item \emph{MMD-based clustering}, which minimizes a statistical discrepancy MMD between clusters, aiming to identify distributional outliers more effectively.
\end{itemize}
These methods are trained to discover typical transaction patterns. Fraud is then inferred by identifying transactions that do not conform to the learned clusters (e.g., assigned to "small" or sparse clusters, or flagged by secondary scoring heuristics).

\paragraph{Results and analysis}. Figure~\ref{fig:580} illustrates confusion matrices for the last scenario of each approach. Both methods successfully identify most non-fraudulent transactions, but differ in detecting fraud: MMD-based clustering achieves fewer false positives while maintaining similar true positive counts compared to $k$-means, suggesting better discrimination of atypical patterns. However, overall recall remains low for both approaches, reflecting the challenge of unsupervised fraud detection in highly imbalanced datasets and highlighting the trade-off between precision and recall in anomaly detection tasks.

\begin{figure}
\centering
\includegraphics[width=1.0\textwidth, keepaspectratio]{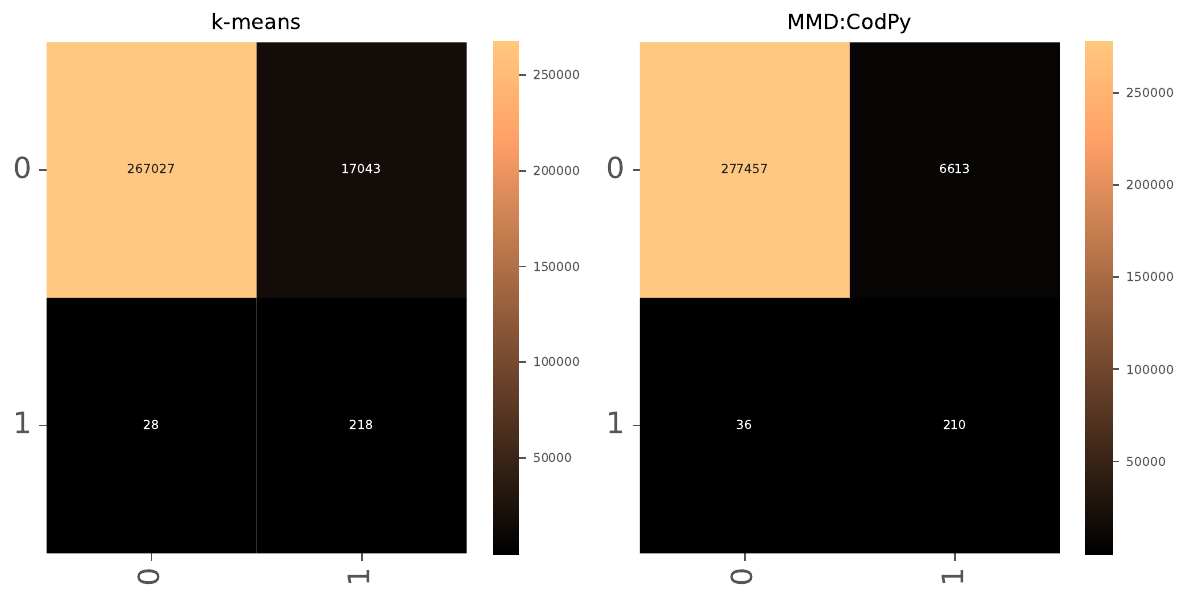}
\caption{\label{fig:580} Confusion matrix obtained using the $k$-means and MMD-based clustering in the final evaluation scenario}
\end{figure}


\subsection{Portfolio of stock clustering}
\label{portfolio-of-stock-clustering}

\paragraph{Objective.} 
Clustering stocks based on their price movement patterns helps identify groups of companies with similar behavior in the market. Such grouping has practical applications in portfolio diversification and risk management. At this stage, we compare traditional $k$-means clustering with a discrepancy-based method using MMD \eqref{dk} minimization.

\paragraph{Database.}  

The data consist of daily stock price differences (closing minus opening prices) from 2010 to 2015 for 86 companies. Each stock is represented by a high-dimensional feature vector capturing daily price movement over the period.

\paragraph{Methods.} 

Two clustering methods are applied.
\begin{itemize}
    \item \emph{k-means clustering}\footnote{{scikit-learn} implementation.}, a classical method that minimizes intra-cluster variance using Euclidean distances.

\item \emph{MMD-based clustering}, implemented via the sharp discrepancy algorithm\footnote{CodPy implementation}, described in Section~\ref{Discrepancy functional and sharp-discrepancy-sequences}. This approach minimizes the MMD, encouraging clusters that are statistically distinct in terms of their kernel mean embeddings. 
\end{itemize}

The number of clusters is fixed at $N=10$ for both methods. Each algorithm is applied to the normalized data\footnote{Using the \texttt{sklearn.preprocessing.Normalizer}} to ensure fair comparison.

\paragraph{Results.} Table~\ref{tab:stock_clust} shows the resulting clusters from both methods. Each cluster contains companies whose stock price dynamics over the 5-year period are considered similar under the respective method:

\begin{longtable}[]{@{}
  >{\centering\arraybackslash}p{(\columnwidth - 4\tabcolsep) * \real{0.0271}}
  >{\centering\arraybackslash}p{(\columnwidth - 4\tabcolsep) * \real{0.6000}}
  >{\centering\arraybackslash}p{(\columnwidth - 4\tabcolsep) * \real{0.3729}}@{}}
\caption{Stock clustering \label{tab:stock_clust}}\tabularnewline
\toprule\noalign{}
\textbf{\#} & \textbf{$k$-means} & \textbf{MMD minimization} \\
\midrule\noalign{}
\endfirsthead
\toprule\noalign{}
\textbf{\#} & \textbf{$k$-means} & \textbf{MMD minimization} \\
\midrule\noalign{}
\endhead
\bottomrule\noalign{}
\endlastfoot
\textbf{0} & Apple, Amazon, Google/Alphabet & ConocoPhillips, Chevron, IBM, Johnson \& Johnson, Pfizer, Schlumberger, Valero Energy, Exxon \\
\hline
\textbf{1} & Boeing, British American Tobacco, GlaxoSmithKline, Home Depot, Lookheed Martin, MasterCard, Northrop Grumman, Novartis, Royal Dutch Shell, SAP, Sanofi-Aventis, Total, Unilever & Intel, Microsoft, Symantec, Taiwan Semiconductor Manufacturing, Texas instruments, Xerox \\
\hline
\textbf{2} & Caterpillar, ConocoPhillips, Chevron, DuPont de Nemours, IBM, 3M, Schlumberger, Valero Energy, Exxon & Dell, HP \\
\hline
\textbf{3} & Intel, Navistar, Symantec, Taiwan Semiconductor Manufacturing, Texas instruments, Yahoo & Coca Cola, McDonalds, Pepsi, Philip Morris \\
\hline
\textbf{4} & Canon, Honda, Mitsubishi, Sony, Toyota, Xerox & Boeing, Lookheed Martin, Northrop Grumman, Walgreen \\
\hline
\textbf{5} & Colgate-Palmolive, Kimberly-Clark, Procter Gamble & AIG, American express, Bank of America, Ford, General Electrics, Goldman Sachs, JPMorgan Chase, Wells Fargo \\
\hline
\textbf{6} & Johnson \& Johnson, Pfizer, Walgreen, Wal-Mart & British American Tobacco, GlaxoSmithKline, Novartis, Royal Dutch Shell, SAP, Sanofi-Aventis, Total, Unilever \\
\hline
\textbf{7} & Coca Cola, McDonalds, Pepsi, Philip Morris & Amazon, Canon, Cisco, Google/Alphabet, Home Depot, Honda, MasterCard, Mitsubishi, Sony, Toyota \\
\hline
\textbf{8} & Cisco, Dell, HP, Microsoft & Apple, Caterpillar, DuPont de Nemours, 3M, Navistar, Yahoo \\
\hline
\textbf{9} & AIG, American express, Bank of America, Ford, General Electrics, Goldman Sachs, JPMorgan Chase, Wells Fargo & Colgate-Palmolive, Kimberly-Clark, Procter Gamble, Wal-Mart \\
\hline
\end{longtable}

The two clustering methods  ---$k$-means and MMD-based clustering--- produce different but meaningful groupings of stocks. $k$-means tends to group stocks with similar Euclidean distance in their normalized return vectors, often clustering by volatility or magnitude of movement. In contrast, the MMD-based method emphasizes differences in distributional features, which can lead to different groupings that may reflect statistical regularities not captured by simple distance metrics.

For example, MMD places several large tech companies into a single group, but distributes others (e.g., Apple, Amazon) differently than $k$-means. Similarly, financial institutions are grouped differently across the two methods, indicating differing perspectives on similarity.

These differences suggest that the clustering methods highlight \textit{complementary structure} in the data, with no single "best" clustering. The choice of method should be guided by the downstream application--for example, whether one prioritizes interpretability, statistical distributional features, or geometric compactness.

 
\section{Application to generative models}
\label{application-to-generative-models}
 
\subsection{Generating complex distributions with CelebA dataset}
\label{generating-complex-distributions}

\paragraph{Objective.} 

Now, we illustrate the practical behavior of the kernel-based generator described in Section~\ref{Encoder-decoders-and-sampling-algorithms} in the context of image generation. Specifically, we evaluate the ability of the generator to produce novel, high-dimensional image samples that approximate the structure of a real-world dataset.

\paragraph{Dataset description.} We use the CelebA dataset\footnote{see \cite{liu2015} and \url{https://mmlab.ie.cuhk.edu.hk/projects/CelebA.html}}, which contains over 200,000 aligned facial images annotated with 40 binary attributes. Each image is represented as a \(218 \times 178 \times 3\) RGB tensor, yielding a flattened vector in \(\mathbb{R}^{116{,}412}\). For this test, we sample \(N_y = 1000\) images from the dataset to serve as the training set \(Y = (y^1, \ldots, y^{N_y}) \subset \mathbb{R}^{116{,}412}\).

\paragraph{Latent space sampling and image generation.}

In image generation tasks, the input is typically a set of real examples, and the goal is to generate novel samples that resemble the training distribution while exhibiting controlled variation.

We aim to illustrate both generative and recognition behavior in latent space. To this end, we define the generator $G_k$ following the encoder--decoder formulation in~\eqref{ED}:
\be \label{generatorA}
  G_{k}(\cdot)  =  K(\cdot,X) \theta, \quad \theta = K(X,X)^{-1} (Y\circ \sigma). 
\ee 
The dataset $Y=(y^1,\ldots,y^N)$ consists of CelebA images, while $X=(x^1,\ldots,x^N)$ denotes a latent representation, typically drawn from a simple distribution (e.g., uniform or Gaussian). The output $G_{k}(x^i)$, $i=1,\ldots,N$ is expected to be the image $y^{\sigma(i)}$ of the dataset, and  $G_{k}(x)$, $x\neq x^i$, is a generated, "fake" image, $x$ being a \textit{latent variable}.

We used \(N_y = 1000\) images of celebrity examples, denoted \(Y = (y^1,\ldots,y^{N_y})\) in the training set. Thus the input distribution dimension is \((1000, 116412)\). We encode this distribution, using an uniform distribution $X$ as latent variable, testing two latent dimensions having size \(D_x = 3,40\), defining a map $Y_k^\sigma(\cdot)$. Once encoded, we generate 16 samples of the uniform distribution $z^i$, $i=1,\ldots,16$, and plot $G_k(z^i)$ in the figures  \eqref{fig:g3},\eqref{fig:g40}, in order to discuss the role of the latent space, similarly to our numerical tests in Section~\ref{Encoder-decoders-and-sampling-algorithms}.

\begin{figure}[!htb]
    \centering
    \begin{subfigure}{0.48\textwidth}
        \centering
        \includegraphics[width=\textwidth]{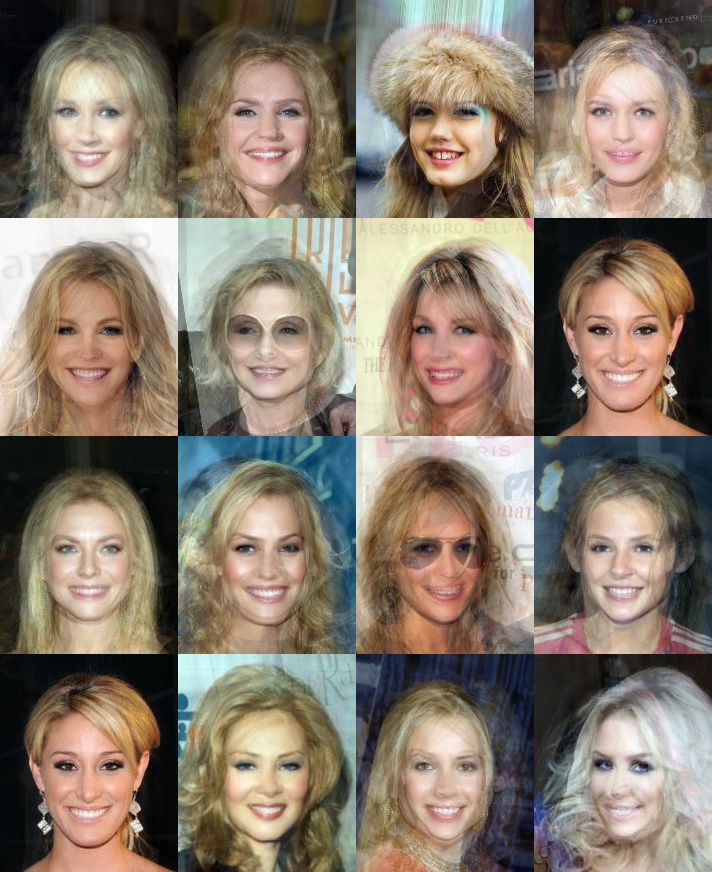}
        \caption{Generated pictures, latent dim 3}
        \label{fig:g3}
    \end{subfigure}
    \hfill
    \begin{subfigure}{0.48\textwidth}
        \centering
        \includegraphics[width=\textwidth]{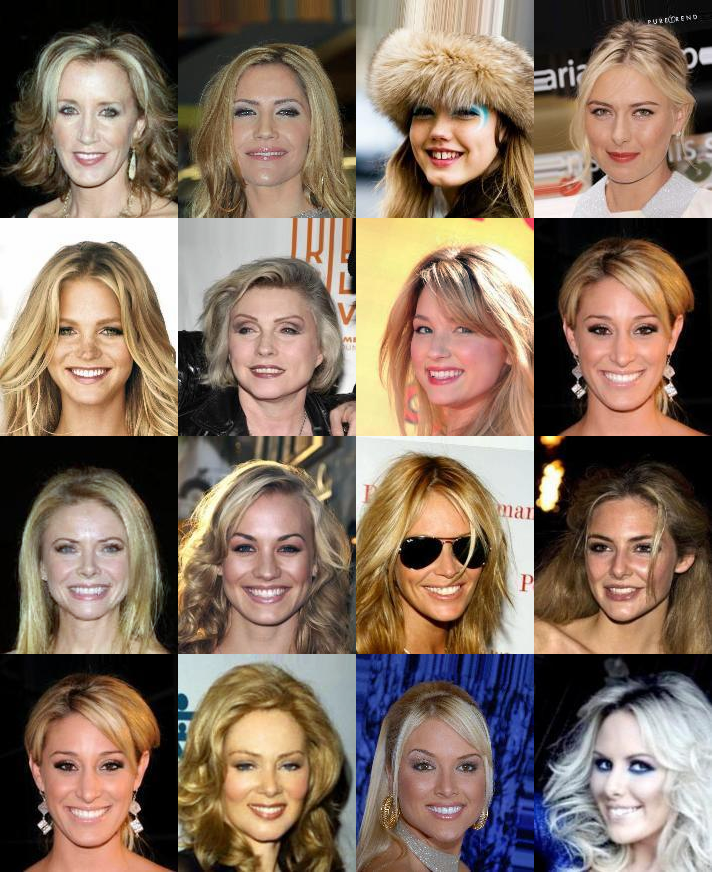}
        \caption{Closest database pictures}
        \label{fig:d3}
    \end{subfigure}
        \caption{
    Qualitative comparison between generated samples and their nearest neighbors in the training dataset, for a latent space of dimension \(D_x = 3\). The generator produces diverse and structured outputs that are visually close to the training data, indicating bad coverage and generalization.}
    \label{fig:latent3}
\end{figure}

\begin{figure}[!htb]
    \centering
    \begin{subfigure}{0.48\textwidth}
        \centering
        \includegraphics[width=\textwidth]{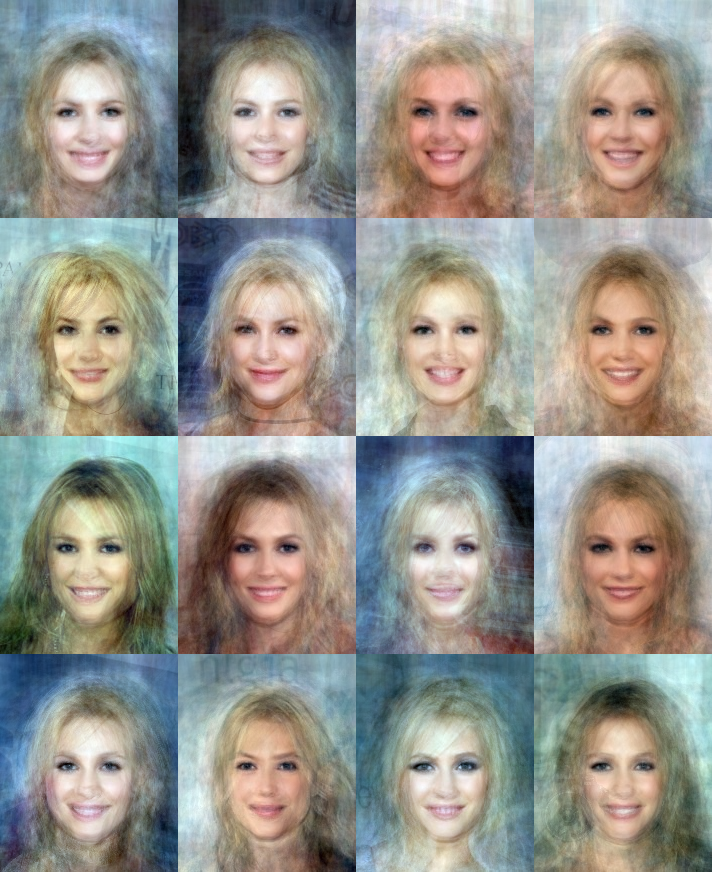}
        \caption{Generated pictures, latent dim 40}
        \label{fig:g40}
    \end{subfigure}
    \hfill
    \begin{subfigure}{0.48\textwidth}
        \centering
        \includegraphics[width=\textwidth]{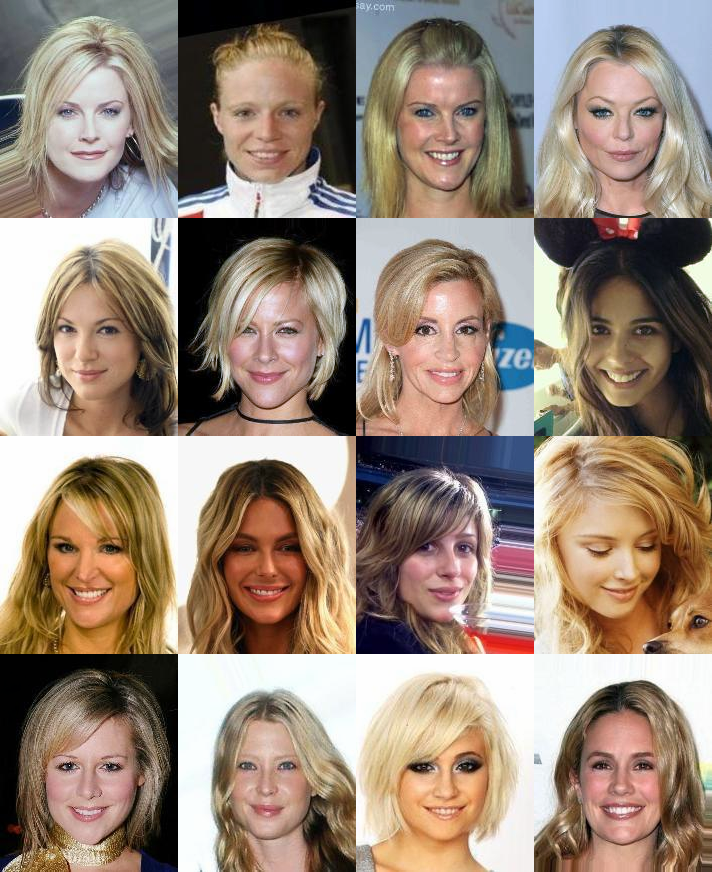}
        \caption{Closest database pictures}
        \label{fig:d40}
    \end{subfigure}
        \caption{
    As the latent dimension increases to \(D_x = 40\), generated samples appear smoother and closer to their nearest training neighbors, with reduced diversity. This reflects the concentration of measure effect in high-dimensional spaces, where random samples tend to cluster near the mean. While higher dimensions offer greater generative flexibility, they may impair coverage and realism.}
    \label{fig:latent40}
\end{figure}

\paragraph{Matching and latent discrepancy.} To evaluate the similarity between generated and training images, we compute latent-space nearest neighbors using a kernel discrepancy measure:
\be\label{791}
y^{i(j)}, \quad i(j) = \arg\min_{j = 1,\ldots,N_y} d_k(x^i, l^j),
\ee
where \(l^j \in \mathbb{R}^{D_x}\) denotes the latent code associated with image \(y^j\). This approach functions as an efficient recognition method in latent space, revealing visual similarity between generated and training images.

\paragraph{Results and analysis.} When \(D_x = 3\), the generator produces a diverse array of samples that are visually distinct, but quite close, from their nearest neighbors. As the latent dimension increases to \(D_x = 40\), the samples become smoother and less varied, suggesting a reduction in generative diversity. This phenomenon is consistent with the \textit{concentration of measure} effect: in higher dimensions, uniformly sampled latent points are less likely to fall near the support of the training distribution.

The observed trade-off between sample diversity and realism is directly influenced by the latent dimensionality. Lower-dimensional spaces promote better coverage and fidelity to the training distribution, while higher dimensions offer more generative flexibility but risk producing blurrier or less realistic outputs.

\paragraph{Outlook and next steps.} 

This test validates the continuity and generalization behavior of the generator $G_k$ and provides a foundation for more advanced refinement techniques. Next, we are going to introduce two such methods aimed at improving perceptual quality and distributional alignment. These methods select latent codes $x \in \mathbb{R}^{D_x}$ such that $G_k(x)$ approximates a separate reference distribution $Y^{\mathrm{ref}} \subset \mathbb{R}^{116{,}412}$, using discrepancy-based selection or projection criteria.

\subsection{Image reconstruction}

\paragraph{Objective.} We next evaluate the reconstruction capability of the generator $G_k$, defined in~\eqref{generatorA}, when used as part of an encoder--decoder framework. Given a set of unseen images $y \in \mathbb{R}^{116412}$ drawn from the CelebA dataset, the objective is to recover an approximate latent representation $x \in \mathbb{R}^{D_x}$ such that the generated output $G_k(x)$ closely matches the input image $y$.

This process consists of three stages, as follows. 
\begin{itemize}
    \item Decoding: The generator $x \mapsto G_k(x)$ maps latent vectors to image space.
    \item Encoding: The inverse mapping $y \mapsto G_k^{-1}(y)$ seeks a latent code that reconstructs a given image. 
    \item Reconstruction: The map $y \mapsto G_k(G_k^{-1}(y))$, is called a reconstruction.   
\end{itemize}

\paragraph{On the non-uniqueness of mappings.}
In our framework, the generator $G_k$ is trained via a kernel ridge regression model \eqref{Pk}, where the pairing between latent codes $x^i \in X$ and images $y^{\sigma(i)} \in Y$ is determined by solving a discrete Gromov--Monge problem. Unlike classical Monge transport on compatible spaces, the Gromov formulation compares intra-domain distances rather than absolute positions. As a result, the optimal permutation $\sigma \in \Sigma$ may not be unique, and structurally similar permutations can yield similar objective values. However, this construction is expected to compute a smooth, one-to-one mapping from two any distributions.

To illustrate this reconstruction process, we select 16 images from the CelebA test set (distinct from those used in training) and compute their reconstructions using the procedure above. figures~\ref{fig:g47} and~\ref{fig:d47} display the reconstructed outputs and their corresponding closest database images (in latent space), respectively.

\paragraph{Observations.} The reconstructed images exhibit improved fidelity relative to raw generative samples (as seen in figures~\ref{fig:g3} and~\ref{fig:g40}), confirming the effectiveness of the encoding procedure. However, some reconstructed images appear visually similar to one another, a phenomenon analogous to \emph{mode collapse} in generative adversarial networks.

In the context of kernel-based generators, this effect can be attributed to the structure of the latent space: the encoder $G_k^{-1}$ projects high-dimensional image data into a lower-dimensional latent domain, where multiple images may share nearly identical latent codes due to the non-injective nature of the inverse mapping. Consequently, their reconstructions via $G_k$ are nearly indistinguishable.

This highlights an inherent trade-off in reconstruction-based approaches: while image fidelity improves through latent-space alignment, diversity may be reduced if the encoding fails to preserve semantic differences between visually distinct inputs. In the next section, we address this issue by introducing a refinement strategy based on Wasserstein adversarial alignment, designed to improve both variety and fidelity of the generated samples.

\begin{figure}[!htb]
    \centering
    \begin{subfigure}{0.48\textwidth}
        \centering
        \includegraphics[width=\textwidth]
        {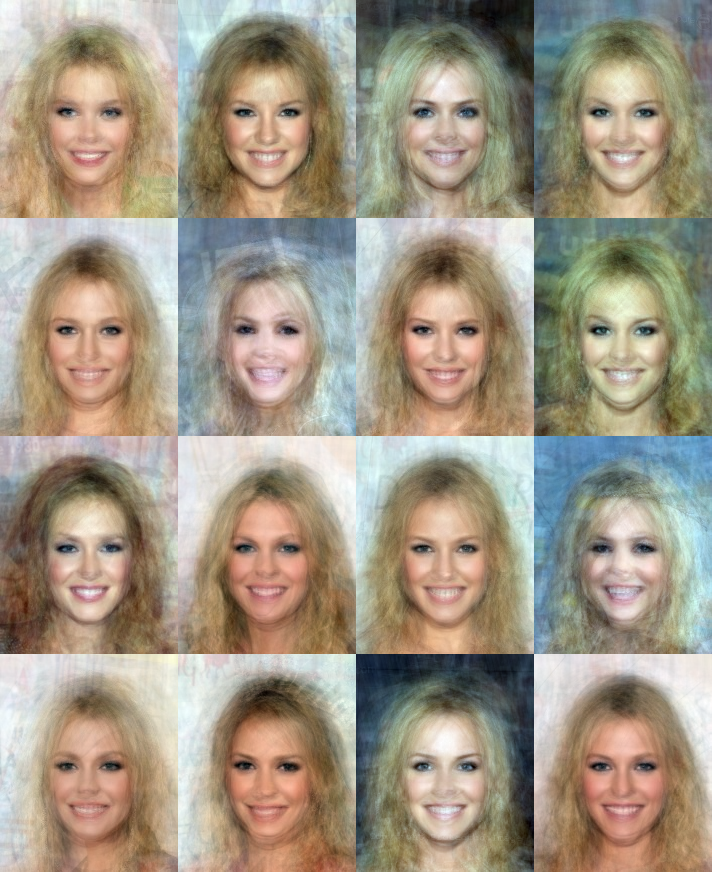}
        \caption{Generated pictures, latent dim 40}
        \label{fig:g47}
    \end{subfigure}
    \hfill
    \begin{subfigure}{0.48\textwidth}
        \centering
        \includegraphics[width=\textwidth]
        {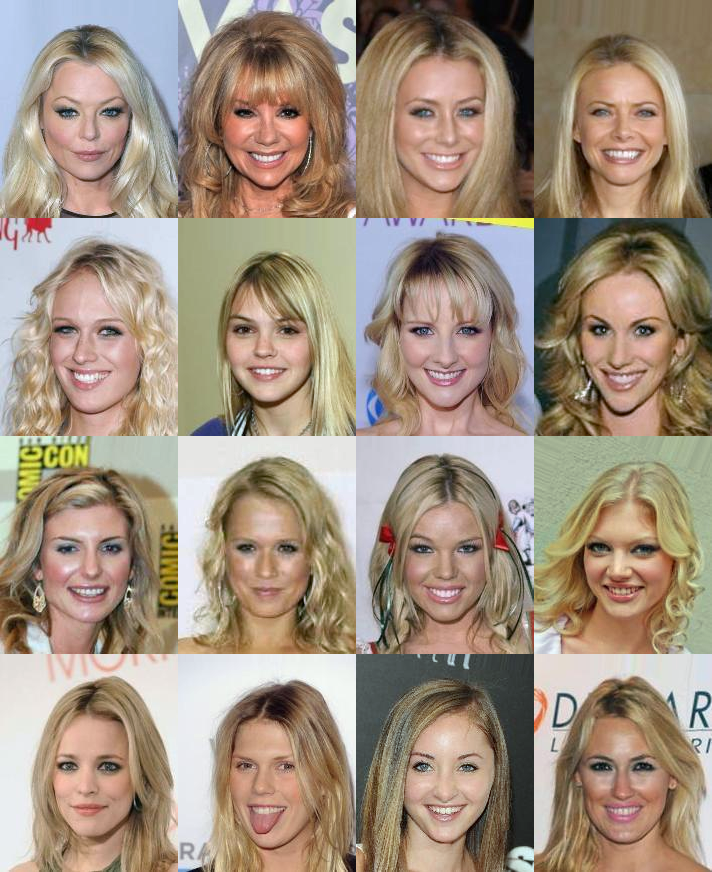}
        \caption{Closest database pictures}
        \label{fig:d47}
    \end{subfigure}
    \caption{
    Image reconstruction example using a kernel-based encoder--decoder pair \( (G_k^{-1}, G_k) \). The left panel shows reconstructed images obtained by encoding and decoding training set samples; the right panel displays the corresponding original images from the dataset. Observe that several reconstructed images on the left appear visually similar, 
indicating a form of mode collapse.}
    \label{fig:reconstruction}
\end{figure}

\subsection{Generative adversarial Wasserstein kernel architectures}

\paragraph{Objective.} 

Now, we refine the generative model introduced previously by incorporating a Wasserstein-based adversarial training procedure. This approach seeks to improve the quality and diversity of generated outputs by minimizing the Wasserstein distance between the generated and real image distributions, in the spirit of generative adversarial networks (GANs) with a focus on the Wasserstein formulation\footnote{see \cite{Goodfellow:2014} and \cite{Arjovsky:2017}}.

\paragraph{Architecture.} 

Consider a set of $N$ distinct real images $Y_D = (y_D^1, \dots, y_D^N)$ and a set of generated images $Y_G$, such that $Y_D \cap Y_G = \emptyset$. The generator distribution $Y_G$ is denoted $Y^\sigma$ in Section~\ref{generating-complex-distributions}. Let $G_k(x)$ denote the generator, where $x$ is a latent variable, and define the set of generated images as
\be
Y_G(X) = \left( y_G(x^1), \dots, y_G(x^N) \right) = \left( G_k(x^1), \dots, G_k(x^N) \right),
\ee
where $X = (x^1, \dots, x^N)$ is the latent code vector. In the adversarial setting, we introduce a \textit{discriminator}, that is, a scalar function $D : y \mapsto \mathbb{R}$ defined on images. The goal is to solve the following saddle-point problem:
\be \label{LGD}
\inf_X \sup_{\left\{ D : D(y^1) - D(y^2) \leq |y^1 - y^2|_p^p \right\}} \mathcal{L}(X, D) = \sum_{n=1}^N D\big( y_G(x^n) \big) - \sum_{n=1}^N D\big( y_D^n \big).
\ee

For fixed $X$, this inner problem is equivalent to the dual formulation of the $p$-Wasserstein distance described in Section~\ref{dual}, that is,
\be
\sup_{\left\{ D : D(y^1) - D(y^2) \leq |y^1 - y^2|_p^p \right\}} \mathcal{L}(X, D) = W_p^p \left( Y_G(X), Y_D \right),
\ee
where $W_p^p$ denotes the $p$-Wasserstein cost. Thus, problem \eqref{LGD} reduces to the following optimization:
\be \label{LGDW}
\inf_{X = (x^1, \dots, x^N)} W_p^p\left( Y_G(X), Y_D \right), \quad W_p^p\left( Y_G(X), Y_D \right) = \sum_n | Y_G(x^n) - Y_D^{\sigma_M(x^n)} |_p^p,
\ee
where $\sigma_M(X)$ is the permutation computed by solving the discrete Monge problem as defined in~\eqref{LSAP}. Observe that $\sigma_M$ is distinct from $\sigma$ used in $Y^\sigma$, and is computed afresh during optimization to match $Y_G(X)$ to $Y_D$.

The formulation \eqref{LGDW} forces the generator to produce images that lie close (in Wasserstein sense) to the target distribution $Y_D$.

\paragraph{Gradient flow.} 

Generative adversarial methods typically use the case $p=1$, but in this work we consider the simpler, yet comparable case $p=2$, for which the gradient is easily computable. Specifically, we compute 
\be
\nabla_X \mathcal{L}(X, D) = \sum_n \nabla_{x^n} \left| Y_D^{\sigma_M(x^n)} - Y_G(x^n) \right|_2^2 = \sum_n \left( Y_G(x^n) - Y_D^{\sigma_M(x^n)} \right) \cdot \nabla Y_G(x^n),
\ee
where $\nabla Y_G(x^n) = \nabla G_k(x^n)$ can be estimated using the gradient formula from Section~\ref{nabla}. This leads to the semi-discrete gradient descent scheme:
\be
\frac{d}{dt} X_t = -\left( Y_G(X_t) - Y_D^{\sigma(X_t)} \right) \cdot \nabla Y_G(X_t),
\ee
which can be integrated numerically using the descent algorithm introduced in Section~\ref{explicit-descent-algorithms}.

\paragraph{Results.} 

To illustrate the effectiveness of this Wasserstein-based method, we solve the minimization problem starting from initial latent variables $X_0$, obtained from the reconstruction shown in Figure~\ref{fig:g47}. The resulting generated images are shown in Figure~\ref{fig:g48}, and the corresponding closest images from the training set $Y_D$ are shown in Figure~\ref{fig:d48}.

\begin{figure}[!htb]
    \centering
    \begin{subfigure}{0.48\textwidth}
        \includegraphics[width=\textwidth]{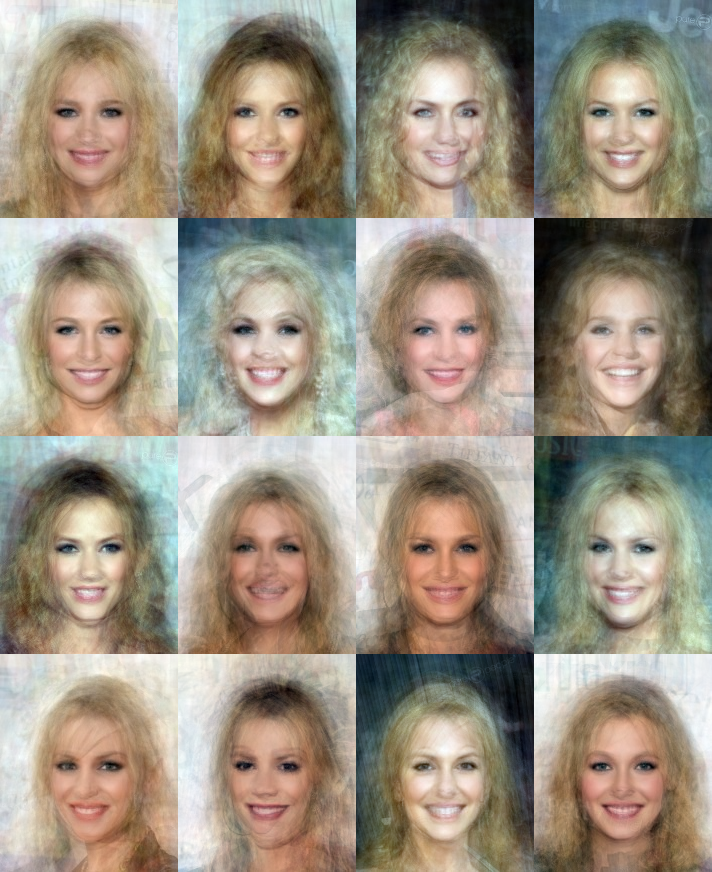}
        \caption{Generated pictures, latent dim 10}
        \label{fig:g48}
    \end{subfigure}
    \hfill
    \begin{subfigure}{0.48\textwidth}
        \includegraphics[width=\textwidth]{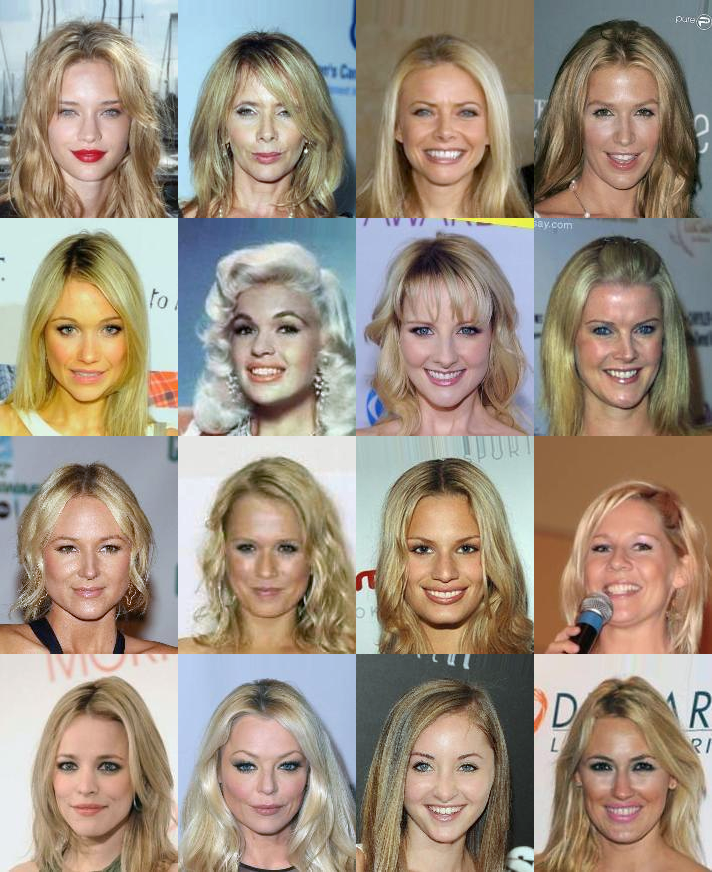}
        \caption{Closest database pictures}
        \label{fig:d48}
    \end{subfigure}
        \caption{Image generation using the Wasserstein-based kernel generative model. 
    The left panel shows images produced by solving the optimal transport problem starting from latent encodings of real images (see Figure~\ref{fig:g47}). 
    The right panel displays their closest counterparts in the training set. 
    Compared to earlier reconstruction results, this method improves both fidelity and diversity, mitigating mode collapse.}
    \label{fig:WassGen}
\end{figure}

This test confirms that Wasserstein-based adversarial training provides improved control over the distribution of generated images. In particular, it reduces the mode collapse phenomenon and results in a generator that balances diversity with realism. The link between the adversarial loss and the optimal transport formulation via the permutation $\sigma_M$ provides a mathematically interpretable mechanism for learning generative models.


\subsection{Conditional image generation and attribute manipulation}

\paragraph{Objective.} We evaluate the capacity of conditional generative models to manipulate semantic image attributes in a controlled fashion. Specifically, we examine the ability of a generator to modify visual attributes--such as presence of glasses or hats--on face images while preserving other characteristics such as identity. This task serves as a proxy to assess the continuity and disentanglement properties of the learned latent representation.

\paragraph{Experimental setup.} We used a filtered subset of the CelebA dataset, selecting 1000 images labeled with the attributes \texttt{[Woman, Light Makeup]} as our base class. We identify a set of images with additional attributes \texttt{[Hat, Glasses]} and manually select four representative samples for manipulation.

The conditional generator is trained on the $1000$ images, conditioned on a two-dimensional attribute vector corresponding to the binary variables \texttt{[Hat, Glasses]}. The latent space comprises a $25$-dimensional standard Gaussian vector. In our test, we fix all components of the latent vector except the two attribute dimensions, which are varied in steps of $0.4$ from $1$ to $0$. The objective is to observe whether smooth interpolation of the attributes leads to semantically coherent image transformations.

\paragraph{Results.} The generated image grid is shown in Figure~\ref{fig:CelebAExample}. The first row corresponds to the original samples with both \texttt{Hat} and \texttt{Glasses} present. Each subsequent row reflects a progressive attenuation of those attributes. The final row (\texttt{[0,0]}) ideally corresponds to versions of the original images without hats or glasses.

\begin{figure}
\centering
\includegraphics[width=0.8\textwidth, keepaspectratio]{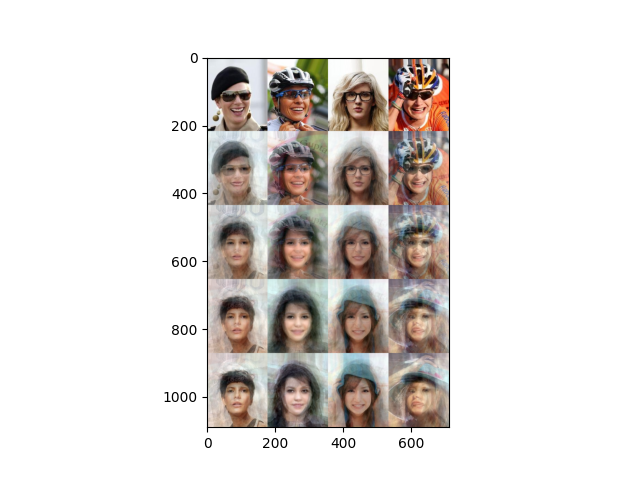}
\caption{\label{fig:CelebAExample} Progressive attribute manipulation from \texttt{[Hat, Glasses] = [1,1]} to \texttt{[0,0]}}
\end{figure}

The results indicate that the model is able to remove the target attributes smoothly. However, image identity is not always preserved, particularly in the final row where the faces sometimes resemble different individuals in the dataset. This behavior highlights a key challenge: while the generator responds to attribute control, over-regularization or latent space collapse may lead to loss of identity fidelity.

We also point out that latent space dimensionality plays a significant role: if it is too low, the generator fails to capture meaningful variation, while if it is too high, the disentanglement deteriorates. In this test, a dimensionality of 25 was found to offer a reasonable compromise via manual tuning.

\subsection{Conditional sampling for data exploration: Iris dataset}
\label{sec:data-exploration}

\paragraph{Objective.} 

This test illustrates how conditional generative modeling can be applied to exploratory data analysis. Using the Iris dataset\footnote{see \cite{fisher1936} and the description in \url{https://archive.ics.uci.edu/dataset/53/iris}}, we study the conditional distribution of selected floral measurements given a fixed feature value. This serves to benchmark various conditional generation methods under a low-data regime, where classical parametric assumptions may fail.

\paragraph{Dataset description.} The Iris dataset consists of 150 samples equally drawn from three species: \textit{Iris setosa}, \textit{Iris versicolor}, and \textit{Iris virginica}. Each observation contains four continuous variables: sepal length, sepal width, petal length, and petal width.

We designate petal width as the conditioning variable $X$, and let the target variable $Y \in \mathbb{R}^3$ represent the remaining three features: petal length, sepal length, and sepal width.

\paragraph{Methodology.} 

Given a specific conditioning value $x_0$ for petal width, we seek to estimate the conditional distribution $Y \mid X = x_0$ and draw samples from it using three generative techniques. 
\begin{itemize}
    \item \textit{Kernel conditional generative model} via the optimal transport approach described in Section~\ref{Conditional-distribution-sampling-model}, with a standard Gaussian prior on the latent space.
    \item \textit{Nadaraya-Watson kernel estimator}, which performs nonparametric regression using the conditioning feature.
    \item \textit{Mixture density networks (MDN)}\footnote{see \cite{bishop1994}}, trained to model conditional density through a mixture of Gaussians.
\end{itemize}

The conditioning value $x_0$ is taken as the empirical mean of petal width over the dataset. Since the dataset contains no sample with exactly this value, we define a reference set for comparison by selecting samples within a small range: those satisfying
\be
|x - \overline{x}| \leq \epsilon \cdot \mathrm{Var}(X),
\ee
with $\epsilon = 0.25$, yielding approximately a dozen reference points.

\paragraph{Results and analysis} 

We generate 500 samples from each conditional model and visualize the resulting univariate marginal cumulative distribution functions (CDFs) against the empirical CDF of the reference set (Figure~\ref{fig:irisexample}). 

\begin{figure}[h]
\centering
\includegraphics[width=0.7\textwidth]{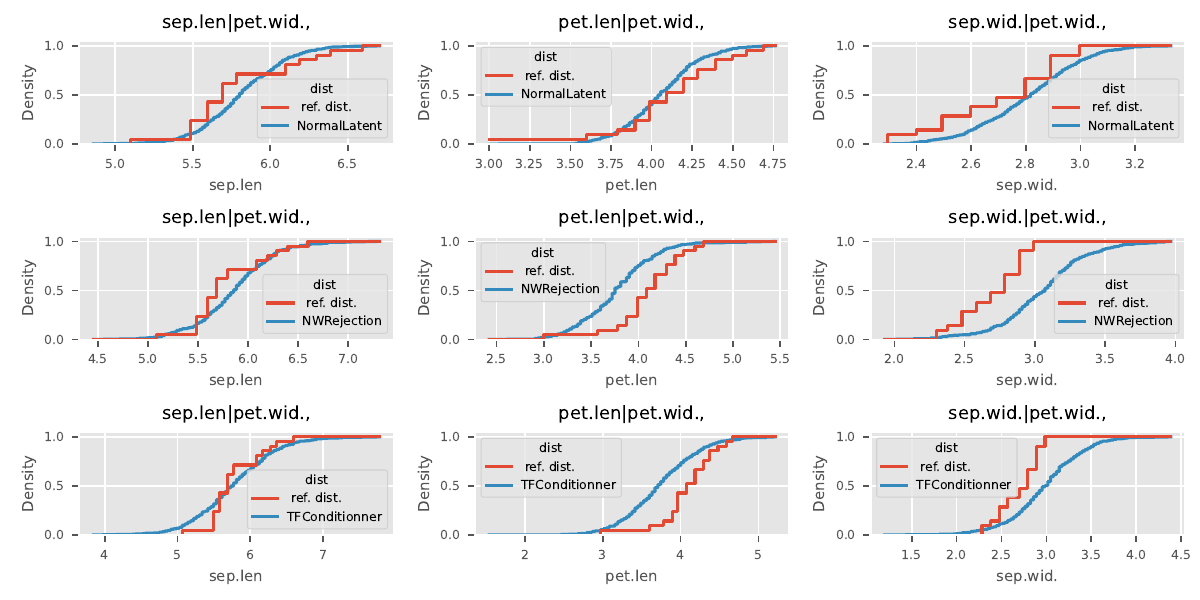}
\caption{Marginal CDFs of generated conditional distributions vs. empirical reference (Iris dataset)}
\label{fig:irisexample}
\end{figure}

To assess higher-order structure, we also plot bivariate marginals and joint sample density visualizations (Figure~\ref{fig:irisbimargexample}) for one of the models.

\begin{figure}[h]
\centering
\includegraphics[width=0.7\textwidth]{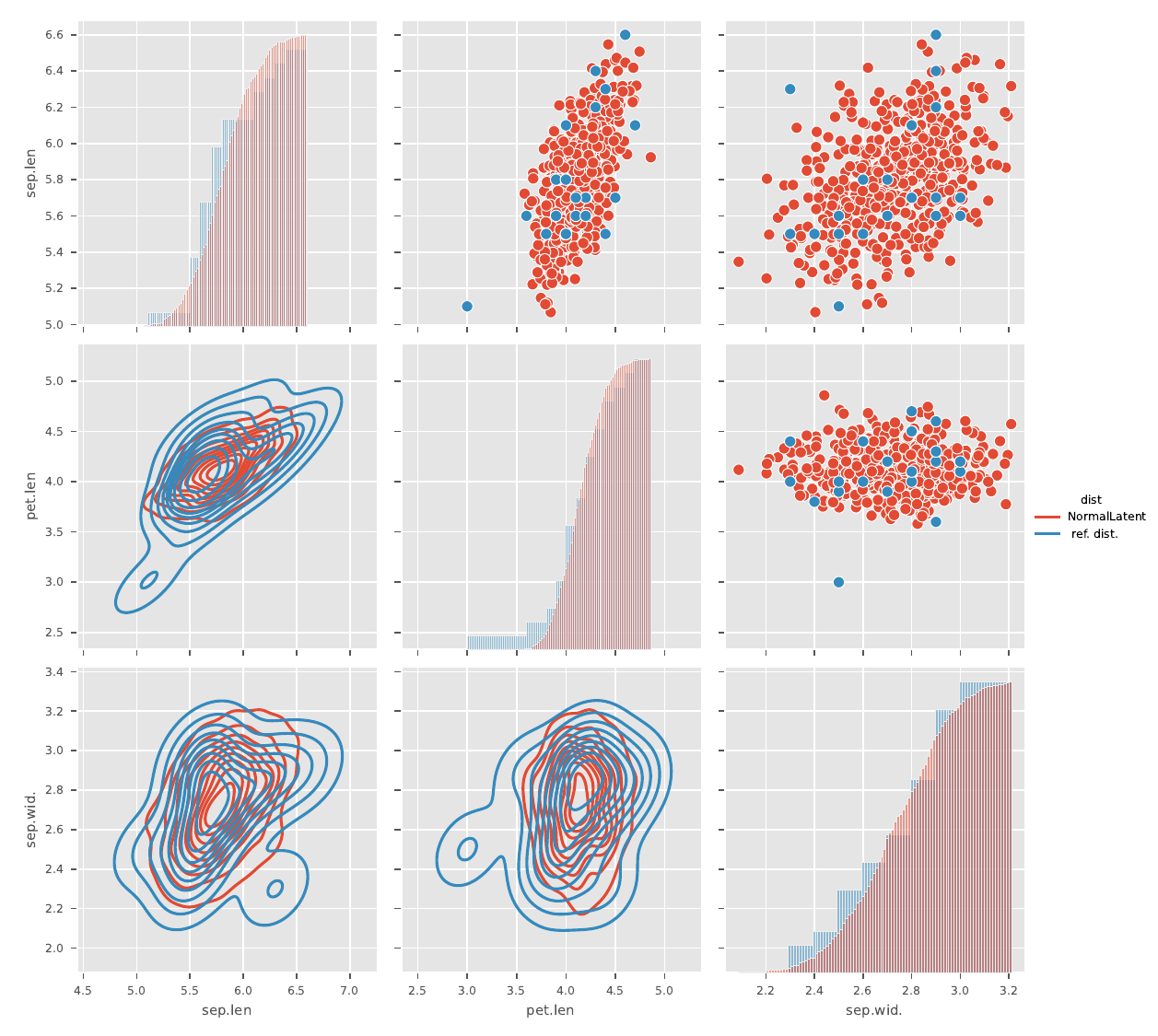}
\caption{Bivariate marginal densities for conditional samples. Center: marginal CDF. Outer: pairwise joint projections.}
\label{fig:irisbimargexample}
\end{figure}

Summary statistics for each generated distribution and its empirical counterpart are presented in Table~\ref{tab:iristab}, including mean, variance, skewness, kurtosis, and the Kolmogorov--Smirnov (KS) distance.

\begin{table}[htbp]
\caption{\label{tab:iristab}Marginal statistics and KS test values for generated conditional distributions. Values in parentheses refer to the reference set.}
\centering
\resizebox{\linewidth}{!}{
\begin{tabular}{l|l|l|l|l|l}
\hline
Method & Mean & Variance & Skewness & Kurtosis & KS Distance \\
\hline
Kernel OT: pet. len & 4.1 (4.1) & --1.2 (0.12) & 0.14 (0.053) & 3 (--0.12) & 0.30 (0.05) \\
Kernel OT: sep. wid & 2.7 (2.8) & --0.5 (--0.13) & 0.048 (0.035) & --0.9 (--0.28) & 0.17 (0.05) \\
Kernel OT: sep. len & 5.8 (5.8) & 0.65 (0.014) & 0.13 (0.073) & 0.15 (0.34) & 0.041 (0.05) \\
NW Estimator: pet. len & 4.1 (3.8) & --1.2 (0.38) & 0.14 (0.17) & 3 (0.72) & 1.7e--6 (0.05) \\
NW Estimator: sep. wid & 2.7 (3.0) & --0.5 (--0.024) & 0.048 (0.12) & --0.9 (--0.23) & 6.8e--6 (0.05) \\
NW Estimator: sep. len & 5.8 (5.8) & 0.65 (0.20) & 0.13 (0.17) & 0.15 (0.68) & 0.27 (0.05) \\
MDN: pet. len & 4.1 (3.7) & --1.2 (--0.29) & 0.14 (0.27) & 3 (0.24) & 6.6e--6 (0.05) \\
MDN: sep. wid & 2.7 (3.1) & --0.5 (--0.024) & 0.048 (0.19) & --0.9 (0.37) & 8.5e--7 (0.05) \\
MDN: sep. len & 5.8 (5.8) & 0.65 (0.19) & 0.13 (0.27) & 0.15 (0.055) & 0.19 (0.05) \\
\hline
\end{tabular}}
\end{table}

The results show that conditional generative methods can approximate local structure in the target distribution even when conditioning on a value not present in the dataset. While none of the methods pass statistical tests perfectly due to the small sample size and arbitrary reference window, the generated samples exhibit plausible empirical behavior. Kernel-based and MDN methods perform comparably, with the Nadaraya--Watson estimator showing slightly lower sample diversity.

\subsection{Data completion via conditional generative modeling}
\label{sec:data-completion}

\paragraph{Objective.} 

In this test, we assess the ability of conditional generative models to perform data completion and synthetic data generation. We focus on reconstructing realistic data samples for a missing class within a labeled dataset. This task serves as a practical benchmark for conditional generative methods under partial observation and class imbalance.

\paragraph{Dataset description.} 

We utilize the Breast Cancer Wisconsin (Diagnostic) dataset\footnote{see \cite{street1993} and the description in \url{https://archive.ics.uci.edu/dataset/17/breast+cancer+wisconsin+diagnostic}}, a classical benchmark in binary classification. The dataset consists of 569 instances, each represented by 30 continuous-valued diagnostic features derived from cell nuclei images. Each sample is labeled as either malignant (212 entries) or benign (357 entries).

For simplicity, we focus on the first four features: \texttt{[mean radius, mean area, mean perimeter, mean texture]}, which form a four-dimensional feature vector $Y \in \mathbb{R}^4$. The binary class label $C \in \{0,1\}$ serves as the conditioning variable, where $C = 1$ denotes malignant.

\paragraph{Experimental setup.} 

To simulate a partial data scenario, we split the malignant class into two equal halves: 106 samples are combined with all benign entries (totaling 463 samples) to form the training set. The remaining 106 malignant entries are held out as a reference for evaluation.

A conditional generator is trained to learn the distribution $Y \mid C$, using the training set. We then sample 500 synthetic instances conditioned on $C = 1$ (malignant), and compare the generated distribution to the withheld malignant class using both visual and statistical tools.

\paragraph{Results and analysis} 

Figure~\ref{fig:breastexample} visualizes the conditional samples using marginal and bivariate projections. The central diagonal displays univariate cumulative distribution functions (CDFs), while off-diagonal plots show joint projections.

\begin{figure}[h]
\centering
\includegraphics[width=0.7\textwidth]{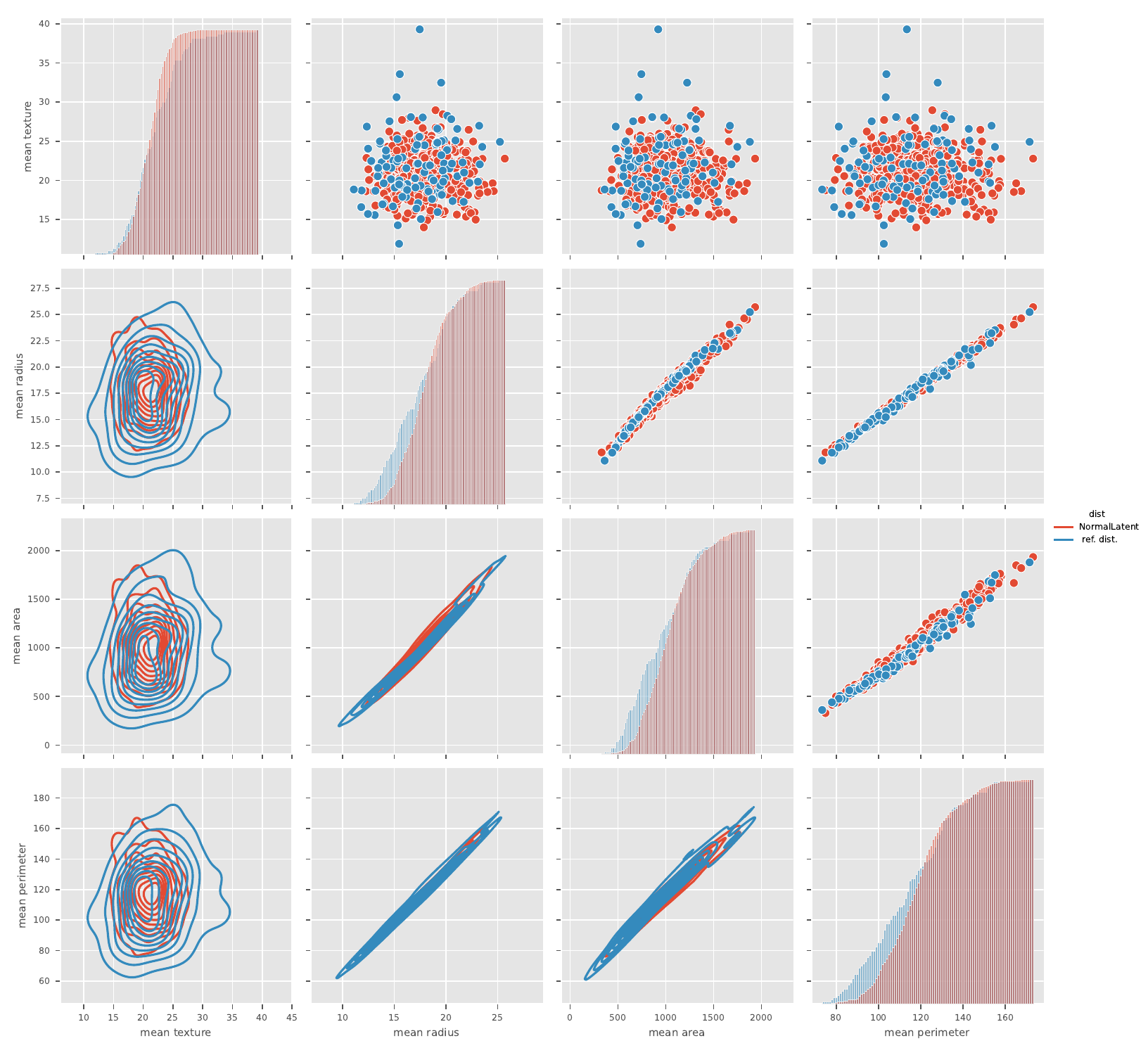}
\caption{\label{fig:breastexample}Conditional generation for the malignant class in the Breast Cancer dataset. Center: CDFs of marginals. Outer: joint projections.}
\end{figure}

Table~\ref{tab:breasttab} shows the mean, variance, skewness, kurtosis, and Kolmogorov--Smirnov (KS) statistics comparing generated samples to the withheld real malignant data.

\begin{table}[htbp]
\caption{\label{tab:breasttab}Marginal statistics and KS test results for synthetic vs. real malignant data. Values in parentheses correspond to the real data.}
\centering
\resizebox{\linewidth}{!}{
\begin{tabular}{l|l|l|l|l|l}
\hline
Feature & Mean & Variance & Skewness & Kurtosis & KS Distance \\
\hline
mean radius & 17 (18) & 0.13 (0.32) & 9.1 (5.3) & --0.56 (0.1) & 0.0042 \\
mean area & $9.6 \times 10^2$ ($1.0 \times 10^3$) & 0.47 (0.49) & $1.1 \times 10^5$ ($7.3 \times 10^4$) & --0.24 (0.14) & 0.0028 \\
mean perimeter & $1.1 \times 10^2$ ($1.2 \times 10^2$) & 0.20 (0.34) & $4.2 \times 10^2$ ($2.6 \times 10^2$) & --0.47 (0.23) & 0.011 \\
mean texture & 22 (21) & 0.92 (0.034) & 18 (6.6) & 2.3 (--0.12) & 0.0019 \\
\hline
\end{tabular}}
\end{table}

The conditional generator demonstrates strong potential for reconstructing the statistical properties of the malignant class, even though it was only partially observed during training. The generated samples exhibit plausible marginal statistics, and KS test distances remain acceptably low, suggesting that the model captures meaningful structural information.

However, we observe sensitivity to the kernel used during training. In particular, we found that a ReLU-based kernel led to better performance in this case. For heavy-tailed distributions, alternative kernels (e.g., Cauchy) may offer improved accuracy.

This test highlights the usefulness of generative models for synthetic data augmentation and imputation, particularly in medical or diagnostic contexts where certain classes may be underrepresented or protected due to privacy constraints.


\section{Large-scale dataset} 
\label{large-scale-dataset}

\subsection{Reproducible kernel ridge regression for large dataset} 
\label{large-scale-dataset-kr}

\paragraph{Problem setup} 

We consider the kernel ridge regression (KRR) formulation \eqref{FIT} and restrict attention to the reproducible setting, where $X = Y$ and $\epsilon = 0$. The algorithmic and memory costs in this setting become prohibitive for large datasets: computing the Gram matrix requires $\mathcal{O}(N_x^2)$ operations, and solving the resulting system involves $\mathcal{O}(N_x^3)$ complexity. To mitigate these costs, one often introduces a smaller representative set $Y \subset X$, at the expense of losing the reproducibility property.

\paragraph{Multiscale strategy} 

In this test, we propose a \textit{divide-and-conquer} strategy that enables scalable and reproducible kernel methods on large datasets. Our approach interprets the parameter $N_y$ in \eqref{FIT} not just as a reduction parameter, but as defining a set of computational units, each associated with a \textit{centroid} $y^n$, such that $Y = (y^1, \dots, y^{N_y})$. These centroids are used to partition the data space via hard clustering: each input point is assigned to exactly one cluster. Let $l : \mathbb{R}^D \to \{1, \dots, N_y\}$ denote the cluster assignment function, as defined in Section~\ref{assignment}.

The overall construction proceeds as follows. Given a clustering $Y$ and assignment map $l(\cdot)$, we first choose a global (coarse) kernel $k_0(\cdot, \cdot)$ and apply the KRR formula \eqref{FIT} using the centroids $Y$ as the training set. This yields a coarse approximation $f_{k_0, \theta}(\cdot)$ and defines a residual error function:
\be
  \epsilon(\cdot) = f(\cdot) - f_{k_0,\theta}(\cdot).
\ee
To refine this approximation, we introduce a second layer of local corrections, and define the full extrapolation operator via:
\be \label{MS}
  f_{k_0, \ldots, k_N, \theta}(\cdot) = f_{k_0,\theta}(\cdot) + \epsilon_{k_{l(\cdot)},\theta}(\cdot),
\ee
where each kernel $k_n$ (for $n = 1, \dots, N_y$) is applied within its corresponding cluster using a local model fitted to the residual error.

This construction yields a hierarchical, two-level model. Each cluster contains approximately $\frac{N_x}{N_y}$ points, allowing the local KRR problems to remain tractable. While the method introduces a new hyperparameter $N_y$, it offers a natural and parallelizable framework for designing scalable kernel methods that retain reproducibility and allow for meaningful error control.

Achieving uniform cluster sizes is important to balance computational load and ensure model consistency. This motivates the use of \textit{balanced clustering} strategies, as described in Section~\ref{balanced-clustering}.

Once the data are partitioned, each computational unit can be run independently--ideally in parallel across $N_y$ threads--but can also be executed sequentially or in a hybrid parallel-sequential setup depending on the available hardware.

From a complexity standpoint, this strategy leads to a reproducible KRR algorithm with overall linear complexity in both input and output sizes, as discussed earlier (albeit with a larger constant factor). Naturally, this gain in efficiency comes at the cost of approximation accuracy. The main source of error in kernel extrapolation can be linked to the distance between input points and the training set. In our multiscale setting, the leading error term becomes
\be
  d_{k_0}(\cdot, Y) \cdot d_{k_{l(\cdot)}}(\cdot, X_{l(\cdot)}),
\ee
where $X_i \subset X$ denotes the data points assigned to the $i$-th cluster. Both components are computable and provide a handle on the total approximation error.

From an engineering perspective, this framework allows for further extensions. The initial kernel $k_0$ may be replaced by simpler, structured regressors--such as polynomial basis functions--to capture global trends or low-order moments. Additional layers can be introduced in a tree-like or hierarchical fashion to refine the approximation across scales. Overlapping clusters can also be used: in this case, one kernel captures a specific feature (e.g., local geometry), while another specializes in orthogonal aspects, effectively acting as a sequence of filters.

More generally, the multiscale construction introduced here lends itself to modular designs based on oriented graphs, enabling adaptive, interpretable architectures for large-scale kernel learning.

\paragraph{Numerical illustration: MNIST classification}\label{multiscale-supervised-learning with MNIST}

We now illustrate the behavior of the multiscale method \eqref{MS} on the MNIST dataset, as introduced in Section~\ref{classification-problem-handwritten-digits}. We present here the classification accuracy and execution times to assess the trade-off between computational efficiency and predictive performance.

In this test, we varied the number of clusters $N_y$ from 5 to 80 and applied several clustering strategies to partition the dataset. The model was evaluated in full extrapolation mode on each cluster. Since balanced clusters are essential for consistent performance, we applied the balanced clustering method described in \eqref{balanced} to compute the cluster assignments. For comparison, we also included a random cluster assignment baseline to assess sensitivity to the clustering scheme.

Results are summarized in Figure~\ref{fig:otmsperf}, where the left panel shows classification accuracy and the right panel presents corresponding execution times.
\begin{itemize}
\tightlist
\item
  In terms of predictive accuracy, no clustering method consistently outperforms the others. The balanced clustering strategy \eqref{balanced} tends to neutralize variability across methods, and even randomly assigned clusters produce respectable scores.
\item
  A smaller number of clusters generally results in better accuracy but comes at a higher computational cost. In particular, using five clusters approaches kernel saturation, leading to improved accuracy but significantly longer runtimes.
\item
  The sharp discrepancy method requires computation of the full Gram matrix at initialization. This step alone took over 100 seconds and dominated the overall execution time for that algorithm.
\end{itemize}

\begin{figure}
\hypertarget{fig:otmsperf}{%
\centering
\includegraphics[width=1.\textwidth]{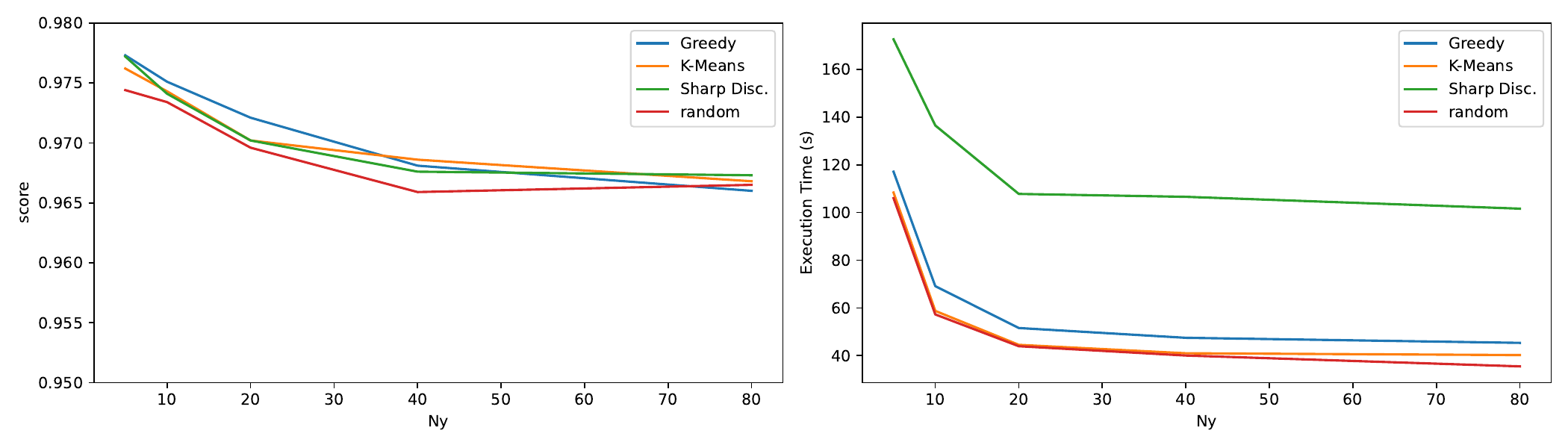}
\caption{Performance of the multiscale method applied to MNIST. Left: classification scores. Right: execution times.}\label{fig:otmsperf}
}
\end{figure}

These results confirm that the proposed multiscale method effectively trades accuracy for computation time, with manageable performance degradation. Importantly, the method demonstrates robustness with respect to the choice of clustering strategy, allowing considerable flexibility in practical implementations. A single computing unit--such as a standard laptop CPU--can handle throughput ranging from approximately 1,000 to 10,000 samples per second, depending on the dataset and kernel parameters.

\subsection{Multiscale strategies for Monge optimal transport on large datasets}\label{multiscale-ot}

\paragraph{Problem setup} The combinatorial approach to optimal transport, as explained in Section~\ref{optimal-transport-and-statistical-kernel-methods}, offers an attractive alternative to entropy-regularized methods such as Sinkhorn--Knopp. However, its reliance on linear sum assignment problem (LSAP) solvers makes it computationally infeasible for large datasets due to their super-linear complexity.

To overcome this limitation, we adapt the multiscale strategy introduced in Section~\ref{large-scale-dataset-kr} to the discrete Monge problem. This adaptation enables approximate solutions to large-scale OT problems via a divide-and-conquer approach, which we now describe and illustrate numerically.

\paragraph{Multiscale Monge problems}

Consider the discrete Monge problem \eqref{LSAP}, where $X, Y \in \mathbb{R}^{N \times D}$ represent empirical distributions of size $N$ in $D$-dimensional space. We begin by selecting two sets of centroids $C_X, C_Y \in \mathbb{R}^{M \times D}$, one for each distribution, using a clustering method. Let $\sigma_X(\cdot)$ and $\sigma_Y(\cdot)$ denote the associated assignment functions that map each point in $X$ and $Y$ to a cluster index in $\{1, \ldots, M\}$, as introduced in Section~\ref{assignment}. Let the resulting clusters be defined as
\be
X_i = \{x \in X : \sigma_X(x) = i\}, \qquad Y_j = \{y \in Y : \sigma_Y(y) = j\}.
\ee
Our goal is to solve a collection of smaller OT problems between pairs $(X_i, Y_{\sigma(i)})$, thereby approximating the global Monge transport. This requires defining a correspondence $\sigma(i)$ between clusters in $C_X$ and $C_Y$, which we now formalize.

Let $M$ denote the number of clusters. We assume that the clusters are of equal size--a requirement inherited from the kernel ridge regression formulation \eqref{FIT}. Ensuring this balance necessitates the use of perfectly balanced clustering methods, as discussed in Section~\ref{balanced-clustering}.

Given equal-cardinality clusters, we define the optimal assignment of clusters via a combinatorial optimization problem. Let $d(\cdot, \cdot)$ be a generic distance function (e.g., Euclidean distance or discrepancy $d_{k_0}$, where $k_0$ is the coarse-level kernel used in \eqref{MS}). Let $\Sigma$ denote the set of permutations on $\{1, \ldots, N\}$. Define the cost matrices:
\be
D_x = d(X, C_X), \quad D_y = d(Y, C_Y), \quad C = d(C_X, C_Y),
\ee
with sizes $(N, M)$, $(N, M)$, and $(M, M)$, respectively. We then seek permutations $\overline{\sigma}_X, \overline{\sigma}_Y \in \Sigma$ minimizing the total transport cost:
\be \label{KM}
  \overline{\sigma}_X, \overline{\sigma}_Y = \mathop{\arg\min}_{\sigma_X, \sigma_Y \in \Sigma} \sum_{n} D_x(n, \sigma_X^{n \bmod M}) + \sum_{n} D_y(n, \sigma_Y^{n \bmod M}) + \sum_n C(\sigma_X^{n \bmod M}, \sigma_Y^{n \bmod M}).
\ee

\noindent
Here, we use the shorthand notation  $\sigma_X^{n \bmod M}  =  \sigma_X(n \bmod M)$, and similarly for $\sigma_Y$, to improve readability and reduce the horizontal length of the expression. Once solved (or approximated), this problem yields the following. 
\begin{itemize}
\item Assignment functions $l_X(\cdot) = \overline{\sigma}_X^{\arg\min_n d(X^n, \cdot) \bmod M}$ and similarly \newline $l_Y(\cdot) = \overline{\sigma}_Y^{\arg\min_n d_{k_0}(Y^n, \cdot) \bmod M}$, associating each point to a centroid in $C_X$ or $C_Y$, respectively.

\item A pairing of clusters $C_X^{\overline{\sigma}_X^n \bmod M} \mapsto C_Y^{\overline{\sigma}_Y^n \bmod M}$, which defines the inter-cluster correspondence.
\end{itemize}

This allows us to define a cluster-level assignment $i \mapsto \sigma(i)$, and approximate the transport map analogously to the supervised multiscale method \eqref{MS}. Specifically, we write:
\be \label{MSOT}
  Y_{k_0, \cdot, k_N, \theta}^{\sigma}(\cdot) = Y_{k_0, \theta}^{\sigma}(\cdot) + \epsilon_{Y^\sigma_{l(\cdot)}, \theta}(\cdot), \quad \text{where } \epsilon(\cdot) = Y^\sigma - Y_{k_0, \theta}^{\sigma}(\cdot).
\ee
Since the dataset is partitioned, the resulting map approximates the global optimal transport while maintaining exact invertibility within each cluster. This forms the basis for an efficient and scalable OT algorithm suitable for large datasets.

\hypertarget{numerical-comparison}{%
\paragraph{Numerical benchmark}\label{numerical-comparison}}

We now present a numerical benchmark which compares our multiscale optimal transport method with the Sinkhorn--Knopp algorithm. Our benchmarking protocol follows the methodology of Pooladian and Niles-Weed\footnote{see \cite{NiPa:2021}; we used the code made publicly available by the authors.}.
The setup is as follows. We generate a distribution $X \in \mathbb{R}^d$, and define $Y = S(X)$, where $S = \nabla h$ is the gradient of a convex function $h$, so that $S$ defines an exact Monge map. We then shuffle $Y$ and provide the input pair $(X, Y)$ to various methods that estimate the transport map $S_{k,\theta}(\cdot)$. Performance is measured using the mean squared error (MSE) between the predicted and exact transport on a test set $Z$ sampled from the same distribution as $X$, computed as $|S_{k,\theta}(Z) - S(Z)|_{\ell^2}$.

For this test, we set $X$ to be uniformly distributed, and define the exact transport as $S(X) = X |X|_2^2$. The results are presented in Figure~\ref{fig:otperf} (MSE vs. dataset size) and Figure~\ref{fig:ottime} (execution time). We test various problem dimensions $d = \{2, 10, 100\}$, and dataset sizes $N = \{256, 512, \ldots, 4096\}$.

Five methods are benchmarked. 
\begin{itemize}
\tightlist
\item
  {COT}, {COT Parallel}, and {COT MS}: These are combinatorial optimal transport (COT) methods based on \eqref{ED} and Algorithm~\ref{alg1}. The base method ({COT}) uses an exact LSAP solver in serial. {COT Parallel} uses a parallelized but sub-optimal LSAP solver. {COT MS} is the multiscale variant using the assignment strategy from \eqref{KM}, with the number of clusters set to $C = N / 256$, chosen so that {COT MS} matches {COT} exactly in the smallest case ($N = 256$).
\item
  {POT} and {OTT}: These use Sinkhorn-regularized optimal transport via the publicly available {POT} and {OTT} libraries. Because the Sinkhorn method is sensitive to the regularization parameter, we manually tuned this value to ensure convergence. {OTT} provides an automatic heuristic, whereas {POT} required trial-and-error tuning.
\end{itemize}

Our main conclusions are as follows.
\begin{itemize}
\item
  \textit{Performance:} As expected, the exact combinatorial methods outperform entropy-regularized approaches ({POT}, {OTT}) in terms of transport accuracy across all problem sizes and dimensions.
\item
  \textit{Scalability:} The main limitation of combinatorial methods lies in computational cost. {POT}, in particular, benefits from GPU parallelization and is highly efficient at scale.
\item
  \textit{Asymptotic efficiency:} Among the combinatorial variants, {COT MS} is the most scalable. Its cost grows linearly with dataset size and remains competitive while delivering exact transport locally. However, this efficiency comes with a trade-off in accuracy, especially in low-dimensional settings.
\end{itemize}

It is worth emphasizing that Sinkhorn-based algorithms are designed to solve regularized optimal transport problems, which are relevant in various applications due to their stability and differentiability. However, for solving the unregularized Monge problem--as considered here--combinatorial methods achieve significantly higher accuracy, especially on small- to medium-scale datasets.

\begin{figure}
\hypertarget{fig:otperf}{%
\centering
\includegraphics[width=0.9\textwidth]{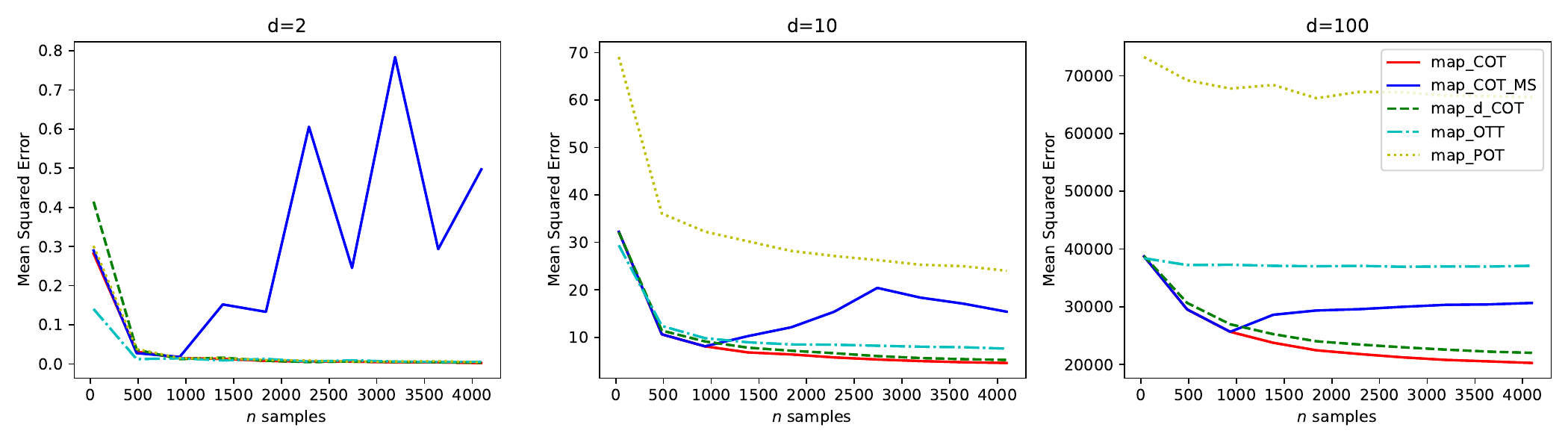}
\caption{Mean squared error (MSE) for various optimal transport methods across increasing dataset sizes and dimensions.}
\label{fig:otperf}
}
\end{figure}

\begin{figure}
\hypertarget{fig:ottime}{%
\centering
\includegraphics[width=1.0\textwidth]{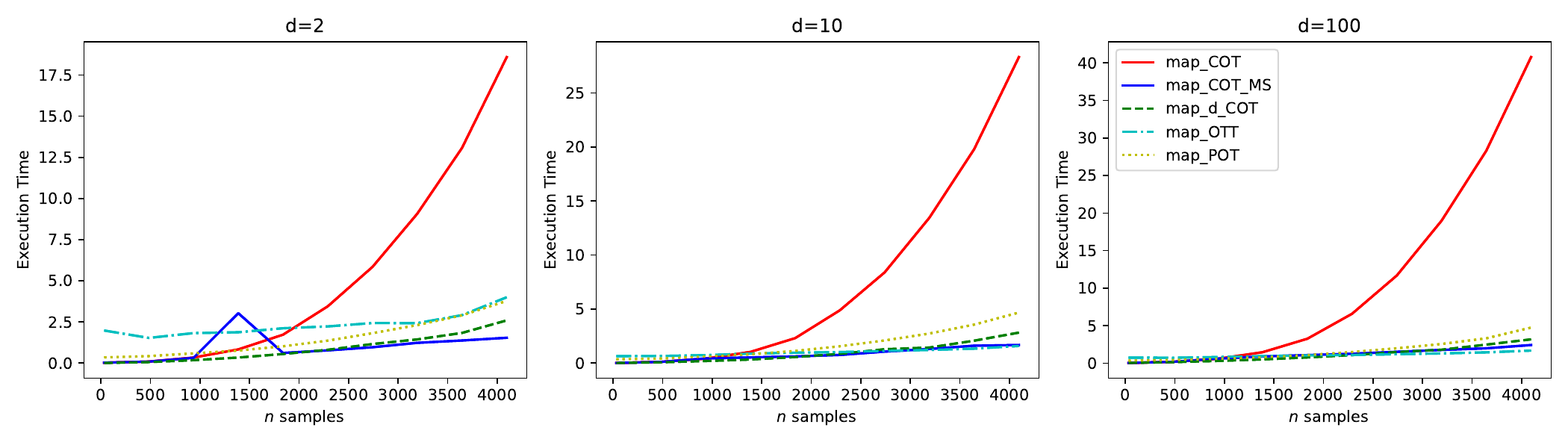}
\caption{Execution time of each method as a function of dataset size. {POT} and {OTT} are regularized approaches using Sinkhorn iterations; {COT} methods are combinatorial.}
\label{fig:ottime}
}
\end{figure}

\paragraph{Qualitative illustration of the multiscale OT approximation}
\label{remark-on-the-multiscale-optimal-transport-approach}

Until now, we benchmarked a multiscale approach for solving the Monge optimal transport problem between a source distribution $X$ and a target distribution $Y$. For large sample sizes $N$, solving the LSAP directly on the full datasets becomes computationally intractable. The multiscale strategy mitigates this by first partitioning both distributions into $C$ clusters via balanced clustering. Each cluster pair is treated as an independent subproblem, allowing for OT computations to be performed efficiently--either in parallel or sequentially--and later recombined into a global transport map.

This approach significantly reduces computation time but introduces an approximation error. To qualitatively assess this error, we conduct the following test.

We generate source and target distributions $X$ and $Y$ as samples from Gaussian mixtures with known means and variances. We use $N = 1024$ samples and choose $C = 2$ clusters. As a reference, we compute the exact OT map using LSAP between $X$ and $Y$ based on pairwise Euclidean distances. We then apply the multiscale OT method using the same distance metric and the balanced clustering procedure described earlier.

Figure~\ref{fig:LSAPOT} shows the results of both approaches. Points from $X$ are shown in red, points from $Y$ in blue, and black lines represent the computed assignments. The left plot shows the exact LSAP solution, where the assignments follow a globally optimal, monotonic structure. The right plot shows the multiscale approximation, where assignments are computed within clusters and then aggregated. Notably, the assignments in the multiscale method exhibit line crossings--something that should not occur under exact Euclidean optimal transport. This highlights the local, cluster-wise nature of the approximation.

\begin{figure}
\centering
\includegraphics[width=0.9\textwidth]{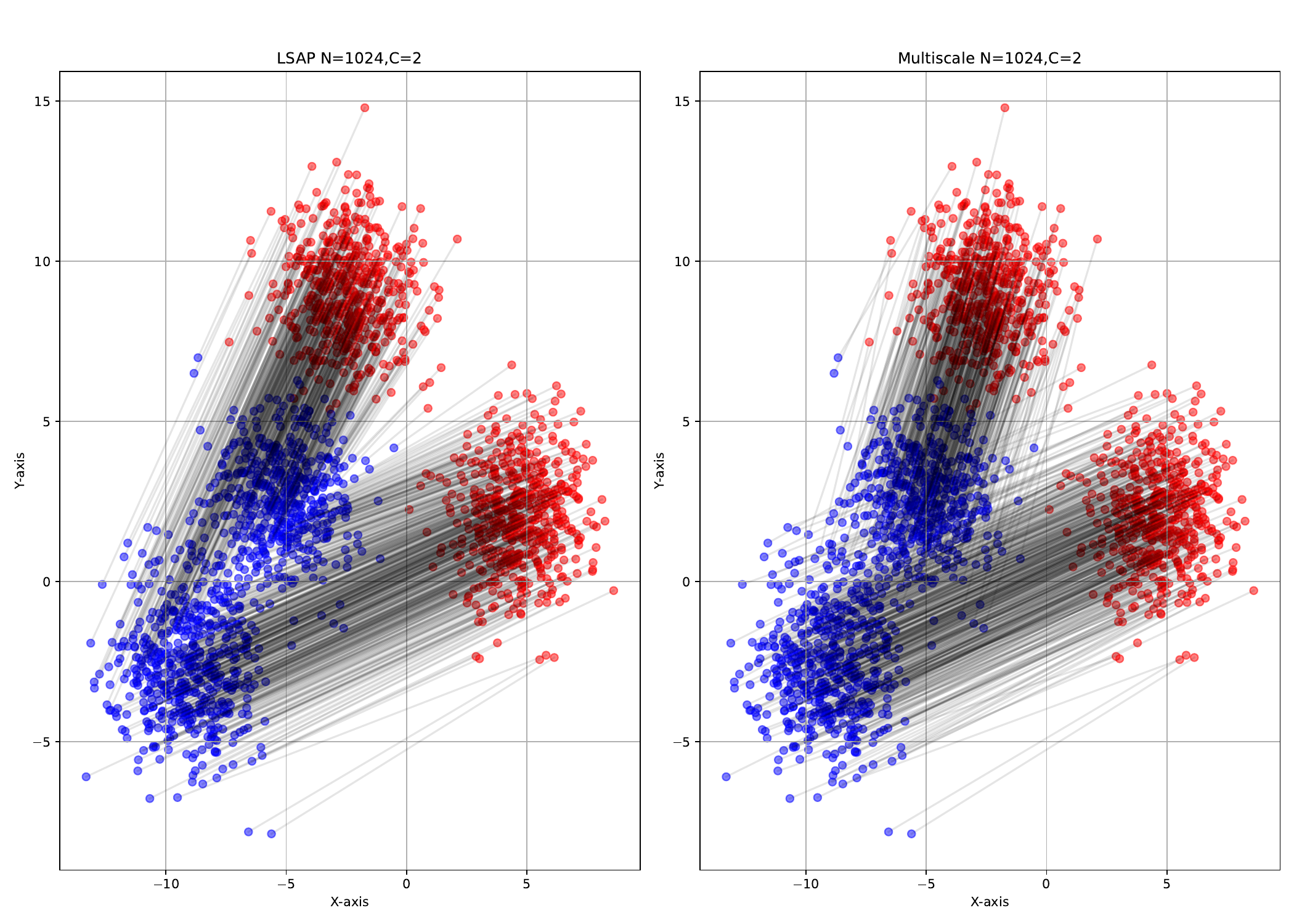}
\caption{\label{fig:LSAPOT}Comparison of exact OT (left) and multiscale OT (right). Red and blue points correspond to source ($X$) and target ($Y$) distributions, respectively. Black lines denote the computed transport assignments. In the multiscale case, assignment lines cross, illustrating the approximation error induced by localizing the transport within clusters.}
\end{figure}


\chapter{Application to physics-informed modeling}
\label{application-to-partial-differential-equations}

\section{Introduction}
\label{PDEs}

We now consider physics-informed systems as maps between data and observables (see Section~\ref{physics-informed-systems} below), where these maps are determined by a physical model. Most of these models are defined through partial differential equations (PDEs). So, the purpose of this chapter is to provide numerical, RKHS-based techniques, in order to approximate solutions to PDEs. We demonstrate here that the approach we propose offers some advantages over traditional numerical methods for PDEs. Moreover, due to the natural bridge between RKHS methods and operators; see Section~\ref{kernel-based-operators}, most of numerical analysis approaches to PDEs can be effortlessly adapted to RKHS frameworks, as for instance the approximation of Green operators or or time-dependent operators, also called generators; see~Section~\ref{time-evolution-operators-based-on-theta-schemes}.

\begin{itemize}
\tightlist
\item
  \textit{Mesh-free methods.} Kernel methods allow for mesh-free (sometimes called mesh-free) formulations to be used. This case corresponds where $X=(x^1,\ldots,x^N)$, a given matrix of data, is called the mesh, and the unknown is $F(X)$. Unlike traditional finite difference or finite element methods, mesh-free methods do not require a predefined mesh, nor to compute connections between nodes of the grid points. Instead, they use a set of nodes or particles to represent the domain. This makes them particularly useful for modeling complex geometric domains.
\item
  \textit{Particle methods.} Kernel methods can be used in the context of particle methods in fluid dynamics, which are Lagrangian methods involving the tracking of the motion of particles. This situation corresponds to cases where $X$ represents a probability $\delta_X$ approximating a probability density $d\mu$, which is the unknown of a physical system. Kernel methods are well-suited for these types of problems because they can easily handle general meshes and boundaries.
\item
  \textit{Boundary conditions.} Kernel methods allow one to express complex boundary conditions, which can be of Dirichlet or Neumann type, or even of more complex mixed-type expressed on a set of points. They also can also encompass free boundary conditions for particle methods, as well as fixed meshes.
  \item \textit{Error analysis.}  Kernel approaches to PDE allows also for standard numerical error analysis, using~\eqref{estim}, which is the only universal regression technic having this property.
\end{itemize}

We are going to provide several illustrations of the flexibility of this approach. The price to pay with mesh-free methods is the computational time, which is greater than the one in more traditional methods such as finite difference, finite element, or finite volume schemes: the RKHS approach usually produces dense matrices, whereas more classical methods on structured grids, due to their localization properties, typically lead to sparse matrix, a property that matrix solvers can benefit on. However, traditional methods have difficulties to tackle high-dimensional data, due to the \textit{curse of dimensionality}\footnote{see curse of dimensionality, \url{https://en.wikipedia.org/wiki/Curse_of_dimensionality}, which inspired us for the CodPy Library}, where the RKHS approach usually excels\footnote{Others universal regressor methods, as neural networks, can also be adapted to this purpose, with a very similar numerical analysis; see for instance the recent work \cite{EF:2025}. Pros, cons, and benchmarks to both approaches remains to be done.}.

In this chapter, we start presenting a series of simple examples, commencing with static models to illustrate our purpose. We then initiate our discussion with some of the technical details relevant to the discretization of partial differential equations via kernel methods, and progress to encompass a spectrum of time evolution equations. Our primary goal is to showcase and the efficacy and broad applicability of mesh-free methods, in the context of, both, structured and unstructured meshes.


\section{Physics-informed modeling}
\label{physics-informed-systems}

A \textit{physics-informed model} is an invertible mapping between a training dataset ($X,f(X)$) and observable parameters  $\varepsilon$:
\be\label{MF}
    \mathcal{L}(X, F(X), \nabla F(X), \ldots) = \varepsilon,
\ee
where
\begin{itemize}
  \item \( \varepsilon \) are parameters associated with the model, as measured temperature or pressure on a precise location, which may differ in shape from the original data, 
  
  \item and \( \mathcal{L}(\ldots)\) is a mapping defined by a physical model, supposed to be invertible, with inverse denoted by 
  \be\label{FM}
      X,F(X),\ldots = \mathcal{L}^{-1}(\varepsilon).
  \ee
\end{itemize}
\noindent This structure can be composed: given maps \( \mathcal{L}_1, \mathcal{L}_2 \), with \( \mathrm{Im}(\mathcal{L}_2) \subseteq \mathrm{Dom}(\mathcal{L}_1) \), their composition \( \mathcal{L}  =  \mathcal{L}_1 \circ \mathcal{L}_2 \) also defines a model with inverse \( \mathcal{L}^{-1} = \mathcal{L}_2^{-1} \circ \mathcal{L}_1^{-1} \), which allows usually to consider elementary bricks to build more sophisticated models. Specific examples can be found in this monograph for time-series analysis; see Section~\ref{free-time-series-modeling}.


\section{Mesh-free methods}
\label{MM}

\subsection{Poisson equation}\label{poisson-equation}

\paragraph{Continuous analysis} 

We illustrate the flexibility of the mesh-free method with two numerical tests.
First, we treat the Laplace-Beltrami operator. For this purpose, we introduce a weighted probability measure on $\mathbb{R}^D$,  $d\mu = \mu(x)dx$, where $dx$ the Lebesgue measure and $\mu \in L^{\infty}$, $\mu(\mathbb{R}^D)=1$, $\mu \ge 0$ on a possibly unbounded, Lipschitz continuous domain $\Omega = \text{supp} \hskip.12cm \mu $. The Poisson equation corresponds to the functional
\be
    J(u) = \frac{1}{2} \int |\nabla u|_2^2 d \mu - \int u \ f d \mu, 
\ee
defined over the weighted Sobolev space $H^1_\mu = \{u : \int (|\nabla u|_2^2 + u^2) d \mu < +\infty\}$. Poisson equation usually considers measures on sets, that is, $d\mu = 1_{\Omega}(x)dx$, where $1_{\Omega}(\cdot)$ is the indicator function of $\Omega$, but this analysis extends to more general measures.

Consider minimizing this functional, which means finding $u = \arg \min_{v \in H^1_\mu } J(v)$. A minimum to this functional yields the following Poisson equation. Namely, minimizing $J$ over $H_\mu^1(\Omega)$ yields the weighted Poisson problem in divergence form:
\begin{align*}
  \text{Find } u\in H_\mu^1(\Omega)\ \text{ such that }\ 
  \int_\Omega \mu\,\nabla u\cdot\nabla \varphi\,dx  =  \int_\Omega \mu\,f\,\varphi\,dx
  \quad \forall\,\varphi\in H_\mu^1(\Omega).
\end{align*}
The corresponding strong form is 
\be
  \nabla \cdot (u \nabla\mu)  = f\ \mu, \quad \text{ supp } u  \subset \Omega, \quad u_{\partial \Omega}=0,
\ee
where \(f\) is sufficient regular and \(\Omega\) is a sufficient regular domain, which is to be considered in the weak sense~\eqref{WFPE}.  We can rewrite it in equivalent form using weighted Laplacian:
\be
\Delta_{\mu} u = \mu^{-1} \nabla \cdot (u \nabla\mu)
\ee
and the PDE is given by
\be
-\Delta_\mu u = f.
\ee

\paragraph{RKHS discretization} 

To compute an approximation of this equation with a kernel method we proceed as follows. 
\begin{itemize}
\tightlist
\item
  Select a mesh \(X \in \mathbb{R}^{N_x,D}\) representing \(\Omega\).
\item
  Choose a kernel \(k\) that generates a discrete RKHS space $\mathcal{H}_{k,X}$ of null trace functions.
\item Consider minimizing
\be
    \inf_{u\in \mathcal{H}_{k,X}}J_k(u), \quad J_k(u) = \frac{1}{2N} \sum_{n=1}^N |\nabla_k u(x^n)|_2^2  - \frac{1}{N}\sum_{n=1}^N (f u)(x^n).  
\ee
\end{itemize}

A kernel approximation of this equation consists in approximating the solution as a function \(u \in \mathcal{H}_k\), that is, in the finite dimensional reproducing kernel Hilbert space generated by the kernel \(k\) and the set of points \(X\), satisfying
\be
  <\nabla_k u ,\nabla_k \varphi>_{\mathcal{H}_k}=-<\Delta_k u , \varphi>_{\mathcal{H}_k} = <f, \varphi>_{\mathcal{H}_k} \quad \text{ for all } \varphi \in \mathcal{H}_k,
\ee
leading to the equation \((\Delta_k u)(X) = f(X)\), \(\Delta_k\) being the approximation of the Laplace-Beltrami operator of Section \eqref{divergence-operator}. A solution to this equation is computed as \(u = (\Delta_k)^{-1} f\), defined in \eqref{Deltainvk}.

Figure~\ref{fig:PDE1} displays a regular mesh for the domain \(\Omega = [0,1]^2\), where \(f\) is plotted in the left-hand side, and the solution \(u\) in the right-hand side. Figure~\ref{fig:PDE2} computes a Poisson equation on an unstructured mesh generated by a bimodal Gaussian random variable, with \(f\) plotted on the left and the solution \(u\) on the right.

\begin{figure}
\centering
\includegraphics[width=0.7\textwidth, keepaspectratio]{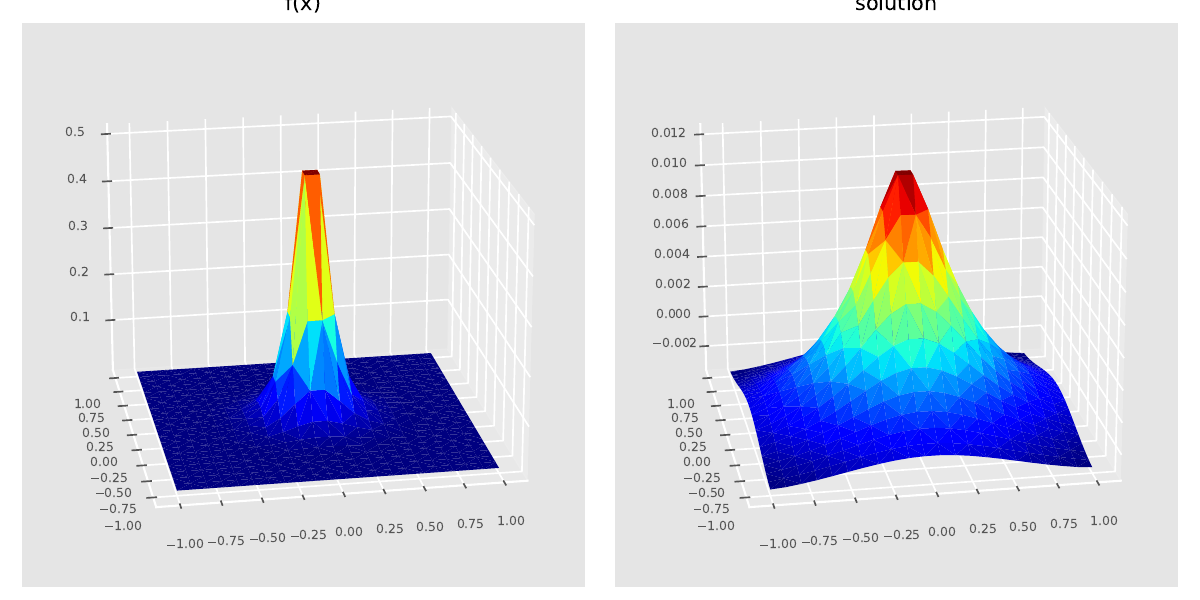}
\caption{\label{fig:PDE1}Computed inverse Laplace operator - regular mesh}
\end{figure}

Both examples (regular vs more complex geometries) illustrates how kernel methods facilitate the use of structured or unstructured meshes, enabling the description of more complex geometries in a unified framework. From a code perspective, there is no difference of treatment or interface to both problems. 

\begin{figure}
\centering
\includegraphics[width=0.7\textwidth, keepaspectratio]{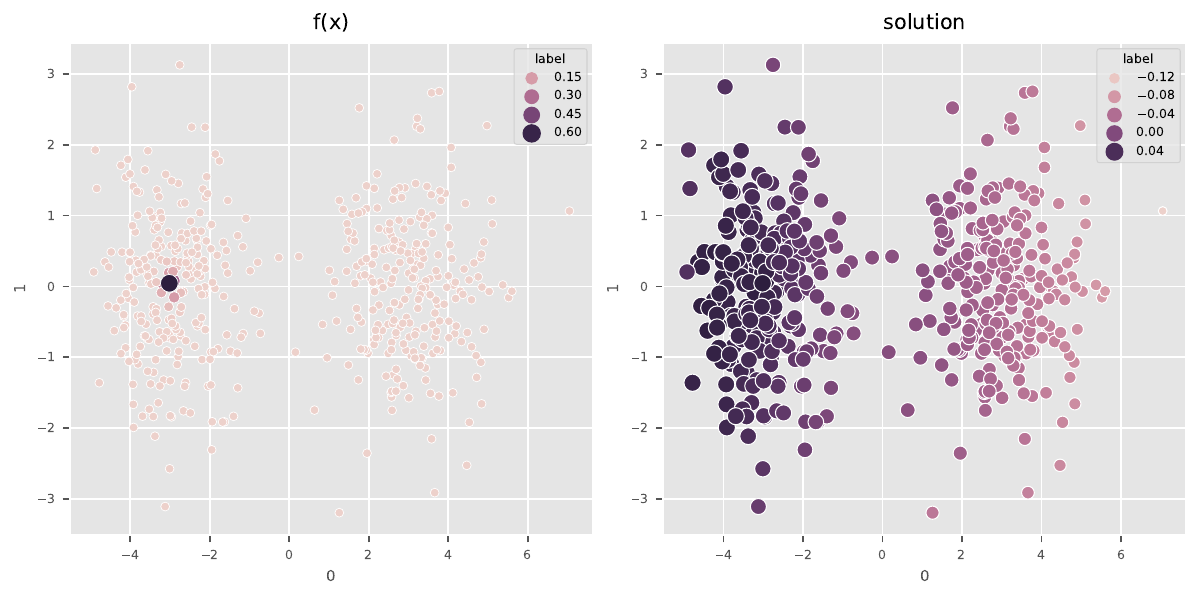}
\caption{\label{fig:PDE2}Computed inverse Laplace operator - irregular mesh}
\end{figure}

\subsection{A denoising problem}\label{a-denoising-problem}

A slight modification of the Poisson equation allows us to define alternative ways to regularize the (optional) \textit{ridge} regularization term in the projection operator \eqref{Pk}, introduced as an additional parameter in the pseudo-inverse formula \eqref{pseudo-inverse}, which we now discuss.
Suppose we want to solve a minimization problem in the form
\be
  \inf_{G \in {\mathcal{H}}_k} J(G), \quad J(G)= \|G-F\|_{\mathcal{H}_k}^2 + \epsilon \|L(G)\|_{L^2}^2
\ee
Here, \(L:\mathcal{H}_k(\Omega) \to L^2(\Omega)\) is a linear operator that serves as a penalty term. A formal solution is given by
\be
  G + \epsilon L^T L G= F
\ee

Numerically, consider \(X \in \mathbb{R}^{N_x,D}\), defining an unstructured mesh \(\mathcal{X}\), together with a kernel \(k\) for defining \(\mathcal{H}(\Omega)\). Denote \(L_k\) the discretized operator. This penalty problem defines a function \(G\)
as follows:
\be
  z \mapsto G(z) = K(X,z) \Big(K(X,X) + \epsilon \Big(L_k^T L_k\Big)(X,X) \Big)^{-1} F(X)
\ee
To compute this function, input \(R = \epsilon L_k^T L_k\) into the pseudo-inverse formula \eqref{pseudo-inverse}.

As an example, consider the denoising procedure, which aims to solve:
\be\label{DNabla}
  \inf_{G \in \mathcal{H}_k} \|G-F\|_{L^2}^2 + \epsilon \|\nabla G\|_{L^2}^2.
\ee
In this case, \(L_k = \nabla_k\), and \(L_k^T L^k\) corresponds to \(\Delta_k\). Figure~\ref{fig:denoising} demonstrates the results of this regularization procedure. The noisy signal (left image) is given by \(F_\eta(x) = F(x) + \eta\), where \(\eta\) is a white noise, and \(f\) is the cosine function \(f(x) = f(x_1, \ldots, x_D) = \prod_{d=1,\ldots,D} \cos (4\pi x_d) + \sum_{d=1,\ldots,D} x_d\). The regularized solution is plotted on the right.

In this case, \(L_k = \nabla_k\), the discrete gradient operator defined at \eqref{nablak}, and \(\nabla_k^T \nabla_k\) is an approximation of the Laplace-Beltrami operator \(\Delta_k\). Figure~\ref{fig:denoising} demonstrates the
results of this regularization procedure. The noisy signal (left image)
is given by \(F_\eta(x) = F(x) + \eta\), where \(\eta  =  \mathcal{N}(0,\epsilon)\) is a white Gaussian noise, \(\epsilon = 0.1\), and \(f(x) = f(x_1, \ldots, x_D) = \prod_{d=1,\ldots,D} \cos (4\pi x_d) + \sum_{d=1,\ldots,D} x_d\) is an example function. The regularized solution is plotted on the right.

\begin{figure}
\centering
\includegraphics[width=0.7\textwidth, keepaspectratio]{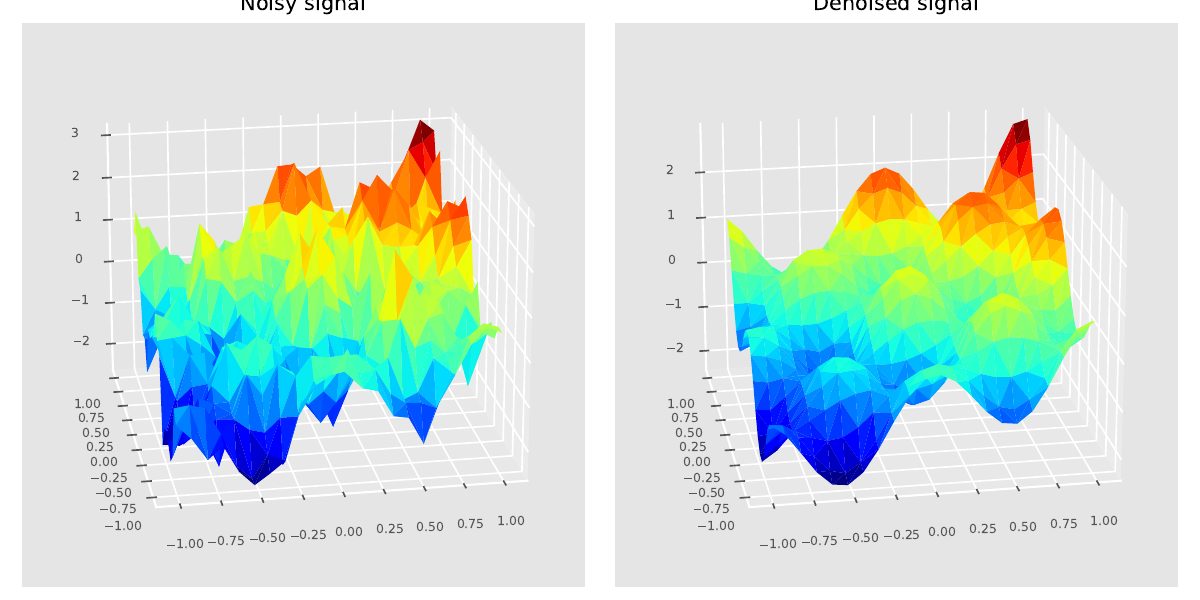}
\caption{\label{fig:unnamed-chunk-62}\label{fig:denoising}Example of denoising signals with a Laplace operator}
\end{figure}


\section{Time-evolution problems}

\subsection{Fokker--Plank and Kolmogorov equations}
\label{time-evolution-operators-based-on-theta-schemes}

\paragraph{Fokker--Plank equation}

Our physics-informed use cases are mainly concerned with Fokker--Planck and Hamilton--Jacobi--Bellman (HJB) equations, which are closely linked. We therefore focus on a dedicated framework for such equations.
Consider first the stochastic differential equation (SDE) describing the dynamics of a Markov-type stochastic process
\(t \mapsto X_t \in \RD\), i.e.~
\be\label{SED}
dX_t = g(X_t)dt+\sigma(X_t)dW_t. 
\ee
Here, \(W_t \in \mathbb{R}^D\) denotes a \(D\)-dimensional,
independent Brownian motion, while \(g \in \mathbb{R}^D\) is a prescribed
vector field and \(\sigma \in \mathbb{R}^{D\times D}\) is a prescribed
matrix-valued field.

Denote by \(d\mu = \mu(t,s,x,y)dx\) (defined for \(t\ge s\)) the
\textit{probability measure} associated with \(X_t\). where $x$ is the spatial integration variable, $y$ is the conditioning point \(X_s = y\) and $\mu$ is the transition density.

We recall that the density
\(\mu\) satisfies the \textit{Fokker-Planck equation}:
\be\label{FP}
    \del_t \mu - \Lcal \mu = 0, \quad \mu(s, \cdot) = \delta_{y},
\ee which is a linear convection-diffusion equation in the variable $x$. Moreover, the
initial data is the Dirac mass \(\delta_y\) at some point \(y\), and we use divergence-form forward operator: 
\be\label{Lcal}
\Lcal \mu =  \nabla \cdot (f (\mu)) + \nabla\cdot( A \nabla \mu), \quad A =  \frac{1}{2} \sigma \sigma^T, \ f(\mu) = \big(-g+(\nabla \cdot A)\big)\mu. 
\ee
Here, \(\nabla\) denotes the gradient operator,
\(\nabla \cdot\) the divergence operator, and
\(\nabla^2 =  (\del_i \del_j)_{1 \leq i,j \leq D}\) is the Hessian
operator. Since the initial data is a Dirac mass in $x$,  the equation~\eqref{Lcal} can be understood in the weak sense: meaning that we look for a measure $t \mapsto d\mu(t,\cdot)$ satisfying, for every smooth function $\varphi \in C_c^{\infty}$,
\be\label{WFP}
 \frac{d}{dt } \int \varphi d \mu = \int  (<f(\mu),\nabla \varphi> + <\nabla \cdot (A \mu),\nabla \varphi>) d\mu.
\ee

\paragraph{Backward Kolmogorov and Feynman-Kac}

The (vector-valued) dual of the Fokker-Planck equation is the
\textit{ Kolmogorov equation}, also known in mathematical finance as the
\textit{Black and Scholes equation}, with a drift given by risk-neutral measure. Given terminal data $P(t,\cdot) = P_t(\cdot)$,  the solution is denoted by $P(s,y) = \mathbb E[\;P_t(X_t)\mid X_s=y\;]$. This equation determines the unknown vector-valued function \(\overline{P} = \overline{P}(s,x)\) as a
solution, with \(s \leq t\), of the reverse-in-time equation: 
\be\label{KE}
\del_t \overline{P} - \Lcal^* \overline{P} =0, \quad \Lcal^* \overline{P} =  - f \cdot \nabla \overline{P} + \nabla \cdot (A \nabla \overline{P}).  
\ee
By the Feynmann--Kac formula, a solution to the Kolmogorov
equation \eqref{KE} can be interpreted as a time-average of an
expectation function, that is, for all $s \le t$ and $x$,
\be
    \overline{P}(s,y) = \mathbb{E}_{s, y}(P(t,\cdot)) = \int P(t,\cdot) d \mu(t,s,\cdot,y).
\ee

We now focus on the Fokker-Planck equation~\eqref{FP}. Reinforcement learning provides tools to treat the Kolmogorov equation in a more general forms, that are Hamilton-Jacobi-Bellman type equations; see~\eqref{HJB}. Let us split the Fokker-Planck equation into two terms. The first is $\partial_t\mu = \nabla \cdot f(\mu)$, called a \textit{scalar hyperbolic conservation law}, and $\partial_t \mu = \nabla \cdot ( A \nabla \mu)$, which is of \textit{diffusion} type. We treat these two components separately.

\paragraph{Lagrangian versus Eulerian semi-discrete schemes}

We provide two RKHS approaches to compute a solution to the Fokker--Planck equation~\eqref{FP}, both derived from the weak formulation~\eqref{WFP} and based on semi-discrete schemes. The first approach is Eulerian in nature and is based on fixed nodes and discrete differential operators, while the second approach is Lagrangian in nature and is based on moving particles and operator splitting. 

\begin{itemize}
\item \textit{Eulerian(fixed mesh)}: Use a fixed set of nodes $X=(x^1,\ldots,x^N)$ to represent the computational domain and a kernel $k$.  Then follow our mesh-free methodology and compute the probability solution $d\mu(t,\cdot) \sim \sum \mu(t,x^n)\delta_{x^n}$ as a solution to
\be \label{EUL}
    \frac{d}{dt} \mu(t,x^n) + \nabla_k \cdot f(\mu)(t,x^n)= \Big( \nabla_k \cdot (A \nabla_k ) \mu\Big)(t,x^n), \quad n=1,\ldots,N.
\ee
This is a semi-discrete linear system of ordinary differential equations, which diffusive part $\nabla_k \cdot (A \nabla_k \mu)$ is positive-definite.
\item \textit{Lagrangian(moving mesh)}: Use a particle representation $t \mapsto X_t=(x^1_t,\ldots,x^N_t)$. We thus define a time dependent kernel $t\mapsto k_t$ to represent the probability solution as $d\mu(t,\cdot) \approx \frac{1}{N}\sum \delta_{x^n_t}$. In this situation, we split the system solving at each time-step first the hyperbolic part with transport methods; see Section~\ref{Conservation-laws}. Then we solve the diffusive part, which leads to the following discrete nonlinear, semi-discrete, system of ordinary differential equations
\be \label{LAG}
    \frac{d}{dt} x^n_t = \Big( \nabla_{k_t} \cdot (A \nabla_{k_t} x^n_t)\Big), \quad n=1,\ldots,N.
\ee

Let us add an observation. Consider the trivial relation $\nabla x = I_D$. Residual kernels~\eqref{RKR} replicate this property as $\nabla_{k_t,X_t} X_t = I_D$, for which the discrete system~\eqref{LAG} can be written more simply as $\frac{d}{dt} x^n_t = \sum_m \nabla_{k_t} \cdot A(t,x^m_t)$. In particular, the heat equation $\partial_t \mu = \Delta \mu$ leads to $\frac{d}{dt} x^n_t = \sum_m (\nabla_{k_t} \cdot I_D)(x^m_t)$. This expression is somehow perturbing, but the divergence, defined in Section~\ref{laplace-beltrami-operator}, is the one of the Laplace-Beltrami operator and $\nabla_k \cdot I_D$ is not trivial.
\end{itemize}

Lagrangian methods are usually more accurate than Eulerian ones, but are more computationally involved. Better accuracy comes from one hand because Lagrangian methods can handle unbounded domains, and from another hand since the system~\eqref{LAG} approximates sharp-discrepancy sequences; see Section~\ref{sharp-discrepancy sequences} below. 

Both approaches, Eulerian versus Lagrangian, end up considering time-dependent systems having form $\frac{d}{dt} u = A(t,\ldots) u$, which are called semi-discrete schemes. Thus methods to go from semi-discrete to fully discrete schemes are required. The next section presents $\theta$-schemes, which is a simple, yet efficient method to fully discrete schemes. The Eulerian approach might require some extra care, so a more elaborate construction, called entropy-satisfying schemes, is also provided in~Section~\ref{entropy-dissipative-schemes}.

\paragraph{Time-dependent generators based on $\theta$-schemes}


When it comes to discretization of time-dependent PDEs, most examples usually resumes to consider the following class of dynamical system with Cauchy initial conditions
\be\label{TES}
\frac{d}{dt} u(t)  = A u(t), \quad u(0) \in \mathbb{R}^{N_x,D}, \quad A \in \mathbb{R}^{N_x,N_x}, 
\ee
where \(A \equiv A(t,x,u,\nabla u)\) can be any matrix valued operator, assumed to be \textit{bounded}, i.e.~satisfying
\be
\langle A u, u\rangle_{\ell^2} \le C \qquad \text{ for all } u \in \mathbb{R}^{N_x,D_x}.
\ee
Thus we describe a classical way to deal with such systems.

Let $\{t^n \}_{n\geq 0}$ be a time grid with steps \(\tau^n=t^{n+1}-t^{n}\). For \(0 \le \theta \le 1\), the \(\theta\)-scheme reads
\be
 \delta_t u(t^n) = \frac{u(t^{n+1})- u(t^{n})}{\tau^{n}} = A \Big(\theta u(t^{n+1}) + (1-\theta) u(t^{n}) \Big) = A u^\theta(t^n).
\ee
A formal solution of this scheme is given by \(u(t^{n+1}) = B(A,\theta,dt) u(t^{n})\), where \(B\) is the \emph{generator} of the equation, defined as
\be\label{generator}
 B(A,\theta,\tau^n) = \Big( I - \tau^n \theta A  \Big)^{-1}\Big( I + \tau^n (1-\theta)A \Big).
\ee
\begin{itemize}
\tightlist
\item
  The value \(\theta = 1\) corresponds to the \emph{implicit} Euler approximation.
\item
  The value \(\theta = 0\) corresponds to the \emph{explicit} Euler approximation.
\item
  The value \(\theta = 0.5\) corresponds to the \emph{Crank--Nicolson}.
\end{itemize}

The Crank--Nicolson scheme is motivated by the following energy estimate, taking the scalar product with \(u^\theta(t^n)\) in the discrete equation, where \(\ell^2\) denoting the standard discrete quadratic norm
\be
 <A u^\theta(t^n), u^\theta(t^n)>_{\ell^2} = \frac{\theta \|u(t^{n+1})\|_{\ell^2}^2 - (1-\theta)\|u(t^{n})\|_{\ell^2}^2 + (1-2\theta)<u(t^{n+1}), u^\theta(t^n)>_{\ell^2}}{\tau^n}.
\ee
For \(\theta \ge 0.5\), an \emph{energy dissipation} \(\|u(t^{n+1})\|_{\ell^2}^2 \le \|u(t^{n})\|_{\ell^2}^2\) is achieved, provided \(A\) is a negative defined operator. Choosing \(\theta \ge 0.5\) leads to \emph{unconditionally} stable numerical schemes. The Crank--Nicolson scheme \(\theta = 0.5\) is a highly versatile choice, which is adapted to \emph{energy conservation}, that is, to operators \(A\) satisfying \(<Au,u>_{\ell^2} = 0\).

\subsection{Hyperbolic conservation laws}
\label{Conservation-laws}
\label{characteristic}
\label{WCL}
\label{entropy-dissipative-schemes}

\paragraph{Purpose.}

We consider conservation laws as measures $t\mapsto d\mu(t,\cdot) = \mu(t,\cdot) dx$, which densities $\mu$ are solutions to the following equation having initial conditions at time $t=0$:
\be\label{CL}
    \partial_t \mu + \nabla \cdot f(\mu) = 0, \qquad \mu(0,\cdot) = \mu_0(\cdot),
\ee
where \(f = (f_d(u))_{1 \leq d \leq D}: \mathbb{R}  \to \mathbb{R}^D\) is a given flux
and \(\nabla \cdot f(u) = \sum_{1 \leq d \leq D} \partial_{x_d} f_d(u)\) denotes its divergence, with \(x=(x_d)_{1 \leq d \leq D}\). Here, a solution to~\eqref{CL} has to be understood in a weak sense, that is, for every smooth function $\varphi$,
\be\label{WFCL}
\frac{d}{dt} \int \varphi d\mu(t,\cdot) = \int <f(\mu), \nabla \varphi>(t,\cdot)dx 
\ee
The conservation law attached to the Fokker-Planck equation~\eqref{FP} is a probability measure, which is our main focus. However,~\eqref{WFCL} holds also for general signed, vector-valued measures, which are Hamilton-Jacobi equations.

\paragraph{Lagrangian approach and the characteristic method}

We now seek solutions to the conservation law~\eqref{CL} determined by the \textit{characteristic} method, which corresponds to weak solutions to~\eqref{WFCL} having form, where $y^{-1}(t,x)$ is the inverse map of $y(t,x)$
\be\label{CM}
    \mu(t,x) = \mu_0(y^{-1}(t,x)), \quad y(t,x) = x + t f'(\mu_0(x)).
\ee
For its derivation, we establish that $\mu = \mu_0\circ y^{-1}$ is a strong solution to the conservation law~\eqref{CL}, which supposes that all quantities below are smooth enough to be derivable. In particular, we observe that the map $t\mapsto y(t,\cdot)$ is invertible for small time $t < T = \inf \{t: \text{det}|\nabla y(t,x)|=0\} $.  We compute pointwise
\be
 \partial_t \mu = <(\nabla \mu_0 )\circ y^{-1},\partial_t y^{-1}>. 
\ee
From~\eqref{CM} and the relation $y^{-1}(t,y(t,x)) = x$, we deduce $\partial_t y^{-1} + \nabla y^{-1}f'(\mu_0)\circ y^{-1}=0$, to get
\be
 \partial_t \mu = -<\nabla y^{-1}(\nabla \mu_0 )\circ y^{-1},f'(\mu_0)\circ y^{-1}> = -\nabla\cdot f(\mu_0 \circ y^{-1})= -\nabla\cdot f(\mu).
\ee

\paragraph{Push-forward interpretation.}

By the definition of the push-forward \eqref{transport}, the characteristic map $y(t,\cdot)$ transports the initial density to time $t$:
\be
  \mu(t,\cdot)\,dx  =  \bigl(y(t,\cdot)\bigr)_{\#}\bigl(\mu_0(\cdot)\,dx\bigr),
  \quad\text{i.e.}\quad
  \int \varphi(z)\,\mu(t,z)\,dz  =  \int \varphi\bigl(y(t,x)\bigr)\,\mu_0(x)\,dx
\ee
for every test function $\varphi$. Equivalently,
\be
  \mu_0(\cdot)\,dx  =  \bigl(y^{-1}(t,\cdot)\bigr)_{\#}\bigl(\mu(t,\cdot)\,dx\bigr).
\ee
This representation holds for $t\in[0,T)$ as long as $y(t,\cdot)$ is a diffeomorphism (i.e., $\det\nabla y(t,\cdot)\neq 0$). Hyperbolic conservation laws are often efficiently modeled in Lagrangian form: a mesh moves along characteristics,
\be
  X_t = X_0 + t\,f'\!\bigl(\mu_0(X_0)\bigr),
\ee
yielding the exact classical solution up to the first shock time (the loss of invertibility of $y$).

In view of the push-forward definition~\eqref{transport}, one can define equivalently $\mu$ as a solution to $y_\#(t,\cdot)d\mu(t,\cdot) = d\mu_0(\cdot)$, defining formally a solution to the conservation law~\eqref{CL}, in the weak sense~\eqref{WFCL}, for any time.

Observe that hyperbolic conservation laws are usually better modeled by Lagrangian methods than Eulerian ones, as they define a mesh, moving accordingly to $X_t = X_0+tf'(\mu_0(X_0))$, which provides an exact solution, as long as this map is invertible of course.

\paragraph{Entropy solutions.}

An entropy function is any convex, scalar-valued, function $U = U (\mu)$, and we denote the entropy variable $v(\mu) = U' (\mu)$. Let us denote the entropy flux $G(u)$ satisfying $G'(u)=v(u)f'(u)$. We verify that any smooth solution to the conservation law~\eqref{CL} satisfies the following entropy relation
\be\label{ECL}
    \partial_t U(\mu) = -U'(\mu)\nabla \cdot f(\mu) = -U'(\mu)f'(\mu)\nabla \mu = -\nabla \cdot G(\mu). 
\ee
This equation is also to be understood in the weak sense, that is, for every smooth $\varphi(\cdot)$,
\be\label{WFECL}
\frac{d}{dt} \int \varphi U(\mu)(t,\cdot) dx + \int <G(\mu), \nabla \varphi>(t,\cdot)dx =0.
\ee
In particular, we deduce that any characteristic solution~\eqref{CM} satisfies the \textit{entropy conservation} $\frac{d}{dt}\int U(\mu)(t,\cdot)dx = 0$. However, there exists others, more \textit{physical} solutions, called \textit{entropy solutions} to~\eqref{CL}, which are solutions satisfying the entropy condition 
\be\label{EC}
   \partial_t \mu + \nabla \cdot f(\mu) = 0, \quad \partial_t U(\mu) + \nabla \cdot G(\mu) \le 0,
\ee
in the weak sense, meaning that the quantity~\eqref{WFECL} is negative, for any positive function $\varphi \ge 0$. 

There exists two families of methods to compute these solutions.
\begin{itemize}
\item \textit{Eulerian:} this method consists in solving in the limiting case \(\epsilon \to 0\) the following viscosity equation version of \eqref{CM}
\be
\label{CMD}
  \partial_t \mu_\epsilon + \nabla \cdot f(\mu_\epsilon) = \epsilon \Delta \mu_\epsilon,
\ee
on a prescribed, fixed mesh $X=(x^1,\ldots)$, and we compute the density $t \mapsto (\mu(t,x^1), \ldots)$ as a function solution to~\eqref{CMD}.
For any \(\epsilon > 0\), the solution \(\mu_\epsilon\) satisfies in a strong sense the \textit{entropy dissipation} property \(\partial_t U(\mu_\epsilon) + \nabla \cdot G(\mu_\epsilon) \le 0\), for any convex entropy - entropy fluxes \(U,G\). In the limiting case \(\epsilon \to 0\), this entropy dissipation holds in the weak sense. 

\item  \textit{Lagrangian:} this method involves direct computations explicitly characterizing the entropy solution, as the Hopf-Lax formula, or the convex hull algorithm\footnote{introduced in \cite{LeFlochMercier-2015}} This latter is computed as
\be
  \mu(t,\cdot) = y^+(t,\cdot)_\# \mu_0(\cdot), \quad y(t,x) = x + t f'(\mu_0(x)),
\ee
where \(y^+(t,\cdot)\) is computed as
\be
  y^+(t,\cdot) = \nabla h^+(t,\cdot), \quad \nabla h(t,\cdot) = y(t,\cdot),
\ee
and \(h^+(t,\cdot)\) is the \textit{convex hull} of \(h\). This approach is usually thought as a  Lagrangian approach, as the mesh $t\mapsto X_t$, representing $t\mapsto d\mu(t,\cdot)$, is an equiweighted representation as $d\mu(t,\cdot) \sim \frac{1}{N} \sum_n \delta_{x^n_t}$, is moving accordingly.

\end{itemize}


\paragraph{Entropy conservative versus entropy dissipative schemes}

We now consider solutions to conservation laws with the Eulerian method (see~\eqref{CMD}). In order to properly compute entropy solutions, it is crucial to ensure the entropy relation $\partial_t U + \nabla \cdot G \le 0$ at a discrete level, for any $\epsilon >0$, to retrieve numerical stability. 

However, the \(\theta\)-scheme framework of Section~\ref{time-evolution-operators-based-on-theta-schemes} cannot guarantee such relations. So we discuss here a family of \textit{entropy-satisfying numerical schemes}\footnote{based on \cite{LeFloch-Mercier:2002, LeFlochMercier-R-2002} and references therein} in the context of finite-difference schemes. These techniques can be extended to the RKHS framework proposed in the present monograph, for Hamilton-Jacobi-type equations~\eqref{CL}, which can be rephrased, using the entropy variable $v(u) = U'(u)$, and assuming $f(u) = g(v)$, $G'(u)=v(u)f'(u)$,
\be\label{ED-2}
\partial_t u  + \nabla \cdot g(v(u)) =0, \quad \partial_t U(u)  + \nabla \cdot G(v(u)) \le 0.
\ee
Consider a time grid $t^n< t^{n+1}< \ldots$, a mesh $X=(x^1,\ldots,x^N)$, $k$ a kernel. This system can be approximated using the viscosity approach~\eqref{CMD}, written as the following semi-discrete scheme:
\be\label{decadix}
  \frac{d}{dt} u_\epsilon(t,x^n) + \nabla_{k} \cdot g(v(u_\epsilon))(t,x^n) = \epsilon \Delta_k v(u_\epsilon(t,x^n)), \quad n = 1,\ldots,N
\ee
In what follows, we denote $u \equiv u_\epsilon$ for concision. Let us now fully discretize this scheme and  denote by \(\tau^n=t^{n+1}-t^{n}\) and \(u_i^n \sim u(t^n,x^i)\) the discrete solution. In the same way, define $U_i^n,v_i^n$ and \(\delta_t f^n = \frac{f^{n+1}-f^{n}}{\tau^n}\) the discrete forward time derivative operator. To approximate such a system numerically, a strategy for building entropy dissipative schemes involves first the choice of a \((q+1)\)-time level interpolation \(u^*(u^q, ..,u^0)\) which satisfy the following conditions. 
\begin{itemize}
\item
  Consistency with the identity \((u^*(u, ..,u) = u)\).
\item
  Invertibility and regularity of the map \(u^q \rightarrow u^*(u^q, ..,u^0)\).
\end{itemize}
This construction already prevailed for the $\theta$-schemes of Section~\ref{time-evolution-operators-based-on-theta-schemes}, where $u^*(u^1,u^0) = \theta u^1+(1-\theta)u^0$. However this settings allows us to consider more sophisticated time integrator, providing higher-order interpolation schemes (see~\eqref{VdMalpha}).

Let \(u^{*,n} = u^*(u^n, ..,u^{n-q})\), and let us choose the entropy variable \(U^*(u^q, ..,u^0)\), with \(U(u^*)\) as a possible choice. We set \(U^{*,n} = U^*(u^n, ..,u^{n-q})\). This variable must enjoy the following properties. 
\begin{itemize}
\tightlist
\item
  Be consistent with the original entropy \(U(u)\) (i.e.~\(U^*(u, ..,u) = U(u)\)).
\item
  Define the \((q+2)\)-time entropy variable \(v^{*,n+1/2}(u^{q+1}, ..,u^0)\), which satisfies
  \be
    \delta_t U^{*,n} = \frac{U^{*,n+1}-U^{*,n}}{\tau^n} = v^{*,n+1/2} \cdot \delta_t u^{*,n}
  \ee\\
  and is consistent with the entropy variable \(v^{*,n+1/2}(u, ..,u) = v(u)\).
\end{itemize}
The semi-discrete system~\eqref{decadix} is then approximated by the \textit{fully discrete} numerical scheme displayed now, where \(u^{n+1}\) is the unknown:
\be\label{e-scheme}
    \delta_t u^{*,n} = \frac{u^{*,n+1}-u^{*,n}}{\tau ^n} = - \nabla_k \cdot g(v^{*,n+1/2})+\epsilon \Delta_k v^{*,n+1/2}.
\ee
These schemes can be fully implicit or explicit with respect to the unknown \(u^{n+1}\), based on the entropy variable choice. They are entropy stable as follows: set \(E^{*,n} = \sum_{n} U(u_i^n)\) and compute
\be
    \delta_t E^{*,n} = \sum_i <\nabla_k\cdot g(v_i^{*,n+1/2}) + \epsilon \Delta_kv_i^{*,n+1/2},v_i^{*,n+1/2}>. 
\ee
The definition of the Laplace-Beltrami operator implies $\sum_i <\Delta_k v_i,v_i>=-\sum_i|\nabla_kv_i|^2 \le 0$, and by definition $(\nabla_k v_i)g(v_i)=-\nabla_k G(v_i)$. Hence we get, for any $\epsilon >0$,
\be
    \delta_t E^{*,n} \le \sum_i \nabla_k\cdot G(v_i^{*,n+1/2}).
\ee
The right-hand side is, by definition, $< G(v^{*,n+1/2}), \nabla_k 1>$, which is zero, provided we consider a kernel satisfying \(\nabla_k 1 \equiv 0\), as the residual kernel~\eqref{RKR}. Hence solutions of these fully discrete schemes enjoys the property \(E^{*,n+1} \le E^{*,n}\), which in turn implies the numerical stability of entropy schemes \eqref{e-scheme}.

\paragraph{Example: Crank--Nicolson for a skew-symmetric operator}

We provide a basic example considering the linear equation 
\be
\partial_t u + Au=0,
\ee
where $A$ is a linear operator (matrix) satisfying $\langle Au,u\rangle = 0$. This equation satisfies the energy conservation $\frac{d}{dt}\int |u(t)|^2dx = 0$. Picking up and the entropy function \(U(u)=u^2\), and with
\be
u^{*,n+1/2} = \frac{u^{n+1}+u^n}{2},
\ee
the update 
\be
\delta_t u^n = A u^{*,n+1/2}
\ee
is the Crank--Nicolson (fully discrete) scheme \(\theta = 1/2\), a second-order accurate scheme in time and preserves the discrete energy exactly.

\paragraph{Numerical illustration with the inviscid Burger equation}

Figure~\ref{fig:CHA1D} illustrates the differences between entropy-conservative and entropy-dissipative solutions for the one-dimensional inviscid Burgers equation
\be
\partial_t \mu + \frac{1}{2}\partial_x \mu^2 =0,
\ee
since Figure~\ref{fig:CHA2D} illustrates the two-dimensional case \(\partial_t \mu + \frac{1}{2}\nabla \cdot(\mu^2,\mu^2)=0\). The left-hand figure is the initial condition at time zero, since the solution at middle represents the conservative solution at time 1, and the entropy solution is plot at right.

\begin{figure}
\centering
\includegraphics[width=0.7\textwidth, keepaspectratio]{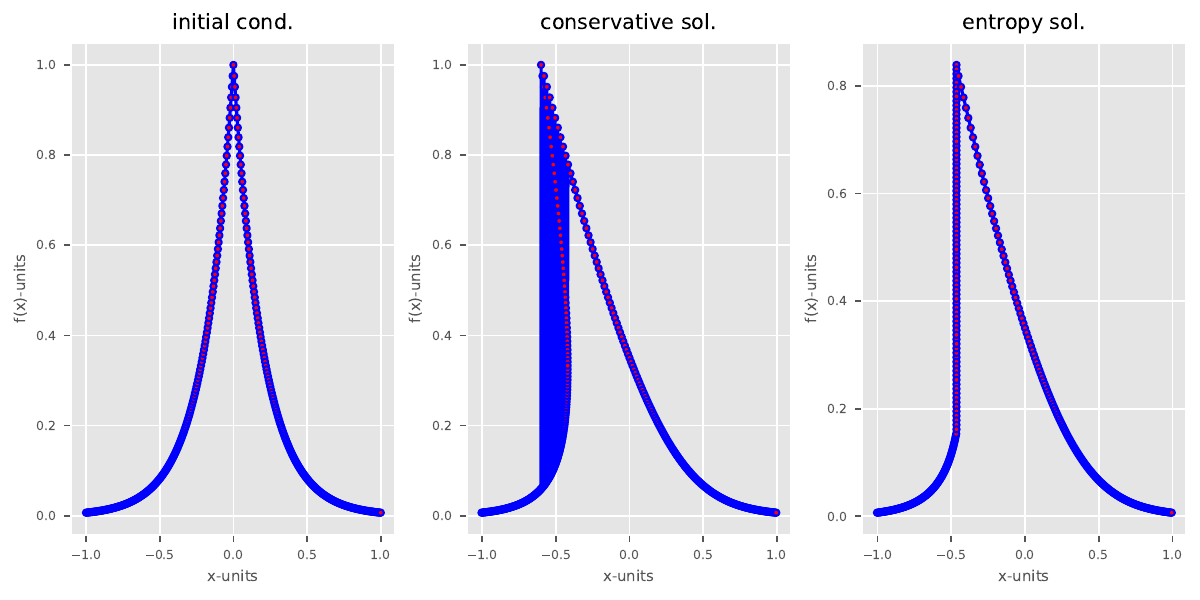}
\caption{\label{fig:unnamed-chunk-63}\label{fig:CHA1D}Convex Hull algorithm: left initial condition $\mu_0$. Middle: conservative solution. Right: entropy solution. Red dots: $y(t,x) = x + t f'(\mu_0)$}
\end{figure}

\begin{figure}
\centering
\includegraphics[width=0.7\textwidth, keepaspectratio]{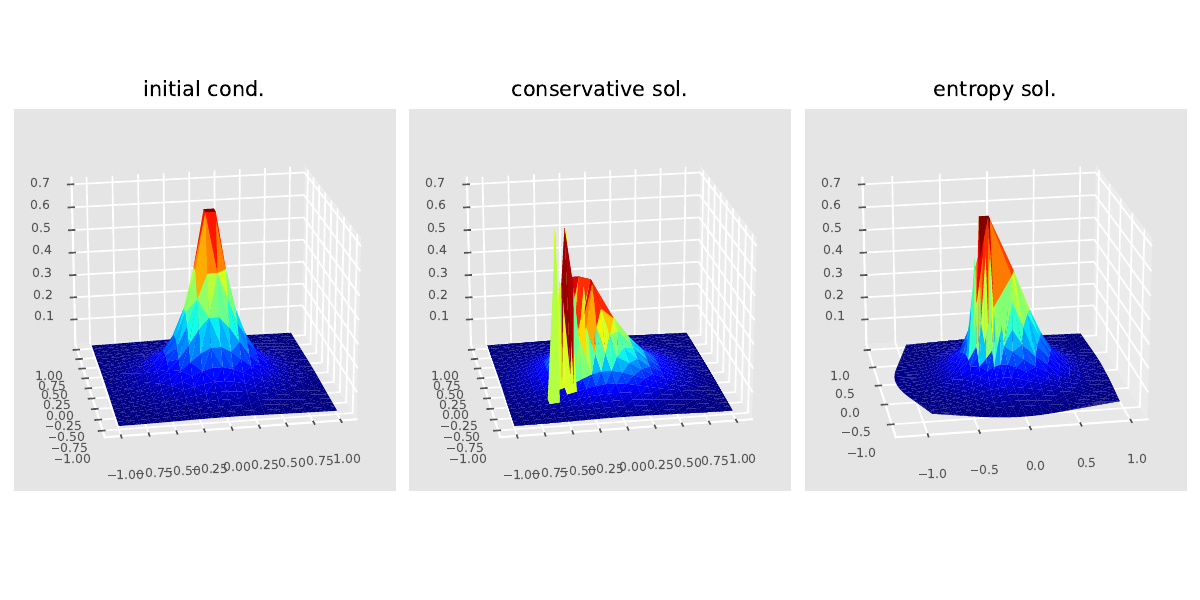}
\caption{\label{fig:unnamed-chunk-64}\label{fig:CHA2D}Convex Hull algorithm}
\end{figure}


\subsection{Diffusion equation}

\paragraph{A mesh-free Eulerian example in a fixed domain}\label{a-mesh-free-eulerian-method-for-a-fixed-domain}

We now illustrate the numerical study of time-dependent PDEs in the context of kernel methods, considering the heat equation equation in a fixed geometry \(\Omega\) with null Dirichlet conditions:
\be
  \partial_t u(t,x) = \Delta u(t,x), \quad u(0,x)=u_0(x),\quad x \in \Omega, \quad u_{\partial \Omega}=0
\ee
To approximate this equation, we follow the following steps. 
\begin{itemize}
\tightlist
\item
  Select a mesh \(X \in \mathbb{R}^{N_x,D}\) for the domain \(\Omega\).
\item
  Pick up a kernel \(k\) generating a space of vanishing trace functions.
\end{itemize}
From~\eqref{EUL}, this equation is discretized as \(\frac{d}{dt} u(t) = \Delta_k u(t)\), and integrated using the evolution operator \(u^{n+1} = B(\Delta_k,u^{n},dt,\theta)\) with \(\theta = 1\), which is the fully implicit case in \eqref{generator}. Figure~\ref{fig:heat2DR} provides a 3-D representation of the initial condition and time evolution of the heat equation on a fixed square.

This approach can be easily adapted to more complex geometries, as demonstrated by Figure~\ref{fig:heat2DIRR}, which shows the heat equation on an irregular mesh generated by a bimodal Gaussian process, as we consider the same setting to the Poisson equation; see Section~\ref{MM}.

\begin{figure}
\centering
\includegraphics[width=0.7\textwidth, keepaspectratio]{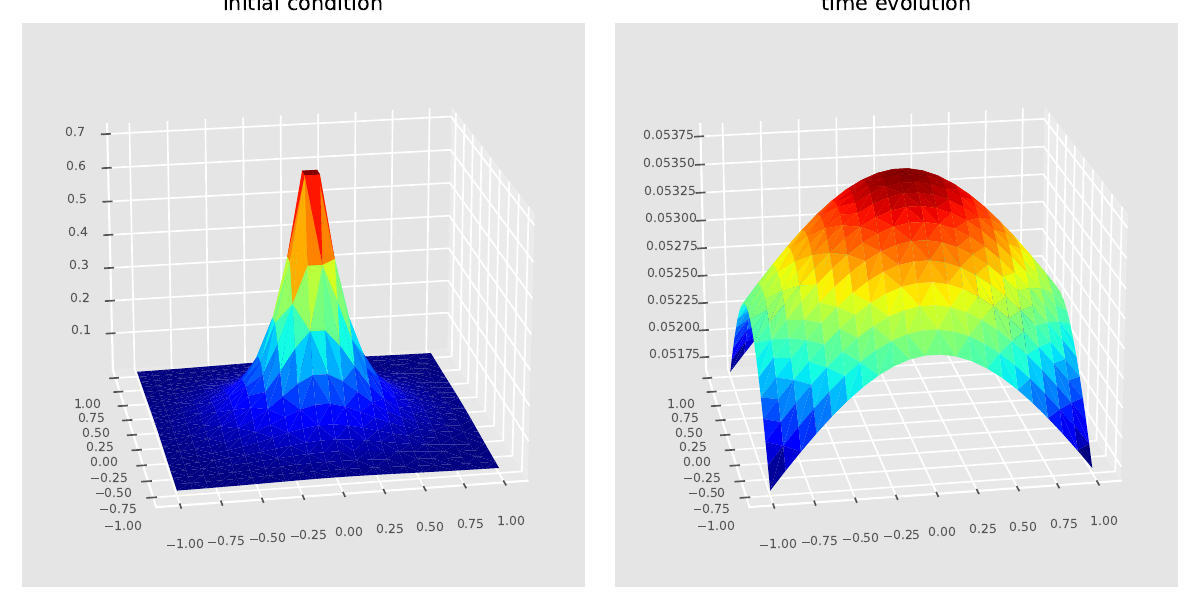}
\caption{\label{fig:heat2DR}A heat equation on a fixed regular mesh}
\end{figure}

\begin{figure}
\centering
\includegraphics[width=0.7\textwidth, keepaspectratio]{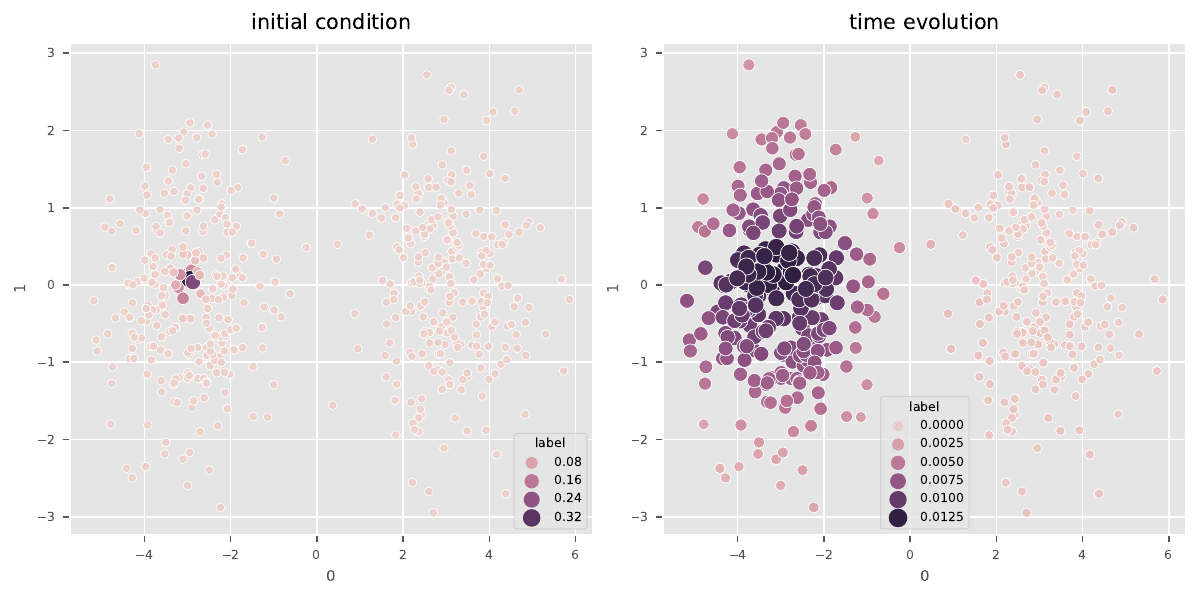}
\caption{\label{fig:heat2DIRR}A heat equation on a irregular mesh}
\end{figure}


\paragraph{A Lagrangian example on unbounded domain}

Next, we consider the heat equation on an unbounded domain, with measure-valued Cauchy initial data, that is:
\be\label{HEU}
  \partial_t \mu = \Delta \mu, \quad \mu(0,x)=\mu_0(x),\quad x \in \mathbb{R}^D.
\ee
The exact Lagrangian (probabilistic) representation is given by the Brownian motion $X_t$ with
\be
dX_t = \sqrt{2}dW_t, \quad X_0 \sim \mu_0,
\ee
so that $\mu(t,\cdot) = \text{Law}(X_t)$.

Following the guidelines for Lagrangian methods~\eqref{LAG}, we discretize this equation as
\be \label{HLE}
    \frac{d}{dt} x^n_t = -(\Delta_{k(t) x}) (x^n_t), \quad n = 1, \dots, N.
\ee
Observe that this equation looks like a heat equation, but with a negative sign. However, a residual kernel~\eqref{RKR} $t\mapsto k_t$  satisfies $<\Delta_{k_t} X_t,X_t>= <\nabla_{k_t} X_t,\nabla_{k_t} X_t> = D$, hence is bounded, and we can integrate using a $\theta$-scheme, which is given is the explicit case $\theta=0$ in the following expression
\be 
    x^n_{t^{k+1}} = x^n_{t^{k}} - \tau^k(\Delta_{k_{t^k}} x)^n_{t^{k}} = x^n_{t^{k}}-\tau^k \Big(\nabla_{k_{t^k}}\cdot I_D\Big)(x^n_{t^{k}}),
\ee
the last coming from the expression $\Delta_k X = \nabla_k^T\nabla_k X = \nabla_k^T I_D$ for residual kernels.
\paragraph{Sharp discrepancy sequences}\label{sharp-discrepancy sequences}

Let $d\mu(\cdot)$ be a probability measure in $\mathbb{R}^D$. Let $X=(x^1,\ldots,x^n)$, consider a kernel $k$, and minimizing the  following discrepancy functional
\be
 \overline{X} = \arg \inf_{X \in \mathbb{R}^{N,D}} d_k(d\mu,\delta_X).
\ee
where $d_k(d\mu,\delta_X)$ is the discrepancy functional expression, which is hybrid between the continuous expression~\eqref{eq:TSP} and the discrete one in~\eqref{dk}: 
\be
 d_k(d\mu,\delta_X) = \frac{1}{N^2}\sum_{n,m} k(x^n,x^m) +\int \int k(x,y) d\mu(x)d\mu(y)-2 \sum_n \int k(x^n,y) d\mu(y). 
\ee
We refer to solutions to the above problems as \textit{sharp-discrepancy} sequences for the measure $d\mu$. Such solutions are quite interesting: as they minimize the error estimate of kernel regression~\eqref{estim}, they compute somehow the best representation of a given probability measure $d\mu$ with equi-weighted Dirac masses.

Following Section~\ref{linear-sum-assignment-problems} and using the discrepancy, a descent algorithm approximates the minimum as, for any $n=1,\ldots,N$ 
\be \label{SDS}
 \frac{d}{dt} x^n_t = \frac{2}{N^2}\sum_m \nabla k_t(x^n_t,x^m_t) -\frac{2}{N} \int \nabla k_t(x^n_t,y) d\mu(y). 
\ee 
This system should converge in infinite time to a solution to the following equation 
\be
\frac{1}{N}\sum_m \nabla k(x^n,x^m) =\int \nabla k(x^n,y) d\mu(y),
\ee 
in the spirit of self-organizing maps\footnote{see wikipedia \url{https://en.wikipedia.org/wiki/Self-organizing_map}}.
As the discrepancy is symmetrical, we can integrate by parts, to get the following equation where $x^1,\ldots,x^N$ are unknowns
\be
 \frac{1}{N}\sum_m \nabla k(x^n,x^m) =-\int  k(x^n,y) (\nabla \mu)(y) dy, \quad n=1,\ldots,N,
\ee
provided $\mu$ is derivable of course. Observe that there are similarities between the two approaches~\eqref{SDS} and~\eqref{HLE}. Indeed, considering the diffusion equation~\eqref{LAG}, we could apply the Lagrangian approach as follows: 
\be
 \overline{X}_t = \arg \inf_{X_t \in \mathbb{R}^{N,D}} d_k(d\mu_t,\delta_{X_t}),
\ee
which leads to a similar system of equations than the semi-discrete scheme~\eqref{LAG}, showing that this scheme allows us to compute sharp-discrepancy sequences.
\begin{figure}
\centering
\includegraphics[width=0.7\textwidth, keepaspectratio]{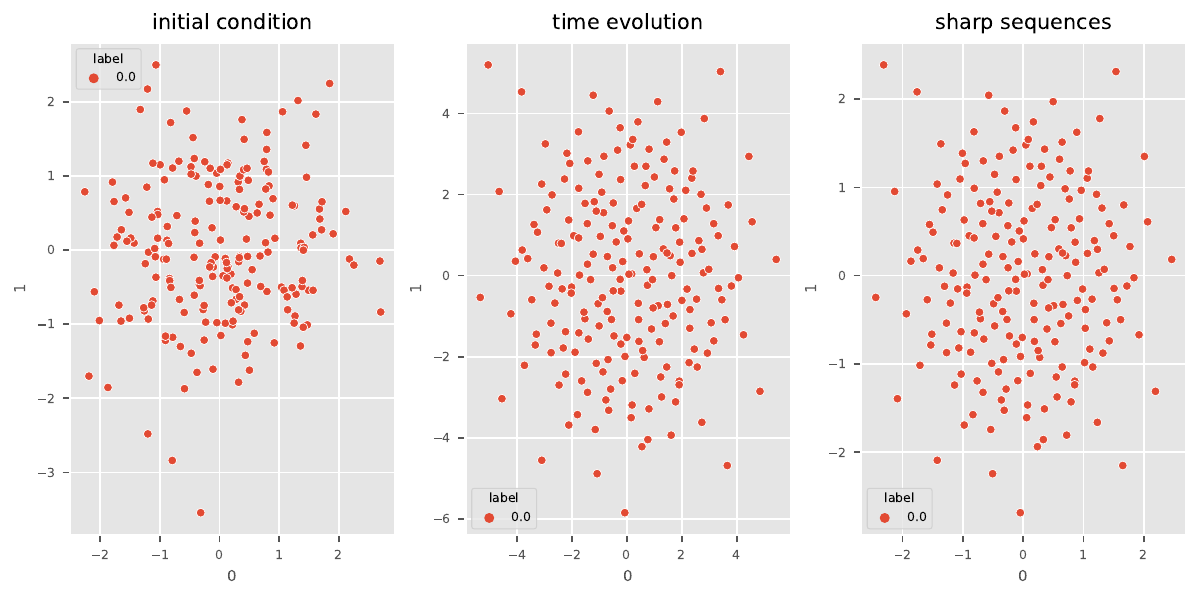}
\caption{\label{fig:heatLAG}A heat equation solved with a Lagrangian method. Left: initial distribution $X_0$. Middle long-time evolution of the system~\eqref{HLE}. Right rescaled to unit-variance (sharp discrepancy sequences for the normal multi-dimensional normal law.}
\end{figure}

We conclude this part by providing a simple example of Lagrangian formulation of the heat equation at Figure~\ref{fig:heatLAG}, corresponding to the computation for a Brownian motion.

\paragraph{Accuracy issue}. 

Observe that such a Lagrangian-based computation hold for any solutions of the Fokker Plank equations~\eqref{FP}, that is, for any stochastic processes having form~\eqref{SED}, and we can check their strong convergence properties. For instance, one\footnote{see \cite{LeFloch-Mercier:2017}} can compute such sequences for the Heston process\footnote{see \url{https://en.wikipedia.org/wiki/Heston_model}}, showing that the convergence rate of such variate is of order
\be
  \Big| \int_{\mathbb{R}^D} \varphi d \mu - \frac{1}{N} \sum_i \varphi (x^i)\Big| \le \frac{\mathcal{O}(1)}{N^2}
\ee
for any sufficiently regular function \(\varphi\). This should be compared to a naive Monte Carlo variate, converging at the statistical rate \(\frac{\mathcal{O}(1)}{\sqrt{N}}\).


\section{Techniques for PDEs}

\subsection{Automatic differentiation}
\label{automatic-differentiation}

\paragraph{Purpose.}

Automatic differentiation (AD) refers to a family of techniques for computing derivatives of functions implemented as programs composed of elementary differentiable operations and control flow. Unlike finite differences, AD yields \emph{exact} derivatives up to floating-point roundoff (i.e., no truncation error), while unlike symbolic differentiation, it operates on the executed program without expression swell\footnote{see, for instance, \cite{GW:2008}}.  

Two computational modes are standard. In \emph{forward-mode} AD, one propagates directional derivatives (Jacobian--vector products, JVPs) through the computational graph, effectively differentiating intermediate variables with respect to (chosen) input directions. In \emph{reverse-mode} AD, one propagates adjoints (vector--Jacobian products, VJPs) backward, accumulating the sensitivity of a scalar output with respect to intermediate variables. For a function $f:\mathbb{R}^D\to\mathbb{R}^m$, forward-mode evaluates columns of the Jacobian at a cost comparable to one function run per direction (favorable when $D$ is small), whereas reverse-mode evaluates gradients of scalar outputs at a cost within a small constant of a single function run (favorable when $m$ is small), with a memory--time tradeoff handled by checkpointing\footnote{see \cite{GW:2008}}. 

\paragraph{Existing libraries.}

Modern libraries provide high-quality implementations of AD, including TensorFlow, PyTorch, Autograd, Zygote (Julia), and JAX.\footnote{Software: TensorFlow \url{https://www.tensorflow.org/}; PyTorch \url{https://pytorch.org/}; Autograd \url{https://github.com/HIPS/autograd}; Zygote \url{https://fluxml.ai/Zygote.jl/latest/}; JAX \url{https://github.com/google/jax}}
Most expose reverse-mode; several additionally expose forward-mode\footnote{E.g., JAX exposes both JVP and VJP; recent PyTorch and TensorFlow releases include forward-mode APIs.}, while 
we prefer to use PyTorch AD as our primary implementation.


\paragraph{Example.}

Figure~\ref{fig:testAAD} illustrates AD-based computation of first- and second-order derivatives for the scalar function $f(x)=x^3/6$. Second derivatives are obtained via nesting (e.g., reverse-over-forward or forward-over-reverse).

\begin{figure}
\centering
\includegraphics[width=0.7\textwidth, keepaspectratio]{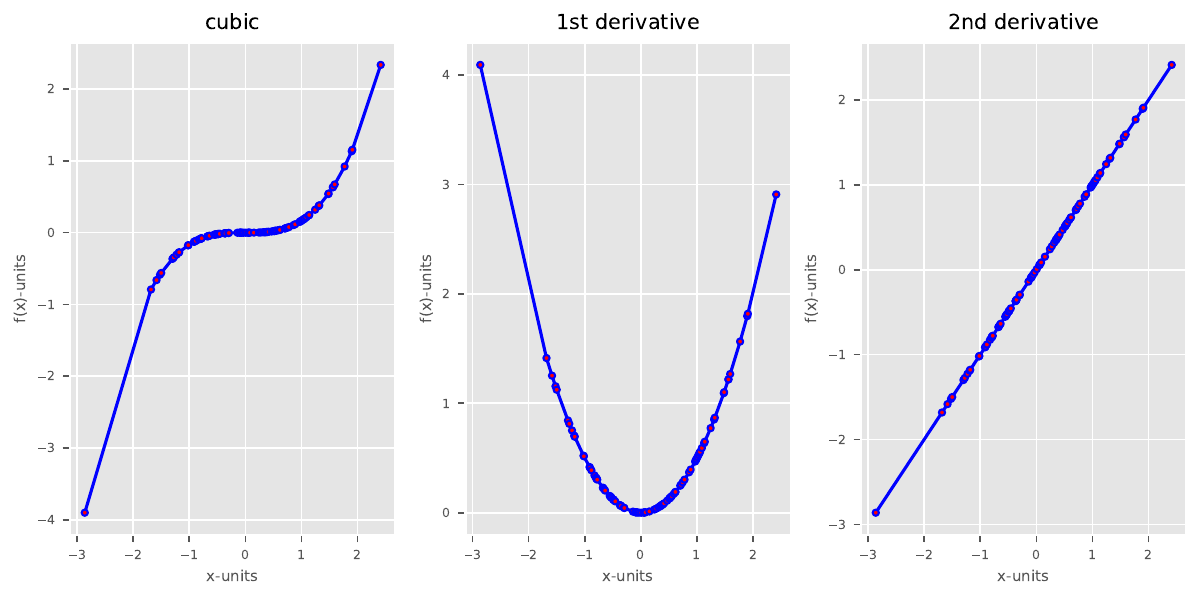}
\caption{\label{fig:testAAD}A cubic function, exact AAD first order and second order derivatives}
\end{figure}

\subsection{Differential machine benchmarks}\label{differential-machine-benchmarks}

AD provides a natural building block for ``differential machines'', i.e., learned models equipped with operators that return derivatives of the learned function. We benchmark two such approaches:
(i) a kernel-based differential operator (gradient/Hessian) in an RKHS (see~\eqref{nablak}), 
 and
(ii) a feed-forward neural network trained from data, with derivatives evaluated by AD.

\paragraph{One-dimensional benchmark.}

Figure~\ref{fig:differentialMlBenchmarks1} shows a one-dimensional test following Chapter~2 protocol. The top row reproduces the function-approximation test. The bottom row shows the derivative estimates on the test set: (left) the ground-truth gradient computed by AD applied to the known $f$; (middle) the kernel gradient operator; (right) two independently trained neural networks evaluated by AD.

\begin{figure}
\centering
\includegraphics[width=0.7\textwidth, keepaspectratio]{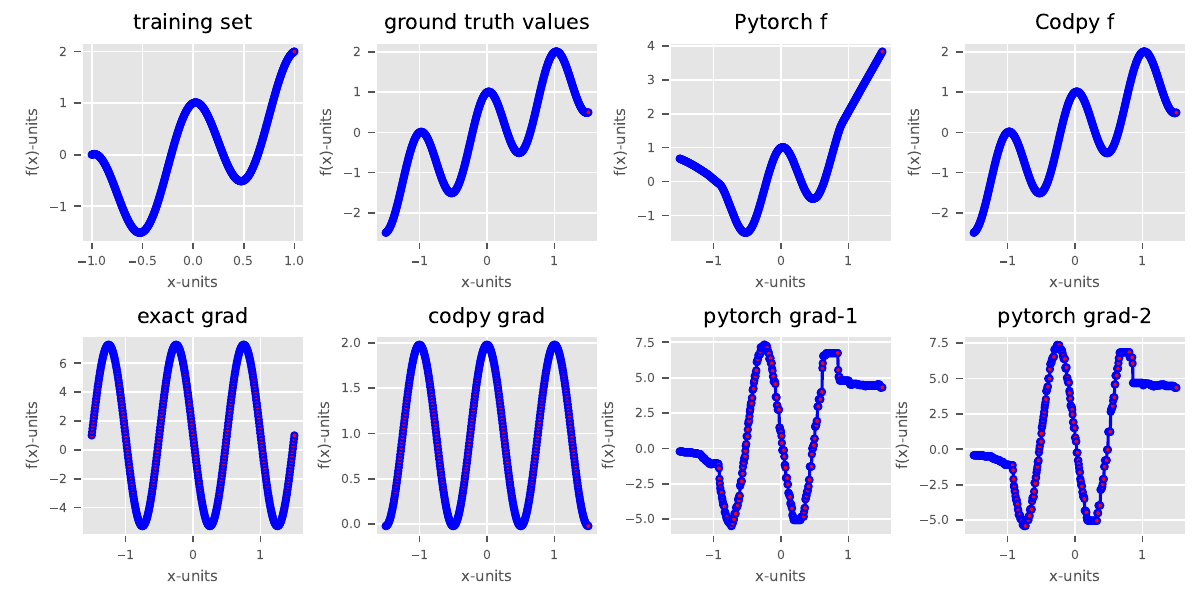}
\caption{\label{fig:differentialMlBenchmarks1}One-dimensional benchmark of differential machines: ground-truth (AD on the known $f$), kernel operator, and two independently trained neural nets (evaluated by AD).}
\end{figure}

\paragraph{Two-dimensional benchmark.}

The same protocol extends to higher dimensions, and Figure~\ref{fig:differentialMlBenchmarks2} shows a two-dimensional test. Concerning these figures, we point out the following facts. 
\begin{itemize}
\tightlist
\item
Neural-network \emph{training} is stochastic (e.g., random initialization, minibatching, Adam), so independently trained models can produce different derivative estimates. \emph{Given a fixed trained model}, the AD evaluation itself is deterministic up to floating-point effects.

\item
In our tests, the kernel-based gradient operator typically attains lower error than the neural-network estimator for the same data budget and tuning. This depends on architecture, regularization, and training protocol, and we therefore present the hyperparameters and random seeds together with the plots.
\end{itemize}

\begin{figure}
\centering
\includegraphics[width=0.7\textwidth, keepaspectratio]{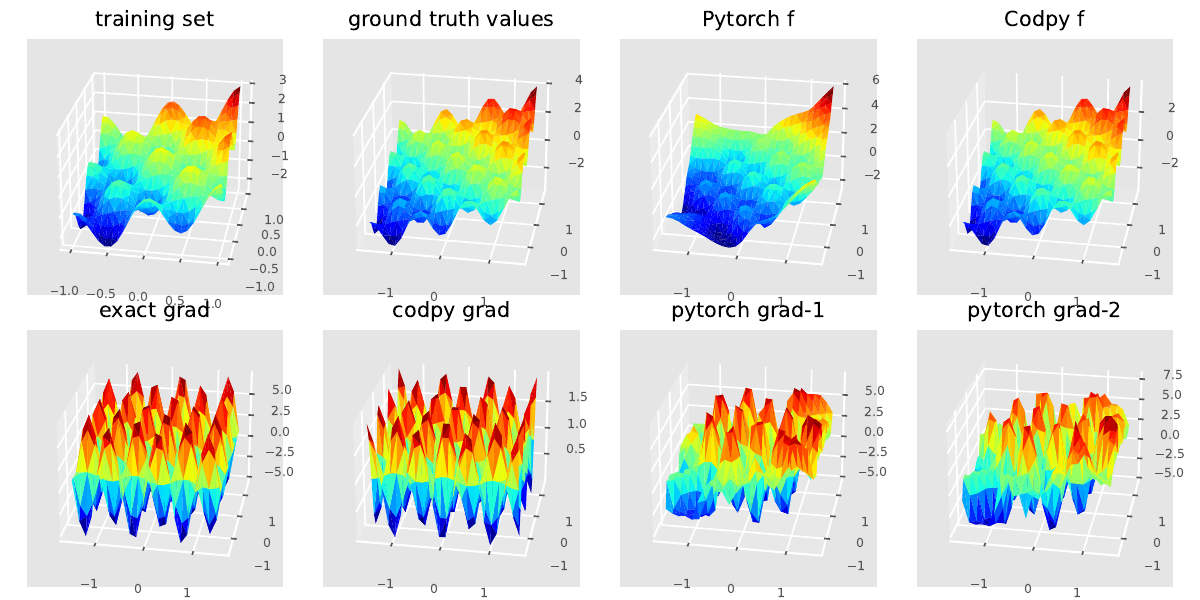}
\caption{\label{fig:differentialMlBenchmarks2}Two-dimensional benchmark of differential machines}
\end{figure}

\subsection{Taylor expansions and differential learning machines}
\label{taylor-expansions-and-differential-learning-machines}

\paragraph{Taylor expansion} 

Taylor expansions using differential learning machines are common for several applications, hence we propose a general function to compute them, that we describe now. We start with the remainder of Taylor expansions.

Let us consider a sufficiently regular, vector-valued map \(f\) defined over \(\mathbb{R}^D\). Considering any sequences of points \(Z,X\) having the same length, the following formula is called a Taylor expansion of order \(p\):
\be\label{eq:taylor}
f(Z) = f(X) + (Z-X)\cdot (\nabla f)(X) + \frac{1}{2}\Big( (Z-X)(Z-X)^T \Big)\cdot (\nabla^2 f)(X) +\ldots+ |Z-X|^p\epsilon(f), 
\ee
where
\be\label{eq:remainder}
R_{p+1}(x,h)
=\frac{1}{p!}\int_{0}^{1} (1-t)^p\, \nabla^{p+1} f(x+th)[\underbrace{h,\ldots,h}_{p+1}]\,dt.
\ee
If $\sup_{t\in[0,1]}\|\nabla^{p+1}f(x+th)\| \le M$ (operator norm), then
\be\label{eq:remainder-bound}
|R_{p+1}(x,h)| \;\le\; \frac{M}{(p+1)!}\,\|h\|^{p+1}.
\ee
(For $f:U\to\mathbb{R}^m$, apply \eqref{eq:taylor} componentwise or interpret $\nabla^k f$ as a tensor-valued multilinear map.)

\paragraph{Learned differential operators}

Given data $X=\{x^j\}_{j=1}^{N_x}\subset U$ and $f(X)=(f(x^1),\ldots,f(x^{N_x}))^\top$, define
\be
\widehat f(\cdot)  =  K(\cdot,X)\,\theta,\qquad
\theta  =  \bigl(K(X,Y)+\alpha\, R(X,Y)\bigr)^{-1} f(X),
\ee
with a kernel Gram matrix $K$ induced by $C^{p+1}$ positive-definite kernel $k$, a symmetric positive semidefinite regularizer $R$, and $\alpha>0$.

Derivatives of $\widehat f$ are obtained by differentiating $K$ in its first argument, that approximates \(\nabla f(x)\),\(\nabla^2 f(x)\) with
\be 
 \nabla f_x = (\nabla K)(\cdot, X)\theta, \quad  \nabla^2 f_x = (\nabla^2 K)(\cdot, X)\theta, 
\ee
where $\theta = (K(Y, X) + \alpha R(Y, X))^{-1}f(X)$, $ \nabla K(\cdot,X) \in \mathbb{R}^{D,N_x}$ and $ \nabla^2 K(\cdot,X) \in \mathbb{R}^{D,D,N_x}$.  

\paragraph{Learned order-2 Taylor approximation.}

Using the learned derivatives, we define the second-order Taylor approximation of $\widehat f$ at $x$ evaluated at $z$ as,
\be\label{eq:learned-taylor}
T_2[\widehat f](z;x)  =  \widehat f(x) \;+\; \nabla \widehat f(x)^\top h \;+\; \tfrac12\, h^\top \bigl(\nabla^2 \widehat f(x)\bigr) h .
\ee

\paragraph{Baselines}

We compare three values. 
\begin{itemize}
\tightlist
\item
  The first one is the reference value for this test. It uses the AAD to compute both \(\nabla f_x,\nabla^2 f_x\).
\item
  The second one, uses a neural network defined with Pytorch together with AAD tools.
\item
  The third one uses the kernel-based Hessian operator as defined earlier in this monograph. 
\end{itemize}
The test is genuinely multi-dimensional, and we illustrate the one-dimensional case in Figure~\ref{fig:taylortest1}.

\begin{figure}
\centering
\includegraphics[width=0.7\textwidth, keepaspectratio]{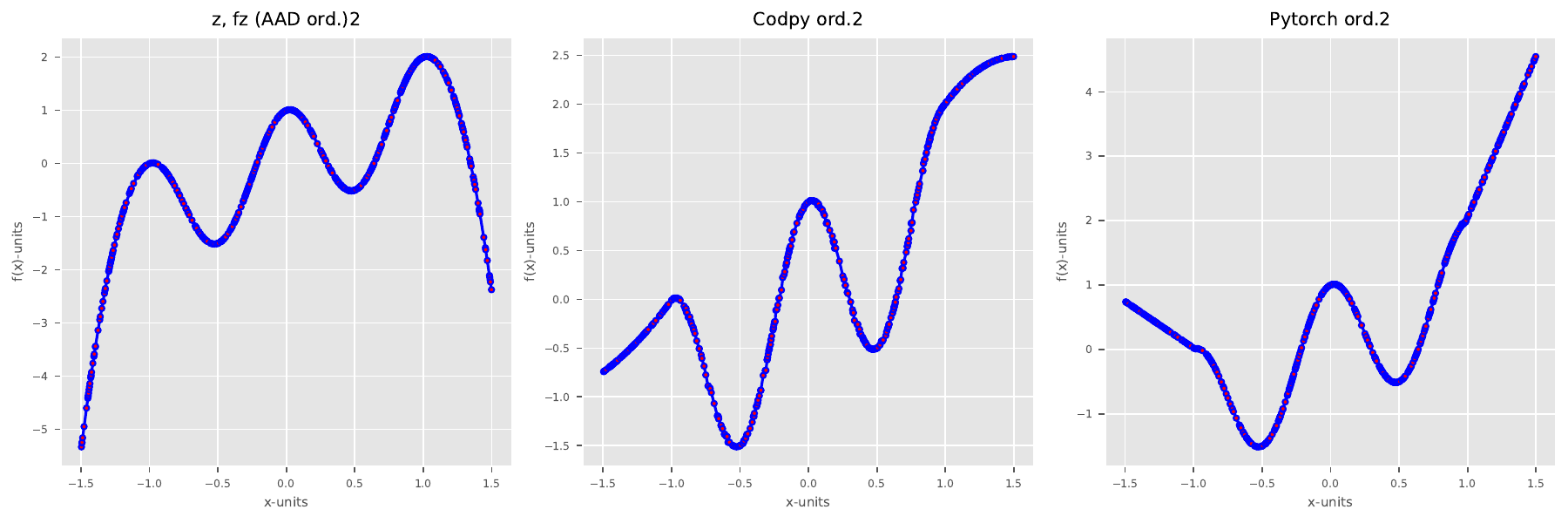}
\caption{\label{fig:taylortest1}A benchmark of one-dimensional learning machine second-order Taylor expansion}
\end{figure}


\section{Discrete high-order approximations}
\label{dix-ddiscrete-high-order-approximations}

Let us denote the Taylor accuracy order \(q>1\). Here, we propose a general \(q\)-point formula in order to approximate any differential operator, accurate at order \(q\). More precisely, consider a sufficiently regular function \(f\), known at \(q\) distinct points \(f(x^k)\), \(x^1 < \ldots<x^q\), and a differential operator \(P^\alpha(\partial) = \sum_{i=0}^{q-1} p_\alpha^i(\partial^{i})\). For any function \(f\), we want to approximate \((P^\alpha(\partial) f)(y) = \sum_{k=1}^{q} f(x^k)\) at some points \(y\). To this aim, consider the Taylor formula
\be
  f(x^k) = f(y)+(x^k-y)\partial f(y) + \ldots +\frac{(x^k-y)^i}{i !}(\partial^i f)(y),\quad k=1,\ldots,q
\ee
with the conventions \(0! = 1\), \(\partial^0f = f\). Multiplying each line by \(\beta_y^k\) and summing leads to
\be
  \sum_{k=1}^{q} \beta_y^k f(x^k) = \sum_{i = 0}^{q-1} (\partial^i f)(y) \sum_{k=1}^{q} \beta_y^k \frac{(x^k-y)^i}{i !}.
\ee
Hence, we rely on a \(q-\)point accurate formula for \(P^\alpha(\partial)\), and we solve the following Vandermonde-type system:
\be\label{VdMalpha}
   \sum_{k=1}^{q} \beta_y^k (x^k-y)^i = (i!) p_\alpha^i,\quad i=0,\ldots,q-1.
\ee
In particular, consider the Vandermonde system
\be\label{VDM}
    A^n \beta^n = (1,0,\ldots,0)^T,\quad A^n = (a_{i,j}^n)_{i,j=0}^q, \quad a_{i,j}^n =  \Big(t^{*n}-t^{n-j}\Big)^j
\ee
for some \(t^{n} \le t^{*} \le t^{n+1}\), and set \(u^*(u^q, ..,u^0)=\sum_{p=0}^q \beta^{n,p}u^{n-p}\). Indeed, there exist\footnote{We refer the reader to~\cite{LeFloch-Mercier:2002}}.
 \(t^{n} \le t^{*} \le t^{n+1}\) such that this operator is of order \(q+2\). The simplest example is the Crank--Nicolson one $u^*(u^1,u^0)=\frac{u^1+u^0}{2}$, $t^* =\frac{t^1+t^0}{2}$.

Conversely, suppose a formula \((Pf)(y^i) = \sum_{k=1}^{q} \beta^k_{y^i}f(x^{i-k})\) is given for distinct points \(y^1<\ldots<y^{N_y}\). To recover \((Pf)f(x^{i})\), \(i=q,\ldots,N_x\), we solve the following linear system:
\be
  (Pf)(x^{i}) = \frac{(Pf)(y^i) - \sum_{k=0}^{q-1} \beta^k_{y^i}f(x^{k})}{\beta^k_{y^i}}, \quad i=q,\ldots,N_x.
\ee


\chapter{Application to reinforcement learning}
\label{application-to-reinforcement-learning}

\section{Introduction}\label{Introduction}

\paragraph{Aim}

This chapter explores the application of kernel methods to reinforcement learning (RL), a framework for sequential decision-making in which an agent interacts with an environment to learn optimal behavior through trial and error\footnote{see \cite{sutton2018}}. The environment is typically modeled as a Markov Decision Process (MDP) $\mathcal{M} = (\mathcal{S}, \mathcal{A}, R, P)$, where the objective is to learn a policy $\pi$ that maximizes expected cumulative reward.

\paragraph{Limitations of deep reinforcement learning.} Classical RL algorithms--such as Q-learning\footnote{see \cite{watkins1992}}, policy gradient methods\footnote{see \cite{silver2014}}, and actor-critic variants--rely on function approximators to estimate value functions. In modern practice, deep neural networks (DNNs) have become the predominant choice due to their expressive capacity. However, DNNs introduce several well-documented challenges, including instability, catastrophic forgetting, sensitivity to initialization, overfitting, and a heavy reliance on hyperparameter tuning\footnote{see \cite{mnih2015, vanhasselt2016}}.

\paragraph{Kernel methods as an alternative.} 

Kernel-based methods offer a principled, nonparametric alternative to DNNs. These methods provide function approximation with well-understood theoretical properties, including generalization bounds and convergence guarantees. Early applications of kernel methods to RL include kernel-based Q-function approximation\footnote{see  \cite{ormoneit2002}} and Gaussian Process Temporal Difference learning (GPTD)\footnote{see \cite{engel2003}}. Subsequent work extended kernel methods to policy gradient estimation in reproducing kernel Hilbert spaces and, more recently, new performance bounds were established for these approaches\footnote{see \cite{Xu2007} and \cite{SF2023}}.

\paragraph{Challenges and motivation.} 

Despite these theoretical strengths, kernel methods are underutilized in RL due to their computational complexity and scalability limitations. Additionally, the lack of software tools tailored to RKHS-based learning has hindered widespread use. Nevertheless, in low- to medium-dimensional problems, kernel methods remain highly competitive--particularly when sample efficiency and interpretability are critical.

\paragraph{Connection to mathematical finance.} 

The structure of RL is closely related to control problems in mathematical finance, where the Hamilton--Jacobi--Bellman (HJB) equation governs decision-making under uncertainty. HJB-based formulations are used to model a wide range of financial problems\footnote{see \cite{LeFloch-Mercier:2020a}}, including option pricing, optimal investment, market making, and algorithmic execution. This connection motivates the application of kernel-based numerical techniques, traditionally used in finance, to reinforcement learning.

\paragraph{Our contribution.} 

The aim of this chapter is to provide a systematic kernel-based framework for reinforcement learning, emphasizing sample-efficient, theoretically grounded algorithms based onthe RKHS methodology. Specifically, our contributions are as follows. 
\begin{itemize}
\item A standardized kernel analysis for value function approximation and Bellman residual estimation, applicable to Q-learning, actor-critic, and HJB-based methods.
\item A practical, sample-efficient kernel RL library,
compatible with Gaussian kernels, clustering-based approximations, and scalable implementations\footnote{see \cite{LeMeMi:2024}}.
\end{itemize}
\noindent The algorithms proposed in this chapter have been implemented in our CodPy Library. Although the focus is on small to moderate datasets, our approach can be extended using clustering, low-rank approximations, or sparsification techniques. Future extensions include incorporating exploration strategies, constructing novel kernels (e.g., convolutional or attention-based), and embedding latent representations for high-dimensional inputs.

\paragraph{Organization.} 

This chapter is structured as follows. 
\begin{itemize}
  \item Section~\ref{Background} reviews foundational RL concepts: MDPs, value functions, Bellman equations, and classical algorithms such as Q-learning, policy gradients, heuristic learning, and HJB-inspired methods.
  \item Section~\ref{Kernel RL Algorithms} introduces our kernel RL framework, including kernel Q-learning, kernel actor-critic, and RKHS-based Bellman residual minimization.
  \item Section~\ref{Experiments} presents empirical results, comparing our approach to deep RL baselines such as DQN and PPO\footnote{see \cite{DQN}  and \cite{PPO}, respectively}.
\end{itemize}


\section{Background}\label{Background}

\subsection{Reinforcement learning}\label{Reinforcement-learning}


\paragraph{General framework for RL.}

 In RL, an \emph{agent} interacts with an
\emph{environment}, over discrete time steps. At each time step \(t\),
the agent observes the current state \(s_t \),
selects an action \(a_t\),
and receives a \emph{reward} signal \(r_t \) along with the
\emph{next state} \(s_{t+1}\) from the environment. 

We adopt the MDP framework which is defined as a tuple
\(\mathcal{M} = (\mathcal{S}, \mathcal{A}, R, P)\), where \(\mathcal{S} \subset \mathbb{R}^D\) is the set of all valid states, \(\mathcal{A}\) is a discrete set of all valid actions, \(R : \mathcal{S} \times \mathcal{A} \rightarrow \mathbb{R}\) is the environment reward function, which provides the reward obtained from taking an action in a given state (\(r_t = R(s_t, a_t)\)) and \(S : \mathcal{S} \times \mathcal{A} \rightarrow \mathcal{S}\) is the environment next state function, which determine the
next states obtained from taking an action in a given state (\(s_{t+1} = S(s_t, a_t)\)).

The rewards and the next state function can be either deterministic or
stochastic, in the latter case there exist joint probability laws
\(\mathbb{P}_R : \mathcal{S} \times \mathcal{A} \rightarrow \mathbb{P}(R)\),
\(\mathbb{P}_S : \mathcal{S} \times \mathcal{A} \rightarrow \mathbb{P}(\mathcal{S})\),
which determine the likelihood of receiving a reward
\(r'\) and transitioning to a state \(s'\) from \(s \in \mathcal{S}\) after taking
action \(a \in \mathcal{A}\): \be\label{tran}
\mathbb{P}_S(s_{t+1} = s' \mid s_t = s, a_t = a), \quad \mathbb{P}_R(r_{t+1} = r' \mid s_t = s, a_t = a)
\ee 
By abuse of notation, we do not distinguish between observed rewards
\(r_t=R(s_t,a_t)\) and their expectations
\(\mathbb{E}[R |s_t, a_t]\). This allows us to simplify the notation by assuming a deterministic reward structure, although we will consider the stochastic case for state transitions while describing HJB equations.

MDP assume the Markov property, which means the future state depends only on the current state and action, not on the sequence of events that preceded it. However, we emphasize that the numerical analysis developed in this chapter holds beyond Markov assumptions.

The goal of RL is to learn a policy \(\pi\) that maximizes the cumulative discounted reward, or \emph{return} \(G_t = \sum_{k=t}^{T_e} \gamma^{k-t} r_k\) at time \(t\),
where \(0 \le \gamma \le 1\) is a discount factor and \(T_e\) is the number of time steps for episode $e$. 

The
\emph{policy} \(\pi(s)=(\pi^1,\dots,\pi^{|\mathcal{A}|})(s)\) is a
probability field describing an agent, determining its actions according to a
probability depending on the states of the environment. This setting
includes both stochastic and deterministic policies, in which later case
$\pi^a(s_t) = \delta^a(a_t)$, where the Kronecker $\delta$ function is defined as $\delta^i(j) = \{1 \text{ if } i=j; 0 \text{ else} \}$.

Our numerical tests deal with learning in the following sense: our agents define a new policy $\pi^{e+1}(\cdot)$ analyzing all past $e$ episodes, and this policy is used in the next game episode $e+1$, covering both online and offline learning.

Let us introduce the buffer, that is, the memory of past $e$ games episodes as a collection of stored
trajectory, denoted \( B^{T^e} = \{(s_t^e, a_t^e, s_{t+1}^e, r_t^e, d_t^e) \}_{t=1}^{T^e}\): states, actions, next states, rewards, dones, the last being a Boolean to indicate if the next state is terminal, corresponding to data terminating a game session. We need this notation for heuristic control, based on episodes, however we consider and note the input data as unordered, possibly unstructured data as buffers \( B^{T} = \{(s_t, a_t, s_{t}', r_t, d_t) \}_{t=1}^{T} \) for kernel bellman error based algorithms (QLearning, Actor Critic, Policy Gradient and HJB).

\paragraph{Bellman equations}
\label{Bellman Equations} 

We now introduce the main definition relative to value and state-action value function and their respective Bellman equations. 

\textit{The state-value function} \(V^{\pi}\) of a MDP is the expected
return starting from a given state \(s\) and following policy \(\pi\),
i.e. \be
V^{\pi}(s) = \mathbb{E}^{\pi} \left[ G_t \mid s_t = s \right] = \mathbb{E}^{\pi} \left[ R(s,\cdot) + \gamma V^{\pi}(s') \mid s_t = s \right], 
\label{valfun_bell}
\ee
where $s'=S(s,\cdot)$.
This is a recursive definition that coincides to the return computation of the last played game ($V^{\pi}(s_t) = G_t$) with a deterministic policy and environment.

\textit{The state-action-value function} \(Q^{\pi}\) measures the
expected return starting from \(s\), taking action \(a\), and following
policy \(\pi\): \be
Q^{\pi}(s,a) = \mathbb{E}^{\pi} \left[ G_t \mid s_t = s, a_t = a \right] = R(s,a) + \gamma \mathbb{E}^{\pi} \left[  Q^{\pi}(s',a) \mid s_t = s, a_t = a \right]
\label{sta_act_valfun}
\ee

\textit{The optimal state-value function} \(V^{\pi^*}(s)\) represents the expected return when starting in state \(s\) and consistently following the optimal policy $\pi^*$ within the environment: \be
 V^{\pi^*}(s) = \mathop{\arg\max}_{\pi} V^{\pi}(s) 
\label{opt_valfun}
\ee \textit{The optimal state-action-value function}
\(Q^{\pi^*}(s,a)\) represents the expected return when starting in state
\(s\) and acting a, and following the optimal policy in the environment, satisfying
\be
Q^{\pi^*}(s,a) = \mathop{\arg\max}_{\pi} Q^{\pi}(s,a)
\label{opt_sta_act_valfun}
\ee

The value functions estimates should satisfy Bellman equation to allow policy improvement. This is the backbone of a class of reinforcement
learning algorithms called \textit{value-based methods}. 
On
the other hand we have \textit{policy-based methods} which goal is to
directly approximate the optimal policy \(\pi^*\).

\subsection{Learning frameworks and control approaches}
\label{Algorithms}

\paragraph{Q-learning}

One of the most prominent value-based algorithms is Q-learning, which
aims to learn the optimal state-action-value function, \(Q^{\pi^*}(s,a)\). 
The core idea
of Q-learning is to solve the optimal Bellman equation \eqref{opt_sta_act_valfun}. For instance, a popular iterative scheme is given by the following iterative rule: 
\be 
Q^{\pi_{n+1}}(s_t,a_t) = Q^{\pi_{n}}(s_t,a_t) + \alpha[R(s_{t},a_t) + \gamma \max_{a} Q^{\pi_{n}}(s_{t+1},a) -Q^{\pi_{n}}(s_t,a_t)],
\label{iter_bell}
\ee where $\alpha$ is a learning rate. Once iterated $N$ times, the so called \textit{greedy} policy corresponds to choosing in-game actions according to $\text{argmax}_{a \in \mathcal{A}} Q^{\pi_{N}}(s,a)$.
However, such deterministic strategies can limit state space exploration, so reinforcement learning algorithms use techniques like the epsilon-greedy strategy to balance exploration and exploitation for effective learning.

Let us have a look at the structure of the scheme \eqref{iter_bell}. This scheme exploits the Markov properties of environments, leveraging the time structure of the buffer data, and can be described as a \textit{backward} time-dependent equation, allowing for numerous optimization techniques. By contrast, in this chapter, we provide a numerical analysis that holds beyond Markov assumptions, trading optimization for generality, and unstructured data, that is, considering $s'_t$ instead of $s_{t+1}$ in \eqref{iter_bell}.

\paragraph{Policy gradient}

Policy gradient methods aim to optimize a policy by directly calculating the gradient of a scalar performance measure with respect to the policy parameters.  These methods fall into two categories: those that directly approximate the policy, and actor-critic methods, which approximate both the policy and the value function.
The general rule can be written as follows: 
\be
 \pi_{n+1}(s) = \pi_n(s) + \lambda A^{\pi_n}(s), 
\label{PG_rule}
\ee
where \( \lambda \) is a learning rate, and \( A^{\pi}\) is a vector field with components \( A^{\pi} = \{A^{1,\pi}, \dots, A^{|\mathcal{A}|,\pi}\} \), under the constraint that the equation \eqref{PG_rule} defines a probability \( \pi_{n+1} \propto \pi_n e^{\lambda A^{\pi_n}} \). The equation \eqref{PG_rule} can be interpreted in the continuous time case as \( \frac{d}{dt} \pi_t(s) = \lambda A^{\pi_t}(s) \). We focus on two cases. 
\begin{itemize}
\item 
The first case is reminiscent of policy gradient methods, where the update aims to improve the value function:
\be \label{eq:advant1}
A^{\pi}(s) = \nabla_{\pi} V^{\pi}(s).
\ee
 
\item 
The second case corresponds to standard actor-critic methods, tailored to minimize the Bellman residuals:
\be \label{eq:advant2}
 A^{\pi^a}(s) = R(s,a) + \gamma V^{\pi}(s') - V^{\pi}(s), \quad s'=S(s,a).
\ee
\end{itemize}
\noindent From a numerical point of view, a common approach is to use a set of parameters \( \theta \) to describe policies \( \pi_n(s) = \pi(s,\theta_n) \), so that \eqref{PG_rule} usually transform into an evolution equation of the parameters \( \theta_{n+1} \leftarrow f(\theta_n)\), where \( f \) depends on the used regression methods.

\paragraph{Hamilton-Jacobi-Bellman equation}

The classical Hamilton--Jacobi--Bellman (HJB) equation is a partial differential equation that arises in continuous-time optimal control. It characterizes the value function of a control problem by capturing how the expected future reward evolves as a function of state dynamics, including both drift and stochastic diffusion terms. While our setting is discrete and data-driven, we draw inspiration from the HJB structure by modeling the transition dynamics via a learned drift and nonparametric stochastic component, leading to a recursive Bellman-type equation. Specifically, we assume the system evolves according to:
\be\label{HJBEq}
    s_{t+1} = s_t + F(s_t, a_t) + \epsilon_t,
\ee
where $F(s_t, a_t)  =  \mathbb{E}[s_{t+1} \mid s_t, a_t] - s_t$ is the drift term, and $\epsilon_t$ is a random perturbation satisfying: $\mathbb{E}[\epsilon_t \mid s_t, a_t] = 0$, with $\epsilon_t \sim \nu(\cdot \mid s_t, a_t)$,
where $\nu(\cdot \mid s, a)$ denotes the unknown conditional distribution of the noise.

This formulation mirrors the structure of controlled stochastic differential equations (SDEs), where $F(s,a)$ plays the role of the drift, and $\epsilon_t$ models the residual stochasticity. The associated recursive value function equation is:
\be
    Q^\pi(s_t, a_t) = R(s_t, a_t) + \gamma \mathbb{E}^\pi \left[ Q^\pi(s', \cdot) \mid s_t, \cdot \right],
\ee
which can be expressed more explicitly in terms of transition probabilities:
\be
    Q^\pi(s_t, a_t) = R(s_t, a_t) + \gamma \int \left[ \sum_{a \in \mathcal{A}} \pi^a(s') Q^\pi(s', a) \right] d\mathbb{P}_S(s' \mid s_t, a_t),
\ee
where $d\mathbb{P}_S(s' \mid s, a)$ denotes the transition probability measure from state $s$ to $s'$ under action $a$, as introduced in~\eqref{tran}. Solving this equation requires estimating both the drift $F(s, a)$ and the transition distribution $d\mathbb{P}_S(s' \mid s, a)$. Traditionally, the latter is approximated by assuming that the noise $\epsilon$ is Gaussian white noise. However, Section~\ref{KHJB} introduces a data-driven approach that estimates these quantities without assuming a parametric form for the conditional noise $\epsilon \sim \nu(\cdot \mid s_t, a_t)$.

\hypertarget{Heuristic-controlled learning}{%
\paragraph{Heuristic-controlled learning}\label{Heuristic-Controlled Learning}}

We define a controller as a deterministic policy parametrized by a set of parameters $\theta$, \( \mathcal{C}^{\theta} : \mathcal{S} \rightarrow \mathcal{A} \),
for \( \theta \in \Theta \), with $\Theta$ being a bounded, closed and convex set, where each $\theta$ defines an agent behavior. At the beginning of each episode $e$, the agent selects $\theta_e$, the parameters for that episode, where $e$ is the episode index. The agent then observes $r_e$, a function of rewards collected during episode $e$, typically defined as the mean of rewards across the episode. The objective is to maximize $\mathbb{E}[r|\theta]$. 

This setting is like a continuous black-box optimization problems, rather than classical reinforcement learning problem, as it is non-associative:  the learned behavior is not tied to specific states of the environment. Only episode-level rewards $r_e$ are observed. 
The action space $\Theta$ is continuous, and there is no concern with the Bellman equations, buffer, and transitions. However we make a strong assumption assuming that the expectation of the rewards $\mathbb{E}[r|\theta]$ is continuous in $\theta$. 

We introduce a model, $R_e(\theta) \sim \mathbb{E}[r|\theta]$, designed to estimate  the conditional expectation of rewards $r$ given $\theta$ based on the first $e$ observations $r_1, \dots, r_e$ and $\overline{\theta_e} = \{ \theta_1, \ldots, \theta_e \}$. 
We consider an optimization function $\mathcal{L}(R_e,\theta)$ and we solve iteratively on each game episode
\be
 \theta_{e+1} = \text{argmax}_{\theta \in \Theta} \mathcal{L}(R_e,\theta)
 \label{HCEq}
\ee

There are multiple choices to pick the functional $\mathcal{L}$, all motivated by providing an exploration incentive. For instance, we experienced that
\be\label{hc_model}
\mathcal{L}(R_e,\theta) = \big( R_e(\theta) - \min_{\overline{\theta_e}}(R_e(\theta)) \big) d(\theta, \overline{\theta_e}),
\ee
where $d(\theta, \overline{\theta_e})$ is the kernel-induced distance gives satisfactory results, but another potential choice is given by $\mathcal{L}(R_e,\theta) = R_e(\theta) - \sigma_e(r|\theta)$ where $\sigma_e(r|\theta)$ is an estimator of the conditional variance of the law $r|\theta$. Observe that this is reminiscent of Upper Confidence Bound (UCB) approach, however the context here is different as we assume continuity in conditional expectations.

This family of problems are paramount for applications: in several situations one can achieve far better results with a good controller than with Bellman-residual approaches on common reinforcement learning tasks, and the same setting applies also to other problems as hyperparameter tuning or PID controllers. Notice however that an expert-knowledge controller is not always available, and no clear way to automatically define one exists to our knowledge. However, due to the importance of these approaches, we propose a kernel method to solve \eqref{HCEq} in Section~\ref{Heuristic-controlled Learning}.


\section{Kernel RL algorithms}
\label{Kernel RL Algorithms}


\subsection{Kernel RL framework}
\label{Kernel RL framework}

\paragraph{Reward regressor.} To estimate the reward function, we use regression. Let \(Z\) be the vector of observed next states-actions
\(Z = [(s_0,a_0), \dots, (s_T,a_T)]\) on the entire buffer. The reward estimator
is determined as
\be \label{eq:update_s_main}
R_{k, \beta}(\cdot) =   K(\cdot, Z)\beta,
\ee where $\beta$ are the parameters fitted through \eqref{FIT} and depend on the buffer size $T$. Similarly, we also define the next state regressor $S_{k, \alpha}(\cdot) = K(\cdot, Z)\alpha$, providing an estimator of the next state function $S(\cdot)$ in the deterministic case, and the conditional expectation $\mathbb{E}[S(\cdot) | s, a]$ in the stochastic case.

\paragraph{Estimating the value functions $V^\pi$ and $Q^\pi$.}

Kernel methods provide efficient methods to estimate $V^\pi$ and $Q^\pi$ values. Kernel methods approximate value functions as
\be \label{eq:update_q_main}
Q_k^\pi(\cdot) = K(\cdot, Z) \theta^\pi, \quad V_k^\pi(\cdot) = K(\cdot, S) \beta^\pi
\ee
To determine
the parameter set \(\theta^\pi\), we solve the Bellman
equation \eqref{sta_act_valfun} on the buffer

\be \label{eq:BE}
Q_k^\pi(Z) = R + \gamma \sum_{a \in \mathcal{A}} \pi^a(S) Q_k^\pi(W^a), \quad W^a = \{(S',a)\},
\ee where we denoted $R = \{r_1, \dots, r_T\} , S = \{s_0, \dots, s_T\}, S' = \{s_1, \dots, s_{T+1}\}$ the rewards, states and next states from the buffer. Observe that the right-hand side requires $|\mathcal{A}| \times T$ evaluations of the regressor of $Q_k^\pi(\cdot)$. 

We plug the expression of the estimator \eqref{eq:update_q_main} into
the previous expression to get

\be \label{eq:BE2}
\Big( K(Z, Z) - \gamma \sum_a \pi^a(S) K(W^a,Z) \Big)\theta^\pi = R.
\ee 
This defines the parameters $\theta$ solving the linear system

\be \label{eq:BE3}
\theta^\pi = \Big( K(Z, Z) - \gamma \sum_a \pi^a(S) K(W^a,Z)\Big)^{-1} R 
\ee
Observe that \(K(W^a,Z)\) is evaluated, or extrapolated, on unseen data \(W^a\), and the formulation \eqref{eq:BE3} inverts a system having size $T \times T$,  which is the square of the buffer size, requiring thus $\mathcal{O}(T^3)$ elementary operations. However, if needed, we can maintain linear computational complexity in the size of the buffer selecting a smaller set $\overline{Z}$, then solving using a least-square inversion matrix $\theta^\pi = \Big( K(Z, \overline{Z}) - \gamma \sum_a \pi^a(S) K(W^a,\overline{Z})\Big)^{-1} R$. To select $\overline{Z}$, one\footnote{see \cite{LeMeMi:2024}} may suggest clustering methods, or a well-chosen subset of $Z$, and also proposes a multiscale method which keeps linear complexity in the buffer size $T$ while considering the whole matrix $K(Z, Z)$; see also Section~\ref{SEK} for an example.

The state action value function is then defined as the regressor
\(Q_k^{\pi}(\cdot) = K(\cdot,Z)\theta^{\pi} \). Similarly, evaluating the Bellman equation
\eqref{valfun_bell} for $V^\pi$ leads to the expression
\be \label{eq:BE4}
 \beta^\pi = \Big(K(S, S) - \gamma K(S', S) \Big)^{-1} \sum_{a \in \mathcal{A}} R_k(S,a)\pi^a (S). 
\ee 
Observe that \(R_k(S,a)\) must therefore be evaluated on
unseen data, hence, a modeling of this function must be provided, as for
instance \eqref{eq:update_s_main}.

\subsection{Kernel Q-learning}
\label{Kernel Q-Learning}

\paragraph{Algorithm.} 

We now describe an iterative algorithm to compute an approximate optimal Q-function using kernel ridge regression. The algorithm iteratively refines its estimate through Bellman updates and interpolation.
This algorithm computes a regressor $Q^{\pi^*}_{k,\theta^{\pi^*}}(\cdot)$ approximating the optimal state-action value function \eqref{opt_sta_act_valfun}. It is an iterative algorithm, each iteration consisting in two main steps. 
\begin{enumerate}
\item Estimation of a first rough set of parameters $\theta^{\pi}_{n+1/2}$.
\item Refining through an interpolation coefficient $\lambda$
\be
\theta_{n+1}^{\pi} = \lambda \theta^{\pi}_{n+1/2} + (1 - \lambda) \theta_{n}^{\pi}.
\ee
\end{enumerate}

\textit{Step 1: Estimating \( \theta^{\pi}_{n+1/2} \)}.
We estimate the parameter \( \theta^{\pi} \) of the \( Q^{\pi}_k \) kernel regressor using \eqref{eq:update_q_main} as follows:
\be \label{qbell} 
 \theta^{\pi}_{n+1/2} = \Big( K(Z, Z) - \gamma \sum_a \pi_{n+1/2}^a(S) K(W^a,Z)\Big)^{-1} R, 
\ee
where the policy is determined from the last iteration $\pi_{n+1/2}^a (\cdot) =  \{1 \text{ if } \mathop{\arg\max}_b Q^{\pi_n}_k(\cdot,b)=a; 0 \text{ else} \}$ and $W, S, R$ are defined in \eqref{eq:BE}.

\textit{Step 2: Computing \( \theta^{\pi}_{n+1} \) via interpolation}.
We further refine at step $n$ interpolating between the previous and newly estimated parameters using a coefficient $\lambda$ aiming to minimize the Bellman residual, defined as
\be \label{eq:bel_res}
 e^{\pi}(Z; \theta) = R + \gamma \max_a Q^\pi_{k, \theta}(W^a) - Q^\pi_{k, \theta}(Z).
\ee
The interpolation coefficient \( \lambda \) is chosen to minimize this residual:
\be
\lambda = \inf_{\beta \in [0,1]} \sum_{z \in Z} e^\pi(z; \beta \theta^{\pi}_{n+1/2} + (1- \beta)\theta_{n}^{\pi}) \big.
\ee
This last quantity is non negative for kernel regressors satisfying \eqref{eq:BE}, as we compute $e^\pi(Z; \theta) = \gamma \big( \max_a Q^\pi_{k, \theta}(W) - \sum_a Q^\pi_{k, \theta}(W) \pi^a \big)  \ge 0$. The iteration stops when the Bellman residual falls below a given threshold, or when further iterations do not yield significant improvement. We observed testally that these steps provide a fast and efficient method to approximate the optimal Bellman equation.

\hypertarget{Kernel-Based Q-value Gradient Estimation}{%
\subsection{Kernel-based Q-value gradient estimation}\label{Kernel-Based Q-Value Gradient Estimation}}

\paragraph{Estimating the derivative of the value function with respect to the policy}

We derive the gradient of the $Q$-function with respect to the policy parameters. Unlike standard policy gradient methods, which optimize the expected return, our approach differentiates the kernel-based Bellman equation to estimate \textit{how the Q-function evolves} with changes in policy parameters. 
We parameterize the policy using a softmax function \eqref{eq:softmax} and denote $\pi(s) = \text{softmax}(y(s))$, with
\(y(s)=\ln \pi(s)\)
where $y(s)$ represents the logits, serving as the underlying parameters of the policy.

To compute the gradient of $Q^\pi$ with respect to the policy parameters, we differentiate the kernel-based Bellman equation \eqref{eq:BE2} with respect to $y^b$:
\be 
\Big( K(Z, Z) - \gamma \sum_a \pi^a(S) K(W,Z) \Big)\partial_{y^b}\theta^\pi = \gamma \sum_a  K(W,Z)\theta^\pi \partial_{y^b} \pi^a(S).
\ee
where $W, S, Z$ are defined in Section~\ref{Kernel RL framework}. Using the kernel-based representation of the Q-function $Q^{\pi}_k(\cdot) = K(\cdot,Z) \theta^\pi$ and the softmax derivative $\partial_{y^b} \pi^a = \pi^a (\delta_{ab} - \pi^b)$, we obtain the following closed-form expression for the gradient of the kernel-based Q-function:
\be \label{gradientBE}
 \nabla_{y} \theta^\pi = \gamma \Big(K(Z,Z) - \gamma \sum_a \pi^a(S) K(W,Z) \Big)^{-1} \sum_a\Big(Q^{\pi}_k(W) \pi^a(\delta_b(a)-\pi^b)\Big)
\ee
Thus, the kernel-based Q-value gradient is given by
\be
\nabla_{y} Q^\pi_k(\cdot) = K(\cdot, Z) \nabla_{y} \theta^\pi.
\ee

It is straightforward to verify $\sum_a \nabla_{y} Q^\pi(W) = 0$. This property ensures that the estimated gradient does not introduce unintended biases.

\paragraph{Gradient of the state value function $V^\pi$}

Similarly, the value function $V^\pi$ is approximated using a kernel-based estimator \eqref{eq:update_q_main}. Differentiating the Bellman equation for $V^\pi$ gives the expression
\be \label{gradientBE2}
  \nabla_y \beta^\pi = \Big(K(S, S) - \gamma K(S', S) \Big)^{-1} \sum_{a \in \mathcal{A}} R_k(S,a)\pi^a(\delta_b(a)-\pi^b)
\ee

\subsection{Kernel Actor-Critic with Bellman residual advantage}
\label{AC}

\paragraph{Advantage function based on Bellman residual error}

In standard actor-critic methods, the advantage function is typically defined as in \eqref{eq:advant2}. However, in our kernel-based approach, we use the Bellman residual error \eqref{eq:bel_res} as the advantage function, which is similar but with the state-action value function instead of the state value function. 

\paragraph{Actor (Policy update using Bellman residual advantage)}

Given our Bellman residual-based advantage function, we define the policy update as
\be \label{UPDATEAC}
\pi_k^a(s_t; \alpha) = \text{softmax} \Big( \ln \pi_k^a(s_t) + \alpha A^{\pi_k}(s_t, a) \Big),
\ee
where $\alpha$ is the learning rate, a hyperparameter, and $A^{\pi_k}(s_t, a)$ is the Bellman residual advantage function. This defines also the probability kernel regressor  $\pi_k(\cdot, \alpha)$ as in \eqref{pi_k}. For policy gradient kernel methods, we observed testally that picking up a learning rate defined as $\alpha = \frac{1}{||A^{\pi_k}||_2^2}$ gives satisfying results. 


\subsection{Kernel non-parametric HJB}\label{KHJB}

\paragraph{Motivation} 
The HJB equation provides a continuous formulation of the optimal control problem, extending the deterministic Bellman equation by incorporating stochastic perturbations. In reinforcement learning, this corresponds to modeling uncertainty in transition dynamics -- particularly when extrapolating value functions to unseen states or actions. From a computational point of view, the HJB viewpoint is a modification of the existing RL algorithms, as Q-Learning or Actor-critic, incorporating stochastic effects.

We introduce a kernel-based, data-driven discretization of the HJB equation, which enables the construction of an explicit \textit{transition operator matrix} $\Gamma(P^a)$, which is a stochastic, or transition probability matrix. This operator is not derived from parametric assumptions or empirical next-state mappings, but is instead learned. In particular, 
we describe the martingale optimal transport (MOT), which, together with the kernel ridge regression \eqref{Pk}, propose a general approach to compute this operator, ensuring compatibility with a learned drift.

The HJB formulation amounts to modify the algorithms introduced for solving the Bellman equation, given by \eqref{eq:BE3} in our RKHS numerical framework, taking into account stochastic effects. Such stochastic effects comes either from uncertainty or distributional variability in off-policy, or limited-datasettings. As a result, it offers a principled alternative for computing value functions in environments where transition dynamics are complex or partially observed.

The MOT approach enforces a martingale structure aligned with a kernel-regressed drift model, providing a structured way to approximate state evolution without relying on conditional densities.

\paragraph{Modeling the drift term}
We begin by modeling the deterministic component of the HJB equation, commonly referred to as the \emph{drift}, denoted by $F$ in \eqref{HJBEq}. The environment dynamics are assumed to follow:
\be
    s_{t+1} = s_t + F(s_t, a_t) + \epsilon_t,
\ee
where $F(s_t, a_t)  =  \mathbb{E}[s_{t+1} \mid s_t, a_t] - s_t$ represents the expected state change, and $\epsilon_t$ is a zero-mean conditional noise: $\mathbb{E}[\epsilon_t \mid s_t, a_t] = 0$. We model $F$ using a kernel regressor:
\be
    F_k(\cdot) = K(\cdot, Z) \theta,
\ee
where $Z = \{z_t = (s_t, a_t)\}_{t=1}^T$ is an empirical dataset of state-action pairs, $K(\cdot, Z)$ is a positive-definite kernel over $\mathcal{S} \times \mathcal{A}$, and $\theta \in \mathbb{R}^{T \times D}$ are fitted coefficients for each state dimension. We define:
\begin{align}
    s_{t+1}^a &= s_t + F_k(s_t, a), \\
    p_{t+1}^a &= (s_{t+1}^a, a),
\end{align}
as the predicted next state and corresponding state-action pair under action $a$.

\paragraph{Martingale optimal transport for transition estimation}
We aim to estimate the transition law:
\be
    \tau_k(s' \mid s, a) \approx \mathbb{P}(s_{t+1} = s' \mid s_t = s, a_t = a).
\ee
Rather than assuming Gaussian noise, we construct an empirical approximation using MOT, consistent with the drift constraint $\mathbb{E}[s' \mid s, a] = s + F_k(s, a)$. Let $P^a = \{(s_t + F_k(s_t, a), a)\}_{t=1}^T$ be the predicted next state-action pairs under action $a$. Let $Z$ denote the empirical support of observed state-action pairs, and consider the corresponding empirical measures $\mu_Z$ and $\mu_{P^a}$ supported on $Z$ and $P^a$, respectively.

The MOT problem seeks stochastic couplings $\Pi(P^a \mid Z), \Pi(Z \mid P^a) \in \mathbb{R}^{T \times T}$ satisfying marginal constraints and a martingale condition. Specifically, for each source $z_i = (s_i, a_i) \in Z$, the forward plan $\Pi(P^a \mid Z)$ satisfies:
\be
    \sum_{j=1}^T \pi_{ij} p_j = z_i + F_k(z_i),
\ee
a martingale constraint for each row $z_i$, where $p_j \in P^a$ is the $j$-th target sample. This enforces that the expected target under the transport plan matches the predicted drift at each source location.

The composition of the forward and backward couplings yields a soft transition operator acting over the empirical dataset. Thus, the conditional transition density $\tau_k(\cdot \mid \cdot)$ is not constructed explicitly as a density function, but instead is represented implicitly through the induced stochastic matrix $\Gamma(P^a)$ described next.

\paragraph{Constructing the transition operator}

For each action $a \in \mathcal{A}$, we compute a separate transition operator $\Gamma(P^a)$, used to model soft transitions in the HJB equation. These operators are later aggregated across all actions using the policy weights $\pi^a(S)$.

Let $P^T  =  \left\{ \left(s_t + F_k(s_t, a_t), a_t\right) \right\}_{t=1}^T$ denote the set of predicted next state-action pairs using the drift model $F_k$ and actions $a_t$ from the dataset. We compute the transition probabilities $\Gamma(P^a) \in \mathbb{R}^{T \times T}$ following the guidelines for the MOT; see Section~\ref{Transition-Probabilities-With-Kernels}, which determine $|\mathcal{A}|$ transition probability matrices. All matrices are of shape $T \times T$ and row-normalized to ensure that $\Gamma(P^a)$ is \textit{row-stochastic}:
\be
\sum_j \Gamma(P^a)_{ij} = 1, \quad \Gamma(P^a)_{ij} \geq 0 \quad \text{for all } i.
\ee

In the full HJB update (see~\eqref{HJBMAT}), we compute $\Gamma(P^a)$ \emph{for each} $a \in \mathcal{A}$, and aggregate them using the policy weights $\pi^a(S) \in \mathbb{R}^T$. This yields a policy-weighted mixture of soft transitions in the value estimation.

Observe that the entries $\Gamma(P^a)_{tu}  =  \tau_k(s_{u+1}^{a} \mid p_{t+1}^{a_t})$ involve \emph{extrapolated} points rather than purely observed data. Specifically, $p_{t+1}^{a_t} = (s_t + F_k(s_t, a_t), a_t)$ is computed using a kernel-based regression model of the drift, not directly from the buffer, the extrapolation being presented in Section~\ref{Perturbative-kernel-regression}.
Similarly, $s_{u+1}^a$ refers to a predicted next state under action $a$, given the estimated drift.

Thus, the transition kernel $\tau_k$ approximates soft transitions between \emph{regressed} state-action pairs, not necessarily those that were observed in the dataset. This generalization allows us to define a smooth operator that supports reasoning about unseen transitions while respecting empirical structure via MOT.

\paragraph{HJB-modified value function approximation} 
We define a kernel-based approximation $Q_{k,\theta}^\pi$ to the value function $Q^\pi$, evaluated on the dataset using a discretized, HJB-inspired fixed-point equation. The formulation integrates over both the uncertainty in future transitions and the stochasticity in the policy:
\be\label{HJB}
Q_{k,\theta}(s_t,a_t) = r(s_t,a_t) + \gamma \sum_{u=1}^T \tau_k(s_{u+1}^{a} \mid p_{t+1}^{a_t}) 
\left[ \sum_{a \in \mathcal{A}} \pi^a(s_{u+1}^a) Q_{k,\theta}(s_{u+1}^a,a) \right].
\ee

Here, $\tau_k(s_{u+1}^a \mid p_{t+1}^{a_t})$ represents a kernel-smoothed estimate of the conditional transition density from query state-action pair $p_{t+1}^{a_t} = (s_t + F_k(s_t, a_t), a_t)$ to target state $s_{u+1}^a$. The inner summation marginalizes over actions, weighting each action-specific value by the policy $\pi^a$ evaluated at the target state. This results in a soft, probabilistic Bellman operator that captures both drift-based transitions and policy-induced uncertainty.

\paragraph{Matrix representation and transition operator} 
To obtain a matrix formulation, define $\Gamma(P^a) \in \mathbb{R}^{T \times T}$ as the transition operator for action $a$, with entries $\Gamma(P^a)_{tu}  =  \tau_k(s_{u+1}^{a} \mid p_{t+1}^{a_t})$.
We also define the policy matrix $\pi^a(S) \in \mathbb{R}^{T}$, representing the probability of taking action $a$ at each training state $s_t$. We lift this into a diagonal matrix $\operatorname{diag}(\pi^a(S)) \in \mathbb{R}^{T \times T}$ to properly weight the contribution of each row of $\Gamma(P^a)$. The resulting operator $\operatorname{diag}(\pi^a(S)) \Gamma(P^a)$ accounts for transitions under both the model dynamics and policy distribution.

Then, the linear system defining the HJB-modified value function approximation becomes:
\be\label{HJBMAT}
\theta = \left( K(Z, Z) - \gamma \sum_{a \in \mathcal{A}} \operatorname{diag}(\pi^a(S)) \Gamma(P^a) K(P^a, Z) \right)^{-1} R,
\ee
where $R \in \mathbb{R}^T$ is the vector of observed rewards and $K(P^a, Z) \in \mathbb{R}^{T \times T}$ denotes the kernel evaluations between the projected points $P^a = \{(s_t + F_k(s_t,a), a)\}_t$ and the training inputs $Z = \{(s_t, a_t)\}_t$. This system generalizes the Bellman kernel system \eqref{eq:BE3}, replacing deterministic transitions with soft, transport-based dynamics that integrate over actions via the policy $\pi$.

\subsection{Heuristic-controlled learning}\label{Heuristic-controlled Learning}

We describe a computationally efficient kernel-based method for solving the heuristic control problem introduced in Section~\ref{Heuristic-Controlled Learning}.

Let $\theta \in \Theta \subset \mathbb{R}^d$ denote the controller parameters, and suppose we aim to estimate the conditional expectation $\mathbb{E}[r \mid \theta]$, where $r$ is the observed reward from a game episode. We approximate this expectation using a kernel ridge regression model of the form:
\be
R_{k,\lambda_e}(\theta) = K(\theta, \overline{\theta}_e) \lambda_e,
\ee
where $K(\cdot, \cdot)$ is a chosen kernel function, $\overline{\theta}_e = \{\theta_1, \dots, \theta_e\}$ is the set of past parameter samples, and $\lambda_e$ are the fitted regression weights based on the past rewards $\{r_1, \dots, r_e\}$ using a standard regularized least squares loss (see~\eqref{FIT}).

To guide exploration, we maximize the acquisition function $\mathcal{L}$ introduced in~\eqref{hc_model}. This function includes a distance term $d(\theta, \overline{\theta}_e)$ that quantifies novelty of a new candidate $\theta$ relative to past samples. A practical choice for $d$ is:
\be
d(\theta, \overline{\theta}_e) = \inf_{i} d_k(\theta, \theta_i),
\ee
where $d_k$ may be defined via the maximum mean discrepancy (MMD) or another kernel-induced metric (see~\eqref{dk}).

We assume the parameter domain $\Theta$ is a known convex compact subset of $\mathbb{R}^d$ and that we can sample uniformly over it. To solve the acquisition maximization problem~\eqref{HCEq}, we employ an adaptive sampling strategy:
\be
\theta_{n+1} = \arg\max_{\theta \in \Theta_n} \mathcal{L}(R_{k,\lambda_e}, \theta), \quad \text{with} \quad \Theta_n  =  \overline{\theta}_e \cup \Theta_{N,n},
\ee
where $\Theta_{N,n}  =  (\theta_n + \alpha^n \Theta_N) \cap \Theta$ is a local neighborhood of candidate points. Here, $\Theta_N$ denotes a fixed-size i.i.d. sample from the prior over $\Theta$, $\alpha \in (0,1)$ is a concentration parameter controlling the shrinking search region, and $n$ is the current iteration.

This procedure balances exploration (through kernel-induced distance) and exploitation (via expected reward), and performs well in high-dimensional, black-box optimization tasks where gradient-based methods are not applicable.


\section{Numerical illustrations}
\label{Experiments}

\subsection{Setup and kernel configuration}

We benchmark five kernel-based reinforcement learning algorithms against two widely used baselines\footnote{see  \cite{PPO} and \cite{DQN}, respectively}: Proximal Policy Optimization (PPO) (denoted \texttt{PPOAgent}) and Deep Q-Network (DQN) (\texttt{DQNAgent}). The kernel-based methods include a heuristic controller-based algorithm and four Bellman residual solvers: Actor-Critic (\texttt{KACAgent}), kernel-based Q-value gradient estimation ({KQGE}), standard kernel Q-learning ({KQLearning}), and its HJB version ({KQLearningHJB}); see Section~\ref{AC} for more details.

All methods are evaluated on two standard benchmark environments from Gymnasium suite\footnote{see \cite{Gymnasium:2024}; Gymnasium, the maintained RL environment API by the Farama Foundation (formerly OpenAI Gym).}: \texttt{CartPole-v1} and \texttt{LunarLander-v3}. Our primary metric of interest is sample efficiency. While this criterion does not favor PPO--an algorithm designed for environments with extensive interactions--it remains a crucial baseline due to its robustness and widespread use.

All tests are conducted on a CPU-based platform\footnote{AMD 7950X3D processor.} We fix the discount factor $\gamma = 0.99$, a value found to yield optimal performance for DQN in preliminary tests.

Each test is repeated independently over ten runs, each with $E = 100$ episodes. We present the mean and standard deviation of two quantities: cumulative rewards and cumulative training time. These are visualized in paired charts for each environment.

We employ the standard Matérn kernel described in Section~\ref{maps-and-kernels} with preprocessing of the inputs via a standard mean map \eqref{standardmap}.
To maintain computational efficiency, the kernel buffer size is capped at 1000. This ensures that learning remains within sub-second runtimes per iteration. For illustration, we focus on the HJB-enhanced version of kernel Q-learning, although similar observations hold for the actor-critic and policy gradient variants.

\subsection{CartPole}\label{Cartpole}
\paragraph{Environment description} The \texttt{CartPole-v1} environment involves balancing a pole on a moving cart by applying discrete left/right forces. The state space is $\mathcal{S} \subset \mathbb{R}^4$, and the action space is binary: $\mathcal{A} = \{0, 1\}$. At each timestep, the agent receives a reward of +1. A trial is considered successful if cumulative rewards reach 1000, in which case the policy is no longer updated.

\paragraph{Training setup.}
To ensure fairness, all algorithms are trained for a maximum of 500 steps per episode and are stopped if cumulative training time exceeds 10 seconds. Exceptions are made for \texttt{PPOAgent} and controller-based agents, which are less time-sensitive.

This environment is fully deterministic, which favors residual Bellman error-based methods. We also include a baseline heuristic controller of the form:

\be
\mathcal{C}_\theta(s) = \text{sign}(\langle \theta, s \rangle),
\ee

with randomly initialized weights $\theta$. In our tests, {KQLearningHJB} achieved higher scores than its standard counterpart, at the cost of increased runtime due to the computation of transition probability matrices. This underscores the tradeoff between sample efficiency and computational cost.

\begin{figure}[!htb]
\centering
\begin{subfigure}{0.48\textwidth}
\centering
\includegraphics[width=\textwidth]{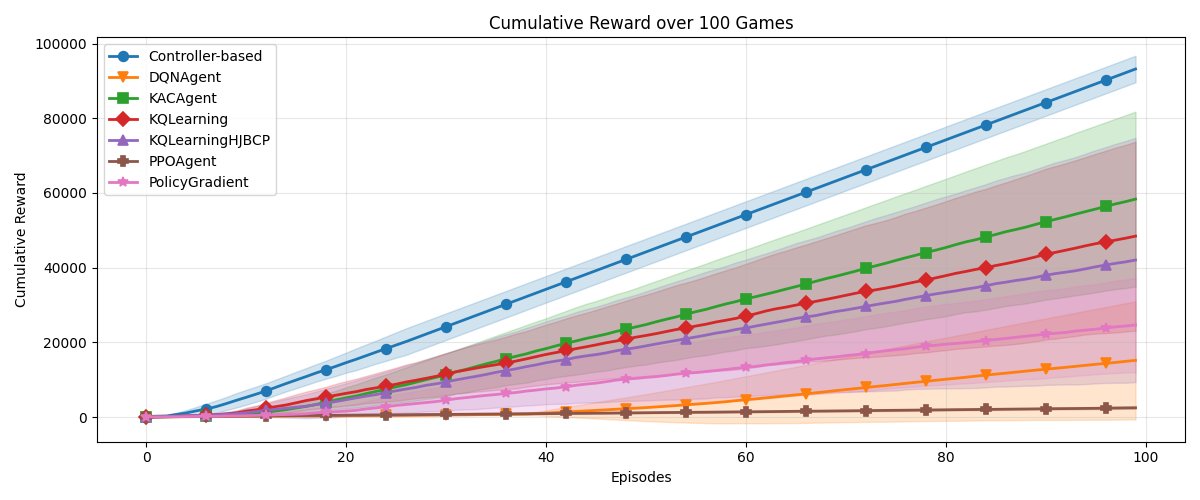}
\caption{Mean cumulative reward per episode}
\label{fig:image1}
\end{subfigure}
\hfill
\begin{subfigure}{0.48\textwidth}
\centering
\includegraphics[width=\textwidth]{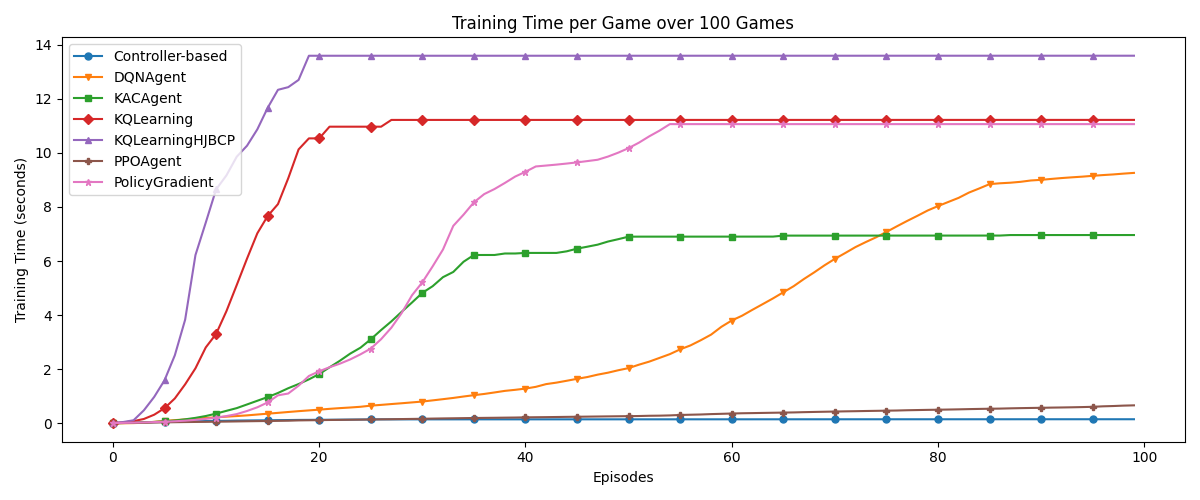}
\caption{Mean cumulative training time per episode}
\label{fig:image2}
\end{subfigure}
\end{figure}

\subsection{LunarLander}\label{Lunar-Lander}

\paragraph{Environment description.} The \texttt{LunarLander-v3} environment involves controlling a lunar module to land safely on a designated pad. The agent operates in an 8-dimensional state space (including positions, velocities, angles, and ground contact flags) with a discrete action space $\mathcal{A} = \{0, 1, 2, 3\}$ corresponding to no-op, left engine, right engine, and main thruster.

\paragraph{Training setup.} In this environment, there is no default limit on the number of timesteps per episode. As a result, some policies may converge to a "hovering" behavior, causing episodes to continue indefinitely. This leads to high variability in the number of timesteps experienced across different algorithms over a fixed number of episodes, resulting in some policies being trained on significantly more data than others and introducing an unfair evaluation. To address this, we impose a hard limit of 2000 timesteps per episode. Furthermore, to ensure fairness across methods, training for all algorithms is stopped once the cumulative training time exceeds 50 seconds.

\paragraph{Short-horizon training (100 episodes).} 

\begin{figure}[!htb]
\centering
\begin{subfigure}{0.48\textwidth}
\centering
\includegraphics[width=\textwidth]{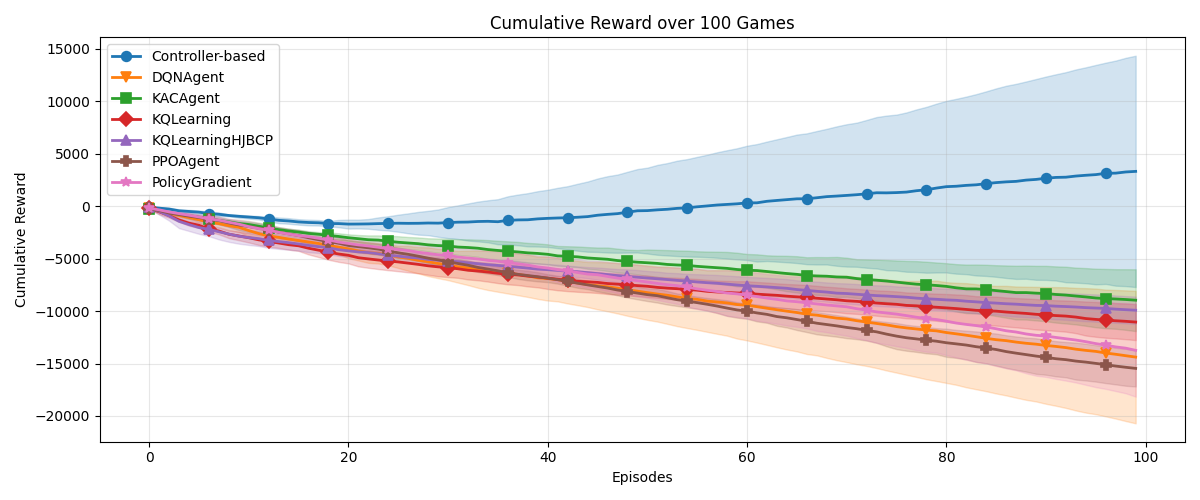}
\caption{Mean cumulative reward per episode}
\label{fig:LNR}
\end{subfigure}
\hfill
\begin{subfigure}{0.48\textwidth}
\centering
\includegraphics[width=\textwidth]{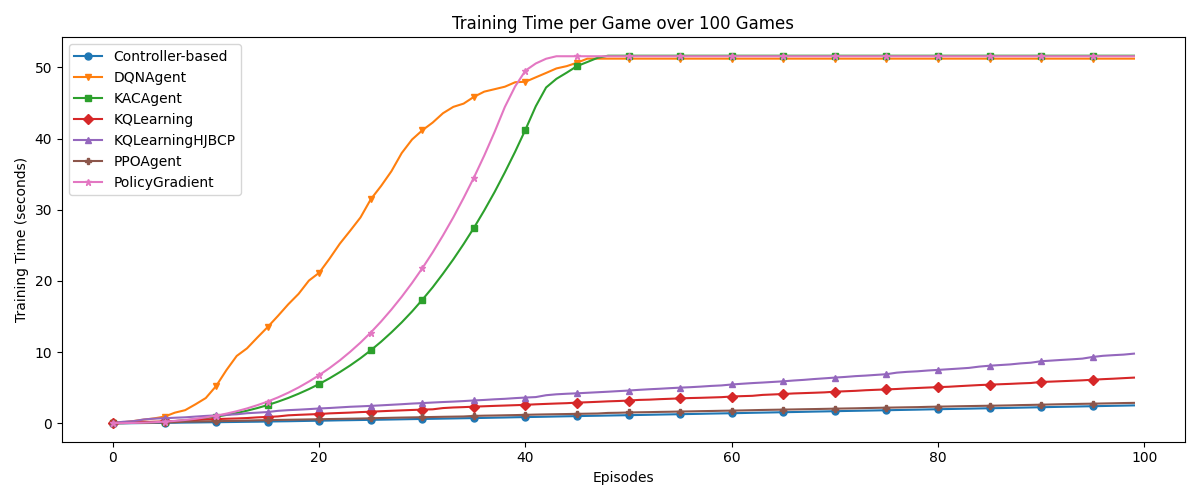}
\caption{Mean cumulative training time}
\label{fig:LNT}
\end{subfigure}
\end{figure}

The best-performing algorithm on this task is the heuristic-controlled learning approach, which achieved strong performance but with high variance. This result relies on a finely-tuned controller $\mathcal{C}_\theta(s)$ with 12 learnable parameters\footnote{see \cite{HeuristicLunar}}.
Unlike CartPole, the LunarLander environment is non-stationary ---the lunar surface changes in each episode. Such variability negatively impacts\footnote{see  \cite{FuMePrNaSh, sutton2018}} Bellman residual-based methods like {KQLearning}, {KACAgent}, and {PolicyGradient}.

Hence, to improve scalability, we introduce an episode-clustered version of {KQLearning} and its HJB variant. This modification enhances training time efficiency, albeit at the cost of slightly lower scores. A detailed description of this variant is provided in Section~\ref{SEK}. The same clustering strategy is applicable to other residual-based algorithms such as {KACAgent} and {KQGE}.

\paragraph{Extended training (600 episodes).} 

In the previous paragraph we ran the algorithms for 100 episodes, to emphasize sample efficiency. However, this setup diverges from standard reinforcement learning benchmarks, which typically allow algorithms to engage in extended environmental interaction.

To provide a more balanced evaluation, we re-ran the tests for a longer number of episodes while maintaining the 50-second time constraint. This offers a fairer comparison with kernel-based agents, while acknowledging that {PPOAgent} is designed for longer training horizons and performs best with large-scale experience collection.

\begin{figure}[!htb]
\centering
\begin{subfigure}{0.48\textwidth}
\includegraphics[width=\textwidth]{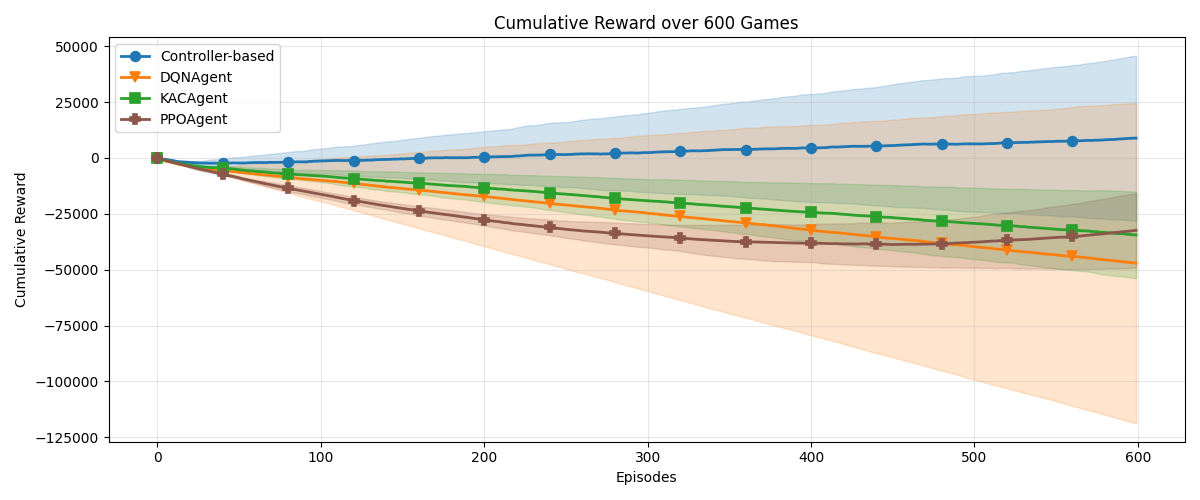}
\caption{Mean cumulative reward per episode (600 episodes)}
\label{fig:PPOLNR}
\end{subfigure}
\hfill
\begin{subfigure}{0.48\textwidth}
\includegraphics[width=\textwidth]{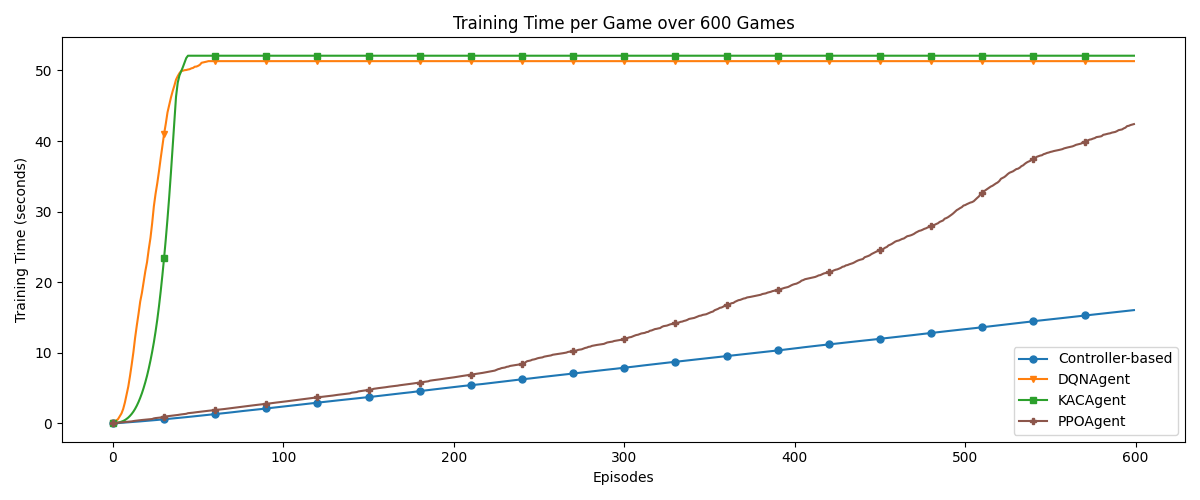}
\caption{Mean cumulative training time per episode}
\label{fig:PPOLNT}
\end{subfigure}
\end{figure}

Figure~\ref{fig:PPOLNR} illustrates cumulative rewards over 600 episodes. Initially, the {PPOAgent} underperforms due to an exploration phase during which the policy collects diverse experience through environment interaction. This behavior is expected and reflects {PPOAgent}’s design, which prioritizes long-term improvement over early performance.

Toward the later episodes, {PPOAgent} shows a clear upward trend in reward accumulation. This progression confirms that given sufficient interaction time, {PPOAgent} can outperform sample-efficient methods. Thus, for the sake of testal transparency, we acknowledge that {PPOAgent} is likely to excel under longer training regimes and should be interpreted accordingly.


\section{Clustering methodology using kernel baseline RL algorithms}
\label{SEK}

Clustering is a natural and efficient idea to lower computation time of kernel methods. The general idea is to define a partition of the buffer $B^{T} = \{(s_t, a_t, s_{t}', r_t, d_t) \}_{t=1}^{T}$ into $I$ clusters $B_i^{T} = \{(s_t^i, a_t^i, s_{t}'^i, r_t^i, d_t^i) \}_{t=1..T^i}$. This defines $I$ different agents $\mathcal{A}_i$, defined as local solver of the Bellman equation \eqref{eq:bel_res}, each using its own kernel $k_i$. Observe that this approach can be compared to the \textit{symphony of experts}\footnote{see \cite{Jon2023}}, but we opted for a different approach that we describe now.

A simple and efficient idea is to define a cluster for each in-game episode, for which we solve the optimal Bellman equation as in Section~\ref{Kernel Q-Learning}, or its HJB version in Section~\ref{KHJB}, and denote its value $q_{k_i}^*(\cdot)$. This is quite adapted to the LunarLander game, for which each episode usually contains a small number of 100/200 steps, implying competitive learning times, as shown in Figure~\ref{fig:LNT}. We experienced that, with this approach, considering local gamma values is paramount, and we took $\gamma^i=\exp(-\frac{-\ln(T^i)}{T^i})$. 

Agents defines a natural distance, which is the closest distance between a given state $s$ and the buffer used by each agent:
\be
d(z,\mathcal{A}_i) = d_i(z,z^{i^*(s)}), \quad i^*(s) = \mathop{\arg\inf}_i d_i(z,Z_i), z=s,a, \ \quad Z_i = \{s_t^i, a_t^i\}_{t= 1 \ldots T^i},
\ee
where $d_i(\cdot,\cdot)$ is a distance. A natural choice for kernel methods is to pick-up the kernel discrepancy as a distance, $d_i(x,y) = k_i(x,x)+k_i(y,y) - 2 k_i(x,y)$, that is, the distance selected for our tests.

During a game, considering a state $s$, we pick-up the action $\mathop{\arg\max}_{a=1,..,|\mathcal{A}|} q^*(s,a)$, where the optimal state value function $q^*(s,a)$ is approximated as follows: we first identify the closest $I \times |\mathcal{A}|$ points to $s,a$, denoted $X = \{s_{i^*(s)},a_{i^*(s)}\}_i$, 
together with their local optimal q-values denoted $Y = \{q_i(s_{i^*(s)},a_{i^*(s)})\}$. These values $(X,Y)$ are then extrapolated in order to provide a value $q_k^*(s,a)$. 


\chapter{Application to mathematical finance}\label{application-to-mathematical-finance}

\section{Aim of this chapter}

This chapter presents practical applications of machine learning techniques within the domain of mathematical finance. The exposition is structured as follows. 
\begin{itemize}
  \item The first part focuses on \emph{time-series modeling and forecasting}. Adopting an econometric perspective, we consider observed time-series data and construct data-driven stochastic processes that are consistent with the empirical observations. These models serve as a foundation for probabilistic forecasting and simulation.

    \item The second part addresses the approximation of \emph{pricing functions}, which, in mathematical finance, are often defined as conditional expectations of future payoffs. We do not describe classical pricing methods--such as Monte Carlo simulations or PDE solvers--to evaluate these quantities. We demonstrate instead how supervised learning techniques, particularly kernel methods, can efficiently approximate these pricing functions from simulated or historical data, where they are given. This allows fast evaluation, differentiation (Greeks), and extrapolation to unseen scenarios, facilitating practical tasks such as hedging, PnL attribution, and risk management. The purpose is to provide  \textit{real-time} methods to risk measurements, as pricing functions are computationally intensive and difficult to use intradays.
    \item Two use cases of interest are also included. The first one concerns reverse-stress test, which is a challenging inversion problem. The second concerns investment strategies with trading signals, which applies the RKHS framework to classical portfolio management.

  Importantly, while this chapter emphasizes data-driven interpolation/extrapolation of pricing functions, the framework also connects to classical PDE-based pricing formulations. In particular, Section~\ref{KHJB} introduces a kernel Hamilton-Jacobi-Bellman (HJB) approach, which provides a data-driven alternative for approximating value functions in stochastic control problems, which is a general setting to pricing methods for mathematical finance.
\end{itemize}


\section{Nonparametric time-series modeling}
\label{free-time-series-modeling}

\paragraph{Setting and notation.}

We consider the problem of modeling and forecasting a stochastic process \( t \mapsto X(t) \in \mathbb{R}^D \), observed on a discrete time grid \( t^1 < \ldots < t^{T_x} \). The data are organized as a three-dimensional tensor:
\be\label{TS}
X = \left(x_d^{n,k}\right)_{d=1,\ldots,D}^{n=1,\ldots,N_x,\ k=1,\ldots,T_x} \in \mathbb{R}^{N_x \times D \times T_x},
\ee
where the data are as follows. 
\begin{itemize}
  \item \( D \) is the dimensionality of the process (e.g., the number of assets),
  \item \( T_x \) is the number of discrete time points,
  \item \( N_x \) is the number of observed trajectories (typically \( N_x = 1 \) in financial applications).
\end{itemize}

\paragraph{Abuse of notation.}

In what follows, we write $X^k$ to denote the entire set of observed data at time index $k$, i.e.,
\be
X^k  =  \left(x_d^{n,k}\right)_{d=1,\ldots,D;\ n=1,\ldots,N_x} \in \mathbb{R}^{N_x \times D}.
\ee
Similarly, we denote by $X^n$ the entire trajectory of the $n$-th sample across all features and time steps:
\be
X^n  =  \left(x_d^{n,k}\right)_{d=1,\ldots,D;\ k=1,\ldots,T_x} \in \mathbb{R}^{D \times T_x}.
\ee
When expressions like $X^{k+1} - X^k$ or $\log X^k$ appear, they are to be understood \emph{componentwise}, i.e., applied to each entry $x_d^{n,k}$. Similarly, increments and transformations on $X$ are performed per trajectory and per feature unless specified otherwise.

This \emph{abuse of notation} allows us to write expressions more compactly, while implicitly handling the full index structure. The notation naturally supports applications with multiple cross-sectional samples (indexed by $n$), such as customer data or multi-asset financial modeling.

\paragraph{Illustrative dataset.}

To ground the methodology, we use historical daily closing prices of three large-cap technology stocks: Google (GOOGL), Apple (AAPL), and Amazon (AMZN), spanning January 1, 2020, to December 31, 2021. The trajectories are visualized in Figure~\ref{fig:plot1}.

This dataset serves as a running example throughout for demonstrating model calibration, trajectory generation, and statistical evaluation. Summary information is provided in Table~\ref{tab:settings}, and Table~\ref{tab:empirical_stats} shows the empirical moments of the log-returns, which we use as a baseline for validating the generative models.

\begin{figure}[h]
\centering
\includegraphics[width=0.7\textwidth]{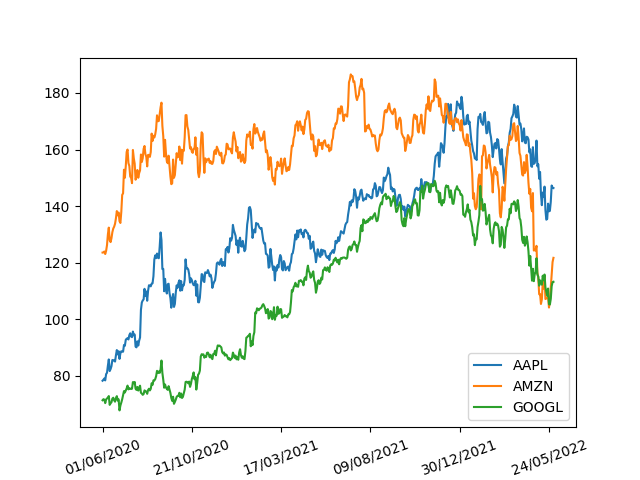}
\caption{\label{fig:plot1}Daily closing prices for AAPL, GOOGL, and AMZN between January 2020 and December 2021}
\end{figure}

\begin{table}[htbp]
\centering
\caption{Configuration for the illustrative dataset}
\label{tab:settings}
\begin{tabular}{l|l|l|l}
\hline
Start Date & End Date & Pricing Date & Symbols \\
\hline
01/06/2020 & 01/06/2022 & 01/06/2022 & AAPL, GOOGL, AMZN \\
\hline
\end{tabular}
\end{table}

\begin{table}[htbp]
\caption{\label{tab:empirical_stats}Empirical descriptive statistics of log-returns for AAPL, AMZN, and GOOGL. These values are computed from historical market data.}
\centering
\begin{tabular}{l|r|r|r}
\hline
Statistic / Stock & AAPL & AMZN & GOOGL \\
\hline
Mean & 0.001241 & -0.000030 & 0.000915 \\
Variance & 0.000402 & 0.000501 & 0.000326 \\
Skewness & -0.069260 & -0.437657 & -0.090484 \\
Kurtosis & 1.954379 & 6.702824 & 1.393783 \\
\hline
\end{tabular}
\end{table}

\subsection{Physics-informed time-series model mappings}\label{free-time-series-models-mappings}

Now, we introduce a \textit{physics-informed mapping} framework for time series, as defined in Section~\ref{physics-informed-systems}. This approach enables the construction of \emph{agnostic models}, in which the observed time series is related to a noise process--referred to as the \emph{innovation process}--via an invertible transformation. Formally, we consider mappings of the form
\be\label{F_map}
    F(X) = \varepsilon,
\ee
where we use the following notation. 
\begin{itemize}
    \item \( \varepsilon \in \mathbb{R}^{N_\varepsilon \times D \times T_\varepsilon} \) denotes the \emph{innovations} or \emph{process-induced noise}, extracted from the observed trajectories through the map \(F\). These innovations represent the residual stochastic component remaining after the dominant structural dynamics have been accounted for. The shape of \(\varepsilon\) need not match that of the original dataset.
    
    \item \( F: \mathbb{R}^{N_x \times D \times T_x} \to \mathbb{R}^{N_\varepsilon \times D \times T_\varepsilon} \) is a continuous and invertible operator. Its inverse enables reconstruction of synthetic trajectories from new innovation samples:
    \be\label{F_inv}
        X = F^{-1}(\varepsilon).
    \ee
\end{itemize}
The mapping \(F\) is said to be \emph{physics-informed} when it is designed to reflect known structural properties of the underlying process dynamics (e.g., conservation laws, symmetries, or stochastic differential equation structures), while still permitting a nonparametric estimation of the residual noise.

This framework can be viewed as a discrete-time analogue of physics-informed modeling approaches used for PDEs. In PDE settings, one typically models the main dynamics through differential operators, while accounting for modeling errors or uncertainties via additional source terms.

Similarly, in our time-series setting, the map \(F\) captures the main structural dynamics of the observed process, while the residual noise term \(\varepsilon\) represents unmodeled effects or stochastic perturbations. This separation mirrors the decomposition of physical systems into deterministic and uncertain components, but applied here in a fully data-driven and nonparametric fashion.

\paragraph{Time-series generation via encoder--decoder framework.}

This model supports a generative implementation via an encoder--decoder approach, leveraging the invertibility of a transformation \( F \) that maps an observed time series \( X \) to a process-induced noise signal \( \varepsilon \). This signal reflects the residual dynamics after removing structured behavior from the data (e.g., trends, seasonality, autoregression). A typical sampling pipeline is outlined in Algorithm~\ref{alg_ts}.

\begin{algorithm}
\label{alg_ts}
\begin{algorithmic}[1]
\REQUIRE Observed time series \( X \in \mathbb{R}^{N_x \times D \times T_x} \); base (latent) noise \( \eta \in \mathbb{R}^{N_\eta \times D_\eta \times T_\eta} \).
\ENSURE Synthetic time series \( \widetilde{X} \in \mathbb{R}^{N_x \times D \times T_x} \)
\STATE Compute process-induced noise via map: \( \varepsilon = F(X) \).
\STATE Generate new noise samples \( \widetilde{\varepsilon} \) by encoder--decoder framework (Section~\ref{Encoder-decoders-and-sampling-algorithms}) from \( \eta \):
\be \label{eq:kergen}
    G_k(\cdot) = K(\cdot, \eta)\theta, \quad \text{where } \theta = K(\eta,\eta)^{-1}(\varepsilon \circ \sigma).
\ee
\STATE Reconstruct synthetic time series via the inverse map: \( \widetilde{X} = F^{-1}(\widetilde{\varepsilon}) \).
\end{algorithmic}
\end{algorithm}

\noindent The generator \( G_k \) can be adapted to produce conditional samples \( \varepsilon \mid \omega \), where \( \omega \in \mathbb{R}^{D_\omega} \) is an exogenous covariate. This is handled via the conditional model introduced in Section~\ref{Conditional-distribution-sampling-model}.

\paragraph{Applications.}

This framework supports several financial modeling tasks. 
\begin{itemize}
  \item Benchmarking: Generate simulated paths on the same time grid as \( X \) to evaluate model fit and forecasting skill.
  \item Monte Carlo forecasting: Generate future paths \( \widetilde{X}^{*k} \) for \( t^{*k} > t^{T_x} \), enabling scenario analysis.
  \item Forward calibration: Generate \( \widetilde{X} \) such that future constraints are met, e.g.,
  \be
  \min_Y d(X, Y) \quad \text{subject to} \quad \mathbb{E}[P(Y^{\cdot, *k})] = c_p,
  \ee
  where \( P \) is a payoff function and \( d \) is a distance metric.

\item PDE pricers: we can use the HJB described in~\eqref{HJB} to get a martingale pricing of financial instruments for the vast majority of quantitative models. This part is already investigated in research papers\footnote{We refer the reader to~\cite{Mercier:2014,LeFloch-Mercier:2020a}} and is omitted here. 
\end{itemize}

\paragraph{Modeling perspective.}

The physics-informed formulation \eqref{MF} is broad enough to encompass traditional models--e.g., Brownian motion, geometric Brownian motion, ARMA models--as special cases where \( \varepsilon \) has a known distribution. By contrast, our nonparametric approach learns \( \varepsilon \) directly from data, avoiding parametric assumptions and enabling more flexible dynamics.

This formulation is particularly well-aligned with modern generative learning methods, such as kernel density estimators and deep generative models, which enable empirical calibration and generation of realistic synthetic trajectories.

\subsection{Brownian motion mappings}
\label{random-walks-and-brownian-motion-mappings}

\paragraph{Brownian motion.}

To illustrate the modeling framework defined by~\eqref{MF}, we begin with a canonical example: the \emph{random walk}. A discrete-time random walk satisfies the recursive relation
\be\label{RW}
X^{k+1} = X^k + \epsilon^k, \quad \text{where } \epsilon^k  =  \left(\epsilon_d^{n,k}\right)_{d,n},
\ee
with $(\epsilon_d^{n,k})$ i.i.d. centered random variables representing innovations or increments.
This fits the general framework by defining the forward map $F = \delta_0$ as the discrete difference operator:
\be
\epsilon^k = \delta_0(X)^k  =  X^{k+1} - X^k.
\ee

The inverse map $F^{-1}$ corresponds to discrete summation (integration):
\be
X^k = X^0 + \sum_{l=0}^{k-1} \epsilon^l, \quad k = 0, \ldots, T_x-1.
\ee
Component-wise, this reads as
\be
x_d^{n,k} = x_d^{n,0} + \sum_{l=0}^{k-1} \epsilon_d^{n,l}.
\ee
If the increments $\epsilon^k$ are identically distributed and centered, then by the \emph{Central Limit Theorem}, the scaled process converges in distribution:
\be
\frac{1}{\sqrt{k}} X^k \xrightarrow[k \to + \infty]{\mathcal{D}} \mathcal{N}(0, \Sigma),
\ee
where $\Sigma = \mathrm{Var}(\epsilon^k) \in \mathbb{R}^{D \times D}$.

\paragraph{Brownian motion increments.}

A similar structure arises for continuous-time \emph{Brownian motion} $W_t$, discretized on a time grid $t^0 < t^1 < \ldots < t^{T_x}$. We define the \emph{normalized increments}
\be
\delta_{\sqrt{t}}^k(W)  =  \frac{W_{t^{k+1}} - W_{t^k}}{\sqrt{t^{k+1} - t^k}}, 
\qquad \delta_{\sqrt{t}}(W)  =  \left(\delta_{\sqrt{t}}^k(W)\right)_{k=0}^{T_x-1}.
\ee
The inverse map reconstructs $W$ via a stochastic sum defined as
\be
W_{t^k} = W_{t^0} + \sum_{l=0}^{k-1} \sqrt{t^{l+1} - t^l} \cdot \epsilon^l, \quad \epsilon^l \sim \mathcal{N}(0, \Sigma).
\ee
We denote this weighted sum operator as $\sum_{(1/2)} \epsilon$, indicating square-root time scaling in the summation. This scheme corresponds to the \emph{Euler--Maruyama method}\footnote{see \cite{kloeden1992}} for simulating Brownian motion:
\be
W_{t^{k+1}} = W_{t^k} + \sqrt{t^{k+1} - t^k} \cdot \epsilon^k, \quad \epsilon^k \sim \mathcal{N}(0, \Sigma).
\ee

\paragraph{Log-normal models.}

Next, we consider \emph{log-normal models}, commonly used in finance applications. They are constructed by composing the logarithm with the difference operator:
\be\label{equa=log=dif}
\epsilon^k = (\delta_0 \circ \log)(X)^k  =  \log X^{k+1} - \log X^k.
\ee
This is understood component-wise across all entries of \(X\), that is, \( \log X^k  =  (\log x_d^{n,k})_{d,n} \), and similarly for \(\exp\) below.

The inverse map reconstructs the time series by exponentiating cumulative sums of increments:
\be
X^k = X^0 \cdot \exp\left(\sum_{l=0}^{k-1} \epsilon^l\right).
\ee
In the continuous-time setting, the scheme becomes
\be\label{PLR}
X_{t^{k+1}} = X_{t^k} \cdot \exp\left(\sqrt{t^{k+1} - t^k} \cdot \epsilon^k\right), \quad \epsilon^k \sim \mathcal{N}(0, \Sigma).
\ee
This can be inverted via
\be
X = (\delta_{\sqrt{t}} \circ \log)^{-1}(\epsilon) = \left(X^0 \cdot \exp\left(\sum_{l=0}^{k-1} \sqrt{t^{l+1} - t^l} \cdot \epsilon^l\right)\right)_k.
\ee
The equation~\eqref{PLR} thus provides a closed-form integral operator of the form:
\be
X = X^0 \cdot \left(\exp \circ \sum_{(1/2)}\right)(\epsilon),
\ee
which fits directly into the general modeling framework described in~\eqref{MF}.

\paragraph{Application to financial data.}

We now apply this scheme to the financial data introduced in Figure~\ref{fig:plot1}. From the observed time series $X$, we compute log-return innovations $\epsilon$ using~\eqref{equa=log=dif}. These empirical innovations are shown in Figure~\ref{fig:plot5} (left panel) for selected assets.

Using the encoder--decoder framework, we sample new innovations $\widetilde{\epsilon}$ via the generator $G_k$ using latent noise $\eta$, and reconstruct synthetic log-returns. The right panel in Figure~\ref{fig:plot5} illustrates samples generated through this approach.

\begin{figure}[h]
\centering
\includegraphics[width=0.9\textwidth]{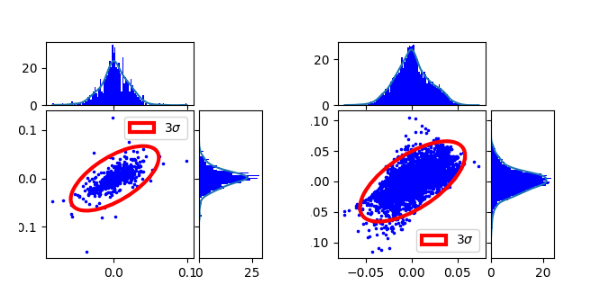}
\caption{\label{fig:plot5} Extracted historical innovations (left) vs. generated innovations via kernel generator~ \eqref{eq:kergen} (right)}
\end{figure}

\begin{figure}
\centering
\includegraphics[width=1.\textwidth, keepaspectratio]{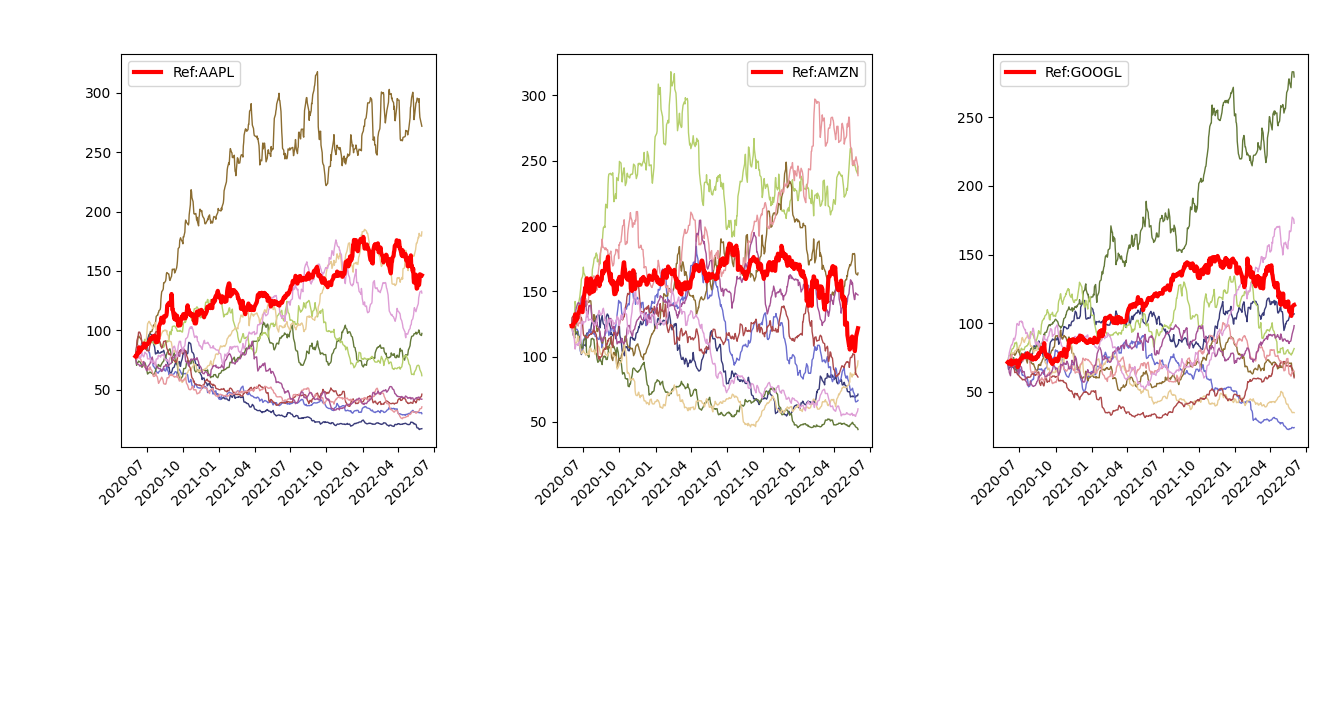}
\caption{\label{fig:plotBM} Ten synthetic paths generated from a log-normal model}
\end{figure}

To assess the fidelity of the generative process, Table~\ref{tab:102} shows descriptive statistics---including the empirical moments and Kolmogorov--Smirnov distances---computed on the log-return \emph{innovations} extracted from historical data and their synthetic counterparts generated via the kernel-based model. This comparison evaluates how well the generative model captures the distributional properties of the noise components underlying asset returns.

\begin{table}[htbp]
\caption{\label{tab:102}Descriptive statistics of log-return innovations for Amazon, Apple, and Google. 
Each cell displays a statistic based on historical data (the corresponding value from generated innovations is shown in parentheses).}
\centering
\begin{tabular}{l|l|l|l}
\hline
Statistic / Stock & AMZN & AAPL & GOOGL \\
\hline
Mean & 0.0012 (0.00054) & -3e-05 (0.00032) & 0.00091 (0.00067) \\
Variance & -0.066 (-0.09) & -0.44 (0.029) & -0.09 (-0.19) \\
Skewness & 0.0004 (0.00034) & 0.0005 (0.00041) & 0.00033 (0.00026) \\
Kurtosis & 2 (0.57) & 6.7 (2) & 1.4 (0.62) \\
KS statistic & 0.48 (0.05) & 0.93 (0.05) & 0.31 (0.05) \\
\hline
\end{tabular}
\end{table}

\subsection{Autoregressive and moving average mappings}\label{auto-regressive-moving-averages-maps}

Autoregressive moving average (ARMA) models\footnote{see \cite{shumway2010}} are standard tools for modeling univariate time series with linear dynamics. An \( \mathrm{ARMA}(p,q) \) process is defined by
\be\label{GARMA}
X^k = \mu + \sum_{i=1}^{p} a_i X^{k-i} + \sum_{j=1}^{q} b_j \epsilon^{k-j},
\ee
where \( \mu \in \mathbb{R} \) is a mean parameter, \( \{a_i\} \) and \( \{b_j\} \) are model coefficients, and \( \{\epsilon^k\} \) is a sequence of i.i.d. white noise variables with zero mean and finite variance \( \sigma^2 \). Classical estimation methods include least squares and maximum likelihood\footnote{see  \cite{BrockwellDavis:2006}}.

In the context of the physics-informed framework \eqref{MF}, this model corresponds to a causal transformation \( F(X) = \epsilon \), whereby the latent noise sequence is inferred directly from the data. When the coefficients \( \{a_i\}, \{b_j\} \) are fixed, the ARMA structure provides a closed-form expression for this mapping.

To recover \( \epsilon^k \) from \( X^k \), we use the infinite-order moving average representation:
\be
\epsilon^k = \mu + \sum_{j=0}^{+\infty} \pi_j X^{k-j},
\ee
where the coefficients \( \{\pi_j\} \) solve the difference equation:
\be
\pi_j + \sum_{\ell=1}^{q} b_\ell \pi_{j-\ell} = -a_j,
\ee
with the conventions \( a_0 = -1 \), \( a_j = 0 \) for \( j > p \), and \( b_j = 0 \) for \( j > q \). This relation may be expressed using the backshift operator \( B \) and polynomials \( \phi(B) \), \( \theta(B) \) as
\be
\phi(B) X^k = \theta(B) \epsilon^k,
\quad \text{where} \quad
\phi(B) = 1 - \sum_{i=1}^p a_i B^i,
\quad \theta(B) = 1 + \sum_{j=1}^q b_j B^j.
\ee

For our simulation, we focus on the autoregressive model \( \mathrm{AR}(p) \), i.e., the special case \( \mathrm{ARMA}(p,0) \). The forward and inverse maps are
\be
F(X^k) = \phi(B) X^k = \epsilon^k, \quad \text{and} \quad X^k = \phi^{-1}(B) \epsilon^k.
\ee
Figure~\ref{fig:plot360} displays synthetic trajectories generated from such a model.

\begin{figure}
\centering
\includegraphics[width=1.\textwidth, keepaspectratio]{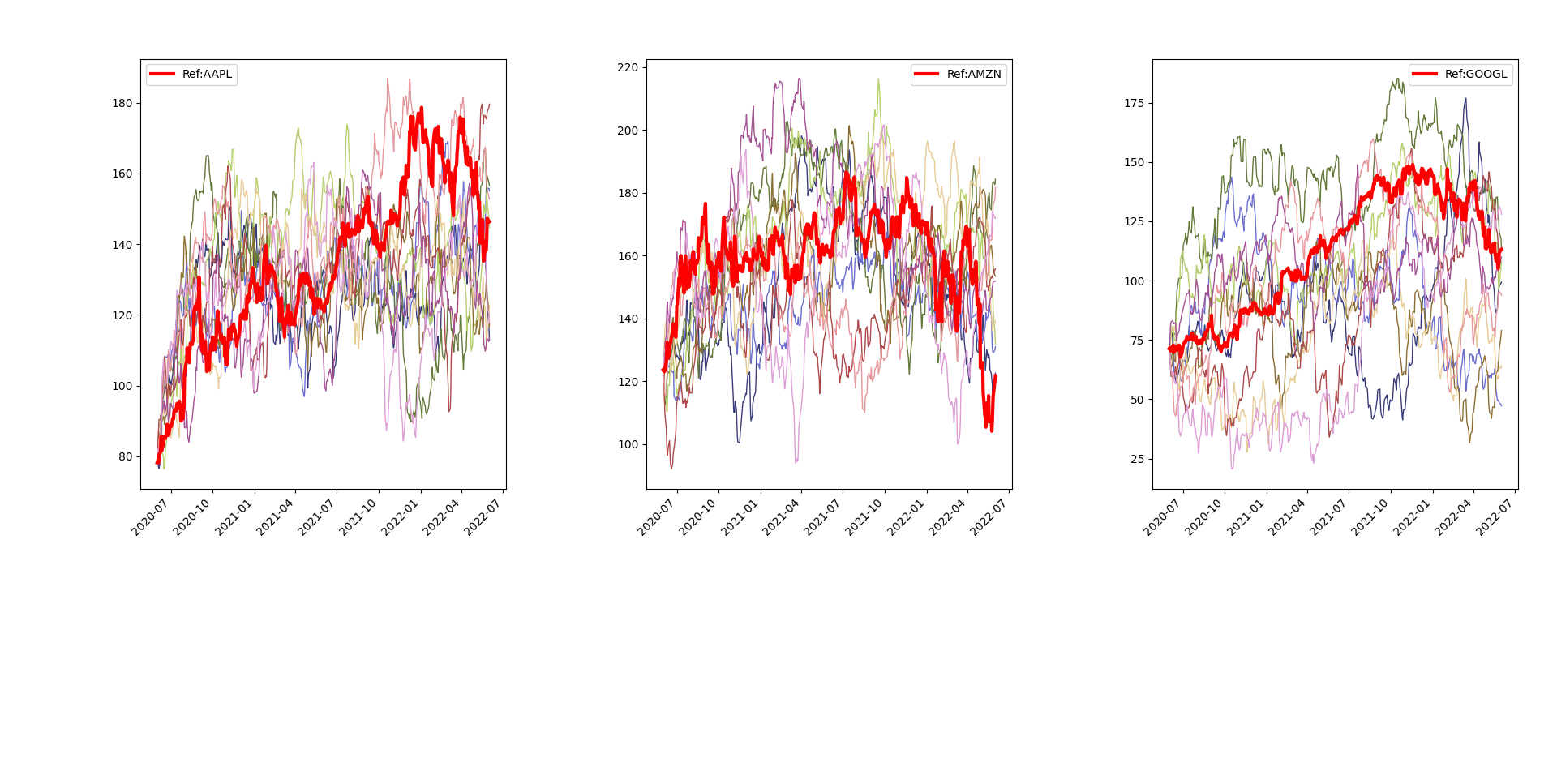}
\caption{\label{fig:plot360} Ten synthetic paths generated from the $\text{AR}(p)$ model}
\end{figure}

\begin{table}[htbp]
\caption{\label{tab:ar_stats}Descriptive statistics of log-returns innovations for Amazon, Apple, and Google under the $\text{AR}(p)$ model. Each cell shows a statistic computed from historical data, and the corresponding value for the AR-map generated data (in parentheses).}
\centering
\begin{tabular}{l|l|l|l}
\hline
Statistic / Stock & AMZN & AAPL & GOOGL \\
\hline
Mean & 0.13 (0.32) & -0.0014 (0.05) & 0.079 (0.26) \\
Variance & 0.012 (-0.11) & -0.3 (-0.18) & 0.095 (-0.27) \\
Skewness & 7.1 (6.2) & 11 (8.4) & 4.1 (3.5) \\
Kurtosis & 1.4 (0.69) & 5.1 (1.8) & 2 (0.79) \\
KS statistic & 0.04 (0.05) & 0.48 (0.05) & 0.0035 (0.05) \\
\hline
\end{tabular}
\end{table}

\subsection{GARCH mappings}\label{garchpq-maps}

Generalized autoregressive conditional heteroskedasticity (GARCH) models\footnote{see \cite{engle1982} and\cite{bollerslev1986}}
 capture time-varying volatility, common in financial time series. The $\text{GARCH}(p,q)$ process takes the form:
\be
\begin{aligned}
X^k &= \mu + \sigma^k \epsilon^k, \\
(\sigma^k)^2 &= \alpha_0 + \sum_{i=1}^{p} \alpha_i (X^{k-i})^2 + \sum_{j=1}^{q} \beta_j (\sigma^{k-j})^2,
\end{aligned}
\ee
where \( \mu \in \mathbb{R} \), \( \{\epsilon^k\} \) is a white noise sequence with unit variance, and \( \sigma^k \) is a stochastic volatility term determined recursively.

To express this in operator form, we set 
\be
\alpha(B) = \sum_{i=1}^{p} \alpha_i B^i, \quad
\beta(B) = \sum_{j=1}^{q} \beta_j B^j.
\ee
Then, we have 
\be
(1 - \beta(B)) (\sigma^k)^2 = \alpha_0 + \alpha(B)(X^k)^2.
\ee
We introduce also the composition \( \pi(B) = [(1 - \beta(B))]^{-1} \alpha(B) \), yielding
\be
(\sigma^k)^2 = \pi(B)(X^k)^2, \quad \text{and} \quad \epsilon^k = \frac{X^k - \mu}{\sigma^k}.
\ee
This defines the forward map
\be
F(X^k) = \frac{X^k - \mu}{\sqrt{\pi(B)(X^k)^2}}.
\ee
This transformation allows the recovery of the noise process from data and can be inverted when \( \pi(B) \) is invertible. The map \( F \) thus defines a GARCH model in the physics-informed model framework.
Figure~\ref{fig:plot362} illustrates synthetic trajectories generated from a $\text{GARCH}(1,1)$ model.

\begin{figure}
\centering
\includegraphics[width=1.\textwidth, keepaspectratio]{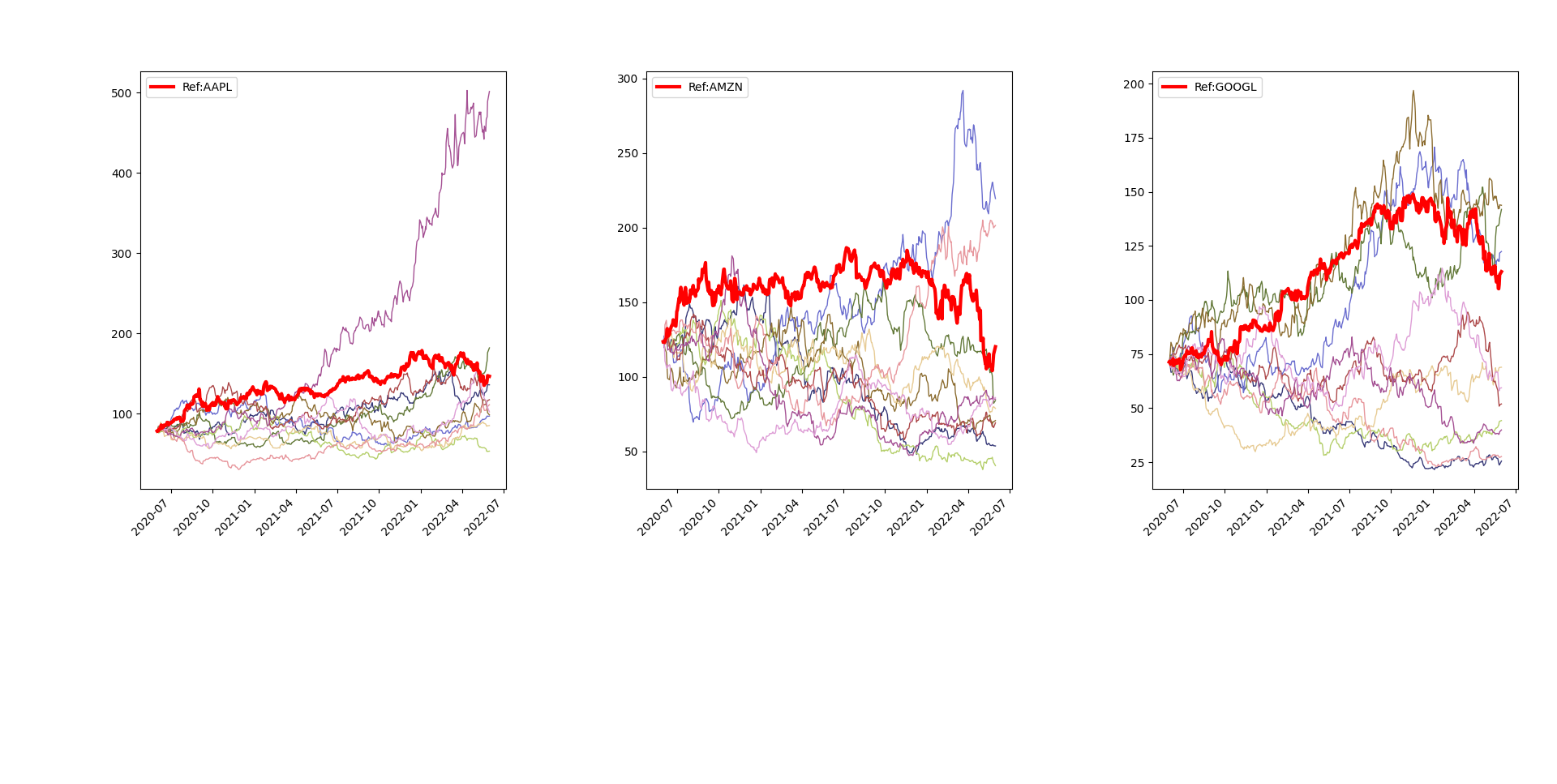}
\caption{\label{fig:plot362} Ten synthetic paths generated from a $\text{GARCH}(1,1)$ model}
\end{figure}

\begin{table}[htbp]
\caption{\label{tab:garch_stats}Descriptive statistics of log-return innovations for Amazon, Apple, and Google under the GARCH model. Each cell shows a statistic from historical data and the corresponding value from the GARCH-map generated data (in parentheses).}
\centering
\begin{tabular}{l|l|l|l}
\hline
Statistic / Stock & AMZN & AAPL & GOOGL \\
\hline
Mean & -6.7e-06 (-0.00062) & -1.2e-06 (-0.00068) & -8.4e-06 (-0.00075) \\
Variance & -0.068 (0.057) & -0.44 (-0.43) & -0.089 (-0.14) \\
Skewness & 0.0004 (0.00038) & 0.0005 (0.00044) & 0.00033 (0.00032) \\
Kurtosis & 2 (0.24) & 6.7 (2.1) & 1.4 (0.34) \\
KS statistic & 0.21 (0.05) & 0.67 (0.05) & 0.23 (0.05) \\
\hline
\end{tabular}
\end{table}

\subsection{Additive noise map}\label{additive-noise-map}

We now introduce a conditional noise transformation of the form \( \eta = \eta_Y(\varepsilon) \), where the process-induced noise \( \varepsilon \) is adjusted based on an exogenous variable \( Y \). This setup corresponds to an \emph{additive noise model}:
\be\label{MFP}
  \eta_Y(\varepsilon) = \varepsilon - f(Y), \qquad \varepsilon = \eta_Y^{-1}(\eta) = \eta + f(Y),
\ee
where our notation is follows. 
\begin{itemize}
  \item \( \eta \in \mathbb{R}^{D_\varepsilon} \) denotes a whitened residual, ideally independent of \( Y \),
  \item \( f : \mathbb{R}^{D_Y} \to \mathbb{R}^{D_\varepsilon} \) is a smooth function modeling the dependence of \( \varepsilon \) on \( Y \). If unknown, \( f \) can be estimated from historical data using the denoising algorithm in~\eqref{DNabla}.
\end{itemize}

This model captures conditional structure in the noise by isolating a predictable component \( f(Y) \). Such conditioning appears naturally in stochastic models like the Vasicek process\footnote{see \cite{vasicek1977}}:
\be
  dr_t = F(r_t)\,dt + \sigma\,dW_t,
\ee
where the drift term \( F(r_t) \) is a function of the state \( r_t \). In discretized form, this becomes:
\be
  \ln X^{*,k+1} = \ln X^{*,k} + f\big(\ln X^{*,k}\big) + \varepsilon^k,
\ee
with \( f \) playing the role of the drift, and \( \varepsilon^k \) representing a process-induced noise. This corresponds to the model:
\be
F  =  \eta_Y \circ \delta_0 \circ L^{2p} \circ \log,
\quad \text{with} \quad Y  =  X^* \circ \log(X).
\ee

Figure~\ref{fig:plot35} shows resampled trajectories generated by this model. The function \( f \) was fitted using the denoising procedure~\eqref{DNabla}, with a regularization parameter \( \lambda = 10^{-3} \), from the historical data pairs \( (X^{*,k}, \varepsilon^{*,k}) \).

\begin{figure}[h]
\centering
\includegraphics[width=1.\textwidth]{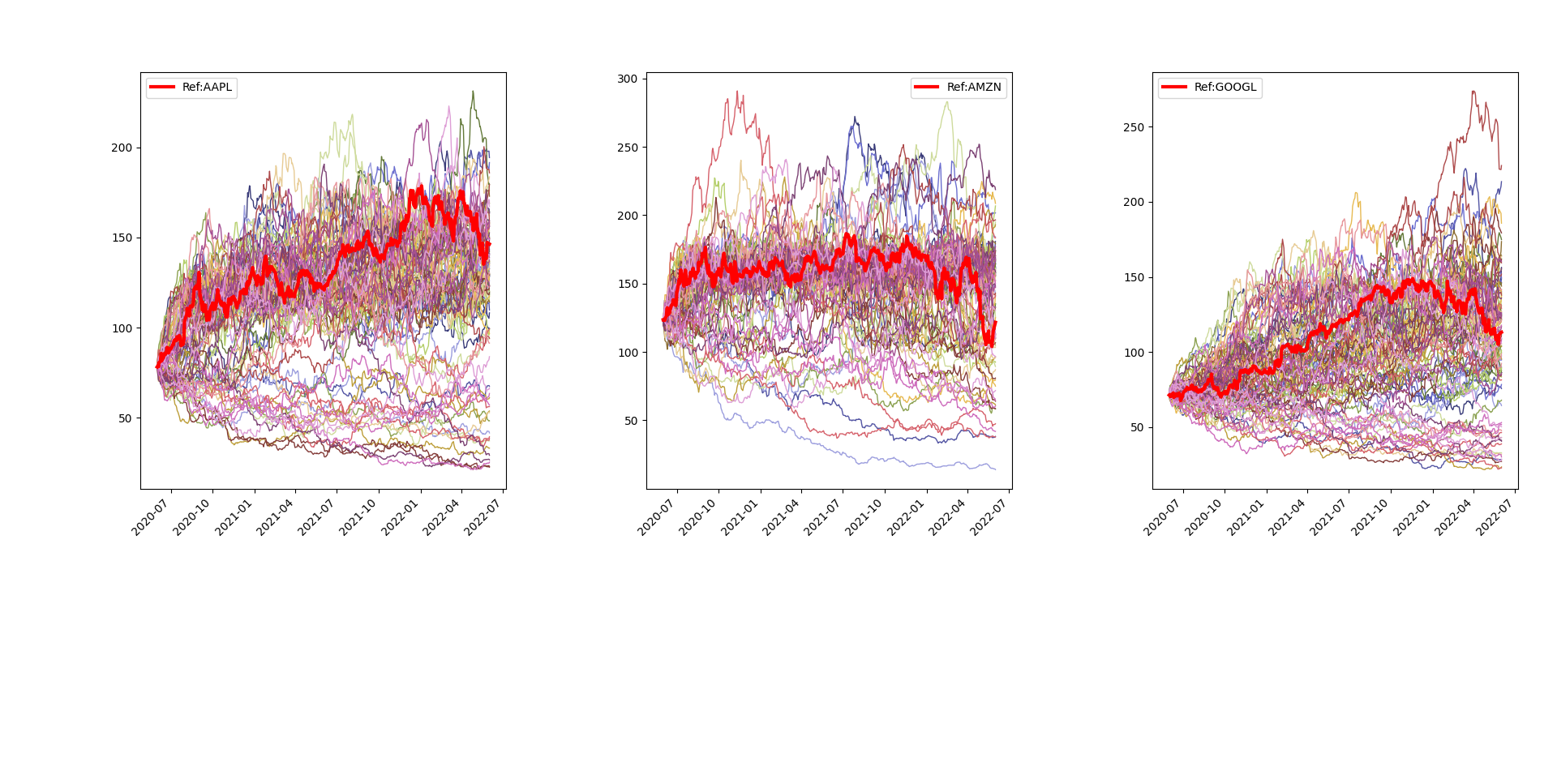}
\caption{\label{fig:plot35} Ten examples of resampled trajectories using the additive noise model}
\end{figure}

\begin{table}[htbp]
\caption{\label{tab:additive_stats}Descriptive statistics of log-return innovations for Amazon, Apple, and Google under the additive noise model. Each cell shows the value computed from historical data and the corresponding value from the additive-map generated data (in parentheses).}
\centering
\begin{tabular}{l|l|l|l}
\hline
Statistic / Stock & AMZN & AAPL & GOOGL \\
\hline
Mean & \(-2.4\times 10^{-6}\) (\(-0.0015\)) & \(-4.8\times 10^{-6}\) (\(-0.00055\)) & \(1.3\times 10^{-6}\) (\(-0.00056\)) \\
Variance & 0.0003 (0.00029) & 0.00034 (0.00035) & 0.00025 (0.00026) \\
Skewness & 0.19 (\(-0.076\)) & \(-0.36\) (\(-0.08\)) & \(-0.1\) (\(-0.049\)) \\
Kurtosis & 2.1 (0.39) & 3.4 (0.43) & 1.3 (0.31) \\
KS statistic & 0.087 (0.05) & 0.45 (0.05) & 0.48 (0.05) \\
\hline
\end{tabular}
\end{table}

\subsection{Conditioned map and data augmentation}\label{conditioned-map-and-data-augmentation}

We now consider conditional generative models, in which the extracted process-induced noise \( \varepsilon \) is conditioned on exogenous or endogenous variables. A basic example is conditioning on the process itself. Specifically, one may model the conditional distribution given current state:
$\varepsilon^k \sim \mathbb{P}(\varepsilon^k \mid X^k)$, where \( \varepsilon = F(X) \) is the noise induced by a physics-informed model \( F \). Numerically, this distribution can be approximated using a conditional generator \( G_k(\cdot \mid X^k) \), as introduced in Section~\ref{Conditional-distribution-sampling-model}.

By composing this generator with a finite-difference transformation, the model defines a discrete-time stochastic scheme of the form:
\be\label{MFCScheme}
  \ln X^{k+1} = \ln X^k + \varepsilon^k \mid \ln X^k,
\ee
which reflects a stochastic process with state-dependent noise. Figure~\ref{fig:plot46} shows resampled trajectories based on this scheme.

\begin{figure}[h]
\centering
\includegraphics[width=1.\textwidth]{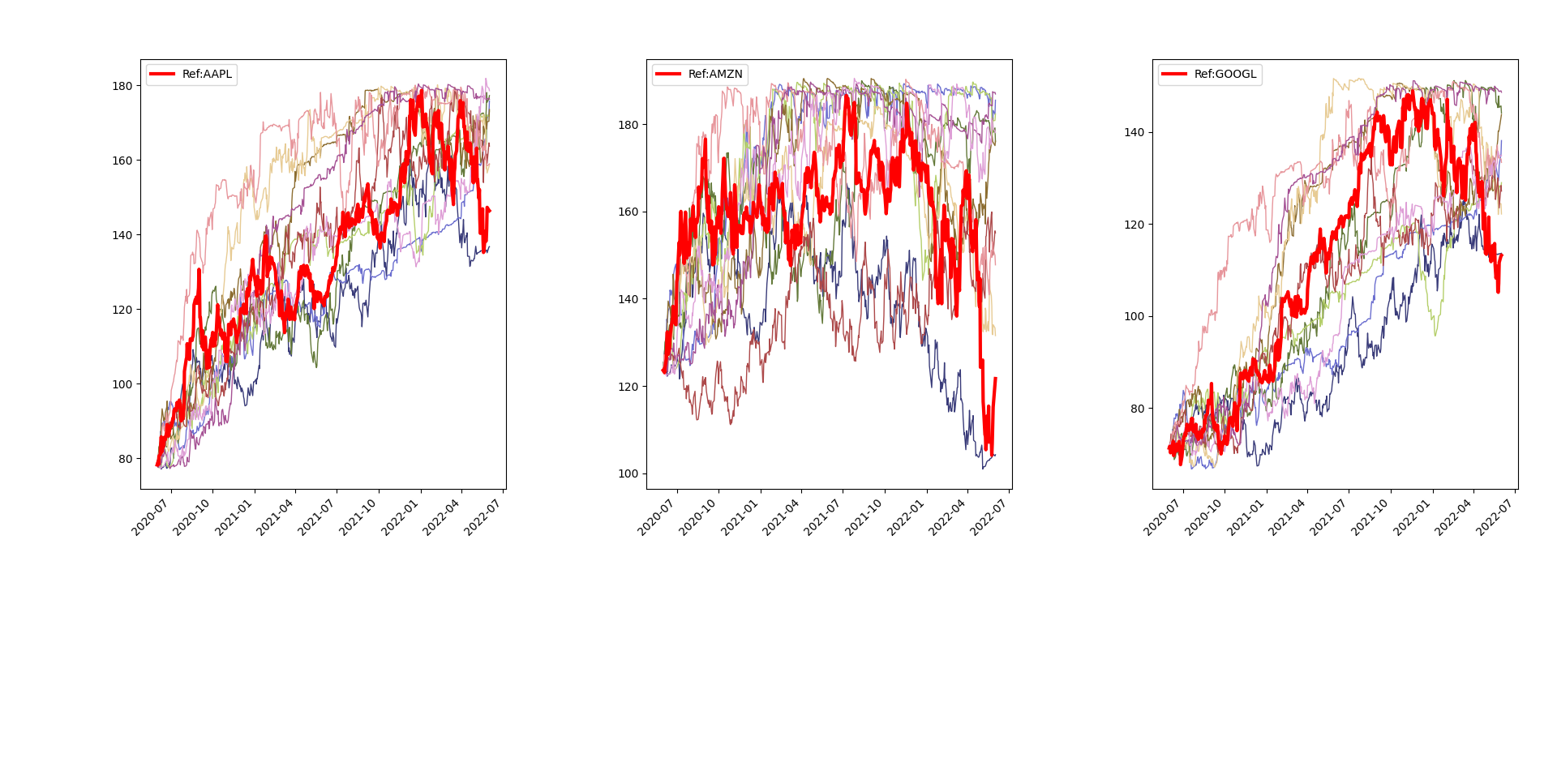}
\caption{\label{fig:plot46} Ten examples of resampled paths using the conditional noise model \eqref{MFCScheme}}
\end{figure}

This model structure is well-suited to approximate weakly stationary processes such as the Cox--Ingersoll--Ross (CIR) model\footnote{see \cite{cox1985}}. More generally, the framework allows for \emph{data augmentation} by introducing additional conditioning features. For instance, one may augment the process by estimating a local volatility proxy:
\be
  \sigma^k(X)  =  \operatorname{Tr} \left( \operatorname{Cov}\left(X^{k-q}, \dots, X^{k+q}\right) \right),
\ee
where \( q \) controls the local window size and \(\operatorname{Tr}(\cdot)\) denotes the trace of the sample covariance matrix.

The resulting model is a stochastic volatility scheme of the form:
\be\label{SVM}
\begin{split}
\ln X^{k+1} &= \ln X^k + \varepsilon_x^k \mid \sigma^k, \\
\sigma^{k+1} &= \sigma^k + \varepsilon_\sigma^k \mid \sigma^k,
\end{split}
\ee
where \( \varepsilon^k = (\varepsilon_x^k, \varepsilon_\sigma^k) \) are the noise components for the process and volatility, respectively. These are sampled using conditional generators \( G_k^X(\cdot \mid \sigma^k) \) and \( G_k^\sigma(\cdot \mid \sigma^k) \).

Figure~\ref{fig:plot56} displays ten generated paths from this stochastic volatility model.

\begin{figure}[h]
\centering
\includegraphics[width=1.\textwidth]{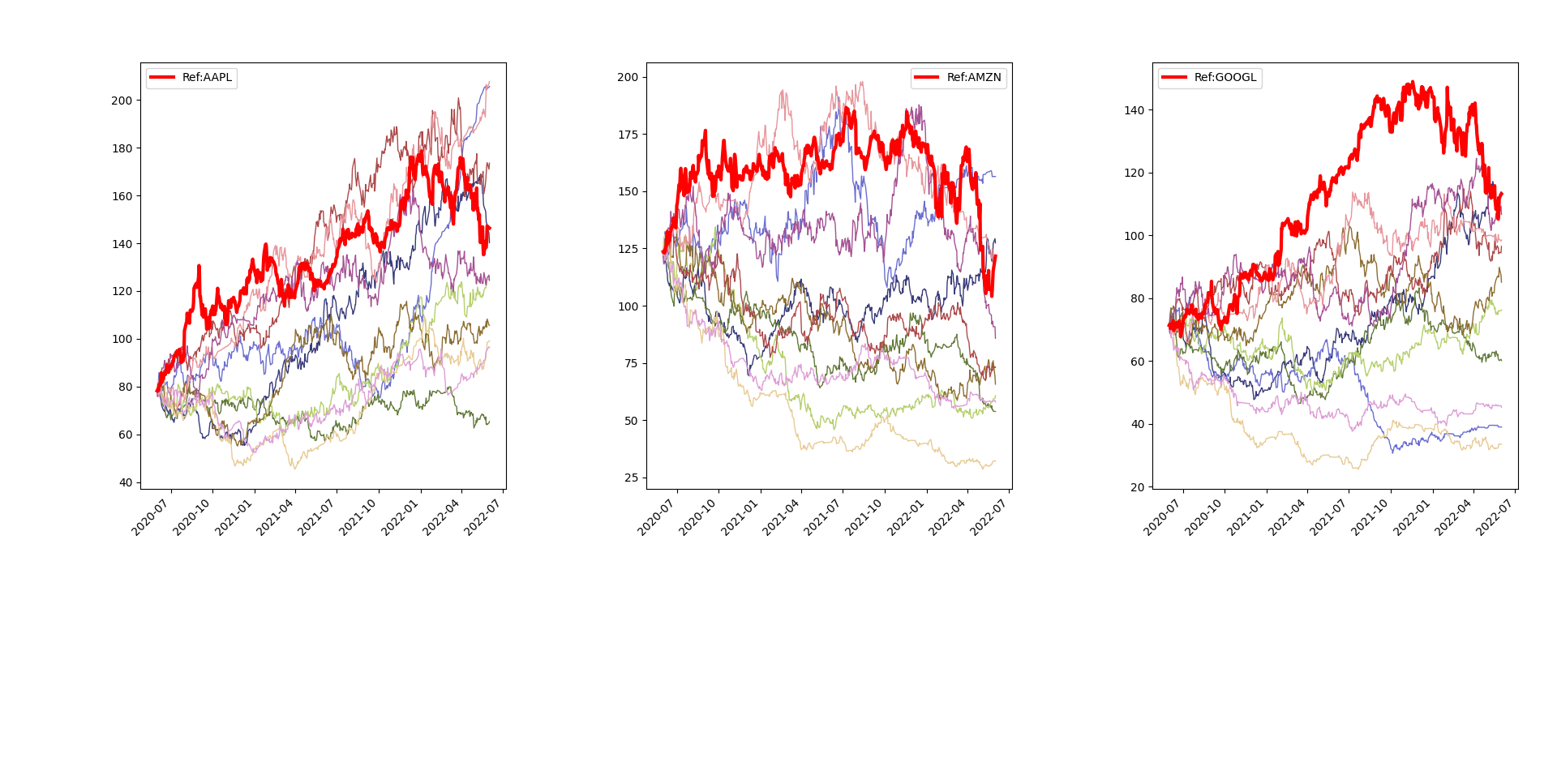}
\caption{\label{fig:plot56} Ten synthetic paths generated using the stochastic volatility model \eqref{SVM}}
\end{figure}

\begin{table}[htbp]
\caption{\label{tab:cond_map_stats}Descriptive statistics of log-return innovations for Amazon, Apple, and Google under the conditional map model. Each cell shows the historical statistic followed by the statistic from the generated data in parentheses.}
\centering
\begin{tabular}{l|l|l|l}
\hline
Statistic / Stock & AMZN & AAPL & GOOGL \\
\hline
Mean & -0.0011 (0.0085) & -0.0028 (-0.025) & 6.8e-06 (0.027) \\
Variance & 0.5 (0.5) & 0.5 (0.46) & 0.5 (0.49) \\
Skewness & 0.033 (0.055) & -0.042 (-0.13) & 0.014 (0.003) \\
Kurtosis & 0.059 (-0.33) & 0.016 (-0.3) & -0.026 (-0.25) \\
KS statistic & 0.99 (0.05) & 0.99 (0.05) & 0.74 (0.05) \\
\hline
\end{tabular}
\end{table}


\section{Benchmarking with synthetic trajectories: a Heston case study}
\label{benchmark-methodology}

\paragraph{Methodology} 

This section introduces a structured methodology for evaluating the generative time-series models described previously. The objective is to assess their ability to reproduce both the statistical properties and functional behavior of known stochastic processes. For this purpose, we use the Heston stochastic volatility model\footnote{see \cite{heston1993}} as a reference. This model features latent variance dynamics and admits closed-form solutions for European option pricing, making it suitable for benchmarking.

Our methodology consists of the following steps. 
\begin{itemize}
  \item Setting. Select a reference stochastic model--in this case, the Heston model--and generate a synthetic path from it using known parameters. This path serves as the observed dataset for calibration.

  \item Calibration. Calibrate both (i) the parameters of the Heston model from the synthetic path, and (ii) a generative model using the same trajectory via the physics-informed modeling framework described in Section~\ref{free-time-series-modeling}.

  \item Reproduction. Verify that the generative model can accurately reproduce the original trajectory from which it was derived. This ensures that the mapping \( F \) in the representation \( \varepsilon = F(X) \) is numerically invertible.

  \item Noise Distribution. Extract the process-induced noise \( \varepsilon \) from both the reference model and the generative model, and compare their statistical properties using descriptive and non-parametric (e.g., Kolmogorov--Smirnov) tests.

  \item Trajectory Generation. Simulate multiple paths from both models and visually assess their similarity. The goal is to confirm that the generative model replicates both first- and higher-order dynamics of the reference process.

  \item Pricing. Use Monte Carlo simulation to price a European call option under each model and compare the resulting estimates to the known analytical price under the Heston model.
\end{itemize}

\paragraph{Benchmark process: Heston model}\label{benchmarks-framework---heston}

The Heston model describes a  \( X_t \) with stochastic volatility \( \nu_t \), governed by the following system of stochastic differential equations (SDEs): 
\begin{equation}
\begin{cases}
dX_t = \mu X_t\,dt + X_t \sqrt{\nu_t}\, dW_t^{(1)}, \\be6pt]
d\nu_t = \kappa(\theta - \nu_t)\,dt + \sigma \sqrt{\nu_t}\, dW_t^{(2)}, \\be6pt]
\langle dW_t^{(1)}, dW_t^{(2)} \rangle = \rho\,dt.
\end{cases}
\end{equation}
We choose Heston parameters satisfying the Feller condition \( 2\kappa\theta > \sigma^2 \), generate one trajectory, and treat it as the observed historical data. From this trajectory, we re-estimate the drift parameter \( \mu \approx \frac{\ln(X_T) - \ln(X_0)}{T} \), and use both the original and generated trajectories for subsequent analysis (see Figure~\ref{fig:comparetrajectoriesHeston}).

\paragraph{Reproducibility test}

To confirm the numerical invertibility of the generative map \( F \), we test whether the model can reproduce the original input trajectory from its associated process-induced noise \( \varepsilon = F(X) \). The result is shown in Figure~\ref{fig:testreproductibilityHeston}.

\begin{figure}
\centering
\includegraphics[width=1.0\textwidth, keepaspectratio]{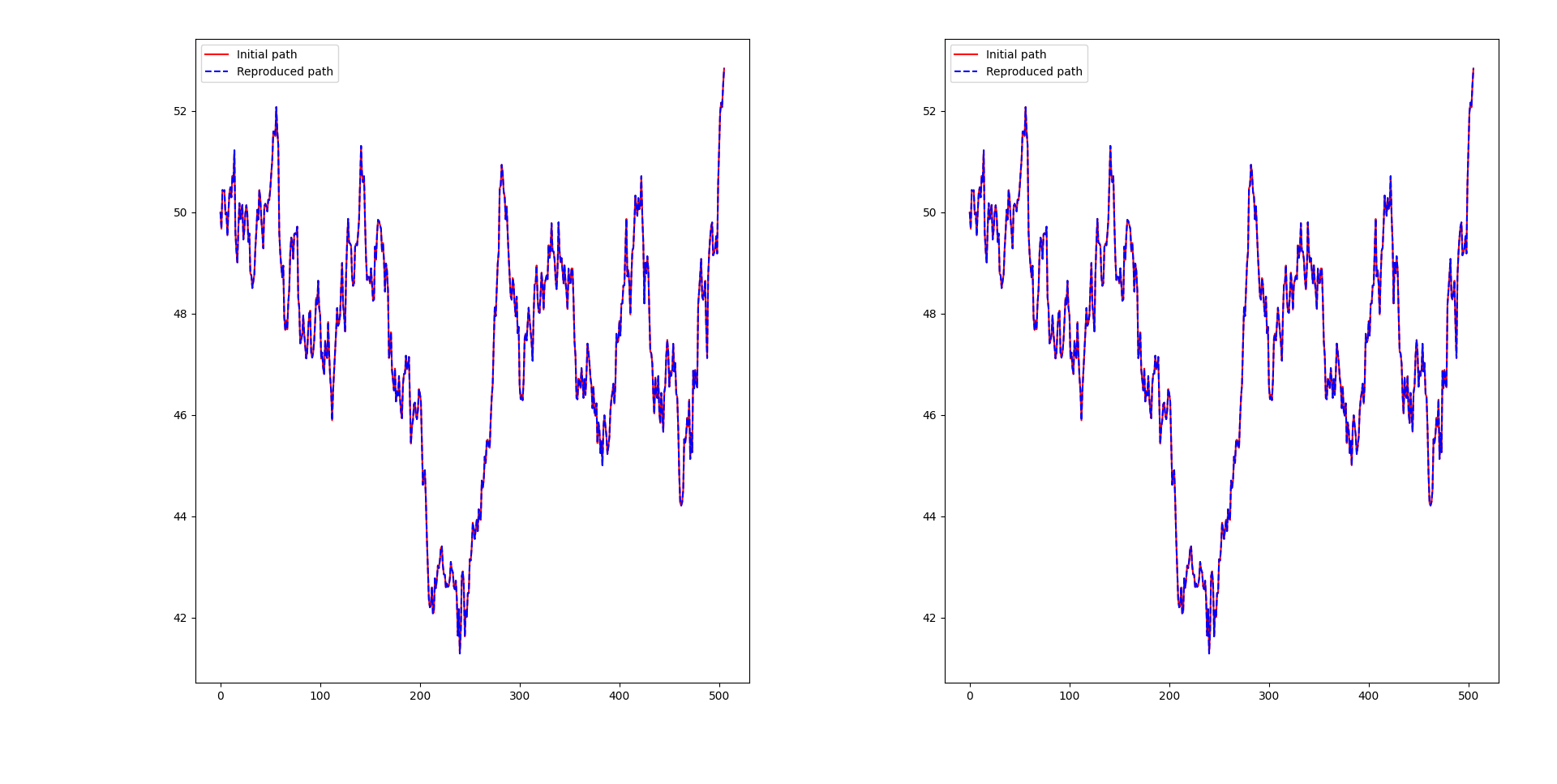}
\caption{\label{fig:unnamed-chunk-113}\label{fig:testreproductibilityHeston} Reproducibility test for a Heston process}
\end{figure}

\paragraph{Noise distribution analysis}\label{benchmarks-distributions}

We compare the distributions of the extracted process-induced noise \( \varepsilon \) between the calibrated Heston model and two generative alternatives: (i) the log-difference model and (ii) a conditional generative model. Visual comparisons are presented in Figure~\ref{fig:comparedistributionsHeston}, with numerical summaries provided in Table~\ref{tab:1903}.

\begin{figure}[h]
\centering
\includegraphics[width=1.0\textwidth]{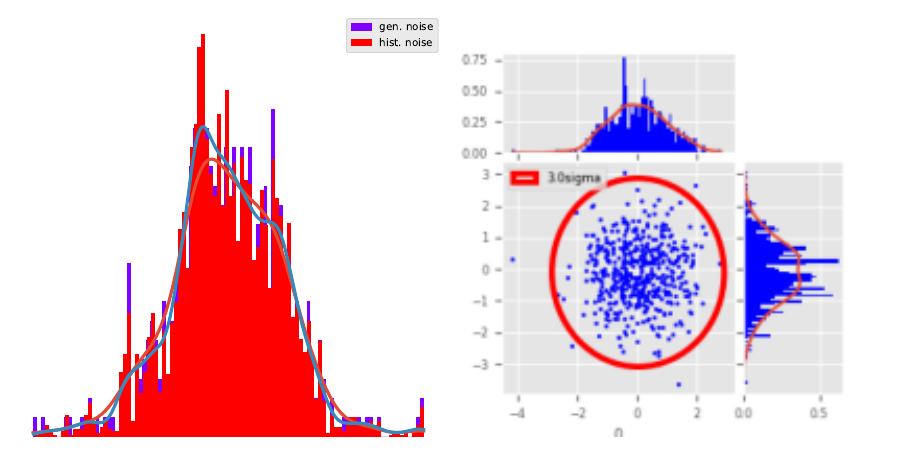}
\caption{\label{fig:comparedistributionsHeston} Empirical distributions of the extracted noise \( \varepsilon \) for the Heston model and two generative variants}
\end{figure}

\begin{table}[htbp]
\caption{\label{tab:1903}Statistical table -- Generative Stats (Calibrated ones)}
\centering
\begin{tabular}[t]{l|l|l|l}
\hline
 & HestonDiffLog lat.:0 & HestonCondMap lat.:0 & HestonCondMap lat.:1 \\
\hline
Mean & 0.00011(0.00013) & -0.0024(0.018) & -0.0032(-0.025) \\
\hline
Variance & -0.045(-0.044) & 0.00053(0.18) & 0.0073(0.15) \\
\hline
Skewness & 9.8e-05(9.8e-05) & 0.9(0.68) & 1(0.97) \\
\hline
Kurtosis & 0.8(0.72) & -0.42(0.056) & -0.49(-0.046) \\
\hline
KS test & 1(0.05) & 0.0088(0.05) & 0.0063(0.05) \\
\hline
\end{tabular}
\end{table}

\paragraph{Trajectory comparison}\label{benchmarks-trajectories}

We simulate 1000 paths from both the Heston model (left) and the generative models (right), using identical random seeds where applicable. The original trajectory is shown in red in both panels of Figure~\ref{fig:comparetrajectoriesHeston}.

\begin{figure}[h]
\centering
\includegraphics[width=1.\textwidth]{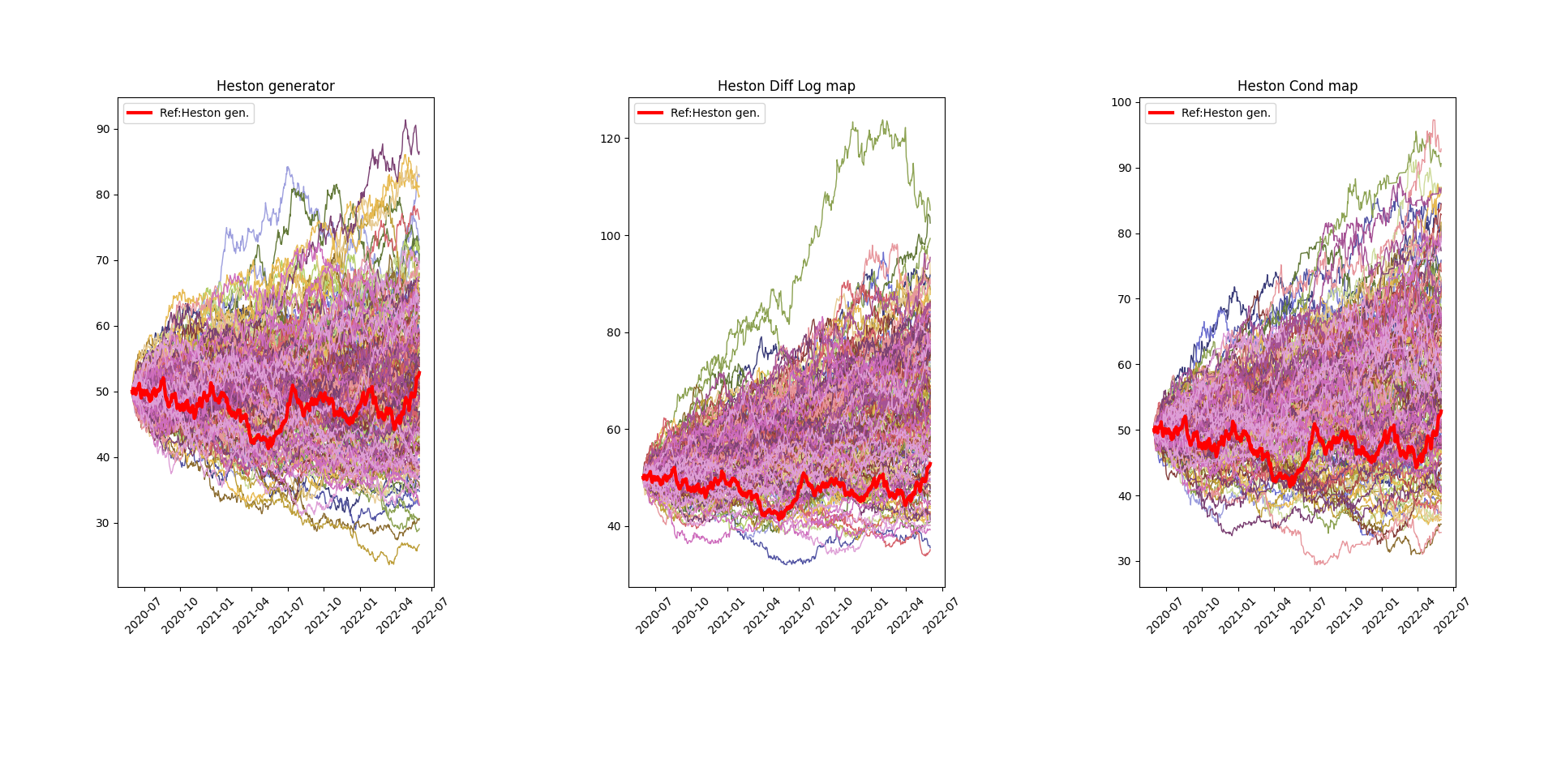}
\caption{\label{fig:comparetrajectoriesHeston} Comparison of paths generated by the Heston model (left) and a generative model calibrated on a single Heston trajectory (right)}
\end{figure}

\paragraph{Option prices benchmark}\label{benchmarks-prices}

We price a European call option with strike \( K = X_T \) and maturity \( T \), using Monte Carlo simulation under each model. The benchmark price from the closed-form Heston formula is also presented when available\footnote{Both Monte Carlo and closed-form prices are computed using the QuantLib library. See \url{https://www.quantlib.org} for more details.}. Results are shown in Table~\ref{tab:1933}.

\begin{table}[htbp]
\caption{\label{tab:1933}European call option prices under the Heston model: Monte Carlo, closed-form solution, and generative models}
\centering
\begin{tabular}[t]{l|r|r|r|r}
\hline
  & MC :PricesDiffLog & Gen :PricesDiffLog & closed pricer & Gen :PricesCondMap\\
\hline
Mean & 7.141976 & 8.552551 & 7.222894 & 7.988374\\
\hline
Var & 79.254114 & 101.864179 & NaN & 58.212532\\
\hline
Lbound & 6.590205 & 7.927006 & NaN & 7.515488\\
\hline
Ubound & 7.693747 & 9.178096 & NaN & 8.461260\\
\hline
\end{tabular}
\end{table}



\section{Extrapolation of pricing functions with generative methods}
\label{pricing-with-generative-methods}

\paragraph{Problem setting.} Let \((X_t)_{t \ge 0} \in \mathbb{R}^D\) be a stochastic process modeling underlying market variables (e.g., asset prices), and let \(V: \mathbb{R}^D \to \mathbb{R}\) denote a payoff function (e.g., option contract, portfolio values). The time-\(s\) value of this instrument under risk-neutral measure \(\mathbb{Q}\) is given by
\be
\overline{V}(s, T, x)  =  \mathbb{E}^{\mathbb{Q}}_{X_s = x} \left[ V(X_T) \right].
\ee
This quantity defines the pricing function (or value surface) for maturity \(T\) and initial condition \(X_s = x\). In practice, we often need to compute \(\overline{V}(s,T,x)\) across many values of \(x\). From an operational perspective, these computations are usually done using overnight batches, as pricing functions are usually too computationally involved to be used on a regular, intraday basis.

\paragraph{Objective.} We aim to approximate the pricing function \((x \mapsto \overline{V}(s,T,x))\) and its derivatives (Greeks) using fast, extrapolation methods of a \textit{given} pricing functions $\overline{V}$. This is especially useful to compute intradays, \textit{real time}. 
\begin{itemize}
    \item Risk measures (e.g., \texttt{Delta}, \texttt{Gamma}, \texttt{Theta})
    \item Scenario-based pricing
    \item PnL attribution over time
\end{itemize}
To achieve this, we train surrogate models using either historical or synthetically generated paths via the physics-informed generative framework introduced in Section~\ref{free-time-series-models-mappings}. We consider two extrapolation methods: 
\begin{itemize}
    \item Taylor approximations of order two. 
    \item Kernel Ridge regression~\eqref{FIT}.
\end{itemize}
We chose Taylor approximations to benchmark against, because this method is used in an industrial context to approximate real time intradays values.

\paragraph{Illustration: basket option.} 

Let \( \omega \in \mathbb{R}^D \) be fixed weights and \(K > 0\) a strike. The payoff of a basket call option is:
\be
V(x)  =  \max(\omega^\top x - K, 0).
\ee
Assume that the asset vector \(X_T\) follows a multivariate log-normal distribution at option maturity $T$, i.e., the log-returns of the underlying basket values are normally distributed and jointly stationary. This naive assumption led to consider the Black-Scholes formula for basket options under log-normality to compute the reference price \( \overline{V}(s,T,x) \) in closed form. Figure ~\ref{fig:plot3} visualizes the pricing surface. Although naive, this analytical formula is only used to provide a tractable ground truth for assessing the accuracy of our surrogate models based on kernel ridge regression ({KRR}) or Taylor approximations.

We approximate \( \overline{V}(s,T,x) \) across a range of values \(x\), using these two extrapolation methods and compare the outputs to the Black-Scholes benchmark. 

\begin{figure}
\centering
\includegraphics[width=0.7\textwidth, keepaspectratio]{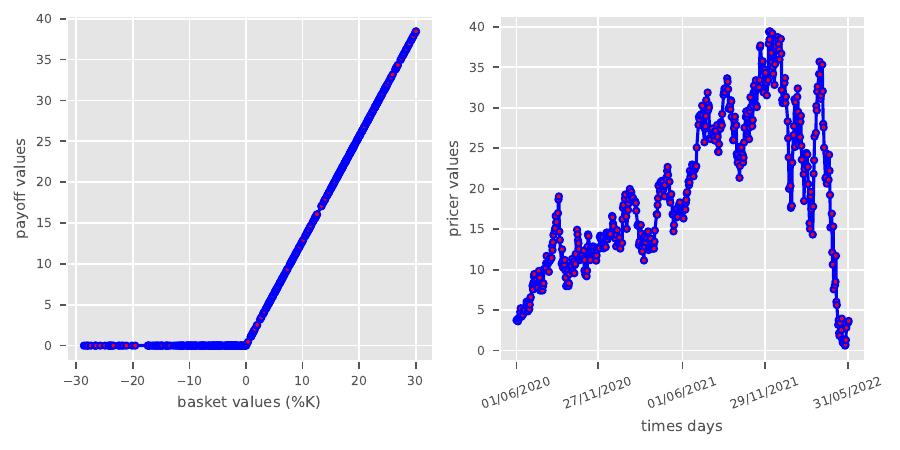}
\caption{\label{fig:unnamed-chunk-118}\label{fig:plot3} Left: payoff $V$ of the basket option. Right: pricing function ($\overline{V}_t = \mathbb{E}^{X_T}(V)$) as a function of time}
\end{figure}

\paragraph{Synthetic scenario generation and PnL attribution}\label{predictive-methods-for-financial-applications.}

The kernel regression~\eqref{FIT} is used to estimate the pricing function \( \overline{V}(t,\cdot) \) at intraday market points \( z \) as \( \overline{V}_{k,\theta,X}(t,z) \). These test points are generated and unseen during the training ($z \notin X$). We now discuss the choice of the training set $X$.

The interpolation quality of \( \overline{V}_{k, \theta,X} \)(z), trained on a dataset \( X \), is governed by the error bound in \eqref{err}, which depends on the maximum mean discrepancy (MMD) \( d_k(z, X) \). To visualize this, Figure~\ref{fig:plot10} shows isocontours of the MMD function \( z \mapsto d_k(z,X) \) for two different training sets $X$ (blue dots), overlaid with test points (red dots):
\begin{itemize}
\item Left and right: The test set (red points) in both panels consists of synthetically generated intraday data with a horizon \(H = 5\). This horizon choice amounts to simulate larger intraday movements of the underlying than expected.
\item Left: Blue corresponds to the  historical basket values, as a function of time. 
\item Right: 
Blue corresponds to synthetically generated basket values, as a function of time, for three time points \(t^0-1\), \(t^0\), \(t^0+1\), with a horizon \(H = 10\) days, corresponding to VaR computations. We added the two extra observation time (\(t^0-1\), \(t^0+1\) to approximate temporal derivatives of the pricing function -- in particular, the \texttt{Theta}: \( \partial_t V \).
\end{itemize}

\begin{figure}
\centering
\includegraphics[width=0.7\textwidth, keepaspectratio]{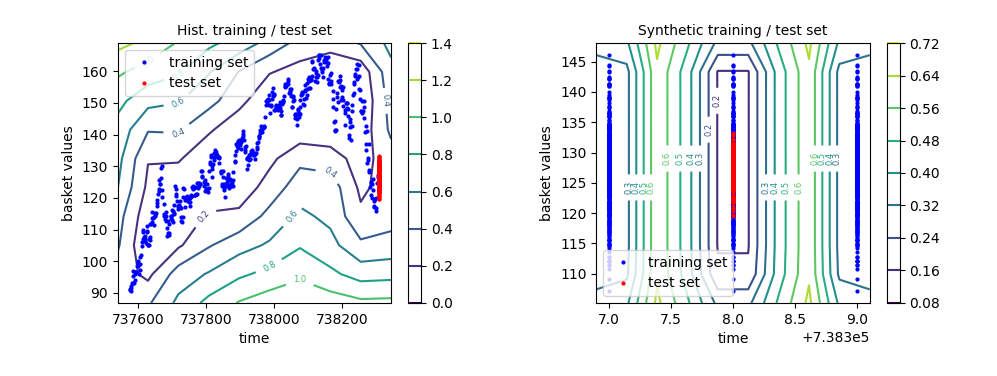}
\caption{\label{fig:unnamed-chunk-119}\label{fig:plot10} Isocontours of interpolation error \( d_k(z, X) \) for two choices of training set $X$ in blue (left observed, historical datas, right VaR data). Red: synthetic generated intradays basket values.}
\end{figure}

This visualization illustrates that synthetic VaR training points (right panel) provide a better coverage of the test domain, leading to improved generalization and more accurate extrapolation. This test thus uses the synthetic VaR data to achieve high accuracy. However, if VaR data are not available, one can extrapolate from historical data, which comes at a price in terms of accuracy.

Figure~\ref{fig:plot11} shows predicted prices at the test points \( z \), compared to ground truth values obtained from the pricing function for both methods.

\begin{figure}
\centering
\includegraphics[width=1.\textwidth, keepaspectratio]{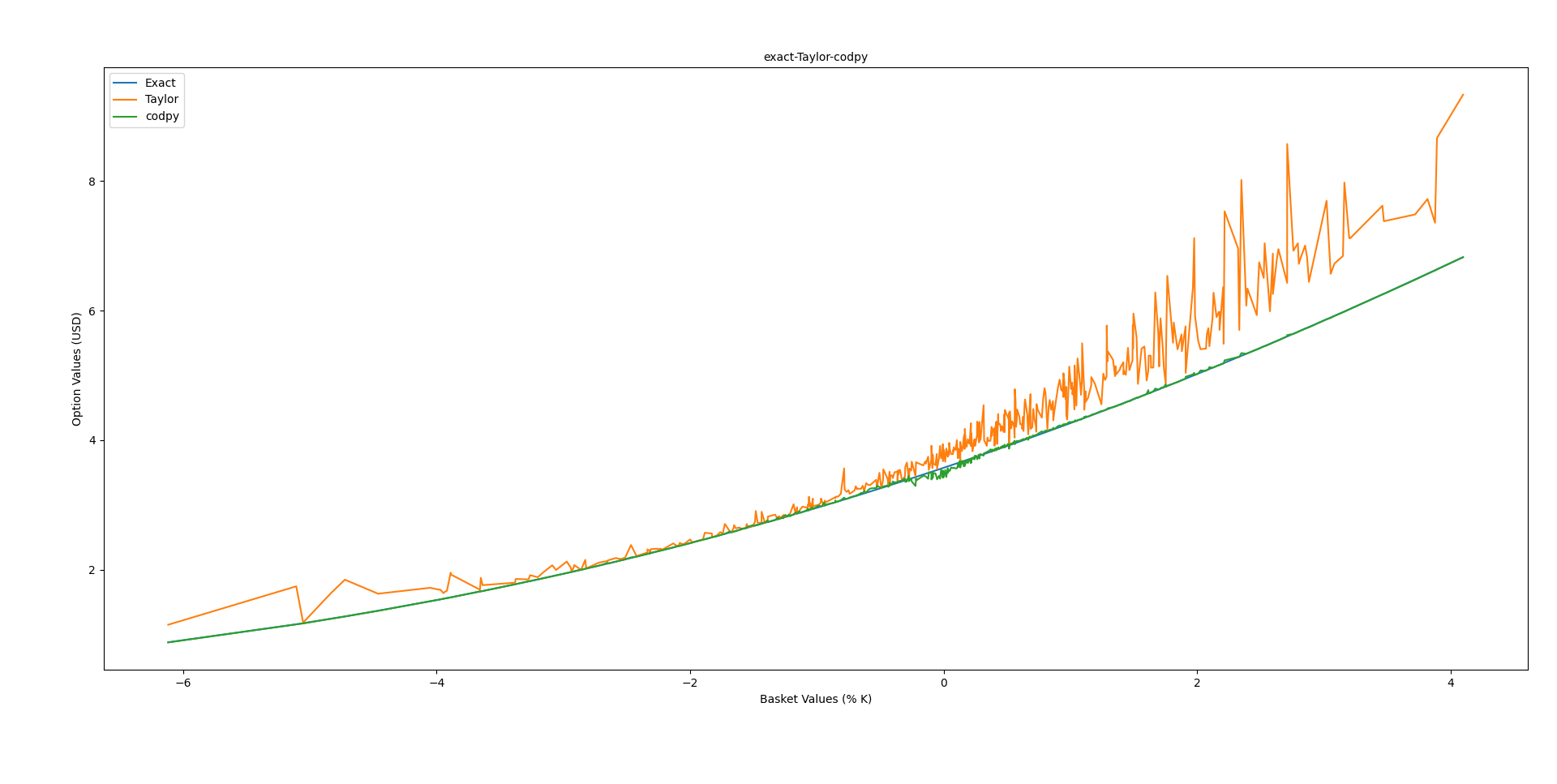}
\caption{\label{fig:unnamed-chunk-120}\label{fig:plot11} Extrapolation of prices at synthetic generated intradays basket values. Green are references Black and Scholes values. Blue (close to green) are kernel ridge regression value $\overline{V}_{k,\theta}$. Orange are extrapolations given by a full second-order Taylor approximation.}
\end{figure}

\paragraph{Greeks.}

To evaluate sensitivities of the pricing function, we compute its partial derivatives--commonly referred to as \emph{Greeks}. Let \( \overline{V}_{k,\theta}(t, \cdot)\) be the output of the {KRR} model approximating the pricing surface.
We denote its gradient with respect to both time and spatial inputs as
\be
\nabla \overline{V}_{k,\theta}(t, \cdot) = \left( \frac{\partial \overline{V}_{k,\theta}}{\partial t}, \frac{\partial \overline{V}_{k,\theta}}{\partial x_1}, \ldots, \frac{\partial \overline{V}_{k,\theta}}{\partial x_D} \right)(t, \cdot),
\ee
evaluated on the test set \( Z \), i.e., \( \nabla \overline{V}_{k,\theta}(t, Z) \in \mathbb{R}^{N_Z \times (D+1)} \).

Both methods (kernel ridge regression and Taylor approximation) provide approximation of greeks, and can be challenged to the reference gradient, which can be computed explicitly with our chosen pricing function. Results are shown in Figure~\ref{fig:plot13}, where each subplot corresponds to one component of the gradient--e.g., \texttt{Theta} (\( \partial_t \)), \texttt{Delta} (\( \partial_{x_1} \)), etc.

Since the training data are sampled i.i.d., raw gradients may exhibit spurious oscillations. We apply the denoising procedure \eqref{DNabla} to stabilize these estimates.

\begin{figure}
\centering
\includegraphics[width=1.\textwidth, keepaspectratio]{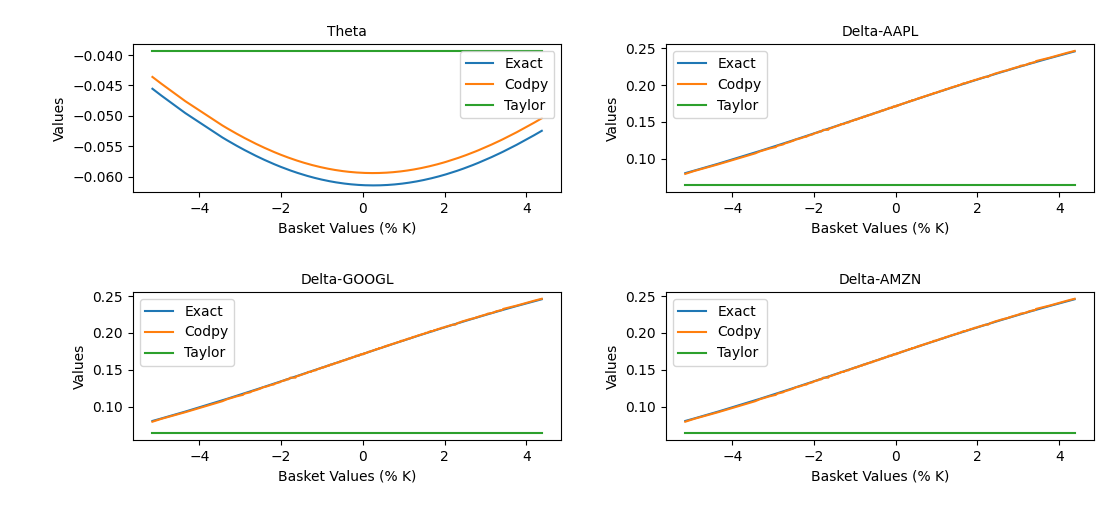}
\caption{\label{fig:unnamed-chunk-121}\label{fig:plot13} Greeks computed using kernel-based derivative estimation and a second order Taylor approximation}
\end{figure}


\section{Application to stress tests and reverse stress tests}
\label{Stress-test-and-reverse-stress-tests}

\paragraph{Description.}
Regulatory authorities impose stress tests and reverse stress tests to analyze extreme tail-risk scenarios. These tools are central to assessing the resilience of investment strategies under severe market conditions.

We focus on \emph{reverse stress testing} (RST), which inverts the traditional question. While standard stress testing asks ``what is the loss given a market scenario?'', RST asks ``which market scenarios are most probable (under a specified reference model) among those that produce a given loss?''. Formally, we start from a target outcome--e.g., portfolio value(s) or PnL, possibly vector-valued--and recover scenarios that map to that outcome. This helps identify portfolio vulnerabilities and informs risk mitigation.

\paragraph{Problem setting.}
We adopt the same setting as in Section~\ref{pricing-with-generative-methods}. Let $(X_t)_{t\ge 0}$ denote the underlying process (see Figure~\ref{fig:plot1}). Consider a basket option with payoff $V(\cdot)$ and pricing function $\overline V(t,\cdot)$ (both plotted in Figure~\ref{fig:unnamed-chunk-118}). Fix a stress-test computation date $T$ and a horizon $H$ (in days). Define the PnL mapping
\be
\mathrm{PnL}_T(x)  =  \overline V(T{+}H, x) - \overline V(T, x_T)\in\RR^{D_p}, 
\ee
where $x\in\RR^{D}$ is a market scenario at $T{+}H$ and $x_T$ is the realized state at time $T$.
Given simulated scenarios $X=(x^1,\ldots,x^{N})$, we obtain PnL samples
\be
P  =  (p^1,\ldots,p^{N})\in(\RR^{D_p})^{N}, \qquad p^n  =  \mathrm{PnL}_T(x^n).
\ee
Our objective is to approximate an \emph{inverse} mapping $X(p)\approx \mathrm{PnL}_T^{-1}(p)$ that returns plausible scenarios $x$ producing a target PnL $p$. Since $\mathrm{PnL}_T:\RR^{D}\to\RR^{D_p}$ can be non-invertible, we rely on the stable inversion framework of Section~\ref{Inversion-of-non-invertible-mappings.}

\paragraph{Methodology.}

Using the free model of Section~\ref{free-time-series-modeling} (here instantiated as the log-normal model; see Section~\ref{random-walks-and-brownian-motion-mappings}, 
 we proceed as follows. 
\begin{enumerate}
\item Generate forward samples $x^n\sim X_{T+H}\mid x_T$ for $n=1,\ldots,N$, with a typical choice $H=10$ days for the VaR horizon.
\item Evaluate PnLs $p^n=\mathrm{PnL}_T(x^n)$ and collect $P=(p^1,\ldots,p^N)$.
\item Fit a kernel ridge regressor (KRR) $f_{k,\theta}:\RR^{D}\to\RR^{D_p}$ to approximate $\mathrm{PnL}_T$, using kernel $k$ and hyperparameters $\theta$: 
      $f_{k,\theta}(x)\approx \mathrm{PnL}_T(x)$.
\item Construct the \emph{stable inverse} $g_{k,\theta}:\RR^{D_p}\to\RR^{D}$ from Section~\ref{Inversion-of-non-invertible-mappings.}, yielding $g_{k,\theta}(p)$ as a plausible scenario that maps to $p$ under $f_{k,\theta}$.
\end{enumerate}

To evaluate robustness on a broad range of PnLs, we draw a \emph{synthetic} target set $P_s=(p_s^1,\ldots,p_s^{N_s})$ using the sampling algorithm of Section~\ref{Encoder-decoders-and-sampling-algorithms}. We then produce reverse-stress scenarios
\be
X_s=(x_s^1,\ldots,x_s^{N_s}), \qquad x_s^n  =  g_{k,\theta}(p_s^n).
\ee
The accuracy is assessed by the reconstruction error
\be
\epsilon^n  =  p_s^n - \mathrm{PnL}_T(x_s^n),
\ee
and by its relative basis-point (bps) version,
\be
\mathrm{bps\_err}^n  =  10^4\!\left(1 - \frac{p_s^n}{\mathrm{PnL}_T(x_s^n)}\right).
\ee

\paragraph{Results.}

Figure~\ref{fig:ReversePrices} (left) shows the reverse-stress scenarios $X_s$ as functions of the sorted synthetic targets $p_s^1<\ldots<p_s^{N_s}$. The recovered scenarios vary smoothly across the PnL spectrum, indicating stability of the inverse. Figure~\ref{fig:ReversePrices} (right) shows $(\mathrm{bps\_err}^1,\ldots,\mathrm{bps\_err}^{N_s})$, demonstrating good out-of-sample accuracy on unseen PnL values.

\begin{figure}[t]
\centering
\includegraphics[width=1\linewidth]{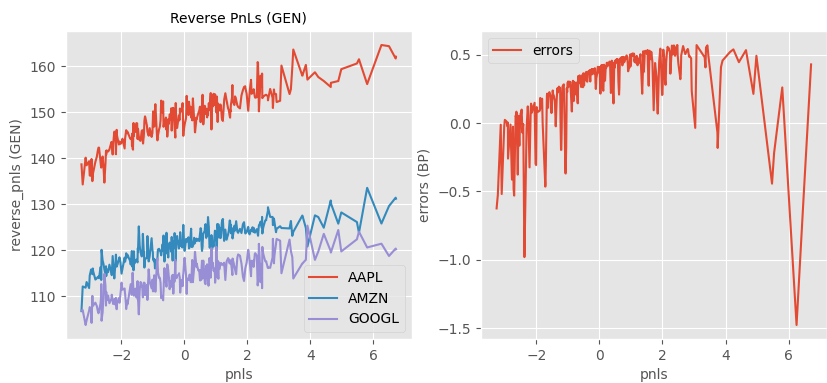}
\caption{\label{fig:ReversePrices}
Reverse-stress scenarios (left) and reconstruction error in bps (right).}
\end{figure}


\section{Application to portfolio management}
\label{extending-quantitative-models-through-conditioning}

\paragraph{Data and set-up.}

We now illustrate that the generative algorithms developed above can be used safely for conditional analysis. We do so on a numerical test in portfolio management with trading signals.

First of all, we consider daily closing prices for a basket of $D=106$ cryptocurrencies from 19/07/2021 to 28/04/2023 ($T_x=649$ trading days). Let $X^n\in\mathbb{R}^D$ denote the price vector at date $t^n$, and $X=(X^0,\ldots,X^{T_x-1})\in\mathbb{R}^{D\times T_x}$. Figure~\ref{fig:cryptos} plots the normalized series $X^n/X^0$.

Following the notation of~\eqref{TS}, we work with simple (componentwise) returns $r^{m}=X^{m}/X^{m-1}-\mathbf{1}$. For a sliding window of length $W$, anchored at index $n$, we collect the $W$-vector of one-step returns
\be
\epsilon^{(n)}  =  \big(r^{n+1},\ldots,r^{n+W}\big)\in\mathbb{R}^{D\times W}.
\ee
The associated inverse map reconstructs a hypothetical price path from a starting level,
\be
F_n^{-1}(\epsilon) = \Big(X^{n}\,,\, X^{n}\!\!\prod_{k=1}^{1}\!\!(\mathbf{1}+\epsilon^{(n)}_k)\,,\, \ldots,\, X^{n}\!\!\prod_{k=1}^{K}\!\!(\mathbf{1}+\epsilon^{(n)}_k)\Big)_{K\ge 0}.
\ee

\paragraph{Objective and portfolios.}

For each anchor $n=0,\ldots,T_x-W$, we form portfolio weights $\omega_n\in\mathbb{R}^{D}$ at rebalancing time $t^{n+W}$. The dollar value is
\be
P^n = \langle \omega_n,\, X^{n+W}\rangle  =  \sum_{d=1}^D \omega_{n,d}\,X^{n+W}_d.
\ee
We use a dollar-neutral long/short specification with leverage bounds,
$ \mathbf{1}^\top \omega_n=0,\quad \|\omega_n\|_\infty\le 1$.

\paragraph{Efficient portfolio (Markowitz).}
A baseline strategy solves a mean--variance problem using expected returns and a covariance estimate over the window:
\begin{equation}
\label{EP}
\overline{\omega}  =  \arg\min_{\omega}
\;\frac{1}{2}\,\omega^\top Q\,\omega \;-\; \lambda\,\omega^\top \overline{\epsilon}
\;+\; \beta\,\|\omega-\omega^0\|_{1}
\quad\text{s.t.}\quad
\mathbf{1}^\top \omega=0,\;\|\omega\|_\infty\le 1,
\end{equation}
where $Q=\mathrm{cov}(\epsilon)\in\mathbb{R}^{D\times D}$ is estimated from the window, $\overline{\epsilon}\in\mathbb{R}^D$ is the corresponding mean return, $\lambda>0$ is a risk--aversion parameter, $\omega^0$ are the current holdings, and the $\ell_1$ penalty $\beta\|\omega-\omega^0\|_1$ models proportional transaction costs.

\paragraph{Methodology (conditioning).}
Our strategies are driven by a random vector of returns $\epsilon$. We compare two otherwise similar settings in which $\epsilon$ is conditioned on different observed variables $\eta$. At each time $t^n$ we form the conditional distribution $\epsilon \,|\, \eta=\eta^n$ using the techniques in Section~\ref{Conditional-distribution-sampling-model}, and compute $\overline{\epsilon}$ and $Q$ from conditional samples\footnote{The same conditional recipe applies to any other quantitative model of the assets.}
 to plug into~\eqref{EP}.

\paragraph{Benchmarks and strategies.}

We present four strategies (see Figure~\ref{fig:plot363}). 
\begin{enumerate}
\item Index (equal weight). An equally weighted portfolio used as a reference.
\item LS (unconditional). Long/short strategy solving~\eqref{EP} using unconditional estimates $(\overline{\epsilon},Q)$.
\item LS$|$CAPM (conditional). Same as (2), but $\epsilon$ is conditioned on a CAPM target
$b^n = r_f + \beta^n\big(R_m^n - r_f\big)$,
where $r_f$ is the per-period risk-free rate, $R_m^n$ is the index return, and $\beta^n$ are per-asset regression coefficients on the index.
\item LS$|$Indicators (conditional). Same as (2), but $\epsilon$ is conditioned on a vector $b^n$ comprising liquidity measures, rolling averages over multiple windows, and their differences (momentum-type signals).
\end{enumerate}
\noindent We evaluate the cumulative performance of strategy $j$ as
\be
P_j^{n}
=
\prod_{k\le n}\frac{\langle \omega_j^{k},\, X^{k+1}\rangle}{\langle \omega_j^{k},\, X^{k}\rangle}\,.
\ee

\begin{figure}[t]
\centering
\includegraphics[width=0.5\textwidth,keepaspectratio]{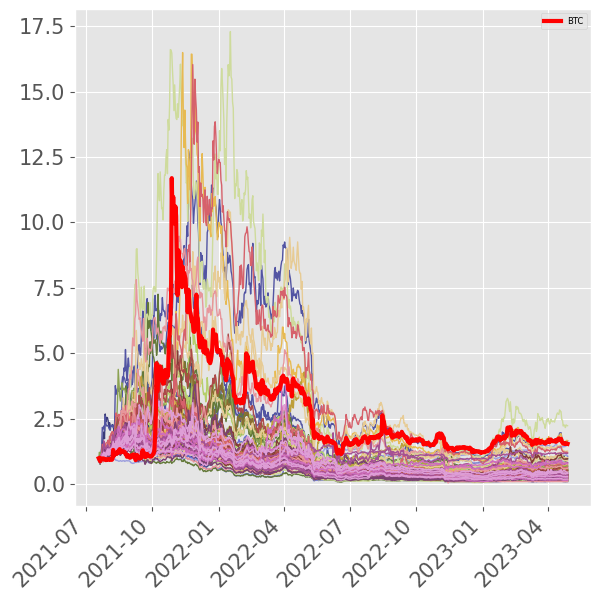}
\caption{\label{fig:cryptos} Normalized daily prices for $D=106$ crypto assets}
\end{figure}

\begin{figure}[t]
\centering
\includegraphics[width=0.7\textwidth,keepaspectratio]{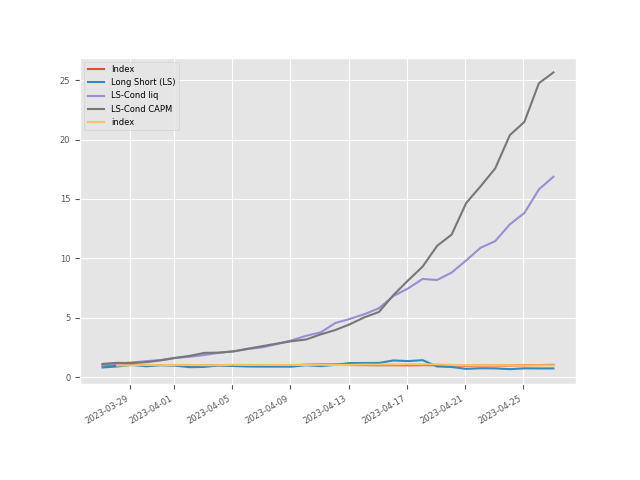}
\caption{\label{fig:plot363} Benchmark performance of the four portfolio strategies}
\end{figure}

In this sample, the conditional strategies (3)--(4) outperform their non-conditioned counterparts. This should be viewed as a validation of the conditional sampling algorithms rather than as an investable result: here the conditioning variables $b^n$ are observed at time $t^n$. In an operational setting, conditioning must be based only on information available at $t^{n-1}$ to avoid look-ahead bias. Determining whether there exists a rebalancing frequency below which the approach remains profitable is an open and practically relevant question.



\cleardoublepage

\printindex

\end{document}